\documentclass[12pt]{amsbook}
\usepackage{mathrsfs}
\usepackage{amsmath,amsthm,amsfonts}
\usepackage{mdframed}
\usepackage{amssymb}
\usepackage{amsfonts}
\usepackage{cite}
  \usepackage{url}
   \usepackage[a4paper, margin=2cm]{geometry}

\usepackage[dvips]{graphicx}
\usepackage{float}
\usepackage{dsfont}

 \usepackage{tikz}
\usetikzlibrary{arrows}

\newtheorem{theorem}{\bf Theorem}[chapter]
\newtheorem{lemma}[theorem]{\bf Lemma}

\theoremstyle{definition}
\newtheorem{definition}[theorem]{\bf Definition}

\theoremstyle{remark}
\newtheorem{remark}[theorem]{\bf Remark}

\newtheorem{corollary}[theorem]{\bf Corollary}

\newtheorem{proposition}[theorem]{\bf Proposition}

\def\upper {^{(i)}}
\def \msp {\mathsf{p}}

\def \bk {\mathbf{k}}
\def \fc {\mathfrak{c}}

\def \fkz {\mathfrak{z}}

\def \M{W}
\def \bh{\mathbb{H}}

\def \inv{^{-1}}

\def \half{\frac{1}{2}}

\def \upper{^{(i)}}

\def \p{\partial}
\def \a{a}

              \def \M{{\mathcal M}}

              \def \Z{{\bf Z}}

              \def \E{{\mathcal E}}

\def \fs {\mathbf s}
\def \ft {\mathbf t}
\def \p{\partial}
\def \E{{\mathcal E}}

\def \v {\vskip 0.1in}
\def \n {\noindent}

\numberwithin{section}{chapter} \numberwithin{equation}{chapter}
\numberwithin{figure}{chapter}

\begin{document}

\frontmatter

\begin{center}
  {\LARGE Lecture Notes on Relative Gromov-Witten Invariants \\}
  \end{center}
\bigskip
  \noindent
  \begin{center}
    {\large   An-Min Li, Li Sheng}\\

      Department of Mathematics, Sichuan University\\
        Chengdu, PRC\\[5pt]
\end{center}

\setcounter{page}{0} \tableofcontents

\mainmatter

\chapter*{Preface}

Almost twenty  years ago, Yongbin Ruan and the first author developed a theory of relative Gromov-Witten invariants and degeneration formulas.
Since then, the formulas has been tested many times. For algebraic case, an algebraic treatment of this theory was developed by Jun Li \cite{L}.

Here is the basic set-up of the theory. Let $(M, \omega)$ be a compact symplectic manifold of dimension $2n+2$, and $\widetilde M = H^{-1}(0)$ for a  local Hamiltonian function
as in the beginning of Section 3 in \cite{LR}.  Under the assumption that the  Hamiltonian vector field $X_H$
 generates a circle action on a neighborhood of $\widetilde M$, there
is a circle bundle  $\pi: \widetilde M \to Z =\widetilde M /S^1$ with a natural symplectic form $\tau_0$ on $Z$. To simplicity, we assume that $\widetilde M$ separates $M$ into two parts   to produce two cylindrical end  symplectic manifold  $M^+$ and $M^-$.  Collapsing the $S^1$-action
 at the infinity, we obtain the symplectic cuts $\overline{M}^+$ and  $\overline{M}^-$, both contain $Z$ as a  codimension two symplectic submanifold.
To obtain and prove the symplectic sum formula, we began with the following strategies.
\begin{enumerate}
\item[(A)]  We introduce the relative moduli spaces for symplectic pairs $(\overline{M}^\pm,Z)$   and  the moduli spaces on $M_\infty$.
 \item[(B)] Then we define the invariants for these moduli spaces, in particular,  including the relative GW invariants of  $(\overline{M}^\pm,Z)$.
\item [(C)]  We relate the Gromov-Witten invariants of $M$ with that of $M_\infty$.
\item[(D)]  Then we  relate the   Gromov-Witten invariants of $M_\infty$ with  the combination of relative invariants of $(\overline{M}^\pm,Z)$.
\end{enumerate}

\v
\v

A core technical issue in \cite{LR}  is to define invariants using virtual techniques. As we know, there had been several different approaches by the time, such as Fukaya-Ono(\cite{FO}), Li-Tian(\cite{LT}), Liu-Tian(\cite{LiuT}),
Ruan(\cite{R2}), Siebert(\cite{S}) and etc.
In \cite{LR}, they used Ruan's virtual neighborhood technique. As for all the other approaches, the smoothness of lower strata of virtual neighborhood is a subtle issue. Li-Ruan provided a much simpler approach by showing that the relevant differential form we try to integrate is in fact decay in certain rate   near lower strata of virtual neighborhood. Therefore, the integrand on top strata is independent of all choices and defines the desired invariants.
Namely, the contribution at lower strata with whatever the possible smooth structure can be ignored. Therefore, we avoided the smoothness problem of lower strata all together. At the time the theory was developed, the above
insight was treated as a technical advance and did not really catch the attention of larger community. With the renew interest on the technical detail of virtual technique during the recent years, Li-Ruan's technique seems to provide the  effcient way to treat the theory as well as many other applications.
\v
This is the draft of lecture  notes for Phd students in Sichuan University.   In this notes we expand  \cite{LR} with much more detailed explanations and calculations.

\chapter{Symplectic manifolds with cylindrical end}
\section{Symplectic cutting}\label{symplectic cutting}

\subsection{Symplectic cutting}\label{symplectic cutting-1}
We recall the construction of symplectic cuts,
\; i.e., a surgery along a
hypersurface which admits a local $S^1$-hamiltonian action ( see \cite{L}). By performing the
symplectic cutting we get two closed symplectic manifolds $\overline{M}^{+}$
and $\overline{M}^{-}$. The symplectic quotient $Z$ is embedded in both
$\overline{M}^{+}$ and $\overline{M}^{-}$ as symplectic submanifolds of
codimension 2.
\v
 Let $(M,\omega)$ be a compact symplectic manifold of dimension
$2n+2$. For simplicity, we assume that $M$ has a global Hamiltonian circle action.
Once we write down the construction, we then observe that a local circle
Hamiltonian action is enough to define a symplectic cutting.
\v
Let   $H:M\rightarrow {\mathbb{R}}$ be a Hamiltonian function such
that there is a small interval $I=(-\ell, \ell)$ of regular
values. Denote $\widetilde{M} =H^{-1}(0)$. Suppose that the
Hamiltonian vector field $X_H$ generates a circle action on
$H^{-1}(I)$. There is a circle bundle $\pi : \widetilde{M}
\rightarrow Z=\widetilde{M}/S^1 $ and a natural symplectic form
$\tau_0$ on $Z$.

  Consider the product manifold $(M\times {\mathbb{C}}, \omega\oplus -idz\wedge
d\bar{z})$. The moment map $F=H-|z|^2$ generates a Hamiltonian circle action
$e^{i\theta}(x, z)=(e^{i\theta}x, e^{-i\theta}z)$ and zero is a regular value of $F$.
We have symplectic reduction
\begin{equation}
\overline{M}^+=\{H=|z|^2\}/S^1,
\end{equation}
and a decomposition
\begin{equation}
\overline{M}^+=\{H=|z|^2\}/S^1=\left(\{H=|z|^2>0\}/S^1\right)\cup
\left(H^{-1}(0)/S^1\right).
\end{equation}
Furthermore,
\begin{equation}
\phi^+: \{H>0\}\rightarrow \{H=|z|^2>0\}/S^1
\end{equation}
defined by
\begin{equation}
\phi^+(x)=(x, \sqrt{H(x)})
\end{equation}
is a symplectomorphism.
Let
\begin{equation}
M^+_b=H^{-1}(\geq 0).
\end{equation}
Then, $M^+_b$ is a manifold with boundary and there is a map
\begin{equation}
M^+_b\rightarrow \overline{M}^+.
\end{equation}
Clearly,  $\overline{M}^+$ is obtained by
collapsing the $S^1$ action of the $H^{-1}(0)$.

To obtain $\overline{M}^-$, we consider circle action $e^{i\theta}(x,z)=(e^{i\theta}x,
e^{i\theta}z)$ with the moment map $H+|z|^2$. $\overline{M}^+, \overline{M}^-$
are called symplectic
cutting of $M$. We define $M^-_b$ similarly. By the construction, $Z=H^{-1}(0)/S^1$ with induced
symplectic structure embedded
symplectically into $\overline{M}^{\pm}$. Moreover, the normal bundles $\mathcal{N}^\pm$ of $Z$ in $\overline{M}^\pm$
satisfy $ \mathcal{N}^+=(\mathcal{N}^-)\inv.$ We call such an intersection pair  a {\em degenerated} symplectic manifold
and denote it by
\begin{equation}\label{eqn_1.6}
\overline{M}^+\cup_Z \overline{M}^-.
\end{equation}
There is a map
\begin{equation}
\pi: M\rightarrow \overline{M}^+\cup_Z \overline{M}^-.\end{equation}
Clearly, we only need a local $S^1$-Hamiltonian action.

\v

\subsection{Symplectic relative pair}
A symplectic relative pair $(X, B)$ is a symplectic manifold
$(X,\omega)$ together with  a symplectic divisor or
codimension two symplectic submanifold $B$ in $X$.
We can standardize the local structure around $B$. The normal bundle $\mathcal N:=\mathcal N_{B|X}$ may be identified with the complementary symplectic bundle. Note that the restriction of $\omega$ to $\mathcal N$ is a symplectic form. Pick a compatible almost complex structure $J$ on $\mathcal N$ such that
$\mathcal N$ is a Hermitian line bundle, we have a metric $<\cdot>$ on $\mathcal N$. Its principal $S^1$-bundle $Y$ is the
 unit circle bundle over $B$ where $S^1$ acts
 as complex multiplication. Then
$
\mathcal N=Y\times_{S^1}\mathbb C.
$

On $Y$, there is a connection 1-form
$\lambda$ which is dual to the vector field $T$ generated by the $S^1$-
action.
Let $\omega_B$ be the symplectic form on $B$, and $\pi:\mathcal N\to B$ be the projection.
\begin{equation}\label{eqn_1.1}
\omega_o:=\pi^\ast\omega_B+\frac{1}{2}d(\rho^2\lambda)
\end{equation}
defines a  form on $\mathcal N\setminus\{B\}$. Here, we take $B$ to be the
 $0$-section,  and $\rho$ to be the radius function on $\mathbb C$.
This form  can be extended over $\mathcal N$ and
 it is a symplectic form over $\mathcal N$. The $S^1$ action is
Hamiltonian in the  sense:
$
i_T\omega_0= -\frac{1}{2}d\rho^2.
$

Let $\mathbb D_\epsilon\subset \mathbb C$ be the  disk of radius $\epsilon$,
$\mathbb D$ be the unit disk and $\mathbb D^*=\mathbb D \setminus \{0\}$. We have the
 following sub-bundles of $\mathcal{N}$:
$$
\mathbb D_\epsilon \mathcal N= Y\times_{S^1} \mathbb D_\epsilon,\;\;\;
\mathcal N^\ast= Y\times_{S^1}\mathbb C^\ast,\;\;\;
\mathbb D^\ast_\epsilon \mathcal N=Y\times_{S^1}\mathbb D^\ast_\epsilon.
$$
The projective completion of $\mathcal{N}$ is
$
Q=Y\times_{S^1} {\mathbb{ CP}}^1.
$
In algebraic situation, $Q=\mathbb P(\mathcal N\oplus \mathbb C)$.
It contains two special sections: the 0-section and the
$\infty$-section, denoted by $B_0$ and $B_\infty$ respectively.
Both of them are identified with $B$.

By the symplectic neighborhood theorem, there
 exists a neighborhood $U\subset X$ of $B$ such that
$
(U,\omega)\cong (\mathbb D_\epsilon \mathcal{N},\omega_o)
$
for some $\epsilon>0$.
Here, $\omega_o$ is given in \eqref{eqn_1.1}.
We normalize the local structure near $B$ such that a neighborhood $U\subset X$ of
$B$ satisfies
\begin{equation}\label{eqn_1.2}
(U,\omega)
\cong (\mathbb D_{\epsilon} \mathcal{N},\omega_o).
\end{equation}
Hence a tubular neighborhood of $B$ is modeled on a neighborhood of $Z$ in $\overline{M}^{+}$ or $\overline{M}^{-}$.
\v
\section{Line bundles over $\overline{M}^{\pm}$ and $\overline{M}^+\cup_Z \overline{M}^-$}\label{line bundle}

Let $(M,\omega)$ be a compact symplectic manifold, $H:M\rightarrow {\mathbb{R}}$  a local $S^1$-Hamiltonian function as in \S\ref{symplectic cutting}. We have symplectic quotient $(Z, \tau_0)$ and normal bundles $\mathcal{N}^{\pm}$ such that
$$
\pi: M\rightarrow \overline{M}^+\cup_Z \overline{M}^-,\;\; \omega^+\mid_{TZ}=\omega^-\mid_{TZ},\;\; \mathcal{N}^+=(\mathcal{N}^-)\inv.
$$
Moreover we have compartible triple $(\tau_0, \tilde{J}, \tilde{g})$ on $Z$.

We can slightly deform $\omega$ to get a rational class $[\omega^*]$ on $M$. By taking multiple, we can assume that $[\omega^*]$ is an integral class on $M$. Therefore, it is the Chern class of a complex line bundle $L$ over $M$.
\v
Similarly, we slightly deform $\omega^{\pm}$ on $\overline{M}^{\pm}$ to get a rational class $[\omega^{\pm*}]$ on $\overline{M}^{\pm}$ such that $\tau_0^{\pm*}:=\omega^{\pm*}\mid_{TZ}$ is nondegenerate, i.e., $(Z, \tau_0^{\pm*})$ is a symplectic submanifold in $(\overline{M}^{\pm}, \omega^{\pm*})$. By taking multiple, we can assume that $[\omega^{\pm*}]$ is an integral class on $\overline{M}^{\pm}$. Therefore, it is the Chern class of a complex line bundle $L^{\pm}$ over $\overline{M}^{\pm}$.

\begin{lemma}  \label{line bundle-1}
We may choose $\omega^{\pm*}$ and the almost complex structure $J^*$ such that
$$\omega^{+*}\mid_{TZ}=\omega^{-*}\mid_{TZ},\;\;\mathcal{N}^+=(\mathcal{N}^-)\inv,\;\;
L^{+}\mid_{Z}=L^{-}\mid_{Z}.$$
\end{lemma}
 \v\n
 {\bf Proof.} The imbedding $i_{\pm}: Z\to \overline{M}^{\pm}$ induces
$$H^2(\overline{M}^{\pm}, \mathbb{R})\to H^2(Z, \mathbb{R}),\;\;\alpha\longmapsto i_{\pm}^*\alpha.$$
The image $i_{\pm}^*(H^2(\overline{M}^{\pm}, \mathbb{R}))$ is a linear subspace in $H^2(Z, \mathbb{R})$. Then $i_{+}^*(H^2(\overline{M}^{+}, \mathbb{R}))\bigcap i_{-}^*(H^2(\overline{M}^{-}, \mathbb{R}))$ is a linear subspace in $H^2(Z, \mathbb{R})$. Since there are $[\omega^+]\in H^2(\overline{M}^{+}, \mathbb{R})$, $[\omega^-]\in H^2(\overline{M}^{-}, \mathbb{R})$ and $[\tau_0]\in H^2(Z, \mathbb{R})$
such that $\tau_0=i_{\pm}^*(\omega^{\pm})\mid_{TZ}$ and $\tau_{0}$ is   a closed  non-degenerate $2$-form on $Z$, we have
$$dim\left(i_{+}^*(H^2(\overline{M}^{+}, \mathbb{R}))\bigcap i_{-}^*(H^2(\overline{M}^{-}, \mathbb{R}))\right)\ne 0.$$
Let $\alpha^{\pm}\in H^2(\overline{M}^{\pm}, \mathbb{Z})$ be a non-degenerate form, by
the local normal form \eqref{eqn_1.1} we conclude that $i_{\pm}^*\alpha^{\pm}\in H^2(Z, \mathbb{Z})$.
\v
We can choose an integral base $\zeta_1,...,\zeta_m$ in $i_{+}^*H^2(\overline{M}^{+}, \mathbb{R}))\bigcap i_{-}^*(H^2(\overline{M}^{-}, \mathbb{R})).$
 We choose integral bases
$$e_{m+1},...,e_k\in H^2(\overline{M}^{+}, \mathbb{Z})  ,\;\;\eta_{m+1},...,\eta_l\in H^2(\overline{M}^{-}, \mathbb{Z})$$
such that
$\zeta_1,...,\zeta_m, i_{+}^*e_{m+1},...,i_{+}^*e_{k}$ (resp. $\zeta_1,...,\zeta_m, i_{-}^*\eta_{m+1},...,i_{-}^*\eta_{l}$) is an integral base in $i_{+}^*(H^2(\overline{M}^{+}, \mathbb{R}))$ (resp. $i_{-}^*(H^2(\overline{M}^{-}, \mathbb{R}))$). Note that $i_{+}^*(H^2(\overline{M}^{+}, \mathbb{R}))\bigcap i_{-}^*(H^2(\overline{M}^{-}, \mathbb{R}))\thickapprox \sum_{i=1}^m \mathbb{R}\zeta_i$. Then we can easily find integral class $[\omega^{\pm*}]$ such that $\omega^{+*}\mid_{TZ}=\omega^{-*}\mid_{TZ}:=\tau_0^*$ and $\tau_0^*$ is a closed  non-degenerate $2$-form on $Z$. There are line bundle $L^{\pm}$ over $\overline{M}^{\pm}$ such that $L^{+}\mid_{Z}=L^{-}\mid_{Z}.$
\v
We choose an almost complex structure $\tilde{J}^*$ on $Z$ such that
$(\tau_0^*, \tilde{J}^*, \tilde{g})$ is a compartible triple on $Z$. $\tilde{J}^*, \tilde{g}$ are lifted to $(J,g)$ on $\mathcal{N}^{\pm}$ in a natural way such that, restricting to the fibre, $J$ is $\sqrt{-1}$ and $g$ is the standard Euclidean metric. Then $\mathcal{N}^+=(\mathcal{N}^-)\inv$.\;\;\;$\Box$.
\v
So there is a symplectic form $\omega^*$ on $M$ and a local circle
Hamiltonian action on $(M, \omega^*)$ such that $\overline{M}^+\cup_Z \overline{M}^-$ can be obtained by symplectic cutting from the local circle
Hamiltonian action on $(M, \omega^*)$. Then there is a
line bundle $L$ over $M$ such that $L=\pi^*L^{\pm}$.
\v

\section{Cylindrical almost complex structures}\label{cylin almost struc}

Let $H:M\rightarrow {\mathbb{R}}$ be a local Hamiltonian function defined on $M$ such
that there is a small interval $I=(-\ell, \ell)$ of regular
values. Denote $\widetilde{M} =H^{-1}(0)$. We may choose a connection 1-form $\lambda $ on
$\widetilde{M}$ such that $\lambda(X_H)=1$ and $d\lambda$
represents the first Chern class for the circle bundle (see
\cite{MS1}). Denote $\xi=ker(\lambda)$. Then $\xi$ is an
$S^1$-invariant distribution and $(\xi,
\pi^{\ast}\tau_0)\rightarrow \widetilde{M}$ is a $2n$-dimensional
symplectic vector bundle. We identify $H^{-1}(I)$ with $I \times
\widetilde{M}$. By a uniqueness theorem on symplectic forms (see
\cite{MS1}) we may assume that the symplectic form on
$\widetilde{M}\times I$ is expressed by
\begin{equation}
 \omega =
\pi^{\ast}(\tau_0 + yd\lambda) - \lambda \wedge dy.
\end{equation}
We assume that the hypersurface $\widetilde{M}=H^{-1}(0)$ devides
$M$ into two parts $M^+$ and $M^-$, which can be written as
$$M^{+}_{0}\bigcup\left\{[-\ell,0)\times \widetilde{M}\right\},$$
$$M^{-}_{0}\bigcup\left\{(0,\ell]\times \widetilde{M}\right\},$$
where $M^{+}_{0}$ and $M^{-}_{0}$ are compact manifolds with
boundary. We mainly discuss $M^+$; the discussion for $M^-$ is
identical. Fix a
function $\phi_{0} :[0, \infty)\rightarrow [- \ell, 0)$ satisfying
\begin{equation}
(\phi_{0})^{\prime}>0,
\;\;\phi_{0}(0)=-\ell, \;\;\phi_{0}(a)\rightarrow 0 \;\;as\;\; a
\rightarrow \infty.
\end{equation} Through $\phi_{0}$ we consider
$M^+$ to be $M^{+} =  M^{+}_{0}\bigcup\{[0, \infty)\times
\widetilde{M}\}$ with symplectic form
$\omega_{\phi_{0}}|_{M^{+}_0}=\omega$, and over the cylinder ${\mathbb
R}\times \widetilde{M}$
\begin{equation}\omega_{\phi_{0}} = \pi^{\ast}(\tau_0 +
\phi_{0}d\lambda) - (\phi_{0})^{\prime}\lambda \wedge d a.
\end{equation}
 Moreover, if we choose the origin of ${\mathbb
R}$
tending to $\infty $, we obtain ${\mathbb
R}\times \widetilde{M}$ in
the limit. Denote by $\mathbb P(\mathcal N\oplus \mathbb C)$ the projective completion of the normal bundle $\mathcal N_{b}={\mathbb
R}\times \widetilde{M}$ with a zero section $Z_{0}$ and an infinity section $Z_{\infty}.$ Topologically, the space $\mathbb P(\mathcal N\oplus \mathbb C)$ can be also obtained from ${\mathbb
R}\times \widetilde{M}$ by collapsing the $S^1$-action on the $\pm\infty$ ends.

\v

We choose a compatible almost complex
structure $\widetilde{J}$ on $Z$ such that
$$g_{\widetilde{J}(x)}(h,k)=\tau_0(x)(h,\widetilde{J}(x)k)\;\;\;\forall\;h,\; k\; \in TZ$$
defines a Riemannian metric. $\widetilde{J}$
and $g_{\widetilde{J}}$ are lifted in a natural way to $\xi $. We
define an almost complex structure $J$  on ${\mathbb R}\times
\widetilde{M}$ as follows:
\begin{eqnarray}\label{complex_structure_I}
J\frac{\partial}{\partial
a}=X_{H},&\;\;\;& JX_{H}=-\frac{\partial}{\partial a},\\
\label{complex_structure_II}
J\xi=\xi, &\;\;\;&J|_{\xi}=\widetilde{J}.\end{eqnarray}
\vskip 0.1in \noindent We denote by $N$
one of $M^{+}$, $M^-$ and ${\mathbb
R}\times\widetilde{M}$,   denote by $\overline{N}$
one of $\overline{M}^{+}$, $\overline{M}^-$ and ${\mathfrak
R}$. We may
choose an almost complex structure $J$ on $N$ such that
\begin{itemize}
  \item[(i)] $J$
is tamed by $\omega_{\phi_0} $ in the usual sense, \par \noindent
 \item[(ii)] Over the tube ${\mathbb R }\times \widetilde {M},$  \eqref{complex_structure_I} and
\eqref{complex_structure_II} hold.
\end{itemize}

\v

Since
$g_{\widetilde{J}}$ is positive, and $d\lambda$ is a 2-form on $Z$
(the curvature form), by choosing $\ell$  small enough we may
assume that $\widetilde{J}$ is tamed by $\tau_0 + yd\lambda $ for
$|y| < \ell $, and there is a constant $C>0$ such that
\begin{equation}\label{tame}
\tau_0(v,\widetilde{J}v)\leq C\left(\tau_0(v,\widetilde{J}v)+
yd\lambda(v, \widetilde{J}v)\right)
\end{equation} for all $v\in TZ,
|y|\leq \ell $. Then $J$ is $\omega_{\phi_0}$-tamed  over the tube.

\vskip 0.1in \noindent
 Then
\begin{equation}\label{omega_forms}
\langle v,w\rangle_{\omega_{\phi_0}} = \frac{1}{2}\left(
\omega_{\phi_0} (v,Jw) + \omega_{\phi_0} (w,Jv)\right) \;\;\;\;\;
\forall \;\; v, w \in TN
\end{equation}
defines a Riemannian metric
on $N$. Note that $\langle \;,\;\rangle_{\omega_{\phi_0}}$ is not
complete.
The length of every orbit of the $S^1$ action on $\widetilde {M}$ with
respect to the metric $\langle\;,\;\rangle_{\omega_{\phi_0}}$ is $\phi_0^{\prime}$,
which converges to zero as $a \rightarrow \pm \infty $. Hence we can view
$\overline{M}^{\pm}$ as the completions of $M^{\pm}$.
\v

We choose another metric $\langle \;,\;\rangle$ on $N$
such that
\begin{equation}\label{omega_forms_1}
\langle \;,\;\rangle = \langle
\;,\;\rangle_{\omega_{\phi_0}} \;\;\;\;on \;\; M^{\pm}_{0}
\end{equation}
 and over
the tubes
\begin{equation} \label{omega_forms_on_tubes}
\langle(a,v),(b,w) \rangle= ab + \lambda (v)\lambda
(w) + g_{\widetilde{J}}(\Pi v, \Pi w),
\end{equation}
 where we
denote by $\Pi:T\widetilde{M}\rightarrow\xi$ the projection along
$X_H$. It is easy to see that $\langle \;,\;\rangle$ is a complete
metric on $N$. \vskip 0.1in \noindent

\begin{remark} To get a line bundle $L$ over $M^+$ we slightly deform $\omega$ on $\overline{M}^{+}$ and take multiple to get integral symplectic forms $\omega^*$ on $\overline{M}^+$ and $\tau_o^*$ on $Z$. Then we choose the metric $g^*$ on $Z$ such that $\widetilde{J}$ is invariant. It follows that $J$ is invariant.
\end{remark}

\vskip 0.1in \noindent
\v We write
$$M^+=M^{+}_{0}\bigcup\left\{[0,\infty)\times \widetilde{M}\right\}.$$
This means we have fixed a coordinate $a$ over the cylinder part.
Similarly, we write
$$M^-=M^{-}_{0}\bigcup\left\{(-\infty,0]\times \widetilde{M}\right\}.$$
\v

\section{$J$-holomorphic maps with finite energy}

Let $(\Sigma,j)$ be a
compact Riemann surface and $P\subset\Sigma$ be a finite
collection of points. Denote $\stackrel{\circ}{\Sigma}
=\Sigma\backslash P.$ Let ${u}:\stackrel{\circ}{\Sigma}
\rightarrow N$ be a $(J,j)$-holomorphic map, i.e., ${u}$ satisfies
\begin{equation}\label{j_holomorphic_maps}
d{u}\circ j={J}\circ d{u}.
\end{equation}  Alternatively
\begin{equation}\label{j_holomorphic_maps-1}
\overline{\partial}_{j,J}u=\half (du+J\circ du \circ j)=0.
\end{equation}

To simplify notations we write $(J,j)$-holomorphic map as
$J$-holomorphic map later.

\v
Denote $$\Phi^+ =\left \{ \phi :[0,
\infty )\rightarrow [-\ell, 0) | \phi^{\prime}
> 0 \right \}.$$
For ${\mathbb
R}\times \widetilde{M}$ let $\ell_1 < \ell_2$ be two real numbers satisfying $ - \ell
< \ell_1 < \ell_2 \leq 0 .$ Denote by $\Phi_{\ell_1,\ell_2}$ the
set of all smooth functions $\phi:{\mathbb
R}\rightarrow
(\ell_1,\ell_2)$ satisfying $$\phi^{\prime}>0,\;\;
\phi(a)\rightarrow \ell_2\;\:\;{\rm as} \;a\rightarrow \infty ,\;\;
\phi(a)\rightarrow \ell_1\;\:\;{\rm as} \;a\rightarrow -\infty.$$ To
simplify notations we use $\Phi $ to denote both $\Phi^+$ and
$\Phi_{\ell_1,\ell_2}$, in case this does not cause confusion.
\v
Following \cite{HWZ1} we
impose an energy condition on $u$.
For any $J$-holomorphic map
$u:\stackrel{\circ}{\Sigma}\rightarrow N$ and any $\phi \in \Phi $
the energy $E_{\phi}(u)$ is defined by
\begin{equation}\label{definition_of_energy}
E_{\phi}(u)=\int_{ {\Sigma}}u^{\ast}\omega_{\phi}.
\end{equation}
Let $z=e^{s+2 \pi it}.$ One
computes over the cylindrical part
\begin{equation}\label{omega_forms_cylinder}
u^{\ast}\omega_{\phi}=
\left(\tau_0 + \phi d\lambda \right)\left((\pi\widetilde{u})_s,
(\pi\widetilde{u})_t\right) + {\phi}^{\prime}(a^{2}_s + a^{2}_t
))ds\wedge dt,
\end{equation} which is a nonnegative integrand.
In fact, by \eqref{tame} we have
\begin{equation}
 \tau_0 \left((\pi\widetilde{u})_s,
 \tilde J(\pi\widetilde{u})_s\right)\leq C\left(\tau_0 +
yd\lambda \right)\left((\pi\widetilde{u})_s,
 \tilde J(\pi\widetilde{u})_s\right)
\end{equation} for all $v\in TZ,
\phi\in \Phi $.

\v
A $J$-holomorphic map $u:\stackrel{\circ}{\Sigma} \rightarrow N $
is called a finite energy $J$-holomorphic map if
\begin{equation}\label{finite_energy_j_holomorphic_maps}
\sup_{\phi \in
\Phi }\left \{\int_{{\Sigma}}u^{*} \omega_{\phi}
\right \}<\infty.
\end{equation} We shall see later that the
condition is natural in view of our surgery. For a $J$-holomorphic
map $u:{\Sigma} \rightarrow {\mathbb
R}\times\widetilde{M}$  we write $u=(a, \widetilde{u})$ and define
\begin{equation}\label{definition_of_energy_on_complex_manifolds}
\widetilde{E}(u)=\int_{{\Sigma}}\widetilde{u}^{\ast}
(\pi^{\ast}\tau_0).
\end{equation}
Denote
$$\widetilde{E}(s)=\int_s^{\infty}\int_{S^1}\widetilde{u}^{\ast}(\pi^{\ast}
\tau_0).$$ Then $$\widetilde{E}(s)=\int_s^{\infty}\int_{S^1}|\Pi
\widetilde{u}_t|^2dsdt,$$
\begin{equation}\label{deriative_of_energy_on_Z}
\frac{d\widetilde{E}(s)}{ds}=-\int_{S^1}|\Pi\widetilde{u}_t|^2dt.
\end{equation}

 \v\n By using the same method as in \cite{H}, one can prove
the following three lemmas: \vskip 0.1in \noindent
\begin{lemma}\label{zero_energy_classification}
 \begin{itemize}
 \item[(1)] Let $u=(a,\widetilde{u}):\mathbb{C} \rightarrow {\mathbb R}
\times \widetilde{M}$ be a $J$-holomorphic map with finite energy.
If $\int_{ \mathbb{C}}\widetilde{u}^{\ast}(\pi^{\ast}\tau_0)=0$, then
$u$ is a constant.
\item[(2)] Let $u=(a,\widetilde{u}):{
\mathbb{R}}\times S^1 \rightarrow {\mathbb R}\times \widetilde{M}$ be a
$J$-holomorphic map with finite energy. If $\int_{\mathbb{R}\times
S^1}\widetilde{u}^{\ast}(\pi^{\ast}\tau_0)=0$, then
$(a,\widetilde{u})=(ks+c, x(kt+d))$, where $x$ is a periodic orbit, $k\in {\mathbb Z}$, $c$ and
$d$ are constants.
\end{itemize}
\end{lemma}
\vskip 0.1in
 \noindent
\begin{lemma}\label{sequence_convergence_of_j_holomorphic_maps} Let
$u=(a,\widetilde{u}):\mathbb{{C}}-D_1 \rightarrow {\mathbb R} \times
\widetilde{M}$ be a nonconstant $J$-holomorphic map with finite
energy. Put $z=e^{s+2\pi it}$. Then for any sequence
$s_i\rightarrow \infty $ , there is a subsequence, still denoted
by $s_i$, such that $$\lim_{i\rightarrow \infty}\widetilde
u(s_i,t)=x(kt)$$ in $C^{\infty}(S^1)$ for some $k$-periodic orbit
$x$.
\end{lemma}
\begin{lemma}\label{partial_a_theta_convergence_zero} Let $u=(a,\widetilde{u}):\mathbb{{C}}-D_1 \rightarrow {\mathbb R}  \times \widetilde{M}$ be a
$J$-holomorphic map with finite energy. Put $z=e^{s+2\pi it}$. Assume that there exists a sequence $s_{i}\to\infty$ such that $\widetilde { u} (s_{i},t)\longrightarrow x(kt)$ in $C^{\infty}(S^1,\widetilde {M})$ as $i\rightarrow \infty$ for some $k\in {\mathbb Z}$.
Then there are constants $c$ and
$d$ such that for all derivatives ${\bf n}=(m_{1},m_{2})\in {\mathbb Z_{\geq 0}^2 }$ with $|{\bf n}|\geq 1$
\begin{equation}
\label{deriatives_a_theta_convergence_zero}
|\partial^{ \bf n }[a-ks-c]|\longrightarrow 0,\;\;\;\;\; |\partial^{\bf n}[\theta-kt-d]|\longrightarrow 0,\;\;\;\;\;\; as \;\;s\longrightarrow \infty,
\end{equation}
uniformly in $t$.
\end{lemma}
 \vskip 0.1in
 \noindent

Let $\pi(x_k(t))=q\in Z$. By Darboux Theorem we can find a local coordinate system $(w_1,...,w_{2n})$ on an open set $\Im\subset Z$ near $q$, $q\in \Im$, such that on $\Im$
\begin{equation}\label{standard_tau}
\tau_0=\sum_{i=1}^n dw_i\wedge dw_{n+i}.
\end{equation}
Then we have a Darboux coordinate
\begin{equation}\label{darboux_coord}
(a,\theta,{\bf w})=(a,\theta,w_1,...,w_{2n}).
\end{equation}
Choose a local trivialization of $\widetilde{M}\rightarrow Z$
 on $\Im$ 
  such that
  \begin{equation} \label{x_Dauboux_coord}
  x_k=\{ 0\leq \theta \leq 1, w= 0 \},
  \end{equation}
   and
\begin{equation}\label{standard_lambda}
\lambda = d\theta + \sum b_i(w) dw_i,
\end{equation}
 where
$b_i(0)= 0$. Obviously, $\xi(\theta,0)$ is spanned by
$\frac{\partial} {\partial w_1},...,\frac{\partial}{\partial
w_{2n}}$. For $y$ small enough we may choose a frame
$e_1,...,e_{2n}$ for $\xi(\theta,w)$ as follows: in terms of the
coordinates $(\theta,w_1,...,w_{2n})$ we assume
 $e_i=( c_i(w),0,..,1,...,0),\;\; i=1,...,2n. $ By $\lambda(e_{i})=0,\;i=1,\cdots,2n,$ we have $$e_i=(-b_i(w),0,..,1,...,0),\;\; i=1,...,2n.$$ We write $$
u(s,t)=(a(s,t),\theta(s,t),{\bf w}(s,t)).$$
  Denote by $L$ the matrix of
the almost complex structure $\widetilde{J}$ on $\xi$ with respect
to the frame $e_1,...,e_{2n}$, and set $\widetilde{J}(s,t)
=L(u(s,t))$. Since $J\xi=\xi$ and $J|_{\xi}=\widetilde{J}$ we can assume that
$Je_{i}=\sum_{j} e_{j}c_{ji}.$ It follows that
\begin{equation}\label{J_action_y}
J(\frac{\partial }{\partial w_{i}})=-b_{i}\frac{\partial }{\partial a}-\sum c_{ji}b_{j}\frac{\partial }{\partial \theta}+\sum c_{ji}\frac{\partial }{\partial w_{j}},
\end{equation}
where we use $J(\frac{\partial }{\partial \theta})=-\frac{\partial }{\partial a}$.
 It follows from \eqref{j_holomorphic_maps} that
\begin{equation}\label{J_holomorphic_equ}
u_{s}+Ju_{t}=0,\;\;\;\; u_{t}-Ju_{s}=0.
\end{equation}
 By \eqref{J_action_y}
  we can write the equation \eqref{J_holomorphic_equ} as follows:
\begin{eqnarray}
\label{coordinates_a}
&&a_s=\lambda(u_t)= \theta_t + \sum b_i({\bf w})(w_{i})_t \\
\label{coordinates_theta}&&a_t=-\lambda(u_s)=- \theta_s - \sum b_i({\bf w})(w_{i})_s \\
\label{coordinates_y}&&w_s + \widetilde{J}(s,t)w_t=0,
\end{eqnarray}
where $(w_{i})_{s}=\frac{\partial w_{i}}{\partial s}$ and $(w_{i})_{t}=\frac{\partial w_{i}}{\partial t}$.

\vskip 0.1in \noindent

\begin{remark}\label{triv_cyli}
Let $(a,\theta,{\bf w})$ be a Darboux coordinates around $x(kt)$ with $\pi(x(kt))=0$. Then $(a,\widetilde{u})$  in  (2) of Lemma \ref{zero_energy_classification} can be rewritten as
$(a,\theta,{\bf w})(s,t)=(ks+c, kt+d,{\bf 0}).$  We call it a trivial cylinder.
\end{remark}

\v\n

\section{Exponential decay of $\widetilde {E}(u)$}

\v
Since $N$ is not compact, to compactify the moduli space of $J$-holomorphic maps with finite energy we have to analyse the
behaviour of holomorphic maps at infinity. In the
case of contact manifolds, a similar analysis for non-degenerare periodic orbit has been done by
Hofer and his collabrators \cite{HWZ1}, \cite{HWZ2} \cite{HWZ3}.
In \cite{LR} the authors adapt the standard $L^2$-moduli space theory, which
has been intensively developed for Chern-Simons theory, to study Bott-type periodic orbits, including the case
of a Hamiltonian $S^1$-action (see the first version \cite{LR1}). According to the suggestion of referee they deleted the part of contact geometry, but kept the $L^2$-theory for the case of a Hamiltonian $S^1$-action in final version. In the case of a Hamiltonian $S^1$-action the contact manifold
$\widetilde M$ is a circle bundle of a complex line bundle, the similar results can be proved
in a rather easy way, using estimates in \cite{MS}. For example, this point of view was employed in \cite{CLSZ}. In this note we choose this point of view.
\v

Denote $\mathbf{D}_{\delta}(0)=\{{\bf w}| \sum (w_i)^2<\delta\}.$

   \begin{proposition}\label{same orbit-1}
Let $u=(a,\widetilde{u}): [0, \infty)\times S^1 \rightarrow {\mathbb R}  \times \widetilde{M}$ be a
$J$-holomorphic map with finite energy, and $s_{i}\to\infty$ be a sequence such that $\widetilde { u} (s_{i},t)\longrightarrow x(kt)$ for some periodic orbit $x(kt)$. Then
for any disk $\mathbf{D}_{\delta}(0)$, there is a $\aleph>0$ such that if $s>\aleph$
then $\pi\circ\widetilde{u}(s,\cdot)\in \mathbf{D}_{\delta}(0).$
\end{proposition}
{\bf Proof.} Take the coordinates transformation  $z=e^{-s-2i\pi t},$ $  \hat{u}(z):=  \hat{u}(s,t)$. Consider the $\tilde{J}$-holomorphic map $\hat{u}=\pi\circ \widetilde{u}:D^{*}_1(0)\rightarrow Z$.
By the Theorem of removal of singularities,
$\hat u$ can be extended to a $\tilde{J}$ holomorphic map
$$\hat{u}:D_{1}(0) \rightarrow Z.$$
Then for any $\delta>0,$ there exists a $\aleph>0$ such that $\hat{u}(D_{e^{-N}}(0))\subset \mathbf{D}_{\delta}(0)$ for all $N>\aleph$,
 where $\mathbf{D}_{\delta}(0)$ is a ball of $\pi x$ in $Z.$  \;\;\;\;$\Box$

\v
For any loop $\gamma:S^1\rightarrow \widetilde{M}$ let $\gamma^*:=\pi\circ\gamma$. Suppose that $\gamma^*(t)$ lies in $\mathbf{D}_{\delta}(0)$.
Set $\gamma^*(t)=(w_1(t),...,w_{2n}(t))$ and
put $\zeta(t):=\gamma^*(t)-0=(w_1(t),...,w_{2n}(t)).$
We define an annulus $W:   [0, 1] \times S^1\rightarrow \widetilde M $ by by $\varpi \zeta(t),$ and define an action functional by
\begin{equation}\label{action-1}
{\mathcal A}(\gamma^*)=-\int_{[0,1]\times S^1 }W^*\tau_0.
\end{equation}
\begin{lemma} \label{inequality-1}
There is a constant $\mathbf{C}_{1}>0$ depending only on $\tilde{J}$ on $Z$ such that for any smooth loop $\gamma^*(t): S^1\rightarrow \mathbf{D}_{\delta}(0)$
\begin{equation}\label{inequality-2}
\mid{\mathcal A}(\gamma^*)\mid \leq \mathbf{C}_{1}\int_{S^1 } |\tfrac{d}{dt}\gamma(t)|^2_{\tilde{J}}.\end{equation}
\end{lemma}
\v\n
The proof is standard (see \cite{MS}).

\v
By Stokes theorem we have $\widetilde{E}(\hat{u},s)={\mathcal A}(\gamma^*_s):={\mathcal A}(s)$ for any $s>\aleph$.
\begin{lemma}  \label{decay of E-1}
There is $s_o$ and constants $\mathfrak{c}_1>0$, $\mathbf{C}_2>0$ such that for any $s >s_o$ we have
\begin{equation}
\widetilde{E}(s)\leq \widetilde{E}(s_{o})e^{-\mathfrak{c}_1( s- s_{o})}.
\end{equation}
\begin{equation}\label{upper_bound_estimates-1}
\int_{s_o}^{s}\|\Pi\widetilde{u}_t\|_{L^2(S^1)}ds
\leq \mathbf{C}_2(\widetilde{E}(s))^{1/2}.
\end{equation}
\end{lemma}
{\bf Proof.}
 Let $s_i$ be a sequence $s_{i}\to\infty$ such that $\tilde { u} (s_{i},t)\longrightarrow x$. Letting $i\rightarrow \infty$ we get
\begin{equation}\label{energy and function}
\widetilde{E}(s):=\widetilde{E}(u;[s,\infty)\times S^1)={\mathcal A}(s).
\end{equation}
Then
\begin{eqnarray}\label{estimates_deriative_of_energy_on_Z_I}
\frac{d\widetilde{E}(s)}{ds}=
-\int_{S^1}|\Pi \widetilde{u}_{t}|^{2}dt\leq -  \mathbf{C}_{1}^{-1}\mathcal{A}(s)= -\mathbf{C}_{1}^{-1}\widetilde{E}(s).\nonumber
\end{eqnarray}
\begin{equation}
\frac{d\widetilde{E}(s)}{ds}\leq -\sqrt{\mathbf{C}_{1}^{-1}} \|\Pi
\widetilde{u}_t\|_{L^2(S^1)} \widetilde{E}(s)^{1\over
2} \nonumber.\end{equation}
It follows that for any $s_1> s \geq s_{o}$
$$\widetilde{E}(s)\leq \widetilde{E}(s_{o})e^{-\mathfrak{c}_1( s- s_{o})},
$$
$$
\int_{s}^{s_1}\|\Pi\widetilde{u}_t\|_{L^2(S^1)}ds\leq {C}_1\left(\widetilde{E}(s)^{1/2}-\widetilde{E}(s_1)^{1/2}\right)
\leq \mathbf{C}_2(\widetilde{E}(s))^{1/2}.
$$
for some constants $\mathfrak{c}_1>0$, $\mathbf{C}_2>0$.\;\;\;\;$\Box$

\v

\section{Convergence to periodic orbits}

We are interested in the
behaviors of the finite energy $J$-holomorphic maps near a
puncture $p$. There are two different types of puncture : the
removable singularities and the non-removable singularities. If
${u}$ is bounded near a puncture, then this puncture is a
removable singularity. In the following, we assume that all
punctures in $P$ are non-removable. Then ${u}$ is unbounded near
the punctures. We mainly discuss $J$-holomorphic maps into ${\mathbb R}\times \tilde{M}$, for $M^{\pm}$ the discusses
are the same.

We fix a Darboux coordinate system $a,\theta,{\bf w}$ as \eqref{darboux_coord}, where ${\bf w}$ is a local coordinates near $\pi (x_k) \in Z$.
We need the following lemma of \cite{HWZ1}.
\begin{lemma}\label{a_lemma_of_hopf}
 Assume $v:[s_o,\infty)\times S^{1}\longrightarrow R^2$ is smooth, bounded and solves the equation
$$
v_{s} +J_{0}v_{t}=g,\;\;\;\;\;\;\;\; where \;\;\; \|g(s)\|\leq C_{o}e^{-\delta s},
$$
for some $\delta>0,$ where the norm is the $L^2(S^1)$-norm. If $v$ satisfies $v_{t}(s,t)\longrightarrow 0$ as $s\longrightarrow \infty$ uniformly in $t$, and moreover has vanishing mean values,
$$
\int_{0}^{1}v(s,t)dt\equiv 0,
$$
then
$$
\int_{s_o}^{s}e^{2\rho s}\|v(s)\|^2ds<\infty,\;\;\;\;\;\int_{s}^{s+1}\|v(s)\|^2ds\leq Ce^{-2\rho s}
$$
for every $0\leq \rho<\delta$ and $\rho<\frac{1}{2},$ where $C>0$ is a constant depending only on $C_{o},\frac{1}{2}-\rho$ and $( v(s_{o}), J_{0}v_{t}(s_{o})).$
\end{lemma}
\vskip 0.1in
 \noindent
 For the reader's convenience we give the proof here.
 \v\n
{\bf Proof. } We first show that $\|g(s)\|\in L^2$ implies $\|v(s)\|\in L^2,$ the norm denoting the $L^2(S^1)$-norm. We make use of the following pointwise identities for a function $w=w(s,t):$
\begin{eqnarray}\label{hwz1}
&& 2\langle w_s, J_{0}w_t\rangle=\frac{d}{ds}\langle w,J_{0}w_t \rangle-\frac{d}{dt}\langle w, J_{0}w_s\rangle \\
&& |w_s|^2+|w_t|^2=|w_s+J_{0}w_t|^2-2\langle w_s, J_{0}w_{t} \rangle.\label{hwz2}
\end{eqnarray}
Since $v$ has mean values zero we can estimate $\|v(s)\|\leq \|v_{t}(s)\|.$ Using  \eqref{hwz1},
integrating by parts, and observing that the integral of the derivative of a periodic function over a period vanishes, and $v$ solves the equation $v_{s}+J_{0}v_{t}=g,$ we obtain
\begin{align}
\int_{s_o}^{s} \|v(s)\|^2ds &\leq \int_{s_o}^{s}(\|v_s(s)\|^2+\|v_t(s)\|^2  )ds \nonumber\\
 & =\int_{s_o}^{s}\|g(s)\|^2ds -( v(s),J_{0}v_{t}(s))+( v(s_{o}),J_{0}v_{t}(s_{o}) ),\nonumber
\end{align}
where $(,)$ denotes the inner product in $L^2(S^1)$.  Since $\|g(s)\|\in L^2$ we conclude for the limit $s\longrightarrow \infty:$
$$
\int_{s_o}^{\infty} \|v(s)\|^2 \leq \int_{s_o}^{\infty} \|g(s)\|^2 +( v(s_{o}), J_{0}v_{t}(s_{o})).
$$
Take now an increasing sequence of monotone increasing functions $\gamma_{n}:\mathbb R \longrightarrow \mathbb R $ satisfying $\gamma_{n}(s)=s$ for $0\leq s\leq n,\;0\leq \gamma'_{n}(s)\leq 1$ for $s\in \mathbb R$, and $\gamma_{n}(s)=const$ for $s\geq n+1.$ Let $\rho>0$ and define the sequence $\hat{v}_{n}=\hat{v}$ as
$$
\hat{v}(s,t) =e^{\rho \gamma_{n}(s)}v(s,t).
$$
Then $\hat{v}$ is smooth, bounded, satisfies $\hat{v}_{t}(s,t)\longrightarrow 0$ as $s\longrightarrow \infty,$ has mean values zero and $\|\hat{v}(s)\|\in L^2$. Differentiating $\hat v$ we obtain
$$
\hat{v}_{s}+J_{0}\hat{v}_{t} =e^{\rho \gamma_{n}(s)}g+ \rho \gamma'_{n}(s)\hat{v}.
$$
If $0<\rho<\delta$ we conclude, in view of the exponential decay of $g$, for $n\geq s_{o}$
\begin{align}
\int_{s_o}^{\infty} \|\hat {v}(s)\|^2ds & \leq \int_{s_o}^{\infty} \|\hat {v}_{s}+J_{0}\hat {v}_{t}\|^2ds+( \hat{v}(s_o),J_{0}\hat {v}_{t}(s_{o}) ) \nonumber\\
& \leq 2\int_{s_o}^{\infty}e^{2\rho s} \|g(s)\|^2ds +2\rho \int_{s_o}^{\infty}\|\hat{v}\|^2 + e^{2\rho s_{o}}( v(s_{o}), J_{0}v_{t}(s_{o})) .\nonumber
\end{align}
Hence
\begin{equation}
\label{condition_of_rho}
 (1-2\rho) \int_{s_o}^{\infty} e^{2\rho \gamma_{n}(s)}\|v(s)\|^2ds\leq C
\end{equation}
with a constant $C$ independent of $n.$ Let $\rho<\frac{1}{2}$; taking the limit as $n\longrightarrow\infty$ we conclude that $e^{\rho s}\|v(s)\|\in L^2$ as claimed.   $\;\;\;\;\;\;\;\;\;\;\;\Box$
\v
By a  similar argument of Lemma \ref{a_lemma_of_hopf} we have
\begin{lemma}\label{cylinder_hofer}
 Assume $v:[-R ,R ]\times S^1\longrightarrow R^2$ is smooth, bounded and solves the equation
$$
v_{s} +J_{0}v_{t}=g,\;\;\;\;\;\;\;\; where \;\;\; \|g(s)|_{|s|\leq B}\|\leq C_{o}e^{-\delta (R-B)},\;\;\;\;\forall 0<B<R ,
$$
for some $\delta>0,$ where the norm is the $L^2(S^1)$-norm. To simplify notations we denote by $v(\pm R)$ the restriction of $v(s,t)$ to $s=\pm R$. If $v$ satisfies
$$ \| v_{t}(\pm R)\|_{L^2(S^{1})}\leq C_{1},\;\;\;
\int_{0}^{1}v(s,t)dt\equiv 0,
$$
then
$$
 \int_{s}^{s+1}\|v(s)\|^2ds\leq Ce^{-2\rho (R-B)},\;\;\;\;\;\forall |s|\leq B
$$
for every $0\leq \rho<\delta$ and $\rho<\frac{1}{2},$ where $C>0$ is a constant depending only on $C_{o},\frac{1}{2}-\rho$ and $\sum |( v(\pm R), J_{0}v_{t}(\pm R))|.$
\end{lemma}
\n{\bf Proof.} As in the proof of Lemma \ref{a_lemma_of_hopf} we have
\begin{align}
\int_{s_1}^{s_{2}} \|v(s)\|^2ds
\leq \int_{s_1}^{s_{2}}\|g(s)\|^2ds + \sum |( v(s_{i}),J_{0}v_{t}(s_{i}))|.
\end{align}
Take  a  even functions $\gamma:[-R,R] \longrightarrow \mathbb R $ satisfying
 $$\gamma (s)=s+R,\;\;\mbox{ for }-R\leq s\leq -1,\; \;\;\;\; \;\;\;0\leq \gamma' (s)\leq 1,\;\;\;\;\forall s\in [-R,0],$$ $$\gamma (s)=const\;\;\mbox{ for }  -\frac{1}{2}<s\leq 0.$$ Let $\rho>0$ and define
$$
\hat{v} (s,t) =e^{\rho \gamma (s)}v(s,t).
$$
Obviously, $
\hat{v} (\pm R,t) = v(\pm R,t).
$  Then $\hat{v}$   has mean values zero over $S^1$ and $\|\hat{v}(s)\|\in L^2$. Differentiating we obtain
$$
\hat{v}_{s}+J_{0}\hat{v}_{t} =e^{\rho \gamma (s)}g+ \rho \gamma' (s)\hat{v}.
$$
If $0<\rho<\delta$ we conclude,
\begin{align}
\int_{-R}^{R} \|\hat {v}(s)\|^2ds & \leq \int_{-R}^{R} \|\hat {v}_{s}+J_{0}\hat {v}_{t}\|^2ds+ \sum|( {v}(\pm R),J_{0} {v}_{t}(\pm R) )| \nonumber\\
& \leq 2\int_{-R}^{R}e^{2\rho \gamma(s)} \|g(s)\|^2ds +2\rho \int_{-R}^{R}\|\hat{v}\|^2 + \sum|( {v}(\pm R),J_{0} {v}_{t}(\pm R) )| .\nonumber
\end{align}
Hence
\begin{equation}
\label{condition_of_rho}
 (1-2\rho) \int_{-R}^{R}\|\hat {v}(s)\|^2ds\leq C_1.
\end{equation}
Let $\rho<\frac{1}{2}$. Then
$$
 \int_{s}^{s+1}\|v(s)\|^2ds\leq Ce^{-2\rho (R-B)},\;\;\;\;\;\forall |s|\leq B.
$$
The lemma follows.\;\;\;\;\;$\Box$
\v

 We need the following lemmas (see \cite{MS}).

\begin{theorem}\label{removal_singular}
Let $(M,\omega)$ be a compact symplectic manifold with $\omega$-tamed complex structure $J.$ Then there exists a constant $\hbar>0$ such that the following holds.
If $r>0$ and $u:B_{r}(0)\rightarrow M$ be a J-holomorphic curve then
\begin{equation}
\int_{B_{r}(0)}|du|^2\leq \hbar\;\;\;\;\;\;\;\; \Longrightarrow \;\;\;\;\;\;\;\;\;\; |du(0)|^2\leq \frac{8}{\pi r^2}\int_{B_{r}(0)}|du|^2.
\end{equation}
\end{theorem}

\begin{lemma}\label{tube_exponential_decay} Let $(M,\omega)$ be a compact symplectic manifold with $\omega$-tamed almost complex structure $J.$ Fix a constant  $\fc\in (0,1)$. There are two positive constant $\mathcal{C}_{1}$ and $\hbar$ depending only on $J$, $\omega$ and $\fc$ such that for any  $J$-holomorphic map  $u:[-R,R ]\times S^1\rightarrow  M$ with
\begin{equation}
 E (u,-R\leq s\leq R) <\hbar,
\end{equation}
we have
 \begin{align}
&E(u,-B\leq s\leq B)\leq \mathcal{C}_{1}e^{-2\fc(R-B)},\;\; \;\;\;\forall\; 0\leq B\leq R,\\
&|\nabla {u}|(s,t)\leq \mathcal{C}_{1}e^{-\fc(R-|s|)},\;\; \;\;\;\forall\; |s|\leq R-1,
\end{align}
 \end{lemma}

 Following \cite{HWZ1} we introduce functions
\begin{equation}
a^{\diamond}(s,t)=a(s,t)-ks,\;\; \theta^{\diamond}(s,t)=\theta(s,t)-kt.
\end{equation}
Denote
\begin{equation}
\pounds=(a^{\diamond},  \theta^{\diamond}).
\end{equation}
We have
\begin{equation}\label{eqn_pounds}
\pounds_s + J_0\pounds_t= h,
\end{equation}
where $h=(\sum b_i({\bf w})(w_{i})_{t}, -\sum b_i({\bf w})(w_{i})_{s})$.
\v
Using Lemma \ref{tube_exponential_decay} we can prove that
\begin{lemma}\label{intexp_L} Let $u:[-R ,R ]\times S^1\rightarrow {\mathbb R}\times \widetilde{M}$ be a $J$-holomorphic maps with finite energy. Suppose that
 $ \widetilde E(u,-R\leq s\leq R)<\hbar$. Then there exists a constant $B>0$ such that
	\begin{align*}
 \int_{S^1} \theta^{\diamond}_{t}(s_{1},t)dt=\int_{S^1} \theta^{\diamond}_{t}(s_{2},t)dt,\;\;\;\forall\;  |s_{1}|,|s_{2}|\leq B.
	\end{align*}
\end{lemma}
\n{\bf Proof.} Consider the $J$-holomorphic map $\hat u=\pi\circ u.$ By Lemma \ref{tube_exponential_decay} there exists a constant $B>0$ such that
$\hat u([-B,B]\times S^{1})$ lies in a local Darboux coordinates system $\mathbf w.$
 Taking derivative $\frac{\p}{\p t}$ of \eqref{eqn_pounds}  we have
 $$
 \label{eqn_pounds}
\theta^{\diamond}_{st}+a^{\diamond}_{tt}=- \left(\sum b_i({\bf w})(w_{i})_{s}\right)_{t}.$$
It is easy to see that $h(s,t)=h(s,t+1)$ and  $a^{\diamond}_{t}(s,t)=a^{\diamond}_{t}(s,t+1)$.
 Integrating this equation   over $S^{1}$ we obtain
\begin{equation}\label{integration_exp}
\left|\frac{d}{ds}\int_{S^{1}} \theta^{\diamond}_{t} dt\right|=\left|\int_{S^{1}} \theta^{\diamond}_{st} dt\right|\leq \left| a^{\diamond}_{t}(s,t+1) - a^{\diamond}_{t}(s,t) \right|+\left|h(s,t+1) - h(s,t) \right| =0.
\end{equation}
Then Lemma follows. $\;\;\;\;\;\Box$
\v\v

Next we prove
\begin{theorem}\label{exponential_estimates_theorem}    Let $u:\mathbb C-D_1\rightarrow {\mathbb R}\times \tilde{M}$ be a J-holomorphic map with finite energy.
Put
$z=e^{s+2 \pi it}$. Then $$\lim_{s\rightarrow \infty}\widetilde
u(s,t)=x(kt) $$ in $C^{\infty}(S^1)$ for some $k$-periodic orbit
$x,$ and there are constants $\ell_0$,
$\theta_0$  such that for any $0<\mathfrak{c}<\min\{\frac{1}{2},\mathfrak{c}_{1} \}$ and for all ${\bf n}=(m_{1},m_{2})\in {\mathbb Z_{\geq 0}^2 }$
\begin{eqnarray}
\label{exponential_decay_a}
|\partial^{\bf n}[a(s,t)-ks-\ell_0]|\leq \mathcal C_{\bf n} e^{-\mathfrak{c}|s|}\\
\label{exponential_decay_theta}
|\partial^{\bf n}[\theta(s,t)-kt-\theta_0]|\leq \mathcal C_{\bf n} e^{-\mathfrak{c}|s|} \\
\label{exponential_decay_y}|\partial^{\bf n}{\bf w}(s,t)|\leq \mathcal C_{\bf n}
e^{-\mathfrak{c} |s|}.
\end{eqnarray}
\end{theorem}
\v
 \v\n
 {\bf Proof.} By Lemma \ref{sequence_convergence_of_j_holomorphic_maps},
  there is a  sequence $s_i\rightarrow \infty $
 such that $$\lim_{i\rightarrow \infty}\widetilde
u(s_i,t)=x(kt)$$ in $C^{\infty}(S^1)$ for some $k$-periodic orbit
$x$.

\v
  For any $(s,t)\in \mathbb C-D_1,$ let $D_{1/2}( s,t )$ be the Euclidean ball centered at $(s,t).$
We have
$$
\tilde E(\tilde u; {D_{1/2}( s,t )})\leq  \widetilde{E}(s-1).
$$
Since $\lim\limits_{s\to\infty} \widetilde E(s)=0,$ we can assume that  $\widetilde E(s_1-1)<\hbar$ for some $s_1\geq s_{0}$. Applying Theorem \ref{removal_singular}  we obtain that
$$
|\nabla {\bf w}(s,t)|\leq  \sqrt\frac{32}{\pi} \sqrt{\widetilde{E}(s-1)}\leq C_{2} e^{- \mathfrak c_{1} s/2 },
$$
where $C_2=\sqrt\frac{32}{\pi}\widetilde{E}(s_0)^{1/2}e^{\mathfrak c_{1} s_{0}/2 }$.
A direct integration give us
\begin{equation} \label{estimate_of_y_energy}
|{\bf w}(s,t)|_{J_0}\leq C_{3} e^{-\mathfrak c_{1} s/2}
\end{equation}
for some constant $C_{3}>0.$
By the standard elliptic estimate we have \eqref{exponential_decay_y}.
\v

\v Integrating the equation $\pounds_{s}+J_{0}\pounds_{t}=h$ over $S^{1}$ we obtain
\begin{equation}\label{integration_exp}
 \left|\frac{d}{ds}\int_{S^{1}} \pounds dt\right|=\left|\int_{S^{1}} \pounds_{s}dt\right|=\left|\int_{S^{1}} hdt\right|\leq C_{4}e^{ -\mathfrak{c}_1s/2},
\end{equation}
for some constant $C_{4}>0.$
Then $\int_{S^{1}}(\theta -kt) dt$ uniformly converges to some constant $\theta_{0}$.
From Lemma \ref{partial_a_theta_convergence_zero} we   conclude that
  $$\theta-kt\longrightarrow \theta_{0},\;\;\; in \;\;C^{\infty}(R).\;\;\;\;\;\;\; $$
Similar, we have
 $$a-ks\longrightarrow \ell_{0},\;\;\; in \;\;C^{\infty}(R),\;\;\;\;\;\;\; $$
 for some constant $\ell_{0}$.
\v
By the same arugment of Lemma \ref{intexp_L} we have
$$
\int_{S^{1}}\pounds_{t}(s,t)dt=\lim_{s\to\infty}\int_{S^{1}}\pounds_{t}(s,t)dt=0.
$$
For any $\bf{ n } = (m_{1},m_{2})\in \mathbb Z^{2}_{\geq 0}$   put $V := {\bf \pounds}_t $ and $g = h_{t}  $,
 we have
$$(\partial^{\bf n}V)_s + J_0(\partial^{\bf n}V)_t = \partial^{\bf n} g , \;\;\;\;\int_{S^1}\partial^{\bf n}V (s,t)dt\equiv 0.$$
 \eqref{exponential_decay_y} gives us
 $$ \| \partial^{\bf n}g\|\leq C_{5}({\bf n},s_{0})e^{-\mathfrak{c}_1 s/2},$$
 where $C_{5}({\bf n},s_{0})>0$ is a constant depending only on $s_{0},$ $\bf n$ and $\tilde E(s_{o})$.
 It follows from \eqref{deriatives_a_theta_convergence_zero} that $(\partial^{\bf n}V)_{t}\longrightarrow 0$ as $s\longrightarrow\infty$ uniformly in $t$.  By Lemma \ref{a_lemma_of_hopf} we have, for any $\mathfrak{c }<\min\{\frac{1}{2},\frac{\mathfrak{c}_1}{2}\}$,
\begin{equation}
\int_{s}^{s+1}\|\partial^{\bf n}V \|_{L^{2}}ds\leq C_{6}({\bf n},s_{0}) e^{-\mathfrak{c } s}
\end{equation}
where $C_{6}({\bf n},s_{0})>0$ is the constant depending only on $ C_{5}({\bf n},s_{0})$ and $(\partial^{\bf n} V(s_{o}),  J_{0} \partial^{\bf n}V_{t}(s_{o})\rangle .$
It follows from the Sobolev imbedding theorem that
 \begin{equation}\label{estimate_deritive_w}
|\pounds_{s}|+|\pounds_{t}|\leq |h|(s,t)+2|V|(s,t)\leq  C_{7}({\bf n},s_{0}) e^{-\mathfrak{c } s}.
 \end{equation}
where $C_{7}({\bf n},s_{0})>0$ is a constant depending only on $s_{0},$ $\tilde E(s_{o})$ and $ \sum\limits_{|\bf n|\leq 2} |(\partial^{\bf n} V(s_{o}),  J_{0} \partial^{\bf n}V_{t}(s_{o})\rangle |.$

Then by a direct integration we can obtain the \eqref{exponential_decay_a} and  \eqref{exponential_decay_theta}.
$\;\;\;\;\;\;\;\;\;\;\;\Box$
\v

Similarly, we have
\begin{theorem}\label{tube_L_decay} Let $u:[-R ,R ]\times S^1\rightarrow {\mathbb R}\times \widetilde{M}$ be a $J$-holomorphic maps with finite energy. Assume that
\begin{enumerate}
\item[(i)] $ \widetilde E(u,-R\leq s\leq R)<\hbar$,
\item[(ii)] $\sum\limits_{n_1,n_{2}\leq 3}\|\nabla^{\bf n} u(-R,\cdot)\|_{L^2(S^{1})}\leq C_{1},\;\;\;\;\sum\limits_{n_1,n_{2}\leq 3}\|\nabla^{\bf n}u(R,\cdot)\|_{L^2(S^{1})}\leq C_{1},$  where ${\bf n}=(m_{1},m_{2}),$
\end{enumerate} Then there exist three constants $k\in \mathbb Z_{>0},$  $C>0$ and   $B>0$ depending only on $\hbar$, $\tilde J$ and $C_1$ such that
 \begin{align}
|\nabla \pounds|\leq Ce^{-\fc(R- |s|)},\;\;\;\forall\; |s|\leq B,
\end{align}
where $\pounds=(a-ks,\theta-kt).$
 \end{theorem}
\n{\bf Proof.}
As in the proof of Lemma \ref{intexp_L} there exists a constant $B>0$ such that
$\hat u([-B,B]\times S^{1})$ in a local Darboux coordinates system $\mathbf w.$ Since $\theta(s,t+1)=\theta(s,t),\mod 1,$ there exists $k\in \mathbb Z_{>0}$ such that $\theta(-B,t)+k=\theta(-B,t+1).$ By Lemma \ref{intexp_L} we have
$$\pounds(s,t)=\pounds(s,t+1),\;\;\;\forall |s|\leq B.$$
Then by a similar argument of the proof of Theorem \ref{exponential_estimates_theorem}  we can prove the lemma.  $\;\;\;\;\;\Box$
 \v

We introduce some terminology.
\v
\begin{definition}\label{terminology on periodic orbit}
{\bf (1).} Let $u:\Sigma
-\{p\}\rightarrow N$ be a $J$-holomorphic map with finite energy,
and $p$ be a nonremovable singularity. If , in terms of local
coordinates $(s,t)$ around $p$, $\lim_{s\rightarrow
\infty}\widetilde u(s,t)=x(kt)$, we say simply that $u(s,t)$
converges to the $k$-periodic orbit $x$. We call $p$ a {\em
positive} (resp. {\em negative} ) end, if $a(z) \rightarrow
\infty$ (resp. $-\infty$) as $z\rightarrow p$.
\v
{\bf (2).} Suppose that $\Sigma_1$ and $\Sigma_2$
join at $p$, and
$u_1:\Sigma_1\to M^+ \; (or\;{\mathbb{R}}\times \widetilde{M}$),
$u_2:\Sigma_2\to {\mathbb{R}}\times \widetilde{M}$ are $J$-holomorphic maps with finite energy. Choose coordinates $(a_1,\theta_1, \mathbf{w})$ on $M^+$ \;(or\;${\mathbb{R}}\times \widetilde{M}$),  $(a_2,\theta_2, \mathbf{w})$ on ${\mathbb{R}}\times \widetilde{M}$ and
choose holomorphic cylindrical
coordinates $(s_1,t_1)$ on $\Sigma_1$ and $(s_2,t_2)$
on $\Sigma_2$ near $p$ respectively. Suppose that
$$\lim_{s_1 \rightarrow \infty} \widetilde{u}_1(s_1,t_1)= x_1(k_1t_1),\;\;\;
\lim_{s_2 \rightarrow -\infty} \widetilde{u}_2(s_2,t_2)= x_2(k_2t_2).$$
We say $u_1$ and $u_2$ converge to a same periodic orbit as the variable tend to $p$, if $k_1=k_2$,
and $\pi(x_1)=\pi(x_2)$, where $\pi$ denotes the projection to $Z$.
\end{definition}

\v

\chapter{The moduli space of stable holomorphic maps}\label{stable holo}

\section{Deligne-Mumford moduli space}\label{Deligne-Mumford moduli space}
\v
First all, we recall some results on the Deligne-Mumford moduli space $\overline{\mathcal{M}}_{g,\mathfrak{n} }$ of stable curves, for detail see \cite{Tromba}, \cite{Wolp-1}, \cite{Wolp-2}.
\v
\subsection{\bf Metrics on $\Sigma$}\label{metric on surfaces} Let $(\Sigma,j,\mathfrak{z})$ be a smooth Riemann surface of genus $g$ with $\mathfrak{n}$ marked points $\mathfrak{z}$. In this paper we assume that $\mathfrak{n}>2-2g$, and $(g,\mathfrak{n})\ne (1,1), (2,0)$. It is well-known that there is a unique complete hyperboloc metric $\mathbf{g}_0$ in $\Sigma\setminus \{\mathfrak{z}\}$ of constant curvature $-1$ of finite volume, in the given conformal class $j$ ( see \cite{Wolp-1}).
Let $\mathbb H=\{\zeta=\lambda+\sqrt{-1}\gamma|\gamma>0\}$ be the half upper plane with the Poincare metric
$$
\mathbf{g}_0(\zeta)=\frac{1}{(Im(\zeta))^2}d\zeta d\bar\zeta.
$$
 Let
$$
\mathbb{D}=\frac{ \{\zeta \in \mathbb H|  Im (\zeta) \geq 1\}}{\zeta \sim \zeta + 1}
$$
be a cylinder, and $\mathbf{g}_0$ induces a metric on $\mathbb{D}$, which is still denoted  by $\mathbf{g}_0$. Let $z=e^{2\pi i\zeta}$, through which we identify $\mathbb{D}$ with $D(e^{-2\pi}):=\{z||z|< e^{-2\pi}\}$.
An important result is that for any marked point $\mathfrak{z}_i$ there exists a
neighborhood $O_i$ of $\mathfrak{z}_i$ in $\Sigma$ such that
$$
(O_i\setminus \{\mathfrak{z}_i\},\mathbf{g}_0)\cong (D(e^{-2\pi})\setminus \{0\},\mathbf{g}_0),
$$
moreover, all $O_i$'s are disjoint with each other. Then we can view $D_{\mathfrak{z}_i}(e^{-2\pi})$ as a neighborhood of $\mathfrak{z}_i$ in $\Sigma$ and $z$ is a local complex coordinate on $D_{\mathfrak{z}_i}(e^{-2\pi})$ with $z(\mathfrak{z}_i)=0$. In terms of the coordinates $z$ the metric $\mathbf{g}_0$ becomes
$$\mathbf{g}_0=\frac{dzd\bar{z}}{|z|^2(log|z|)^2}.$$
The distinguished coordinates $z$ is unique modulo a unimodular factor. We call $z$ the cusp coordinate. For any $c>0$ denote $$\mathbf{D}(c)= \bigcup_{l=1}^{\mathfrak{n}}  D_{\mathfrak{z}_l}(c) ,\;\;\;\Sigma(c)=\Sigma\setminus \mathbf{D}(c).$$ Let $\mathbf{g}'=dzd\bar{z}$ be the standard Euclidean metric on each $D_{\mathfrak{z}_i}(e^{-2\pi})$. We fix a smooth cut-off function $\chi(|z|)$ to glue $\mathbf{g}_0$ and $\mathbf{g}'$, we get a smooth metric $\mathbf{g}$  in the given conformal class $j$ on $\Sigma$ such that
\[
\mathbf{g}=\left\{
\begin{array}{ll}
\mathbf{g}_0 \;\;\;\;\; on \;\;\Sigma\setminus \mathbf{D}(e^{-2\pi}),    \\  \\
\mathbf{g}'\;\;\;on\; \mathbf{D}(\frac{1}{2}e^{-2\pi})\;\;.
\end{array}
\right.
\]
\v
\v
\v
Put $z=e^{s+2\pi\sqrt{-1}t}$. We call $(s,t)$ the cusp cylinder coordinates.
\v
Let $\mathbf{g}^{c}=ds^2+dt^2$ be the  cylinder metric on each $D_{\mathfrak{z}_i}^{*}(e^{-2\pi})$.
We also define another metric $\mathbf{g}^\diamond$ on $\Sigma$ as above by glue $\mathbf{g}_0$ and $\mathbf{g}^c$, such that
\[
\mathbf{g}^{\diamond}=\left\{
\begin{array}{ll}
\mathbf{g}_0 \;\;\;\;\; on \;\;\Sigma\setminus \mathbf{D}(e^{-2\pi}),    \\  \\
\mathbf{g}^{c}\;\;\;on\; \mathbf{D}(\frac{1}{2}e^{-2\pi})\;\;.
\end{array}
\right.
\]

\v
The metric $\mathbf{g}$ (resp. $\mathbf{g}^{\diamond}$) can be generalized to marked nodal surfaces in a natural way. Let $(\Sigma, j,\mathfrak{z})$ be a
marked nodal surfaces with $\mathfrak{e}$ nodal points $\mathbf{q}=(q_{1},\cdots,q_{\mathfrak{e}})$. Let $\sigma:\tilde{\Sigma}=\sum_{\nu=1}^{\mathfrak{r}}\Sigma_{\nu}\to \Sigma$ be the normalization. For every node $q_i$ we have a pair $\{\mathsf{a}_i, \mathsf{b}_i\}$. We view $\mathsf{a}_i$, $\mathsf{b}_i$ as marked points on $\tilde{\Sigma}$ and define the metric $\mathbf{g}_{\nu}$ (resp. $\mathbf{g}^{\diamond}_{\nu}$) for each $\Sigma_{\nu}$. Then we define
$$\mathbf{g}:=\bigoplus_{1}^{\nu} \mathbf{g}_{\nu},\;\;\;\;\;\mathbf{g}^{\diamond}:=\bigoplus_{1}^{\nu} \mathbf{g}^{\diamond}_{\nu}.$$

\subsection{Teichm\"uller space}\label{sub_sect_Teich}

Denote by $\mathcal{J}(\Sigma)\subset End(T\Sigma)$ the manifold of all $C^{\infty}$ complex structures on $\Sigma$, let $\mathcal{G}$ denote the manifold of $C^{\infty}$ Riemannian metrics with constant scalar curvature $-1$ on $\Sigma$. Denote by $Diff^+(\Sigma)$ the group of orientation preserving $C^{\infty}$ diffeomorphisms of $\Sigma$, by $Diff^{+}_{0}(\Sigma)$ the identity component of $Diff^+(\Sigma).$ $Diff^+(\Sigma)$ acts on $\mathcal{J}(\Sigma)$ and $\mathcal{G}$ by
$$(\phi^*J)_x:=(d\phi_x)^{-1}J_{\phi(x)}d\phi_x,\;\;\;
(\phi^*g)(x)(w,v):=g(\phi(x))(d\phi(x)w, d\phi(x)v)$$
for all $\phi\in Diff^+(\Sigma)$, $x\in \Sigma$, $w,v\in T_x\Sigma.$ There is a bijective,  $Diff^+(\Sigma)$-equivariant correspondence between $\mathcal{J}(\Sigma)$ and $\mathcal{G}$:
$$\mathcal{J}(\Sigma)\cong \mathcal{G}.$$
Put
$$\mathbf{P}:=\mathcal J(\Sigma)\times (\Sigma^{\mathfrak{n}} \setminus \Delta),$$
where $\Delta\subset \Sigma^{\mathfrak{n}}$ denotes the fat diagonal.
The orbit spaces are
\begin{align*}
\mathcal{M}_{g,\mathfrak{n}}=\left(\mathcal J(\Sigma)\times (\Sigma^{\mathfrak{n}} \setminus \Delta)\right)/ Diff^+(\Sigma),\;\;\;\mathbf{T}_{g,\mathfrak{n}}=\left(\mathcal J(\Sigma)\times (\Sigma^{\mathfrak{n}} \setminus \Delta)\right)/ Diff^+_{0}(\Sigma).
\end{align*}
$\mathcal{M}_{g,\mathfrak{n}}$ is called the Deligne-Mumford space, $\mathbf{T}_{g,\mathfrak{n}}$ is called the Teichm\"uller space. The mapping class group of $\Sigma$ is
$$Mod_{g,\mathfrak{n}}=Diff^{+}(\Sigma)/Diff^{+}_0(\Sigma).$$
It is well-known that $Mod_{g,\mathfrak{n}}$ acts properly discontinuously on $\mathbf{T}_{g,\mathfrak{n}}$ and $$\mathcal{M}_{g,\mathfrak{n}}=\mathbf{T}_{g,\mathfrak{n}}/Mod_{g,\mathfrak{n}}$$ is a complex orbifold of dimension $3g-3 +\mathfrak{n}$. Let $\pi_{\mathcal{M}}:\mathbf{T}_{g,\mathfrak{n}}\to \mathcal{M}_{g,\mathfrak{n}}$ be the projection.

\v

\v
Consider the principal fiber bundle
$$
Diff^+_0(\Sigma) \to \mathbf{P}   \to \mathbf{T}_{g,\mathfrak{n}}
$$
and the associated fiber bundle
$$
\pi_{\mathbf{T}}:   \mathcal{Q}:= \mathbf{P}\times_{Diff^+_0(\Sigma)}\Sigma\to \mathbf{T}_{g,\mathfrak{n}},
$$
which has fibers isomorphic to $\Sigma$ and is equipped with $\mathfrak{n}$ disjoint sections
$$\mathcal Y_i := \left\{[j,\mathfrak{z}_1,\dots,\mathfrak{z}_{\mathfrak{n}},z]\in \mathcal{Q}\,:\,z=\mathfrak{z}_i\right\},\qquad
i=1,\dots,\mathfrak{m}.
$$
It is commonly called the universal curve over $\mathbf{T}_{g,\mathfrak{n}}$. The following result is well-known (cf \cite{RS}):
\begin{lemma}\label{slice} Suppose that $\mathfrak{n}+ 2g \geq 3$. Then for any $\gamma_o=[(j_o,\mathfrak{z}_o)]\in \mathbf{T}_{g,\mathfrak{n}}$, and any $(j_o,\mathfrak{z}_o)\in \mathbf{P}$ with $\pi_{\mathbf{T}}(j_o,\mathfrak{z}_o)=\gamma_o$ there is an open neighborhood $\mathbf{A}$ of zero in $\mathbb{C}^{3g-3+\mathfrak{n}}$ and a local holomorphic slice
$\iota=(\iota_{0},\cdots,\iota_{\mathfrak{n}}):\mathbf{A}\to \mathbf{P} $ such that
\begin{equation}\label{slice-1}
\iota_0(o)=j_o,\qquad \iota_i(o)=\mathfrak{z}_{io} \qquad i=1,\dots,\mathfrak{n},
\end{equation}
and the map
$$
\mathbf{A}\times Diff^{+}_{0}(\Sigma) \to \mathbf{P}:(a ,\phi)\mapsto (\phi^*\iota_{0}(a),\phi^{-1}(\iota_1(a)),
\cdots,\phi^{-1}(\iota_{\mathfrak{n}}(a))
$$
is a diffeomorphism onto a neighborhood of the orbit of $(j_o,\mathfrak{z}_o)$.
\end{lemma}
\v\n

From the local slice we have a local coordinate chart on $U$ and a local trivialization on $\pi_{\mathbf{T}}^{-1}(U)$:
\begin{equation}\label{local coordinates}
\psi: U\rightarrow \mathbf{A},\;\;\;\Psi:\pi_{\mathbf{T}}^{-1}(U)\rightarrow \mathbf{A}\times \Sigma,\end{equation}
 where $U\subset \mathbf{T}_{g,\mathfrak{n}}$ is a open set. We call $(\psi, \Psi)$ in \eqref{local coordinates} a local coordinate system for $\mathcal{Q}$. Suppose that we have two local coordinate systems
\begin{equation}\label{local coordinates-1}
(\psi, \Psi): (O, \pi_{\mathbf{T}}^{-1}(O))\rightarrow (\mathbf{A}, \mathbf{A}\times \Sigma),\end{equation}
\begin{equation}\label{local coordinates-2}
(\psi', \Psi'): (O', \pi_{\mathbf{T}}^{-1}(O'))\rightarrow (\mathbf{A}', \mathbf{A}'\times \Sigma).\end{equation}
Suppose that $O\bigcap O'\neq \emptyset.$ Let $W$ be a open set with $W\subset O\bigcap O'$.
Denote $V=\psi(W)$ and $V'=\psi'(W)$. Then ( see \cite{RS})
\begin{lemma}\label{local coordinates-3}
$\psi'\circ \psi^{-1}|_{V}:V\to V'$ and $\Psi'\circ \Psi^{-1}|_{V}: {V}\times \Sigma\to {V'}\times \Sigma$ are holomorphic.
\end{lemma}
The diffeomorphism group $ Diff^+(\Sigma)$ acts on $\Sigma^{n}\setminus \Delta$
 by
\begin{equation}\label{group action-1}
\varphi^*(j,\mathfrak{z}_1,\dots,\mathfrak{z}_{\mathfrak{n}})
:=(\varphi^*j,\varphi^{-1}(\mathfrak{z}_1),\dots,\varphi^{-1}(\mathfrak{z}_{\mathfrak{n}})).
\end{equation}
It is easy to see that $\mathbf{g}$ is $Diff^+(\Sigma)$-invariant.
\v
 Let $\overline{\mathcal{M}}_{g,\mathfrak{n}}$ be the Deligne-Mumford compactification space, $g_{\mathsf{wp}}$ be the Weil-Petersson metric on $\overline{\mathcal{M}}_{g,\mathfrak{n}}$. Denote by $\overline{\mathbb{B}}_{g,\mathfrak{n}}$ the groupoid whose objects are stable marked nodal Riemann surfaces of type $(g,\mathfrak{n})$ and whose morphisms are isomorphisms of marked nodal Riemann surfaces. J. Robbin, D. Salamon \cite{RS} used the universal marked nodal family to give an orbifold groupoid structure on $\overline{\mathbb{B}}_{g,\mathfrak{n}}$. Then $\overline{\mathcal{M}}_{g,\mathfrak{n}}$ has the structure of a complex orbifold, and $\mathcal{M}_{g,\mathfrak{n}}$ is an effective orbifold. It is possible that $(g_i,\mathfrak{n}_i)= (1,1)$ for some smooth component $\Sigma_i$, in this case we consider the reduced effective orbifold structure.

\v

\v
\section{Weighted sobolev norms}\label{weight_norms}

\v
We mainly discuss holomorphic maps into $M^+$. For holomorphic maps into $M^-$ and ${\mathbb R}  \times \widetilde{M}$ the discussions are the same.
There is a Riemannian metric
\begin{equation}\label{definition_of_metrics}
G_J(v,w):=<v,w>_J:=\frac{1}{2}\left(\omega(v,Jw)+\omega(w,Jv)\right)
\end{equation}
for any $v, w\in TM^+$. Following \cite{MS} we choose the complex linear connection
 $$
 \widetilde {\nabla}_{X}Y=\nabla_{X}Y-\tfrac{1}{2}J\left(\nabla_{X}J\right)Y $$
 induced by the Levi-Civita connection $\nabla$ of the metric $G_{J}.$
\v
\subsection{\bf Norms for maps from smooth Riemann surfaces}\label{norms for smooth surfaces}

Let $(\Sigma, \mathbf{y}, \mathbf{p})$ be a stable smooth marked surface
of genus $g$ with $m$ distinct marked points ${\bf y}=(y_1,...,y_m)$, $\mu$  distinct puncture points ${\bf p}=(p_1,...,p_{\mu})$. Put $\stackrel{\circ}{\Sigma}= \Sigma - \{\mathbf{y,p}\}$.
Let $u:\stackrel{\circ}{\Sigma} \rightarrow M^+$ be a $(j,J)$-holomorphic map. We choose cusp cylinder coordinates $(s,t)$ near each puncture point $p_j$.
Over each tube the linearized operator
$D_{u}$ takes the following form (see \cite{LS-1})
\begin{equation}
D_{u}=\frac{\partial}{\partial
s}+J_0\frac{\partial} {\partial t}+ F_{u}^{1}+F_{u}^{2} \frac{\partial }{\partial t}.
\end{equation}
Then \eqref{exponential_decay_a}, \eqref{exponential_decay_theta} and \eqref{exponential_decay_y} hold.
\v
 We introduce some notations. There is a bundle $\mathbb{H}\to Z$, whose fibre at $p\in Z$ is $T_pZ\oplus span\{\frac{\p}{\p a}, \frac{\p}{\p \theta}\}$, and a bundle $\mathbb{H}^*\to Z$, whose fibre at $p\in Z$ is $T_pZ$. Let $\{\frac{\p}{\p a}, \frac{\p}{\p \theta}, \frac{\p}{\p w_1},...,\frac{\p}{\p w_n}\}$ be a base of $\bh_p$.
There is a projection $\pi:\bh_p \to \bh^*_p$ given by
$$c_1\frac{\p}{\p a} + c_2\frac{\p}{\p \theta} + \sum b_i\frac{\p}{\p w_i}\longmapsto \sum b_i\frac{\p}{\p w_i}.$$
To simplify notations we identify $\mathbb{H}_{u(p)}$ ( resp. $\mathbb{H}^*_{u(p)}$ ) with $(u^*\mathbb{H})_p$ ( resp. $u^*\mathbb{H}^*_p$ ) and denote it by $\mathbb{H}_p$ ( resp. $\mathbb{H}^*_p$ ).

 By the elliptic regularity  we have, for any $k>0,$
\begin{equation}\label{expontent_decay_of_s}
 \sum_{i+j=k} \left|\frac{\p^{k}F^{i}_{u}}{\partial^ i s \partial ^j t}\right|\leq
C_{k}e^{-\fc s},\;\;\;\;i=1,2 .
\end{equation}
for some constants $C_{k}>0$, $\mathfrak{c}>0$.
Therefore, the operator $H_s=J_0\frac{d}{dt}+F_{u}^{1}+F_{u}^{2} \frac{\partial }{\partial t}$ converges
to $H_{\infty}=J_0\frac{d}{dt}$. Obviously, the operator $D_u$ is
not a Fredholm operator because over each nodal end the
operator $H_{\infty}=J_0\frac{d}{dt}$ has zero eigenvalue. For each puncture $p_i$ the $\ker
H^i_{\infty}$ consists of constant vectors in $\mathbb{H}_{p_i}$. To recover a Fredholm theory we use
weighted function spaces. We choose a weight $\alpha$ for each
end. Fix a positive function $W$ on $\Sigma$ which has order equal
to $e^{\alpha |s|} $ on each end, where $\alpha$ is a small
constant such that $0<\alpha<\fc $ and over each end
$H_{\infty}- \alpha = J_0\frac{d}{dt}- \alpha $ is invertible. We
will write the weight function simply as $e^{\alpha |s|}.$ For  given integer $k>4$  and for any
section $h\in C^{\infty}(\Sigma;u^{\ast}TM^+)$ and section $\eta \in
C^\infty(\Sigma, u^{*}TM^+\otimes \wedge^{0,1}_jT^{*}\Sigma)$ we define the norms
\begin{eqnarray}\label{norm_k_alpha}
&&\|h\|_{k,2,\alpha}=\left(\int_{\Sigma}e^{2\alpha|s|} \sum_{i=0}^{k} |\nabla^i h|^2
dvol_{\Sigma}\right)^{1/2},\\
\label{norm_k-1_alpha}
&&\|\eta\|_{k-1,2,\alpha}=\left(\int_{\Sigma}e^{2\alpha|s|}\sum_{i=0}^{k-1} |\nabla^i \eta|^2
dvol_{\Sigma}\right)^{1/2}.
\end{eqnarray}
Here all norms and
covariant derivatives are taken with respect to the  metric $G_J$ on $u^{\ast}TM^+$ and
the metric $\mathbf{g}$ on $(\Sigma, j, \mathbf{y}, \mathbf{p})$, $dvol_{\Sigma}$ denotes the volume form with respect to $\mathbf{g}$.
Denote by $W^{k,2,\alpha}(\Sigma;u^{\ast}TM^+)$ and
$W^{k-1,2,\alpha}(\Sigma, u^{*}TM^+\otimes\wedge^{0,1}_jT^{*}\Sigma)$ the complete spaces with respect to the norms \eqref{norm_k_alpha} and \eqref{norm_k-1_alpha} respectively.
\v
 The operator $D_u:   W^{k,2,\alpha}\rightarrow W^{k-1,2,\alpha}$
is a Fredholm operator so long as $\alpha$ does not lie in the spectrum
of the operator $H^i_{\infty}$ for all $i=1,\cdots,\mu$.
 \v
\begin{remark} The index $ind (D_u,\alpha)$ does not change if
$\alpha$ is varied in such a way that $\alpha$ avoids the spectrum of
$H^i_{\infty}.$ Conversely, the index will change if $\alpha$ is moved
across an eigenvalue. We will choose $\alpha$ slightly larger than zero such
that at each end it does not across the first positive eigenvalue.
\end{remark}
\v
For each point $p_i \in \{p_1,...,p_{\mu}\}$, $i=1,...,\mu$, let $h^i_{0}\in \ker H^i_{\infty}$. Put
$H_{\infty}=(H^1_{\infty},...,H^\mu_{\infty})$, $h_{0}=(h^1_{0},...,h^\mu_{0}).$
We choose coordinates $(a, \theta)$ over the cylinder end of $M^+$. For each $p_i$ we choose a local Darboux coordinate $\mathbf{w}_{i}$ near $\pi\circ \tilde{u}(p_{i})\in  Z. $
$h_0$ may be considered as a vector field in the coordinate neighborhood.
We fix a cutoff function $\varrho$:
\[
\varrho(s)=\left\{
\begin{array}{ll}
1, & if\ |s|\geq d, \\
0, & if\ |s|\leq \frac{d}{2}
\end{array}
\right.
\]
where $d$ is a large positive number. Put
$$\hat{h}_0=\varrho h_0.$$
Then for $d$ big enough $\hat{h}_0$ is a section in $C^{\infty}(\Sigma; u^{\ast}TM^+)$
supported in the tube $\{(s,t)||s|\geq \frac{d}{2}, t \in {S} ^1\}$.
Denote
$${\mathcal W}^{k,2,\alpha}=\{h+\hat{h}_0 | h \in
W^{k,2,\alpha},h_0 \in \ker H_{\infty}\}.$$
\vskip 0.1in
\noindent
 We define weighted Sobolev  norm  on ${\mathcal W}^{k,2,\alpha}$ by $$\| h+\hat{h}_{0}\|_{\mathcal{W},k,2,\alpha}=
 \|h\|_{ k,2,\alpha} + |h_{0}| .$$

\v
 Denote by $\widetilde{\mathcal{B}}$ the space of $\mathcal{W}^{k,2,\alpha}$-maps $u:\Sigma \longrightarrow M^+$
 with a fixed homology class $A\in H_{2}(\overline{M^+} ,\mathbb Z),$ that is
 $$
\widetilde{ \mathcal{B}}=\{u\in \mathcal{W}^{k,2,\alpha}(\Sigma,M^+)|\;\bar{u}_{\star}([\Sigma])=A\},
 $$
 where $k>4$ is an even integer.  Then $\widetilde{\mathcal{B}}$ is an infinite dimensional Banach manifold. For any $u\in \widetilde{\mathcal{B}},$ the tangent space at $u$ is $T_u \widetilde{\mathcal{B}}=\mathcal{W}^{k,2,\alpha}(\Sigma,u^{\star}TM^+). $
The exponential map for $(M^+,G_J)$ provides a coordinate chart at $u$.

 The map $u$ is called a $(j,J)$-holomorphic map if $du\circ j=J\circ du$. Alternatively
\begin{equation}\label{holo}
\bar\p_{j,J}(u):=\half\left(du + J(u) du\circ j\right)=0.
\end{equation}
Let $\widetilde{\mathcal E}$ be the infinite dimensional Bananch bundle over $\widetilde{\mathcal B}$ whose fiber at $b=(j,\mathbf{y},\mathbf{p},u)$ is
$$W^{k-1,2,\alpha}(\Sigma, u^{*}TM^+\otimes\wedge^{0,1}_jT^{*}\Sigma).$$
The Cauchy-Riemann operator defines a Fredholm section $\overline{\partial}_{j,J}:\widetilde{\mathcal B}\longrightarrow \widetilde{\mathcal E}.$
\v
The diffeomorphism group $ Diff^+(\Sigma)$ acts on $(\Sigma^{m+\mu}\setminus \Delta)\times \widetilde{\mathcal B}$ and $(\Sigma^{m+\mu}\setminus \Delta)\times \widetilde{\mathcal E}$ by
\begin{equation}\label{group action-1}
\varphi^*(j,\mathbf{y},\mathbf{p},u)
:=(\varphi^*j,\varphi^{-1}(\mathbf{y}),\varphi^{-1}(\mathbf{p}), u\circ \varphi)
\end{equation}
\begin{equation}\label{group action-2}
\varphi^*\kappa= \kappa\cdot d\varphi\;\;\;\forall\; \kappa\in W^{k-1,2,\alpha}(\Sigma, u^{*}TM\otimes\wedge^{0,1}_jT^{*}\Sigma)
\end{equation}
for $\varphi\in Diff^+(\Sigma)$. Put
$$Aut(j,\mathbf{y},\mathbf{p},u)=
\{\phi\in Diff^+(\Sigma)|\varphi^*(j,\mathbf{y},\mathbf{p},u)=(j,\mathbf{y},\mathbf{p},u)\}.$$
We call it the automorphism group at $(j,\mathbf{y},\mathbf{p},u)$.
\v
Our moduli space ${\mathcal{M}}_{A}(M^{+},g,m+\mu,{\bf k})$ is the quotient space
$${\mathcal{M}}_{A}(M^{+},g,m+\mu,{\bf k})=((\Sigma^{m+\mu} \setminus \Delta)\times \overline{\partial }_{J}^{-1}(0) )/Diff^+(\Sigma).$$
\v
For any $[b_o]=[(q_o,u)]\in  {\mathcal{M}}_{A}(M^{+},g,m+\mu,{\bf k})$ with $[q_{o}]\in \mathcal M_{g,m+\mu}$
let $\gamma_o=[(j_o,\mathfrak{z}_o)]\in \mathbf{T}_{g,m+\mu}$, $(j_o,\mathfrak{z}_o)\in \mathbf{P}$ with $\pi_{\mathcal M}(\gamma_o)=[q_{o}]$ and $\pi_{\mathbf{T}}(j_o,\mathfrak{z}_o)=\gamma_o$. Choose a local coordinate system $(\psi,\Psi)$ on $U$ with $\psi(\gamma_o)=a_o$ for $\mathcal{Q}$ as in \eqref{local coordinates}, we have a local coordinate chart on $U$ and a local trivialization on $\pi_{\mathbf{T}}^{-1}(U)$:
\begin{equation}\label{local coordinates-1}
\psi: U\rightarrow \mathbf{A},\;\;\;\Psi:\pi_{\mathbf{T}}^{-1}(U)\rightarrow \mathbf{A}\times \Sigma,\end{equation}
 where $U\subset \mathbf{T}_{g,m+\mu}$ is an open set. We can view $a=(j,\mathfrak{z})$ as parameters, and the domain $\Sigma$ is a fixed smooth surface. Denote by $j_a$ the complex structure on $\Sigma$ associated with $a=(j,\mathfrak{z})$ and put $j_{a_o}:=j_o$. The Weil-Pertersson  metric induces a $Diff^+(\Sigma)$-invariant distance $d_{\mathbf{A}}(a_o,a)$ on $\mathbf{A}$ such that $d^{2}_\mathbf{A}(a):=d^{2}_\mathbf{A}(a_o,a)$ is a smooth function on $\mathbf{A}$. Denote by $\mathbf G_{a}$ the isotropy group at $a$, that is
  $$
  \mathbf G_{a}=\{\phi\in Diff^{+}(\Sigma)\;|\; \phi^{*}(j,\mathfrak{z})=(j,\mathfrak{z})\}.
  $$
  Since $\mathcal{M}_{g,m+\mu}$ is an effective orbifold, we can choose $\delta$ small such that $\mathbf{G}_{a}$ can be imbedded into $\mathbf{G}_{a_o}$ as a subgroup for any $a$ with $d_{\mathbf{A}}(a_o,a)<\delta$. Denote by $im(\mathbf{G}_{a})$ the imbedding.
\v
Let $b_{o}=(a_{o},u)=(j_o,\mathfrak{z}_{o},u)$ be the expression of $[(\gamma_o, u)]$ in this local coordinates.
Set
$$
\widetilde{\mathbf{O}}_{b_{o}}(\delta,\rho):=\{(a,v)\in \mathbf{A}\times \widetilde{\mathcal{B}} \;|\; d_{\mathbf{A}}(a_o,a)<\delta, \|h\|_{j_a,k,2}<\rho \},
$$
$$\mathbf{O}_{[b_{o}]}(\delta,\rho)=
\widetilde{\mathbf{O}}_{b_{o}}(\delta,\rho)/G_{b_o},$$
where $v=\exp_{u}(h),$ $G_{b_o}$ is the isotropy group at $b_o$, that is
$$G_{b_{o}}=\{\phi\in Diff^{+}(\Sigma)\;|\; \phi^{*}(j_o,\mathfrak{z}_o,u)=(j_o,\mathfrak{z}_o,u)\}.$$
 Obviously, $G_{b_o}$ is a subgroup of $\mathbf{G}_{a_o}$. Note that both $d_{\mathbf{A}}$ and $\|h\|_{j_a,k,2,\alpha}$ are $Diff^+(\Sigma)$-invariant, we may identified $\mathbf{O}_{[b_{o}]}(\delta,\rho)$ with a neighborhood of $[b_o]\in \mathcal{M}_{g,m+\mu}(A)$ in $\mathcal{B}_{g,m+\mu}(A)$.
\v
\subsection{\bf Pregluing}\label{gluing}

\v
Let $(\Sigma, j,{\bf y}, \mathbf{p}, q)$ be a marked nodal Riemann surface of genus $g$ with $m$ marked points ${\bf y}$,  $\mu$ puncture points ${\bf p}$ and one nodal point $q$.  We write the marked nodal Riemann surface as
$$
\left(\Sigma=\Sigma_{1}\wedge\Sigma_{2},j=(j_{1},j_{2}),{\bf y}=({\bf y}_{1},{\bf y}_{2}), {\bf p}=({\bf p}_{1},{\bf p}_{2}),q=(p_1,p_2)\right),$$
where $(\Sigma_{i},j_{i},{\bf y}_{i}, {\bf p}_{i}, q_i)$, $i=1,2$, are smooth Riemann surfaces. We say that $q_1,q_2$ are paired to form $q$. Assume that both $(\Sigma_{i},j_{i},{\bf y}_{i},{\bf p}_{i},q_{i})$ are stable. We choose metric $\mathbf{g}_i$ on each $\Sigma_i$ as in \S\ref{metric on surfaces}.
We choose the cusp cylindrical coordinates $(s_i, t_i)$ near $q_{i}$.  In terms of the cusp cylindrical coordinates we write
$$\Sigma_1\setminus\{q_1\}\cong\Sigma_{10}\cup\{[0,\infty)\times S^1\},\;\;\;\Sigma_2\setminus\{q_2\}\cong\Sigma_{20}\cup\{(-\infty,0]\times S^1\}.$$
Here $\Sigma_{i0}\subset \Sigma_i$, $i=1,2$, are compact surfaces with boundary. We introduce the notations
$$
\Sigma_i(R_0)=\Sigma_{i0}\cup \{(s_i,t_i)|\;|s_i| \leq R_0\},\;\; \;\;\;\;\;\Sigma(R_0)=\Sigma_1(R_0)\cup \Sigma_2(R_0).$$
For any gluing parameter $(r,\tau)$ with $r\geq R_{0}$ and $\tau\in S^1$ we construct a surface $\Sigma_{(r)}$
with the gluing formulas:
\begin{equation}\label{gluing surface}
s_1=s_2 + 2r,\;\;\;t_1=t_2 + \tau.
\end{equation}
where we use $(r)$ to denote gluing parameters.
\v
We will use the cylinder coordinates to describe the construction of $u_{(r)}:\Sigma_{(r)}\to M^+$. Write
$$u=(u_1,u_2),\;\;u_{i}:\Sigma_i\rightarrow M^+\; with \;\;u_{1}(q)=u_{2}(q).$$
We choose local normal coordinates $(x^{1},\cdots,x^{2m})$ in a neighborhood  $O_{u(q)}$ of $u(q)$ and choose $R_0$ so large that $u(\{|s_i|\geq \frac{r}{2}\})$ lie in $O_{u(q)}$ for any $r>R_0$. We glue the map $(u_1,u_2)$ to get a pregluing maps $u_{(r)}$ as follows. Set\\
\[
u_{(r)}=\left\{
\begin{array}{ll}
u_1 \;\;\;\;\; on \;\;\Sigma_{10}\bigcup\{(s_1,\theta_1)|0\leq s_1 \leq
\frac{r}{2}, \theta_1 \in S^1 \}    \\  \\
u_1(q)=u_2(q) \;\;
on \;\{(s_1,\theta_1)| \frac{3r}{4}\leq s_1 \leq
\frac{5r}{4}, \theta_1 \in S^1 \}  \\   \\
u_2 \;\;\;\;\; on \;\;\Sigma_{20}\bigcup\{(s_2,\theta_2)|0\geq s_2
\geq - \frac{r}{2}, \theta_2 \in S^1 \}     \\
\end{array}.
\right.
\]
To define the map $u_{(r)}$ in the remaining part we fix a smooth cutoff
function $\beta : {\mathbb{R}}\rightarrow [0,1]$ such that
\begin{equation}\label{def_beta}
\beta (s)=\left\{
\begin{array}{ll}
1 & if\;\; s \geq 1 \\
0 & if\;\; s \leq 0
\end{array}
\right.
\end{equation}
and $\sqrt{1-\beta^2}$ is a smooth function,  $0\leq \beta^{\prime}(s)\leq 4$ and $\beta^2(\frac{1}{2})=\frac{1}{2}.$
 We define\\
$$u_{(r)}= u_1(q)+ \left(\beta\left(3-\frac{4s_1}{r}\right)(u_1(s_1,\theta_1)-u_1(q)) +\beta\left(\frac{4s_1}{r}-5\right)(u_2(s_1-2r,\theta_1-\tau)- u_2(q))\right).$$
\vskip 0.1in
\noindent
\v

Denote
$$\beta_{1;R}(s_1)=\beta\left(\frac{1}{2}+\frac{r-s_1}{R}\right),\;\;\;
\beta_{2;R}(s_{2})=\sqrt{1-\beta^2\left(\frac{1}{2}-\frac{s_{2}+r}{R}\right)}, $$
where $\beta$ is the cut-off function defined in \eqref{def_beta}.
For any $\eta \in
C^{\infty}(\Sigma_{(r)};u_{(r)}^{\ast}TM^+\otimes \wedge_{j}^{0,1}T\Sigma_{(r)})$,
let
$$
\eta_{i}(p) =\left\{
\begin{array}{ll}
\eta  & if\;\; p\in \Sigma_{i0}\cup\{|s_{i}|\leq r-1\}\\
\beta_{i;2}(s_{i})\eta(s_{i},t_{i}) & if\;\; p\in \{r-1\leq |s_{i}|\leq r+1\} \\
0 & otherwise.
\end{array}
\right..
$$
If no danger of confusion  we will simply write $\eta_{i}=\beta_{i;2}\eta.$ Then $\eta_{i}$  can be considered as
a section over $\Sigma_i$. Define
\begin{equation}
\|\eta\|_{r,k-1,2,\alpha}=\|\eta_{1} \|_{\Sigma_1,j_1,k-1,2,\alpha} +
\|\eta_{2} \|_{\Sigma_2,j_2,k-1,2,\alpha}.
\end{equation}
\v
We now define a norm $\|\cdot\|_{r,k,2,\alpha}$ on
$C^{\infty}(\Sigma_{(r)};u_{(r)}^{\ast}TM^+).$ For any section
$h\in C^{\infty}(\Sigma_{(r)};u_{(r)}^{\ast}TM^+)$ denote
$$h_0=\int_{ {S}^1}h(r,t)dt,$$
$$h_1(s_1,t_1) = (h-\hat h_{0})(s_1,t_1)\beta_{1;2}(s_1),\;\;\;h_2(s_2,t_2)= (h-\hat h_{0})(s_2,t_2)\beta_{2;2}(s_{2}).$$
We define
\begin{equation}
\|h\|_{r,k,2,\alpha}=\|h_1\|_{\Sigma_1,j_1,k,2,\alpha} +
\|h_2\|_{\Sigma_2,j_2,k,2,\alpha}+|h_{0}|.
\end{equation}
Denote the resulting completed spaces by $W^{k-1,2,\alpha}(\Sigma_{(r)};u_{(r)}^{\ast}TM^+\otimes \wedge_{j_r}^{0,1}T\Sigma_{(r)})$ and $W^{k,2,\alpha}(\Sigma_{(r)};u_{(r)}^{\ast}TM^+)$  respectively.
\v
This pregluing procedure can be generalized to pregluing several nodes.

\subsection{\bf Norms for maps from marked nodal Riemann surfaces}\label{norms for smooth surfaces}

One can generalize the norms $\|h\|_{k,2,\alpha}$,
$\|\eta\|_{k-1,2,\alpha}$ and $\widetilde{\mathbf{O}}_{b_{o}}(\delta,\rho)$, $\mathbf{O}_{[b_{o}]}(\delta,\rho)$ to marked nodal Riemann surfaces. Let $(\Sigma,j,{\bf y}, {\bf p},\nu)$ be a marked nodal Riemann surface with nodal structure $\nu$
of genus $g$ with $m$ distinct marked points ${\bf y}=(y_1,...,y_m)$, $\mu$  distinct puncture points ${\bf p}=(p_1,...,p_{\mu})$. Denote by $\mathbf{q}=(q_1,...,q_{\mathbf{e}})$ the set of nodal points of $\Sigma$. Put $\stackrel{\circ}{\Sigma}= \Sigma - \{\mathbf{y,p,q}\}$.
Then $\stackrel{\circ}\Sigma$ is a Riemann surface with additional punctures $a_{j}, b_{j}$ in the place of the $j$th node of $\Sigma,$ $j=1,\cdots, \mathfrak{e}$. Let $u:\Sigma\to M^+$
is a smooth map satisfying the nodal conditions
$$
\{a_j,b_j\}\in\nu\implies u(a_j)=u(b_j).
$$
For each node $q_{j},j=1,\cdots,\ell$ there is a neighborhood   isomorphic to
$$\{(z_{j},w_{j})\in \mathbb C^{2}| |z_{j}|<1,|w_{j}|<1,z_{j}w_{j}=0\}.$$
Denote by $\Sigma_{i}$ the connected components  of $\stackrel{\circ}\Sigma$, $i=1,\cdots, \iota$. Suppose that $\Sigma_i$ has $n_{i}$ marked points, $q_{i}$ punctures and has genus $g_{i}$.
\v
We can parameterize a neighborhood of $\stackrel{\circ}\Sigma$ in the deformation space  by Beltrami differentials.  Let $z_i$ (resp. $w_{i}$) be a local   coordinate around $a_{i}$ (resp. $b_{i}$), $z_i(a_{i})=0,w_i(b_{i})=0$, $i=1,\cdots, \mathfrak{e}$.
Let $\mathbb{U}_{j}=\{p\in\Sigma||z_{j}|(p)<1\}$ and  $\mathbb{V}_{j}=\{p\in\Sigma||w_{j}|(p)<1\}$ be
disjoint neighborhoods of the punctures $a_{j}$ and $b_{j},j=1,\cdots,\mathbf{e}$.
We pick an open set $\mathbb{U}_{o}\subset \stackrel{\circ}\Sigma$ such that each component of $\stackrel{\circ}\Sigma$ intersects $\mathbb{U}_o$ in a nonempty relatively compact set and the intersection $\mathbb{U}_o\bigcap(\mathbb{U}_j\cup \mathbb{V}_j)$ is empty for all $j$. Denote $N=\sum\limits_{i=1}^{\fc} (3g_{i}-3+n_{i}+q_{i})$.
Choose Beltrami differentials $\nu_{j}, j=1,\cdots,N$ which are supported in $\mathbb{U}_o$ and
form a basis of the deformation space at $\Sigma $.
Let $\mathbf{s}=(\fs_{1},\cdots,\fs_{N})\in \mathbb C^{N}$,
$\nu=\sum\limits_{i=1}^{N} \fs_{i}\nu_{i}.$ Assume $|\mathbf{s}|$ small enough such that $|\nu|<1$. The nodal surface $\Sigma_{\mathbf{s},0}$ is obtained by solving the Beltrami equation $\bar{\p}w=\nu(\mathbf{s})w$.
\v
	We recall the plumbing construction for $\Sigma$ with a pair of punctures $a_{j},b_{j}$.  	
Let $z_{j,\mathbf{s}},$ $w_{j,\mathbf{s}}$ be the canonical coordinates in
	$\mathbb{U}_{j}, \mathbb{V}_{j}$ near $a_{j},b_{j}$ respectively, thus
$$
ds^2_{\mathbf{s},0} (z_{j,\mathbf{s}})=\frac{|dz_{j,\fs}|^2}{|z_{j,\fs}|^2 \log^2 |dz_{j,\fs}|},\;\;\;ds^2_{\fs,0} (w_{j,\fs})=\frac{|dw_{j,\fs}|^2}{|w_{j,\fs}|^2 \log^2 |dw_{j,\fs}|}.
$$
where   $ds^2_{\fs,0} $ be the normalized hyperbolic metric on $\Sigma_{\fs,0}$  of curvature $-1$.
 As \cite{DasMes} denote
$$
	F_{j,\mathbf{s}}=z_{j}\circ z_{j,\mathbf{s}}^{-1},\;\;\;\;	G_{j,\mathbf{s}}=w_{j}\circ w_{j,\mathbf{s}}^{-1}.
$$
  By the removalbe singularity theorem and setting $\widetilde F_{j,\mathbf{s}}=F_{j,\mathbf{s}}/F_{j,\mathbf{s}}'(0)  $ and  $\widetilde G_{j,\mathbf{s}}=G_{j,\mathbf{s}}/G_{j,\mathbf{s}}'(0)$, if necessary , we can assume that
  $$
  F_{j,\mathbf{s}}(0)=0,\;\;\;  F'_{j,\mathbf{s}}(0)=1,\;\;\;G_{j,\mathbf{s}}(0)=0,\;\;\;  G'_{j,\mathbf{s}}(0)=1.\;\;\;  $$
Since $\mathbb{U}_{o}$ is disjoint from the $\mathbb{U}_{j}, \mathbb{V}_{j}$, the $F_{j,\mathbf{s}},G_{j,\mathbf{s}}$ are also holomorphic onto their image. For any $  \mathbf{t}=(\ft_{1},\cdot,\cdot,\ft_{\mathbf{e}})$
with $0<|\ft_{j}|<1$, remove the discs $|z_{j}|< |\ft_{j}|$ and $|w_{j}|< |\ft_{j}|$ when $|\ft_{j}|$ small, and identify  $z_{j}$ via the plumbing equation  $$w_{j}=\frac{\ft_{j}}{z_{j}}.$$ We can rewrite the equation as
	$$
	(F_{j,\mathbf{s}}\circ z_{j,\mathbf{s}})\cdot (G_{j,\mathbf{s}}\circ w_{j,\mathbf{s}})=\ft_{j}.$$
Then we form a new Riemann surface $\Sigma_{\mathbf{s},\ft}.$ We call $(\ft_{1},\cdot\cdot\cdot,\ft_{\mathbf{e}})$ plumbing coordinate. We obtain a family of Riemann surfaces over $\Delta_{\mathbf{s}}\times \Delta_{\ft}$,
whose fiber over $(\fs,\ft)$ is the Riemann surface $\Sigma_{\fs,\ft}$, where $\Delta_{\mathbf{s}}=(\Delta)^N\subset \mathbb{C}^N$, $\Delta_{\ft}=(\Delta)^{\mathbf{e}}\subset \mathbb{C}^{\mathbf{e}}$ are  polydiscs.
\v
In the coordinate system $(\fs,\ft)$ the $g_{\mathsf{wp}}$ metric induces a $Diff^+(\Sigma)$-invariant distance $d_{\mathbf{s},\ft}(\cdot,\cdot)$ on $\Delta_{\mathbf{s}}\times \Delta_{\ft}$.
Put
$$O(\delta)=\{(\fs,\ft)\mid d_{\mathbf{s}, \ft }((0,0),(\fs,\ft))<\delta \}.$$
We can choose $\delta$ small such that $\mathbf{G}_{(\fs,\ft)}$ can be imbedded into $\mathbf{G}_{(0,0)}$ as a subgroup for any $(\fs,\ft)\in O(\delta)$. Denote by $im(\mathbf{G}_{(\fs,\ft)})$ the imbedding.
\v
Let $u_{\mathbf{s},0}:\Sigma_{\mathbf{s},0}\to M^+$ be a $W^{k,2,\alpha}$-map. We can construct $u_{\mathbf{s},\ft}:\Sigma_{\mathbf{s},\ft}\to M^+$.
For any $h\in C^{\infty}_{c}(\Sigma_{\fs,0};u_{\fs,0}^{\ast}TM^+)$ and
any section $\eta \in
C^{\infty}_{c}(\Sigma_{\fs,0}, u_{\fs,0}^{\ast}TM^+\otimes \wedge^{0,1}_jT^{*}\Sigma_{\fs,0})$
we define the norms $\|h\|_{\fs,k,2,\alpha}$ and $\|\eta\|_{\fs,k-1,2,\alpha}$. For any section
$h\in C^{\infty}(\Sigma_{\fs,\ft};u_{\fs,\ft}^{\ast}TM^+)$
and any $\eta \in
C^{\infty}(\Sigma_{\fs,\ft};u_{\fs,\ft}^{\ast}TM^+\otimes \wedge_{j}^{0,1}T\Sigma_{\fs,\ft})$, we define the norms $\| h\|_{\fs,\ft,k,2,\alpha}$ and $\|\eta\|_{\fs,\ft,k-1,2,\alpha}$. We do these in the same way as for one node case.
\v
Let $b_o=(\Sigma,0,0,u)$.
Set
$$\widetilde{\mathbf{O}}_{b_{o}}(\delta_{b_o},\rho_{b_o}):=\left\{((\fs,\ft),v_{\fs,\ft})\;|\; d_{\fs,\ft}((0,0),(\fs,\ft))<\delta_{b_o}, \|h\|_{\fs,\ft,k,2,\alpha}<\rho_{b_o} \right\},$$
$$\mathbf{O}_{[b_{o}]}(\delta_{b_{o}},\rho_{b_{o}})=
\widetilde{\mathbf{O}}_{b_{o}}(\delta_{b_{o}},\rho_{b_{o}})/G_{b_o},$$
where $v_{\fs,\ft}=\exp_{u_{\fs,\ft}}(h).$

\v

\section{Holomorphic cascades in $M^{\pm}$}\label{holo block}

\v
 We discuss $M^+$, for $M^-$ the discusses are the same. A {configuration} in $M^+$ is a tuple $(\Sigma,j,{\bf y}, {\bf p},\nu, u)$ where $(\Sigma,j,{\bf y}, {\bf p},\nu)$ is a marked nodal Riemann surface (see~\cite[\S3]{RS})  of genus $g$ with $m$ distinct marked points ${\bf y}=(y_1,...,y_m)$, $\mu$  distinct puncture points ${\bf p}=(p_1,...,p_{\mu})$, and $u:\Sigma\to M^+$
is a smooth map satisfying the nodal conditions.
 Suppose that $u(z)$ converges to a $k_i$-periodic
orbit $x_{i}$ as $z$ tends to $p_i$. We call $(\Sigma,j,{\bf y}, {\bf p},\nu, u)$
 a relative {configuration} in $M^{+}$. Let $\Sigma =\bigcup_{i=1}^d\Sigma_{i}$, where $(\Sigma_i, j_i)$ is a smooth Riemann surface.
\v
The configuration $(\Sigma,j,{\bf y},{\bf p},\nu, u)$ is called {holomorphic} if the restriction to every $\Sigma_i$ of $u$ satisfies \eqref{j_holomorphic_maps-1}.
\v
Let $\mathcal{J}(\Sigma)\subset End(T\Sigma)$ denote the manifold of complex structures on $\Sigma$. Denote by $Diff^+(\Sigma,\nu)$ the group of orientation preserving diffemorphisms of $\Sigma$, that preserve the nodal structure. Denote by $Diff^+_{0}(\Sigma,\nu)$ the identity component of $Diff^+(\Sigma,\nu).$
For any $\varphi\in Diff^+(\Sigma,\nu)$, $\varphi$ acts on $\mathcal J(\Sigma)\times (\Sigma^{m+\mu} \setminus \Delta)$ by
the holomorphic diffeomorphisms
$$(j,{\bf y},{\bf p})\longmapsto (\varphi^*j, \varphi^{-1}(\mathbf{y}), \varphi^{-1}({\bf p})),$$
where $\Delta\subset\Sigma^{m+\mu}$ denotes the fat diagonal, i.e. set of all
$(m+\mu)$-tuples of points in $\Sigma^{m+\mu}$ where at least two components are equal.
\v

\v
\v

\begin{definition}\label{holomorphic_block_map_equiv_M}
Two relative $(j,J)$-holomorphic configurations $b=(j ,{\bf y},\nu,{\bf p},u)$
 and $\check{b} =(\check{\bf j} ,{\bf \check{y}},\check \nu,{\bf \check{p}}, \check{u})$  in $M^+$
 are called equivalent if there exists a $\varphi\in Diff^+(\Sigma,\nu)$ such that
\begin{itemize}
\item[{\bf(1)}] $\varphi(j ,{\bf y},\nu,{\bf p})=(\check{\bf j} ,{\bf \check{y}},\check \nu,{\bf \check{p}})$,
\item[{\bf(2)}]  $\check{u}= u\circ \varphi$.
\end{itemize}
\end{definition}

\v
\begin{definition}  We put \begin{align*} Aut(b)=
\{\varphi\in Diff^+(\Sigma,\nu)| &\varphi \mbox{ is an automorphism satisfying (1), (2)} \; in \; Definition \; \ref{holomorphic_block_map_equiv_M}  \} .\end{align*}
We call it the automorphism group of $b$.
\end{definition}

\begin{definition}
A relative $(j,J)$-holomorphic configuration $b$ in $M^{+}$ is called stable if $Aut(b)$ is a finite group.
\end{definition}

We collapse the $S^1$-action on $\widetilde{M}=H^{-1}(0)$
to obtain symplectic manifolds $\overline{M}^{+}$ and $\overline{M}^{-}$. The
reduced space $Z$ is a codimension 2 symplectic submanifold of both
$\overline{M}^{+}$ and $\overline{M}^{-}$. By using the removable singularities theorem we get a
$(j,J)$-holomorphic map $\bar{u}$ from $\Sigma $ into $\overline{M}^{\pm}$.
Therefore, we have a natural identification of finite energy $J$-holomorphic maps
into $M^{\pm}$ and $(j,J)$-holomorphic maps into the closed symplectic manifolds
$\overline{M}^{\pm}$. Under this identification, the condition
that $u$ converges
to a $k$-multiple periodic orbit at $p$ is naturally interpreted as
$\bar{u}$ being tangent to $Z$ at $p$ with order $k$. Let $A=[\bar{u}(\Sigma)]$.
It is obvious that
\begin{equation}\label{energy_bound}
E_{\phi}(u)=\omega_{\phi}(A)
\end{equation}
which is independent of $\phi $.
We fix a homology class $A\in{H_{2}(\overline{M}^{+},\Z)}$ and a set
$\{k_1,...,k_{\mu}\}$.
We have
$$\sum_{i=1}^{\mu} k_i= \#(A\cap Z).$$

\v

Denote by ${\mathcal{M}}_{A}(M^{+},g,m+\mu,{\bf k},\nu)$ the space of equivalence classes
of all relative  stable $(j,J)$-holomorphic configurations in $M^{+}$ representing the homology class $A$ and converging
to a $k_i$-periodic orbit as $z$ tends to $p_i$.
The moduli space ${\mathcal{M}}_{A}(\overline{M}^{+},g,m+\mu,{\bf k},\nu)$ can be identified with ${\mathcal{M}}_{A}(M^{+},g,m+\mu,{\bf k},\nu)$. We call ${\mathcal{M}}_{A}(M^+,g,m+\mu,{\bf k},\nu) $ a holomorphic  cascade in $M^+$.

\v
\begin{remark}\label{holo block-1}
It is possible that there are several holomorphic cascades ${\mathcal{M}}_{A_i}(\overline{M}^{+},g_i,m_i+\mu_i,{\bf k}_i,\nu_{i})$, $i=1,...,l$, in $M^+$.
\end{remark}

\section{Holomorphic cascades in ${\mathbb{R}}\times \widetilde{M}$}

Note that the space ${\mathcal{M}}_{A}(M^{+},g,m+\mu,{\bf k},\nu)$ is not large enough to compactify the Moduli space of all relative stable $(j,J)$-holomorphic configurations in $M^{+}$, we need consider ${\mathcal{M}}_{A}({\mathbb{R}}\times \widetilde{M},{\bf k}^+, {\bf k}^-,\nu)$, which will be studied in this section.

Let $(\Sigma, j, {\bf y},\nu, {\bf p}^{+},{\bf p}^{-})$ be a marked nodal Riemann surface of genus $g$ with $m$ marked points
${\bf y}=(y_1,...,y_m)$ and $\mu^{\pm}$ puncture points ${\bf p}^{+}=(p_1^{+},...,p_{\mu^{+}}^{+})$, ${\bf p}^{-}=(p_1^{-},...,p_{\mu^{-}}^{-})$.
Let $u:{\Sigma} \rightarrow {\mathbb{R}}\times \widetilde{M}$ be a $(j,J)$-holomorphic map. Suppose that $u(z)$ converges to a $k_i^{+}$(resp. $k_j^{-}$)-periodic orbit $x_{k_i^{+}}$(resp. $x_{k_j^{-}}$) as $z$ tends to $p_i^{+}$(resp. $p_j^{-}$). We call $(\Sigma, j, {\bf y},\nu, {\bf p}^{+},{\bf p}^{-}, u)$ a relative $(j,J)$-holomorphic configuration in
${\mathbb{R}}\times \widetilde{M}$.

\begin{definition}\label{equivalent_R}
Two relative holomorphic configurations $b =( j, {\bf y},\nu, {\bf p}^{+},{\bf p}^{-}, u)$
and $\check{b} =(\check{\bf j},{\bf \check{y}}, \nu, \check{\bf p}^{+},\check{\bf p}^{-}, \check{u})$  in ${\mathbb{R}}\times \widetilde{M}$
 are called equivalent if there exists a $\varphi\in Diff^+(\Sigma,\nu)$ such that
\begin{itemize}
\item[{\bf(1)}] $\varphi( j ,{\bf y},\nu,{\bf p}^{+},{\bf p}^{-})=(\check{\bf j} ,{\bf \check{y}},\check \nu,\check{\bf p}^{+},\check{\bf p}^{-})$,
\item[{\bf(2)}]  $\check{u}= u\circ \varphi$.
\end{itemize}
\end{definition}

\begin{definition}  Put \begin{align*}Aut(b)=\{\varphi\in Diff^+(\Sigma,\nu)| &\varphi \mbox{ is an automorphism satisfying (1), (2)}  \\ &\mbox{ in Definition \ref{equivalent_R}}  \} .\end{align*} We call it the automorphism group of
$b$.
\end{definition}

\begin{definition}
A relative $J$-holomorphic configuration $b$ in ${\mathbb{R}}\times \widetilde{M}$ is stable if $Aut(b)$ is a finite group.
\end{definition}

We collapse the $S^1$-action at $\pm \infty$
to obtain the symplectic manifold $\mathbb P(\mathcal N\oplus \mathbb C)$. The
reduced space $Z$ is a codimension 2 symplectic submanifold of $\mathbb P(\mathcal N\oplus \mathbb C)$. By using the removable singularities theorem we get a $(j,J)$-holomorphic map $\bar{u}$ from $\Sigma $ into $\mathbb P(\mathcal N\oplus \mathbb C)$.
Therefore, we have a natural identification of finite energy $(j,J)$-holomorphic maps
into ${\mathbb{R}}\times \widetilde{M}$ and $(j,J)$-holomorphic maps into the closed symplectic manifold
$\mathbb P(\mathcal N\oplus \mathbb C)$.  Let $A=[\bar{u}(\Sigma)]$.
It is obvious that
\begin{equation}E_{\phi}(u)=\omega_{\phi}(A)
\end{equation}
which is independent of $\phi $.
We fix a homology class $A\in{H_{2}(\mathbb P(\mathcal N\oplus \mathbb C),\mathbb{Z})}$ and   fixed sets
$\{k_1^{+},...,k_{\mu^{+}}^{+}\}$, $\{k_1^{-},...,k_{\mu^{-}}^{-}\}$ .
We have
$$\#(A\cap Z_{\infty})=\sum_{i=1}^{\mu^{+}} k_i^{+},\;\;\;\;\;\;\;  \#(A\cap Z_{0})=-\sum_{i=1}^{\mu^{-}} k_i^{-}.$$
Then  we define
${\mathcal{M}}_{A}({\mathbb{R}}\times \widetilde{M},g,m+\mu^{+}+\mu^{-},{\bf k}^+, {\bf k}^-,\nu)$ to be the space of equivalence classes
of all relative stable $(j,J)$-holomorphic configurations in ${\mathbb{R}}\times \widetilde{M}$ representing the homology class $A$ and converging
converges to a $k_i^{\pm}$-periodic orbit as $z$ tends to $p_i^{\pm}$.
The moduli space ${\mathcal{M}}_{A}({\mathbb{R}}\times \widetilde{M},g,m+\mu^{+}+\mu^{-} ,{\bf k}^+, {\bf k}^-,\nu)$ can be identified with ${\mathcal{M}}_{A}(\mathbb P(\mathcal N\oplus \mathbb C),g,m+\mu^{+}+\mu^{-} ,{\bf k}^+, {\bf k}^-,\nu)$.

There
is a uniform bound on $\widetilde{E}$ for any $(j,J)$-holomorphic configuration in
${\mathcal{M}}_{A}({\mathbb{R}}\times \widetilde{M},g,m+\mu^{+}+\mu^{-},{\bf k}^+, {\bf k}^-,\nu)$.
\v
We call ${\mathcal{M}}_{A}({\mathbb{R}}\times \widetilde{M},g,m+\mu^{+}+\mu^{-},{\bf k}^+, {\bf k}^-,\nu) $ a holomorphic cascade in ${\mathbb{R}}\times \widetilde{M}$.
\v
There is a natural $\mathbb{C}^*$ action on ${\mathbb{R}}\times \widetilde{M}$, which is given in terms of the coordinates $(a,\theta,\mathbf{w})$ by
\begin{equation}
a\longmapsto a + r,\;\;\theta\longmapsto \theta + \vartheta,\;\;\mathbf{w}\longmapsto \mathbf{w}
\end{equation}
for any $z=e^{r+2\pi\sqrt{-1}\vartheta}$. This action induces an action of $\mathbb{C}^*$ on ${\mathcal{M}}_{A}({\mathbb{R}}\times \widetilde{M},g,m+\mu^{+}+\mu^{-},{\bf k}^+, {\bf k}^-,\nu) $ in a natural way:
$$e^{\ell+2\pi \sqrt{-1} \theta_{0}}\cdot(a(s,t),\theta(s,t),\mathbf{w}(s,t))=
(a(s,t)+\ell,\theta(s,t)+\theta_{0},\mathbf{w}(s,t)).$$
 Denote
$${\mathcal{M}^*}_{A}({\mathbb{R}}\times \widetilde{M},g,m+\mu^{+}+\mu^{-},{\bf k}^+, {\bf k}^-,\nu)={\mathcal{M}}_{A}({\mathbb{R}}\times \widetilde{M},g,m+\mu^{+}+\mu^{-},{\bf k}^+, {\bf k}^-,\nu)/\mathbb{C}^*.$$
\v
\begin{remark}\label{holo block-1}
It is possible that there are several holomorphic cascades  ${\mathcal{M}}_{A_i}({\mathbb{R}}\times \widetilde{M},g_i,m_i+\mu_i,{\bf k}_i,\nu_{i})$, $i=1,...,l$, in one copy of ${\mathbb{R}}\times \widetilde{M}$. The $\mathbb{C}^*$ action on ${\mathbb{R}}\times \widetilde{M}$ induces a $\mathbb{C}^*$ action on every holomorphic cascade. We call the quotient holomorphic cascades
, together with the copy of ${\mathbb{R}}\times \widetilde{M}$, a holomorphic block, denoted by ${\mathcal{M}^*}({\mathbb{R}}\times \widetilde{M},\cdot)$.
\end{remark}

\begin{remark}\label{explain the termilonogy cascades} We explain why we use the terminology "holomorphic cascade". To define relative Gromov-Witten invariants we need a compactified moduli space ${\overline{\mathcal{M}}}_{A}(M^{+};g,m+\mu,{\bf k},\nu)$.
Roughly speaking, each element in ${\overline{\mathcal{M}}}_{A}(M^{+};g,m+\mu,{\bf k},\nu)$ is one obtained by gluing several holomorphic cascades in $N$ ( see \S\ref{Weighted dual graph with an oriented decomposition} and the Chapter \S\ref{compact_theorem} ). There is a natural partial order in the set of the holomorphic cascades, which looks like
"the water flows from a higher level to a lower one". So we use the terminology "holomorphic cascades".
\end{remark}

\v
The following lemma is well-known (see \cite{MS})
 \begin{lemma}\label{energy_lower_bound} Let $(M,\omega)$ be a compact symplectic manifold with $\omega$-tamed almost complex structure $J.$ Then there exists $\hbar>0$ such that if $u:\Sigma\longrightarrow M$ is a nonconstant J-holomorphic map from a closed Riemann surface $\Sigma$ to $M,$ then
$$
\int_{\Sigma} u^{*} \omega\geq  \hbar .
$$
\end{lemma}
\v\n

\v
\begin{corollary}\label{energy_lower_bound-1}
There is a constant $\hbar>0$ such that
for every finite energy $J$-holomorphic map $u=(a,\widetilde{u}):
{\Sigma}\rightarrow {\mathbb{R}}\times \widetilde{M} $ with
$\widetilde{E}(u)\neq 0 $ we have $\widetilde{E}(u)\geq \hbar$.
\end{corollary}
\noindent
{\bf Proof:} Consider the $\widetilde{J}$-holomorphic map
$$\hat{u}=\pi\circ \widetilde{u}:\stackrel{\circ}{\Sigma} \rightarrow Z .$$
$\hat{u}$ extends to a $\widetilde{J}$-holomorphic curve
from $\Sigma $ to $Z$. Then the assertion follows from Lemma \ref{energy_lower_bound}. $\;\;\;\Box$
\vskip 0.1in
\noindent

Since $Z$ is compact, there exists a constant
  $C>0$ such that
\begin{equation*}
\tau_0(v,\widetilde{J}v)\leq
C d\lambda(v, \widetilde{J}v) )
\end{equation*} for all $v\in TZ.$  Let $b=(j,\mathbf {y},\mathbf {p},u)\in {\mathcal{M}}_{A}({\mathbb{R}}\times \widetilde{M},g,m+\mu^{+}+\mu^{-},{\bf k}^+, {\bf k}^-,\nu) $ and $u=(a,\widetilde{u}):
{\Sigma}\rightarrow {\mathbb{R}}\times \widetilde{M} $ be a finite energy $J$-holomorphic map. If  $\widetilde{E}(u)= 0$, we have $2g+m-2>0.$
Assume that $\widetilde{E}(u)\neq 0 $.
It follows from Corollary \ref{energy_lower_bound-1} that
\begin{equation}
\hbar\leq \widetilde E(u)\leq C\int_{\Sigma} u^{*}d\lambda =C(\sum_{i=1}^{\mu^{+}} k_i^{+} -\sum_{i=1}^{\mu^{-}} k_i^{-})
\end{equation}
By  \eqref{tame} and Corollary \ref{energy_lower_bound-1}  we have
$$
\omega_{\phi}(A)=E_{\phi}(u)\geq C\widetilde E(u)\geq C\hbar.
$$
 Then we have
 \begin{lemma}\label{lem_com_Proj}
 For any holomorphic cascade ${\mathcal{M}}_{A}({\mathbb{R}}\times \widetilde{M},g,m+\mu^{+}+\mu^{-},{\bf k}^+, {\bf k}^-,\nu) $  one of following holds
\begin{itemize}
\item[(1)] $2g+m-2>0$,
\item[(2)] $A\neq 0$ and $\sum_{i=1}^{\mu^{+}} k_i^{+} -\sum_{i=1}^{\mu^{-}} k_i^{-}\geq 1.$
\end{itemize}
\end{lemma}

\section{Homology class}\label{homology_class} Let $b=(u_1,u_2; \Sigma_1 \wedge \Sigma_2,j_1,j_2)\in {\mathcal{M}}_{A_1}(\overline{M}^{+},g_1,m_1+1,k,\nu_{1})\times_Z {\mathcal{M}}_{A_2}(\mathbb P(\mathcal N\oplus \mathbb C),g_2,m_2+\mu + 1,{\bf k}^+, k,\nu_{2})$,
where $(\Sigma_1,j_1)$ and $(\Sigma_2,j_2)$ are smooth Riemann surfaces of genus $g_1$ and $g_2$ joining at $p$ and
$u_1: \Sigma_1 \rightarrow M^+$, $u_2: \Sigma_2 \rightarrow {\mathbb{R}}\times \widetilde{M}$ are J-holomorphic maps such that
$u_i(z)$ converge to the same $k$-periodic orbit $x$ as $z\rightarrow p$.

Denote by $Z^{(1)}_{\infty}\in \overline{M}^{+}$, $Z^{(2)}_{0}, Z^{(2)}_{\infty}\in \mathbb P(\mathcal N\oplus \mathbb C)$ the divisors respectively.
We have
\begin{equation}
[u_1(\Sigma_1)]=A_1,\;\;[u_2(\Sigma_2)]=A_2,\end{equation}
\begin{equation}\label{homology-0}
\#(A_1\cap Z^{(1)}_{\infty})=k,\;\;\#(A_2\cap Z^{(2)}_{0})=-k, \;\;\#(A_2\cap Z^{(2)}_{\infty})=\sum_{i=1}^{\mu} k_i.\end{equation}

For any parameter $(r):=(r,\tau)$, we glue $M^+$ and ${\mathbb{R}}\times \widetilde{M}$
 to get again $M^+$, glue $\Sigma_1$ and $\Sigma_2$ to get $\Sigma_{(r)}$, and construct a pre-gluing map $u_{(r)}:\Sigma_{(r)}\rightarrow M^+$.
 It is easy to see that
\begin{equation}\label{homology-1}
A:=[u_{(r)}(\Sigma_{(r)})]= A_1 + A_2,\;\;\;\#(A\cap Z^{(2)}_{\infty})=\sum_{i=1}^{\mu} k_i.\end{equation}
In general the homology class of $u_{(r)}$ depends on the $J$-holomorphic curve representatives
$u_1, u_2$ instead of the homology classes $A_1, A_2$. One can understand it as follows. Recall that there is a map
$$\pi: \overline{M}^+\rightarrow \overline{M}^+\cup_Z \mathbb P(\mathcal N\oplus \mathbb C).$$
$\pi$ induces a homomorphism
$$\pi_*:  H_2(\overline{M}^+, \mathbb Z)\rightarrow H_2(\overline{M}^+\cup_Z \overline{M}^-, \mathbb Z).$$
If $b'=(u'_1,u'_2; \Sigma_1 \wedge \Sigma_2,j_1,j_2)\in {\mathcal{M}}_{A_1}(\overline{M}^{+},g_1,m_1+1,k,\nu_{1})\times_Z {\mathcal{M}}_{A_2}(\mathbb P(\mathcal N\oplus \mathbb C),g_2,m_2+\mu + 1,{\bf k}^+, k,\nu_{2})$ is another element and glued to $u_{(r)}'$.
When $\ker \pi_*\neq 0$, $[u_{(r)}(\Sigma_{(r)})]$ and $[u_{(r)'}(\Sigma_{(r)})]$ could be different from a
vanishing 2-cycle $T^2$ in $\ker \pi_*$, i.e., $A'= A + T^2$. We have
\begin{equation}\label{homology_condition_2}
\#([u_{(r)}(\Sigma_{(r)})]\cap Z^{(2)}_{\infty})=\sum_{i=1}^{\mu} k_i= \#([u'_{(r)}(\Sigma_{(r)})]\cap Z^{(2)}_{\infty}).
\end{equation}

On the other hand, as $u_1$, $u_2$ are $J$-holomorphic maps, and out of the gluing part, $u_{(r)}$ is also $J$-holomorphic map, we have
$$E(u_1)=\omega(A_1),\;\;E(u_2)=\omega(A_2),$$
$$E(u'_1)=\omega(A_1),\;\;E(u'_2)=\omega(A_2),$$
and
\begin{equation}\label{energy_condition_1}
\mid E(u_{(r)})- E(u'_{(r)})\mid \rightarrow 0\;\;\;as\;\; r\rightarrow \infty.\end{equation}
When we compactify  our moduli space of relative stable
$J$-holomorphic maps we need only the properties
\eqref{energy_condition_1} and \eqref{homology_condition_2}, so we write
$A\circeq A_1 + A_2$, and say that the elements in ${\mathcal{M}}_{A_1}(\overline{M}^{+},g_1,m_1+1,k,\nu_{1})\times_Z {\mathcal{M}}_{A_2}(\mathbb P(\mathcal N\oplus \mathbb C),g_2,m_2+\mu + 1,{\bf k}^+, k,\nu_{2})$ have the homology class $A \circeq A_1+A_2$.
\v
This can be immediately generalize to the case that two holomorphic cascades have several common nodal points.

\section{Weighted dual graph}
\v
\subsection{Weighted dual graph for holomorphic cascades}\label{Weighted dual graph-1}
It is well-known that the moduli space of stable maps in a compact symplectic manifold has a stratification indexed by the combinatorial type
of its decorated dual graph. In this section we generalizes this construction
to holomorphic cascades in $N$, where $N$ is one of $M^+$, $M^-$ and
${\mathbb{R}}\times \widetilde{M}$.

Let $G$ be a connected graph. Denote $G=(V(G), E(G))$, where $V (G)$ is a
finite nonempty set of vertices and $E(G)$ is a finite set of edges.

\begin{definition} Let $g$, $m$ and $\mu$ be nonnegative
integers, $A\in H_{2}(\overline{M}^+,\mathbb Z)$.
A $(g,m+\mu, A, \mathbf{k})$-weighted dual graph $G$
consists of $(V(G),E(G))$
together with 4 weights:
\begin{itemize}
\item[(1)] $\mathfrak{g}:V(G) \rightarrow \mathbb Z_{\geq 0}$
assigning a nonnegative integer $g_{v}$ to each vertex $v$ such
 that $$g=\sum_{v\in V(G)}g_{v}+b_{1}(G),$$
where $b_1(G)$ is the first Betti number of the graph $G$;
\item[(2)] $\mathfrak{m}$: assign $m$ ordered tails
$\mathfrak{m}=(t_{1},\cdots,t_{m})$ to $V(G):$  attach $m_v$
tails to $v$ for each $v\in V(G)$ such that $m=\sum m_v$,
\item[(3)] $\mathfrak{h}:V(G) \rightarrow H_{2}(\overline{M}^+,\mathbb Z)$
assigning a $A_{v}\in  H_{2}(\overline{M}^+,\mathbb Z)$ to each vertex
$v$ such that $A=\sum A_{v}.$
\item[(4)] $\mathfrak{l}^{\bf k}$: assign  $\mu$ ordered
weighted half edges $(k_{1}e_{1},\cdots,k_{\mu}e_{\mu})$
to $V(G).$ First we assign half edges
$\mathfrak{l}=(e_{1},\cdots,e_{\mu})$ to $V(G)$: attach $\mu_v$
half edges to $v$ for each $v\in V(G)$ such that $\mu=\sum \mu_v$.
Then we assign  $\mu$ ordered weights
${\bf k}=(k_1,...,k_{\mu})\in (\mathbb{Z}^+)^{\mu}$ to the
half edges $\mathfrak{l}=(e_{1},\cdots,e_{\mu})$ such that
$$\sum_{i=1}^{\mu} k_i=Z^*(A),$$ where $Z^*$ is the Poincare dual of $Z$.
Denote the weighted half edges by
$$\mathfrak{l}^{\bf k}=(k_{1}e_{1},\cdots,k_{\mu}e_{\mu}).$$

 \end{itemize}

\v

\end{definition}

We denote the weighted dual graph $G$ by
$(V(G),E(G),A,\mathbf{k},\mathfrak{g}, \mathfrak{m}, \mathfrak{l}^{\bf k},
\mathfrak{h})$, or simply by $(V(G),E(G),g,m+\mu,A,\mathbf{k})$.
\v
By a leg of $G$ we mean either a tail or a
half-edge.

\begin{definition}
Let $G$ be a weighted dual graph. A vertex $v$ is called stable if one of the following holds:
\begin{itemize}
\item[(1)]  $2g_{v}+ val(v)\geq 3,$ where $val(v)$ denotes the sum of the number of legs attached to $v$;
\item[(2)] $A_v\not=0$.
\end{itemize}
  $G$ is called stable if all vertices are stable.
 \end{definition}

\subsection{Holomorphic cascades in $M^+$ of type $G$}\label{Holomorphic block of type G}
Let $G$ be a stable weighted dual graph with $N$ vertices $(v_1,...,v_N)$, $m$ tails and $\mu $
 half edges. We associate $G$ with a holomorphic cascade in $M^+$ as follows. Let
 $(\Sigma,{\bf y}, {\bf p})$ be a nodal Riemamm surface with $m$ marked points and $\mu$ puncture points.
 Let $A\in H_{2}(\overline{M}^{+},\mathbb Z)$.  A  stable
 $J$-holomorphic map of type $G$ is a quadruple $$(\Sigma, {\bf y},{\bf p}; u)$$ where
 $u:\Sigma\rightarrow M^{+}$ is a continuous map
satisfying the following conditions:
\begin{itemize}
\item [{\bf [A-1]}] $\Sigma = \bigcup_{v=1}^N \Sigma_v$, where  each $v\in V(G)$  represents a
 smooth component $\Sigma_{v}$ of $\Sigma$.
\item [{\bf [A-2]}]for the $i$-th tail attached to the vertex $v$
there exists the $i$-th marked  point $y_i\in \Sigma_{v}$, $m_v$
is equal to the number of the marked points on $\Sigma_v$,
\item [{\bf [A-3]}]for the $i$-th half edge attached to
the vertex $v$ there exists the $i$-th puncture point
$p_i\in \Sigma_{v}$, $\mu_v$ is equal to the number of the
puncture points on $\Sigma_v$,
\item [{\bf [A-4]}] if there is an edge connected the vertices $v$ and $w$,
then there exists a node between $\Sigma_{v}$ and $\Sigma_{w},$  the number  of edges between $v$ and $w$ is equal to the number of node points between $\Sigma_{v}$ and $\Sigma_{w}$;
\item[{\bf[A-5]}] the restriction of $u$ to each component $\Sigma_v$ is $J$-holomorphic.
\item[{\bf[A-6]}]  $u$ converges exponentially to
$(k_{1},\cdots, k_{\mu})$ periodic orbits
$(x_{k_1},...,x_{k_{\mu}})$ as the variable tends to the puncture
$(p_1,...,p_{\mu})$; more precisely, $u$ satisfies \eqref{exponential_decay_a}-\eqref{exponential_decay_y};
\end{itemize}
\v
Similarly, we can define $(g,m+\mu^{-}+\mu^{+},A,\mathbf{k}^-,\mathbf{k}^+)$-weighted
dual graph and associate it with a holomorphic cascade in
${\mathbb{R}}\times \widetilde{M}$. For every holomorphic cascade in ${\mathbb{R}}\times \widetilde{M}$ we have a $\mathbb{C}^*$ action, so we take the quotient.

\begin{remark}
For $(g,m+\mu^-+\mu^+,A,\mathbf{k}^-, \mathbf{k}^+)$-weighted dual graph,
there are weights
$$\mathfrak{l}^{\bf k^-}:
(k^-_{1}e^-_{1},\cdots,k^-_{\mu^-}e_{\mu^-})\to V(G)\;\;
and\;\;
\mathfrak{l}^{\bf k^+}: (k^+_{1}e^+_{1},\cdots,k^+_{\mu^+}e_{\mu^+})\to
V(G),$$ where
$$(k^-_{1},\cdots,k^-_{\mu^-})\in (\mathbb{Z}^+)^{\mu^-},\;\;\;
(k^+_{1},\cdots,k^+_{\mu^+})\in (\mathbb{Z}^+)^{\mu^+}.$$
The weight $\mathfrak{h}:V(G) \rightarrow H_{2}(\mathbb P(\mathcal N\oplus \mathbb C),\Z)$ and the
weights $\mathfrak{l}^{\bf k^{\pm}}$
satisfy
$$\sum_{i=1}^{\mu^{+}} k_i^{+}-\sum_{i=1}^{\mu^{-}} k_i^{-}=\#(A\cap Z_{\infty})+\#(A\cap Z_{0}).$$
\end{remark}

\subsection{Weighted dual graph with an oriented decomposition}\label{Weighted dual graph with an oriented decomposition}
Let $G$ be a connected graph with $V=\{v_1,...,v_N\}$ and let $\mathfrak{d}$
be a partition  of $\{1,2,...,N\}$,
which induces a decomposition of $V$
$$\mathfrak{d}: V= \mathfrak{A}_0\cup \mathfrak{A}_1\cup\mathfrak{A}_2\cup...\cup\mathfrak{A}_k$$
such that
\begin{description}
\item[(1)] $\mathfrak{A}_0= \bigcup_{\alpha=1} G_{0\alpha}$ with $G_{0\alpha}\cap G_{0\beta}=\phi $ for $\alpha\ne \beta$. Each $G_{0\alpha}$ is a connected subgraph of $G$, and it is a
 $(g_{\alpha},m_{\alpha}+\mu_{\alpha},A_{\alpha}, \mathbf{k}_{\alpha})$-weighted dual graph associated with a holomorphic cascade ${\mathcal{M}}_{G_{0\alpha}}$ in $M^+$. We call each ${\mathcal{M}}_{G_{0\alpha}}$ a cascade of level $0$.
\item[(2)] $\mathfrak{A}_i=\bigcup^{i_c}_{a=1} B_{ia}$, $1\leq i\leq k$, where $B_{ia}=\cup_\beta G_{ia\beta}$ with $G_{ia \beta}\cap G_{ia \gamma}=\phi $ for $\beta\ne \gamma$. For any fixed $a$, $1\leq a\leq i_c$, each $G_{ia \alpha}$ is a connected subgraph of $G$, and it is a $(g_{\alpha},m_{\alpha}+\mu^{-}_{\alpha}+\mu^{+}_{\alpha}, \mathbf{k}_{\alpha}^-,\mathbf{k}_{\alpha}^+)$-weighted dual graph associated with a holomorphic cascade ${\mathcal{M}}_{G_{ia \alpha}}$ in a copy of ${\mathbb{R}}\times \widetilde{M}$, and $B_{ia}$ corresponding to a holomorphic block, denoted by ${\mathcal{M}^*}({\mathbb{R}}\times \widetilde{M}, \cdot)_{ia}$. We call each ${\mathcal{M}}_{G_{ia \alpha}}$ a holomorphic cascade of level $-i$ and call $G_{ia\beta}$ and $G_{ia\alpha}$ lie in the same level. For $a\ne b$, ${\mathcal{M}^*}({\mathbb{R}}\times \widetilde{M}, \cdot)_{ia}$ and ${\mathcal{M}^*}({\mathbb{R}}\times \widetilde{M}, \cdot)_{ib}$ lie in different copy of ${\mathbb{R}}\times \widetilde{M}$.
\item[(3)] For each half edge $e^-$ attached to
a vertex $v$ in some $G_{ia \beta}$ there is a unique half edge $e^+$ attached to
a vertex $ v'$ in $G_{j b \alpha}$ with $j<i$ such that $k^+=k^-$. Then there is an edge $\ell\in R(G)$ connecting
$v'\in G_{j b \alpha}$ and $v\in G_{ia \beta},$ and
$\ell$ is the composition of a half edge $e^+$ and
a half edge $e^-$. There is a natural orientation
$\overrightarrow{\ell}:v'\to v.$
We denote simply by $v'\xrightarrow{\ell} v$.
\v
Let ${\mathcal{M}}_{A}({\mathbb{R}}\times \widetilde{M},g,m+\mu^{+}+\mu^{-},{\bf k}^+, {\bf k}^-,\nu)$ and ${\mathcal{M}}_{A'}({\mathbb{R}}\times \widetilde{M},g',m'+\mu^{+'}+\mu^{-'},{\bf k}^{+'}, {\bf k}^{-'},\nu')$ be the holomorphic cascades corresponding to $G_{ia \beta}$ and $G_{j b \alpha}$ respectively. Suppose that $(\Sigma,u)\in {\mathcal{M}}_{A}({\mathbb{R}}\times \widetilde{M},g,m+\mu^{+}+\mu^{-},{\bf k}^+, {\bf k}^-,\nu)$, $(\Sigma',u')\in {\mathcal{M}}_{A'}({\mathbb{R}}\times \widetilde{M},g',m'+\mu^{+'}+\mu^{-'},{\bf k}^{+'}, {\bf k}^{-'},\nu')$, $\Sigma$ and $\Sigma'$ joint at $q$. Then $u(z)$ and $u'(z')$ converge to the same $k=k^+=k^-$
periodic orbit on $\widetilde{M}$ as the variables tend to $q$.

 \item[(4)] For any  $G_{0\alpha}$ and $G_{ia\beta}$, denote
 $$A_{0\alpha}=\sum_{v\in V(G_{0\alpha})}A_{v}\in H_{2}(\overline{M}^+,\mathbb Z),
 \;\;\;\;\;\;\;\;\; A_{ia\beta}=\sum_{v\in V(G_{ia\beta})}A_{v}\in H_{2}(\mathbb P(\mathcal N\oplus \mathbb C),\mathbb Z).$$
We have
$$ \sum_{e_{j}\in G_{0\alpha}}k_{j}+\sum_{\ell^{+}\in G_{0\alpha}}k_{\ell^{+}}=\#(A_{0\alpha}\cap Z_{\infty}) \mbox{ for any } G_{0\alpha}  $$ and
$$ \sum_{e_{j}\in G_{ia\beta}}k_{j}+\sum_{\ell^{+}\in G_{ia\beta}}k_{\ell^{+}}-\sum_{\ell^{-}\in G_{ia\beta}}k_{\ell^{-}}=\#(A_{ia\beta}\cap Z_{\infty}) +\#(A_{ia\beta}\cap Z_{0}) ,\;\;\;\; \mbox{ for any } G_{ia\beta}.$$ \end{description}
Denote
$$A\circeq\sum_{\alpha}A_{0\alpha}+\sum_{i,a,\beta}A_{ia\beta}.$$

\v\v
Let $g$, $m$ and $\mu$ be nonnegative
integers, $A\in H_{2}(\overline{M}^+,\mathbb Z)$.
Let $(V(G),E(G))$ be a connected graph and $\mathfrak{d}$ be
an oriented decomposition satisfying {\bf (1), (2), (3),(4)}.
We call the graph $G$ a $(g,m+\mu,A,\mathbf{k},\mathfrak{d})$-
weighted dual graph  with an oriented decomposition $\mathfrak{d}$.
Denote the $(g,m+\mu,A,\mathbf{k},\mathfrak{d})$ weighted dual graph
by
$$(V(G),E(G),A,\mathfrak{g}, \mathfrak{m}, \mathfrak{l}^{\bf k},\mathfrak{h},\mathfrak{d})$$ or
denoted simply by $G(\mathfrak{d})$.
\v
Two $(g,m+\mu,A,\mathbf{k},\mathfrak{d})$ weighted dual graphs  $G_1(\mathfrak{d}_1)$ and $G_2(\mathfrak{d}_2)$ are called isomorphic
if there exists a bijection $T$ between their vertices and edges keeping  oriented decomposition and all weights.
Let $S_{\mathfrak{g},\mathfrak{m},\mathfrak{l}^{\bf k},\mathfrak{k},\mathfrak{h},\mathfrak{d}}$
 be the set of isomorphic classes of  $(g,m+\mu,A,\mathbf{k},\mathfrak{d})$ weighted dual graphs.
Given $g$, $m$, $\mu$, $A\in H_{2}(\overline{M}^{+},\mathbb Z)$ and the
weight  $\mathbf{k}=(k_1,...,k_{\mu})$, denote by $S(g,m+\mu,A,\mathfrak{l}^{\bf k})$ the union of all possible $S_{\mathfrak{g},\mathfrak{m},\mathfrak{l}^{\bf k},\mathfrak{h},\mathfrak{d}}$.
\v
For every $(g,m+\mu,A,\mathbf{k},\mathfrak{d})$ weighted dual
graph $G(\mathfrak{d})$ we can associate a space $\mathcal{M}_{G(\mathfrak{d})}$ of  the equivalence class of
stable $J$-holomorphic maps of type $G(\mathfrak{d})$ as in \S\ref{Holomorphic block of type G}.
We call $\mathcal{M}_{G(\mathfrak{d})}$ a holomorphic cascade system.
\v\v

\begin{remark} Suppose that $v_{a}\xrightarrow{\ell} v_{b}$ for some
$v_{a}\in G_A$, $v_{b} \in G_B$, Then $\Sigma_{v_a}$ and $\Sigma_{v_b}$
have a node $q$ and $u\mid_{\Sigma_{v_a}}$, $u\mid_{\Sigma_{v_b}}$
converge to the same $k_{\ell}$
periodic orbit $x_{k_{\ell}}$ on $\widetilde{M}$ as the variables tend to the nodal point $q$.
\end{remark}
\v
Given $g$, $m$, $\mu$, $A\in H_{2}(\overline{M}^{+},\mathbb Z)$, and the weight  $\mathbf{k}=(k_1,...,k_{\mu})$
we define
 $${\overline{\mathcal{M}}}_{A}(M^{+};g,m+\mu,{\bf k},\nu)=\bigcup_{G(\mathfrak{d})\in S(g,m+\mu,A,\mathfrak{l}^{\bf k})} \mathcal M_{G(\mathfrak{d})} .$$

\v\v

Denote by ${\mathcal D}^{J,A}_{g,m+\mu,{\bf k}}$ the number of
all possible $S_{\mathfrak{g},\mathfrak{m},
\mathfrak{l}^{{\bk}},\mathfrak{h},\mathfrak{d}}$.
By Lemma \ref{lem_com_Proj}  and the compactness of $\overline{\mathcal{M}}_{g,m+\mu}$ we get
\begin{lemma} ${\mathcal
D}^{J,A}_{g,m+\mu,{\bf k}}$ is finite.
\end{lemma}

\v

\v\v

\chapter{Compactness Theorems}\label{compact_theorem}

\section{Bubble phenomenon }\label{bubble_phenomenon}

\vskip 0.1in
\noindent

\v
\subsection{Bound of the number of singular points }\label{number of singular points}

 Following McDuff and Salamon \cite{MS} we have the notion of
singular points for a sequence $u^{(i)}$ and the notion of mass of singular
points. Suppose that $(\Sigma\upper; {\bf y}\upper,{\bf p}\upper)$ is
stable for every $i$ and
converges to $(\Sigma; {\bf y}, {\bf p})$ in
${\overline{\mathcal{M}}}_{g,m+\mu }$. We view nodal points,
marked points as puncture points. For each $({\Sigma}\upper;{\bf y}\upper, {\bf p}\upper)$ we choose metric $\mathbf{g}^{\diamond}$. We show that there is a constant $\hbar >0$ such that the mass of every singular
point is large than $\hbar$. Let $q\in \Sigma$ be a singular point and $q^{(i)}\in \Sigma\upper$, $q^{(i)}\rightarrow q$. In case $u^{(i)}(q^{(i)})\in M^+_0$
the argument is standard (see \cite{MS}). We only consider $\widetilde{E}$ over the cylinder end.
Without loss of generality we assume that $q^{(i)}$ is not a
nodal point of $\Sigma\upper$.  In term of  the cylinder coordinates, we have
$D_{1/2}(q^{(i)})\subset \Sigma^{(i)}-\{nodal\; points\},$ where $D_{1/2}(q^{(i)})=\{(s\upper, t\upper)\;|\; (s\upper -s\upper(q\upper))^2+(t\upper -t\upper(q\upper))^2 \leq 1/4\}$. We identify $q\upper$ with $0$ and consider $J$-holomorphic maps $u^{(i)}: D_{1/2}(0)\to N.$
\v
The proof of the following  lemma is similar to Theorem 4.6.1 in \cite{MS}.
\v\n
\begin{lemma}\label{lower_bound_of_singular_points}
 Let $u^{(i)}:D_{1/2}(0)\to \mathbb{R}\times \widetilde{M}$ be a sequence of J-holomorphic maps with finite energy such that
  $$\sup_{i} {E}_{\phi}(u\upper)<\infty,\;\;\;\;\; |du^{(i)}(0)|\longrightarrow \infty ,\;\; as \;\;i\to \infty.$$
Then there is a constant $\hbar >0$ independent of $u^{(i)}$ such that,  for every $\epsilon>0$
\begin{equation}
 \liminf\limits_{i\to\infty}\widetilde{E}(u^{(i)};D_{\epsilon}(0)) \geq \hbar.
\end{equation}
\end{lemma}
\noindent

\v
By Lemma \ref{lower_bound_of_singular_points} we conclude that the   singular points are isolated and the limit
$$m_{\epsilon}(q)=\lim_{i \rightarrow \infty} \widetilde{E}(u^{(i)};D_{q^{(i)}}
(\epsilon,h^{(i)}))$$ exists for every sufficiently small $\epsilon>0$.
The mass of the singular point q is defined to be
$$m(q)=\lim_{\epsilon\rightarrow 0}{m_{\epsilon}(q)}.$$
\vskip 0.1in

\v
Denote by $P\subset \Sigma $ the set of  singular points for $u^{(i)}$,
the nodal points and the puncture points. By Lemma \ref{lower_bound_of_singular_points} and \eqref{energy_bound}, $P$ is
a finite set. By definition, $|du^{(i)}|$ is uniformly bounded on every
compact subset of $\Sigma - P$.
 We call a translations along $\mathbb R$ a  $\mathcal{T}$- rescaling.
 By a    possible   $\mathcal{T}$- rescalings and passing to a  subsequence we may
assume that $u^{(i)}$ converges uniformly  with all derivatives on every compact
subset of $\Sigma - P $ to a $J$-holomorphic map $u:\Sigma - P
\rightarrow N.$
 Obviously, $u$ is a finite energy $J$-holomorphic map.
\v\n

We need to study the behaviour of the sequence $u\upper$ near each
singular point for $u\upper$. Let $q \in \Sigma $ be a   singular point for
$u\upper$. We have three cases.
\v\n
{\bf (a)} $q\in \Sigma -\{nodal\; points,\;puncture\; points,\; marked\; points\}$. We consider $J$-holomorphic maps $u^{(i)}: D_{1}(0)\to N$.
\begin{enumerate}
\item[(a-1)] there are $\epsilon > 0$ and a compact set $K\subset N$ such that $u^{(i)}(D_{\epsilon}(q))\subset
K$.
\item[(a-2)] $q$ is a nonremovable singularity.
\end{enumerate}
\v\n
{\bf (b)} $q\in \{nodal\; points, puncture\; points\}$. We discuss only the nodal points, the discussions for puncture points are similar. A neighborhood of  a nodal pint $q$ is two discs $D_{1}(0)$ joint at $0$, where $D_1(0)=\{|z|^2\leq 1\}$.
\begin{enumerate}
\item[(b-1)] there is a compact set $K\subset N$ such that $u^{(i)}(q^{(i)})\subset
K$.
\item[(b-2)] $q$ is a nonremovable singularity.
\end{enumerate}
\v\n
{\bf (c)} $q\in \{marked\; points\}$. A neighborhood of $q$ is $D_{1}(0)$ with $q=0$. We consider $J$-holomorphic maps $u^{(i)}: D_{1}(0)-\{0\}\to N$.
\v
For (a-1), (b-1), {\bf (c)} we  construct bubbles as usual for a compact symplectic manifold (see \cite{RT}, \cite{PW},
\cite{MS}). We call this type of bubbles (resp. bubble tree) the normal bubbles (resp. normal bubble tree).
\v\n
 \subsection{Construction of the bubble tree for (a-2)}\label{bubble_tree_a}
\v\n
We use cylindrical coordinates $(s,t)$
and write
$$u\upper(s,t)=(a\upper(s,t), \widetilde{u}\upper(s,t))=(a\upper(s,t), \theta\upper(s,t), {\bf w}\upper(s,t))$$
$$u(s,t)=(a(s,t),\widetilde{u}(s,t))=(a(s,t), \theta(s,t), {\bf w}(s,t)).$$

\v
\v
Note that the gradient $|du^{(i)}|$ depends not only on the metric $<, >$ on $N$ but also depends on the metric on $\Sigma^{(i)}$.  The energy don't depend on the metric on $\Sigma^{(i)}$. To construct bubble tree in present case it is more
convenient to take the family of metrics $\mathbf{g}$ on each $({\Sigma}\upper;{\bf y}\upper, {\bf p}\upper)$ in a neighborhood of $(\Sigma;{\bf y}, {\bf p})$ in ${\overline{\mathcal{M}}}_{g,m+\mu}$.
\v

By Theorem \ref{exponential_estimates_theorem}, we have
$$\lim_{s \rightarrow \infty} \widetilde{u}(s,t)=x(kt)$$
in $C^{\infty}(S^1)$, where $x(\;,\; )$ is a $k$-periodic orbit on $\widetilde{M}$.
Choosing $\epsilon $ small enough we have $$|{m_{\epsilon}(q)}-m(q)|\leq \frac{1}{10} \hbar.$$
For every $i$ there exists $\delta_i>0$ such that
\begin{equation}\label{choose_delta}
\widetilde{E} (u\upper; D_{\delta_i}(0))  =m(q) - \frac{1}{2} \hbar.
 \end{equation}
Since $u^{(i)}$ converges uniformly with all derivatives to $u$ on any compact set
of $D_{\epsilon}(0)-\{0\}$, $\delta_i $ must converge to 0.
Put
\begin{equation}\label{eqn_cpt_i}
\hat{s}\upper= s + log \delta_i,\;\;\;\hat{t}\upper= t,\end{equation}
\begin{equation}\label{a_rescaling_DT}
\hat{a}\upper= a +k\log\delta_i,\;\;\;\hat{\theta}\upper=\hat{\theta}=\theta.
\end{equation}
 Define the $J$-holomorphic curve $v\upper(\hat{s}\upper,t)$ by
\begin{equation}\label{bubble construction}
v\upper(\hat{s}\upper,t)=(\hat{a}\upper(\hat{s}\upper,t),\widetilde{v}\upper(\hat{s}\upper,t))=
\end{equation}
$$\left(a\upper(-\log\delta_{i}+ \hat{s}\upper,t) +k\log\delta_i ,
\widetilde{u}\upper(-\log\delta_{i}+ \hat{s}\upper,t)\right).$$

\vskip 0.1in
\noindent

\begin{lemma}\label{bubble_tree_a-1}
Suppose that $0$ is a nonremovable singular point of $u$.
Define the $J$-holomorphic map $v\upper$ as above. Then there exists a
subsequence (still denoted by $v\upper$) such that
\begin{description}
\item[(1)] The set of singular points $\{Q_1,\cdot \cdot \cdot,Q_d\}$ for
$v\upper$ is
finite and tame, and is contained in the disc $D_{1}(0)=\{z \mid \mid z
\mid \leq 1\};$
\item[(2)] The subsequence $v^{(i)}$ converges with all derivatives uniformly on
every compact subset of  ${\mathbb C}\backslash\{Q_1,\cdot \cdot \cdot,Q_d\}$
to a nonconstant
J-holomorphic map $v:{\mathbb C}\backslash\{Q_1,\cdot \cdot \cdot,Q_d\}
\rightarrow {\mathbb R}\times \widetilde{M};$
\item[(3)] $\widetilde{E}(v)+\sum\limits_{1}^{d}m(Q_i)=m(0). $
\item[(4)] $\lim\limits_{s \rightarrow \infty} \widetilde{u}(s,t)=
\lim\limits_{\hat{s} \rightarrow -\infty} \widetilde{v}(\hat{s},t).$
Furthermore, we choose  the Dauboux coordinates $(\theta,{\bf w})$ near $x$ on $\widetilde{M}$ and write
$$u(s,t)=(a(s,t),\theta(s,t),{\bf w}(s,t)),\;\;\;v(\hat s, \hat t)=(\hat a(\hat s, \hat t),\hat \theta(\hat s, \hat t),\hat{\bf w}(\hat s, \hat t)).$$ Then there are constants $0<\mathfrak{c}<\frac{1}{2}$, $C_{\bf n}>0$
 such that for all ${\bf n}=(m_{1},m_{2})\in \mathbb Z_{\geq 0}^2 $
\begin{eqnarray}
\label{uv_decay_a}
&|\partial^{\bf n}[a(s,t)-ks-\ell_0]|\leq C_{\bf n} e^{-\mathfrak{c}|s|},&
 |\partial^{\bf n}[\theta(s,t)-kt-\theta_0]|\leq C_{\bf n} e^{-\mathfrak{c}|s|}; \\
\label{uv_decay_theta}
&|\partial^{\bf n}[\hat a(\hat s,\hat t)-k\hat s-\hat \ell_0]|\leq C_{\bf n} e^{-\mathfrak{c}|\hat s|},&
|\partial^{\bf n}[\hat \theta(\hat s,\hat t)-k\hat t-\hat \theta_0]|\leq C_{\bf n} e^{-\mathfrak{c}|\hat s|},
\end{eqnarray}
where $\hat \ell_{0}$, $\ell_{0}$, $\theta_{0}$, and $\hat \theta_{0}$ are constants ( see Theorem \ref{exponential_estimates_theorem} ) and
\begin{equation}\label{connect condition}
\hat \ell_{0}=\ell_{0},\;\;\;\theta_{0}=\hat \theta_{0}.\end{equation}
\item[(5)]  $\widetilde{E}(v)>\frac{1}{4} \hbar$;
\end{description}
\end{lemma}
\vskip 0.1in
\noindent
{\bf Proof: } The proofs of
{\bf (1)}, {\bf (2)} and {\bf (3)} are standard (see \cite{MS}), we omit them here. We only prove {\bf (4)} and {\bf (5)}.
 \vskip
0.1in
{\bf (4) } Consider the $\widetilde{J}$-holomorphic map
$$\hat{u}\upper=\pi \circ \widetilde{u}\upper:\Sigma \rightarrow Z .$$
Write $A(r,R)=D_{R}(0)-D_{r}(0)$. Since $E(\hat{u}\upper, A(R\delta_i, \epsilon ))=\widetilde{E}(u\upper; A(R\delta_i ,\epsilon )),$ we have $$E(\hat{u}\upper, A(R\delta_i, \epsilon )) \leq
\frac{2}{3}\hbar$$ as $i$ big enough.  By Lemma \ref{tube_exponential_decay}, there exists a $T_0>0$ such that for
$T>T_0$
\begin{equation}\label{bubble_energy_1}
E(\hat{u}\upper; A(R\delta_i e^T, \epsilon e^{- T})) \leq
{Ce^{-2\fc T}}E(\hat{u}\upper; A(R\delta_i, \epsilon))
\end{equation}
and
\begin{equation}\label{bubble_energy_2}
\int\limits_{S^1}d(\hat{u}\upper(\epsilon e^{- T+it}),\hat{u}\upper(R\delta_ie^{T+it}))dt\leq Ce^{- \fc T}\sqrt{E(\hat{u}\upper;A(R\delta_i,\epsilon))}.
\end{equation}
We choose $T$ large such that $Ce^{-2\fc T}<\frac{1}{2}$. It follows from \eqref{bubble_energy_1} that
\begin{equation}\label{bubble_energy_3}
E(\hat u\upper; A(R\delta_i ,\epsilon )) \leq
\frac {1}{1 - Ce^{-2\fc T}}\left(E(\hat u^{(i)}; A(\epsilon e^{- T}, \epsilon))
+ E(\hat v\upper;A(R , Re^T))\right).
\end{equation}
Since $u\upper \rightarrow u $ and $v\upper \rightarrow v $ uniformly on any compact sets, by \eqref{bubble_energy_3} we have
\begin{equation}\label{bubble_energy_4}
\lim_{\epsilon \rightarrow 0, R\rightarrow \infty}\lim_{i\rightarrow
\infty}E(\hat u\upper; A( R\delta_i,\epsilon))=0.
\end{equation}
Suppose that $v$ converges to a $k^{\prime}$-periodic orbit $x^{\prime}$.  From \eqref{bubble_energy_2} and \eqref{bubble_energy_4}, we obtain
$$\int\limits_{S^1}d(\pi\circ x(kt),\pi\circ x^{\prime}(k^{\prime}t))dt=0.$$
Therefore  $x=x'$.
\v
Since $\widetilde{E}(u\upper, A(R\delta_i, \epsilon )) \leq
\frac{2}{3}\hbar$, by Lemma \ref{intexp_L} we have
$$
\theta\upper(s_{1},t+1)-\theta\upper(s_{1},t)=\theta\upper(s_{2},t+1)-\theta\upper(s_{2},t),\;\;\;\;\forall\;  -\log \epsilon+B\leq s_{1},s_{2}\leq -\log (R\delta_{i})-B.
$$
Since $u$ (resp. $v$) converges to a $k$ (resp. $k^{\prime}$)-periodic orbit,
there exist a constant $R_{0}>0$ such that
$$
\theta(s,t+1)-\theta(s,t)=k,\;\;\;\hat{\theta}(\hat s,\hat t+1)-\hat{\theta}( \hat s,\hat t)=k',\;\;\;\;\forall\; |s|,|\hat s|>R_{0}.
$$
Then by the locally uniformly convergence of $u\upper$ and $v\upper $,
and  \eqref{bubble construction} we have $k=k'.$

\v
Let $N>0$ be a constant  such that
\begin{equation}\label{bound_energy_rubber1}
\lim_{i\rightarrow \infty} \sup \widetilde {E}(u\upper; N \leq s\leq -\log\delta_{i}-N ) \leq
\frac{1}{2} \hbar,
\end{equation}
Then \eqref{uv_decay_a} and \eqref{uv_decay_theta} follows from Theorem \ref{exponential_estimates_theorem} and Theorem \ref{tube_L_decay}.
Denote $\pounds\upper=(a\upper-ks,\theta\upper -kt)$.
By Theorem \ref{tube_L_decay} and by integrating we have for any fixed $B>N$ and for any fix $t^{*}$,
\begin{equation}\label{pounds_z_1}
|\pounds\upper(B,t^{*}) -\pounds\upper(-B-\log\delta_{i},t^{*})|\leq  {\mathcal C_{1}}{\fc}^{-1}e^{-\mathfrak{c} (B-N)}.
\end{equation}

Note that
$$a\upper(-B-\log\delta_{i},t )-k(-B-\log\delta_{i})=  \hat a\upper (-B,t)-k(-B)$$
and
 $$\theta\upper(s,t)-kt=\hat\theta\upper(\hat s,\hat t) -k\hat t.$$
 Then \eqref{pounds_z_1} can be re-written as
$$\left|a\upper(B,t^{*})-kB -\left[\hat a\upper(- B,t^*)-k(-B)\right]\right|\leq  Ce^{-\mathfrak{c} (B-N)},$$
$$\;\;\left|\theta \upper(B,t^*)-kt^{*}-\left[\hat \theta\upper(-B,t^*)-kt^{*}\right]\right|\leq  Ce^{-\mathfrak{c} (B-N)}.$$    Letting $i\rightarrow \infty$, $B\rightarrow \infty$,  \eqref{connect condition} follows.
\v
 {\bf(5)}  For any fixed $\epsilon>0,$ we have $\widetilde{E} (u\upper; A(\delta_{i},\epsilon))\geq \tfrac{ \hbar}{3}$ as $i$ big enough. By the same argument of  \eqref{bubble_energy_3}, we have
 \begin{equation}
 \widetilde{E}(u^{(i)}; A(\epsilon e^{- T}, \epsilon))
+ \widetilde{E}(v\upper;A(1  ,  e^T))   \geq  (1 - Ce^{-2\fc T})\tfrac{ \hbar}{3}.
\end{equation}
Since $u\upper \rightarrow u $ and $v\upper \rightarrow v $ uniformly on compact sets, let $i\rightarrow \infty$ we have
  \begin{equation}
  \widetilde{E}(u,  -\log \epsilon  \leq s\leq -\log \epsilon+T)+\widetilde{E}(v,-T\leq \hat s\leq 0)\geq (1 - Ce^{-2\fc T})\tfrac{ \hbar}{3}.
  \end{equation}
  Choose $\epsilon$ small enough such that $ \widetilde{E}(u,s\geq -\log \epsilon)\leq  \tfrac{\hbar}{12}.$ Letting $T\rightarrow \infty,$ then (5) follows. $\Box$
\v

We can repeat this again to construct bubble tree.

\v\n
\begin{remark}\label{bubble_tree_a-2}
Note that the coordinates $(s,t)$, $(\hat s, \hat t)$ in ${\bf (4)}$ maybe not the cusp cylinder coordinates. If we choose the cusp cylinder coordinates then \eqref{connect condition} does not hold in general.
\end{remark}

In the case (b-2), we use the same method of as the case (a-2) to construct a bubble $S^2$ with
  $\widetilde{E}(v)\mid_{S^2}> \frac{1}{3} \hbar$, inserted between $\Sigma_1$ and $\Sigma_2$. The same results as Lemma \ref{bubble_tree_a-1} still hold. We can repeat this again to construct bubble tree.
\vskip 0.1in
\noindent

\v

\section{Procedure of re-scaling}\label{procedure_of_re-scaling}
\v
Let $\Gamma\upper=(\Sigma\upper;{\bf y}\upper,{\bf p}\upper, u\upper) \in \overline{\mathcal{M}}_{A}(M^{+},g,m+\mu,{\bf y},{\bf p},{\bk})$ be a sequence
with $\Sigma\upper = \bigcup_{v=1}^N \Sigma\upper_v$. Assume
that there is one component $\Sigma\upper_v$ that has
genus $0$ and is unstable. Let $\Sigma_1\upper$ be such a component.
We identify $\Sigma_1\upper$ with a sphere $S^2$, and consider
$u\upper: S^2\rightarrow N$. We discuss several cases:
\vskip 0.1in
\noindent
{\bf 1).} $u\upper\mid_{S^2}$ has no singular point. Then $\|\nabla u\upper\mid_{S^2}\|$ are uniformly bounded above. As
$(u\upper,S^2;{\bf y}\upper,{\bf p}\upper)$ is stable,
$E(u\upper)\mid_{S^2}\geq \hbar$ or
$\widetilde E(u\upper)\mid_{S^2}\geq \hbar$. Then
$u\upper\mid_{S^2}$ locally uniformly converges to $u: S^2\to N$
with $E(u)\mid_{S^2}\geq \hbar$ or
$\widetilde E(u)\mid_{S^2}\geq \hbar$, so $(u; S^2)$ is stable.
\v\n
{\bf 2).} There are some singular points. For the cases
(a-1), (b-1), {\bf (c)} in \S\ref{number of singular points}
we construct bubbles as usual for a compact symplectic manifold.
For the cases (a-2), (b-2) we construct bubbles as in
\S\ref{bubble_tree_a}. By $({\bf 3})$ and $({\bf 5})$
of Lemma \ref{bubble_tree_a-1} we get a stable map $(v,S^{2})$.
It is possible that $u\mid_{S^2}$ is a point in $N$ and
$\#\{marked \;points, nodal \;points, singlar \;points\}<3$. In this case
we contract $S^2$ as a point.
\v
We can repeat the procedure to construct bubble tree.
\v
In the following we
assume that
$(\Sigma\upper;{\bf y}\upper,{\bf p}\upper)$ is stable and converges to $(\Sigma;{\bf y}, {\bf p})$ in
${\overline{\mathcal{M}}}_{g,m+\mu}$. For simplicity, we consider the case as in Figure 3.1.
The other cases are similar.

\begin{figure}[!htbp] \label{Figure2}
\centering\includegraphics[height=6.5cm]{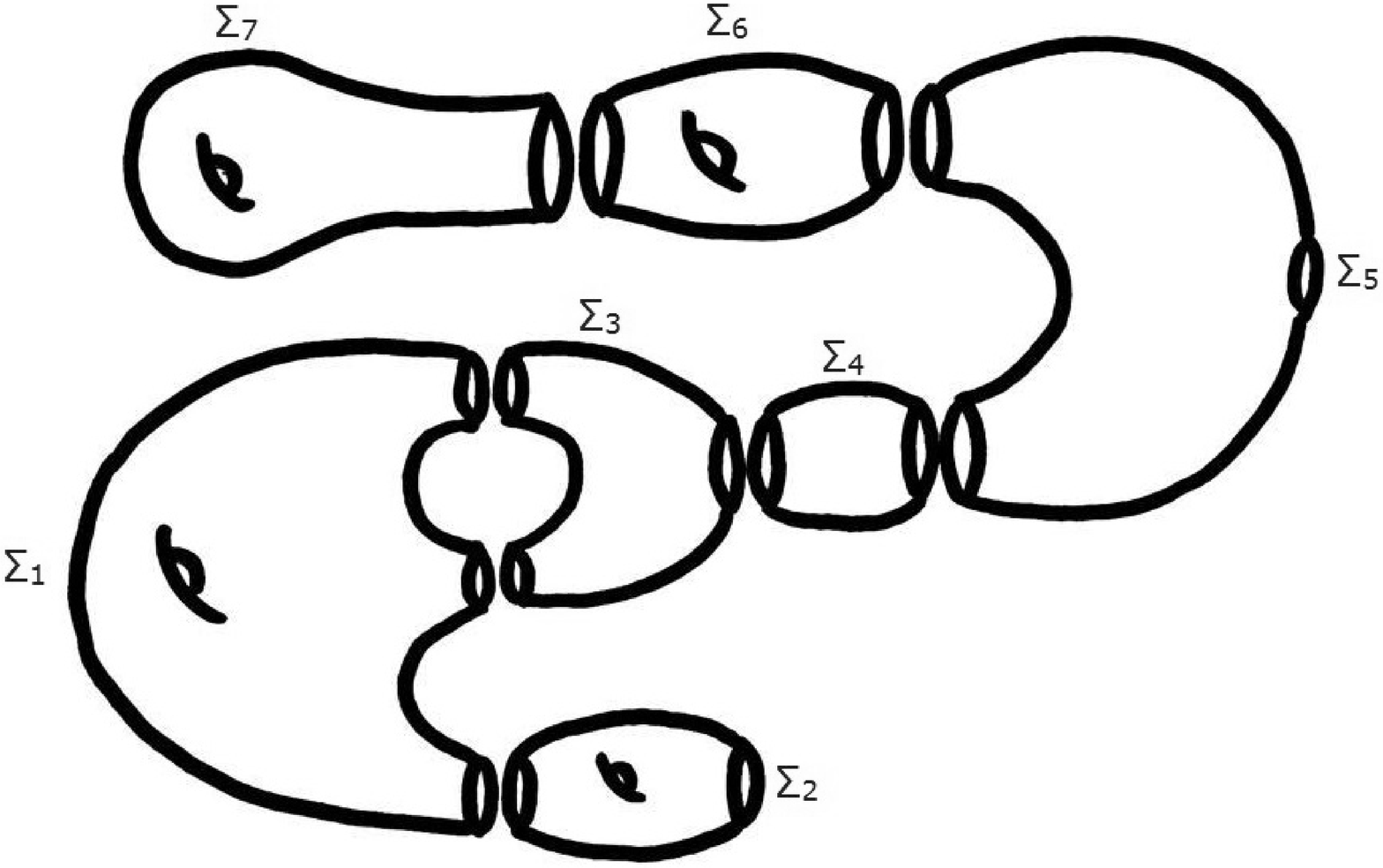}
\caption{}
\end{figure}

\v
\begin{enumerate}
\item[(\bf 1)] By Lemma \ref{lower_bound_of_singular_points} the number of singular points of $\Sigma$ is finite. Denote by $P\subset \Sigma $ the set of  singular points for $u\upper$,
the nodal points and the puncture points.
\item[(\bf 2)] We first find a component $\Sigma_k$ of $\Sigma$, for example $\Sigma_1$ in Figure 3.1, such that
$$\left|a_{1}\upper(z_{1}\upper)\right|\leq \min_{j\in I\setminus\{1\}}\left|a_{j}\upper(z_{j}\upper)\right|,\;\;\forall i.$$
Here $z_{1}\in K_1\subset \Sigma_1\setminus P$, $z_{j}\in K_j\subset \Sigma_j\setminus P$, $K_j$ being some compact set, and we identify $K_i$ with a compact set in $\Sigma_i$, $i\in I$.
We assume that $\sup\limits_{i}\left|a_{1}\upper(z_{1}\upper)\right|<\infty$, that is, $u\upper(z\upper_1)\subset M^+$.
  Find a set $J\subset I$ such that $j\in J$ if and only if
  $u(\Sigma_j)\subset M^+$, for example $ \Sigma_1$ and $\Sigma_7$ in Figure 3.1, i.e., $J=\{1,7\}.$
  Let $\Sigma^{1\bigstar}=\Sigma - \bigcup_{j\in J} \Sigma_j$. $\Sigma^{1\bigstar}$ may have several connected components. For example $\Sigma^{1\bigstar}$ in Figure 3.1 has two connected components: $\Sigma_2$ and $\Sigma_3\cup \Sigma_4 \cup \Sigma_5 \cup \Sigma_6$.
\item[(\bf 3)] For every connected component of $\Sigma^{1\bigstar}$ we repeat the discussion in {(\bf 2)}. For example $u(\Sigma_2)\subset one\;copy \;of\; {\mathbb{R}}\times \widetilde{M}$,  $u(\Sigma_3), u(\Sigma_6)\subset another \;copy \;of\; {\mathbb{R}}\times \widetilde{M}$. We repeat this procedure. We will stop after finite steps.
\item[(\bf 4)] Then we construct bubble tree for every singular point independently to get $\Sigma^{\prime}$, where
${\Sigma}^{\prime}$ is obtained by joining chains of ${\bf P^1}s$ at some double points of $\Sigma$ to separate the
two components, and then attaching some trees of ${\bf P^1}$'s.
For example in Figure 3.1 we have

\end{enumerate}
\v

For every sequence $\Gamma\upper=(u\upper,\Sigma\upper;{\bf y}\upper,{\bf p}\upper)
\in \overline{\mathcal{M}}_{A}(M^{+},g,m+\mu,{\bk})$, using our
procedure we get an element $\Gamma=(u, \Sigma',{\bf y},{\bf p})$
of $\mathcal{M}_{G(\mathfrak{d})}$ for some weighted dual graph with an
oriented decomposition $G(\mathfrak{d})$.

\v
We obtain \vskip 0.1in
\noindent
\begin{theorem}\label{compact_moduli_space}
 ${\overline{\mathcal{M}}}_{A}(M^{+};g,m+\mu,{\bk},\nu)$ is compact.
\end{theorem}
\v\n

\chapter{Local regulization for each holomorphic cascade
}\label{regulization}

The local regulization for each holomorphic cascade
 is very similar to the local regulization for the moduli space of stable holomorphic maps in a closed $C^{\infty}$ symplectic manifold ( see \cite{LS-2}).

\section{Local regularization for ${\mathcal{M}}_{A}(M^+,g,m+\mu,{\bf k},\nu) $
}\label{local regularization-1}
\subsection{Local regularization-Top strata
}\label{local regularization-top}

It is well-known that if $D_{u}$ is surjective for any
$b = (j,{\bf y},{\bf p}, u) \in \mathcal{M}_A(M^+,g,m+\mu)$, then $\mathcal{M}_A(M^+,g,m+\mu)$ is a smooth manifold. When the transversality fails we need to take the regularization.
Suppose that $m+\mu + 2g \geq 3.$
\vskip 0.1in
\noindent

Let   $[b_o]=[(p_o,u)]\in \mathcal{M}_A(M^+,g,m+\mu)$ and let $\gamma_o\in \mathbf{T}_{g,m+\mu}$ such that $\pi(\gamma_o)=p_o$, where $\pi: \mathbf{T}_{g,m+\mu}\to \mathcal{M}_{g,m+\mu}$ is the projection.
We choose a local coordinate system $(\psi,\Psi)$ on $U$ with $\psi(\gamma_o)=a_o$ for $\mathcal{Q}$.
We view $a=(j,{\bf y,p})$ as a family of parameters defined on a fixed $\Sigma$. Denote
$$\widetilde{\mathcal B}(a)=\left\{u\in W^{k,2,\alpha}(\Sigma,M^{+})|\;u_{*}([\Sigma])=A\right\}.$$
Let $\widetilde{\mathcal E}(a)$ be the infinite dimensional Banach bundle over $\widetilde{\mathcal B}(a)$ whose fiber at $v$ is $$W^{k-1,2,\alpha}(\Sigma,v^{*}TM^{+}\otimes\wedge_{j_a}^{0,1}T^{*} \Sigma),$$
where we denote by $j_a$ the complex structure on $\Sigma$ associated with $a=(j,{\bf y,p})$. We will denote $j_{a_o}:=j_o$. We have a continuous family of Fredholm system
$$\left(\widetilde{\mathcal B}(a),\; \widetilde{\mathcal{E}}(a),\;
\bar{\p}_{j_a,J}\right)$$
parameterized by $a\in \mathbf{A}$ with $d_{\mathbf{A}}(a_o,a)<\delta$. For any $v\in \widetilde{\mathcal B}(a)$ let $b=(a,v)$ and denote $\widetilde{\E}(a)|_{v}:=\widetilde{\E}|_{b}.$ Let $b_o=(a_o,u)$. Choose $\widetilde{K}_{b_o} \subset \widetilde{\E}|_{b_o}$ to be a finite dimensional subspace such that every member of $\widetilde{K}_{b_o}$ is in $C^{\infty}(\Sigma, u^{*}TM^{+}\otimes\wedge^{0,1}_{j_o}T^{*}\Sigma)$ and
\begin{equation}\label{regularization_operator}
\widetilde{K}_{b_o} + image D_{b_o} = \widetilde{\E}|_{b_o},
\end{equation}
where $D_{b_o}=D\bar{\p}_{j_o,J}$ is the vertical differential of $\bar{\p}_{j_o,J}$ at $u$.

 Let $G_{b_o}$ be the isotropy group at $b_o$.
In case the isotropy group $G_{b_{o}}$ is non-trivial, we must construct a $G_{b_{o}}$-equivariant regularization.
Note that $G_{b_{o}}$ acts on $W^{k-1,2,\alpha}(\Sigma, u^{*}TM^{+}\otimes\wedge^{0,1}_{j_o}T^{*}\Sigma)$ in a natural way: for any $\kappa\in W^{k-1,2,\alpha}(\Sigma, u^{*}TM^{+}\otimes\wedge^{0,1}_{j_o}T^{*}\Sigma)$ and any $g\in G_{b_{o}}$
$$g\cdot \kappa=\kappa\circ dg\in W^{k-1,2,\alpha}\left(\Sigma, u^{*}TM^{+}\otimes\wedge^{0,1}_{j_o}T^{*}\Sigma\right).$$
Set
$$\bar{K}_{b_{o}}=\bigoplus_{g\in G_{b_{o}}}
g \widetilde{K}_{b_{o}}.$$
Then $\bar{K}_{b_{o}}$ is $G_{b_{o}}$-invariant. To simplify notations we assume that $\widetilde{K}_{b_o}$ is already $G_{b_{o}}$-invariant.
As in \cite{LS-2} one can prove that there are constants $\delta>0$, $\rho>0$ depending on $b_o$ such that there is an isomorphism
$$P_{b_{o}, b}:\widetilde{\E}_{b_o}\rightarrow \widetilde{\E}_{b}\;\;\;\forall \;b\in \widetilde{\mathbf{O}}_{b_{o}}(\delta,\rho).$$

 Now we define a thickned Fredholm system $(\widetilde{K}_{b_o}\times \widetilde{\mathbf{O}}_{b_{o}}(\delta,\rho), \widetilde{K}_{b_o}\times \widetilde{\E}|_{\widetilde{\mathbf{O}}_{b_{o}}(\delta,\rho)}, S)$.
Let $(\kappa, b)\in \widetilde{K}_{b_o}\times \widetilde{\mathbf{O}}_{b_{o}}(\delta,\rho)$, $b=(a,v)\in \widetilde{\mathbf{O}}_{b_{o}}(\delta,\rho)$. Define
\begin{equation}\label{local regu}
\mathcal{S}(\kappa,b) = \bar{\partial}_{j_a,J}v + P_{b_o,b}\kappa.
\end{equation}
We can choose $(\delta, \rho)$ small such that
the linearized operator $D\mathcal{S}_{(\kappa,b)}$ is surjective for any $b\in \widetilde{\mathbf{O}}_{b_{o}}(\delta,\rho)$.
\v

If we fix the complex structure $j_o$ and $\mathbf{y,p}$, then $W^{k,2,\alpha}(\Sigma;u^{\ast}TM^{+})$ is a Hilbert space. It is well-known that $\|h\|_{j_o,k,2,\alpha}^{2}$ is a smooth function ( see \cite{MVV}). Now the $\|h\|_{j_a,k,2,\alpha}^{2}$ is a family of norms, so the following lemma is important. The proof can be found in \cite{LS-2}.
\begin{lemma}\label{smoothness of norms}
For any $[b_o]=[(p_o,u)]\in  \mathcal{M}_A(M^+,g,m+\mu)$ and any local coordinates $(\psi,\Psi)$ on $U$ with $\psi: U\rightarrow \mathbf{A}\ni a_o$
the norm
$\|h\|_{j_a,k,2,\alpha}^{2}$ is a smooth function in $\widetilde{\mathbf{O}}_{b_{o}}(\delta,\rho)$.
\end{lemma}

  As in \cite{LS-2} we have
\begin{lemma}\label{lem_orbi_T} There exist two constants $\delta_{o},\rho_{o}>0$ depend only on $b_{o}$ such that for any $\delta<\delta_{o},\rho<\rho_{o}$ the following hold.
 \begin{itemize}
\item[(1)] For any $p\in \widetilde{\mathbf{O}}_{b_{o}}(\delta,\rho)$, let $G_{p}$ be the isotropy group at $p$, then $im(G_{p})$ is a subgroup of $G_{b_o}$.
\item[(2)] Let $p\in \widetilde{\mathbf{O}}_{b_{o}}(\delta,\rho)$
be an arbitrary point with isotropy group $G_p$, then there is a $G_p$-invariant neighborhood $O(p)\subset \widetilde{\mathbf{O}}_{b_{o}}(\delta,\rho)$ such that for any $q\in O(p)$, $im(G_{q})$ is a subgroup of $G_{p}$, where $G_p$, $G_q$ denotes the isotropy groups at $p$ and $q$ respectively.
\end{itemize}
\end{lemma}

\subsection{Local regularization for lower stratum : without bubble tree}\label{without bubble tree}

Let $\Sigma$, $\stackrel{\circ}\Sigma$, $\Sigma_i$ be as in \S\ref{norms for smooth surfaces}. We choose local plumbing coordinates $(\fs,\ft)$ and construct $\Sigma_{\fs,\ft}\to \Delta_{\fs}\times \Delta_{\ft}$.
Consider the family of Bananch manifold
$$
 \widetilde{\mathcal B}(\fs,\ft)=\{u\in   W^{k,2,\alpha}( \Sigma_{\fs,\ft},M^{+})|\;u_{*}([\Sigma])=A\}.
 $$
Let $\widetilde{\mathcal E}(\fs,\ft)$ be the infinite dimensional Banach bundle over $\widetilde{\mathcal B}(\fs,\ft)$ whose fiber at $b=(\fs,\ft,u)$ is
$W^{k-1,2,\alpha}(\Sigma_{\fs,\ft}, u^{*}TM^{+}\otimes\wedge^{0,1}_{j_{\fs,\ft}}T^{*}\Sigma_{\fs,\ft}).$
We have a continuous family of Fredholm system
$$\left(\widetilde{\mathcal{B}}(\fs,\ft),\widetilde{\mathcal E}(\fs,\ft),\bar{\p}_{j_{\fs,\ft},J}\right)$$
parameterized by $(\fs,\ft)\in \Delta_{\fs}\times \Delta_{\ft}$. Let $b_o=(0,0,u)$, $b=(\fs,\ft,v)$. We use the same method as in \S\ref{local regularization-top} to choose $\widetilde{K}_{b_{o}}=\bigoplus_{i=1}^{\iota}\widetilde{K}_{b_{oi}} \subset \widetilde{\E}|_{b_{o}}=\bigoplus_{i=1}^{\iota}\widetilde{\E}_{b_{oi}}$ to be a finite dimensional subspace such that
\begin{itemize}
\item[(1)]  Every member of $\widetilde{K}_{b_{oi}}$ is in
$C^{\infty}\left(\Sigma_{i,0}, u^{*}_{i}TM^{+}\otimes\wedge^{0,1}_{j_{oi}}T^{*} \Sigma_{i, 0}\right)$ and supports in the compact subset $\Sigma_{0,0}(R_0) $ of $\Sigma_{0,0}$.
\item[(2)] $\widetilde{K}_{b_{oi}} + image D_{b_{oi}} = \widetilde{\E}|_{b_{oi}},\;\;\;\forall i=1,2,...,\iota.$
\item[(3)] $\widetilde{K}_{b_{oi}}$ is $G_{b_{oi}}$-invariant.
\end{itemize}
where we denote by $j_{oi}$ the complex structure on $\Sigma_i$ associated with $(0, 0)$, and
\begin{equation}\label{eqn_sig_R}
W(R_{0}):=\cup_{l=1}^{\mathbf{e}} (\{  |z_{l}|< e^{-R_{0}}\}\cup \{ |w_{l}|< e^{-R_{0}}\})\cup \mathbf D(e^{-R_{0}}),\;\;\;\;\;\;\Sigma_{\fs,\ft}(R_{0})=\Sigma_{\fs,\ft}\setminus W(R_{0}) .
\end{equation}
for a constant $R_{0}>1$. We identify each $\Sigma_{\fs,\ft}(R_0)$ with $\Sigma_{0,0}(R_0):=\Sigma(R_0)$ for $|\fs|,|\ft|$ small. Denote by $j_{\fs,\ft}$ the family of complex structure on $\Sigma(R_0)$. Denote $j_o:=j_{0,0}$.
Then when $|H|$ small
$$\Psi_{j_{o},j_{\fs,\ft}}: W^{k-1,2,\alpha}(\Sigma(R_0), u^{*}TM^{+}\otimes \wedge^{0,1}_{j_{o}}T^{*}\Sigma(R_0))\to W^{k-1,2,\alpha}(\Sigma(R_0), u^{*}TM^{+}\otimes \wedge^{0,1}_{j_{\fs,\ft}}T^{*}\Sigma(R_0))$$
is an isomorphism. Let $P_{b_{o}, b}=\Phi\circ\Psi_{j_o,j_{\fs,\ft}}$.
We fix a smooth cutoff
function $\beta_{R_{0}}: {\mathbb{R}}\rightarrow [0,1]$ such that
\begin{equation}\label{def_beta1}
\beta_{R_{0}} (s)=\left\{
\begin{array}{ll}
0 & if\;\; |s| \geq R_0 \\
1 & if\;\; |s| \leq R_0-1.
\end{array}
\right.
\end{equation}
 As in \cite{LS-2} we have
\begin{lemma}\label{isomorphism-1} Let
$\bar{\mathcal E}(\fs,\ft)$ be the infinite dimensional Banach bundle over $\widetilde{\mathcal B}(\fs,\ft)$ whose fiber at $b=(\fs,\ft,u)$ is
$$\bar{\mathcal E}_{(\fs,\ft,u)}:=\{\beta_{R_{0}}(s)\eta\;| \;\eta\in \widetilde{\mathcal E}_{(\fs,\ft,u)}\}.$$
Then there are constants $\delta>0$, $\rho>0$ depending on $b_o$ such that there is an isomorphism
$$P_{b_{o}, b}:\bar{\E}_{b_o}\rightarrow \bar{\E}_{b}\;\;\;\forall \;b\in \widetilde{\mathbf{O}}_ {b_{o}}(\delta,\rho).$$
\end{lemma}

\v
Now we define a thickned Fredholm system $(\widetilde{K}_{b_o}\times \widetilde{\mathbf{O}}_{b_{o}}(\delta,\rho), \widetilde{K}_{b_o}\times \widetilde{\E}|_{\widetilde{\mathbf{O}}_{b_{o}}(\delta,\rho)}, S)$.
Let $(\kappa, b)\in \widetilde{K}_{b_o}\times \widetilde{\mathbf{O}}_{b_{o}}(\delta,\rho)$, $b=(a,v)\in \widetilde{\mathbf{O}}_{b_{o}}(\delta,\rho)$. Define
\begin{equation}\label{local regu}
\mathcal{S}(\kappa,b) = \bar{\partial}_{j,J}u + P_{b_o,b}\kappa.
\end{equation}
We can choose $(\delta, \rho)$ small such that
the linearized operator $D\mathcal{S}_{(\kappa,b)}$ is surjective for any $b\in \widetilde{\mathbf{O}}_{b_{o}}(\delta,\rho)$.

\subsection{Local regularization for lower stratum : with bubble tree}\label{with bubble tree}
\v

{\bf A-G-F procedure.} We introduce the A-G-F procedure.
\v
Consider a strata $\mathcal{M}^{\Gamma}$ of $\overline{\mathcal{M}}_{A}(M^{+},g,m+\mu,\mathbf{k},\nu)$. Let
$b_{o}=[(\Sigma,j,\mathbf y,u)]\in \mathcal{M}^{\Gamma}$.
Then $(\Sigma, j,\mathbf y)$ is a marked nodal Riemann surface. Suppose that $\Sigma$
has a principal part $\Sigma^P$ and some bubble tree $\Sigma^B$ attaching to $\Sigma^{P}$ at $q$. \v
Let $u=(u_{1},u_{2})$ where $u_{1}:\Sigma^{P}\to M^{+}$ and $u_{2}:S^{2}\to M^{+}$ are $J$-holomorphic maps with $u_{1}(q)=u_{2}(q)$.

 We consider the simple case $\Sigma^{B}=(S^2,q)$ with $[u_2(S^2)]\ne 0,$ the general cases are similar. Denote $b_{oo}:=(S^{2},q,u_2)$,
$$
\widetilde O_{b_{oo}}(\rho_{o})=\left\{ v\in \mathcal W^{k,2,\alpha}((S^2,q),u_2^*TM^{+})| \|h\|_{k,2,\alpha} \leq \rho_{o},\mbox{ where } v=\exp_{u_2}(h)\right\}.
$$
$$
O_{b_{oo}}(\rho_{o})=\widetilde O_{b_{oo}}(\rho_{o})/G_{b_{oo}}
$$
where $G_{b_{oo}}=\{\phi\in Diff^{+}(S^2)\mid \phi^{-1}(q)=q, \; u_2\circ\phi=u_2\}$ is the isotropy group at $b_{oo}$.
\v

We can choose a local smooth codimension-two submanifold $Y$ such that $u_2(S^2)$ and $Y$ transversally intersects, and $u_2^{-1}(Y)=\mathbf{x}=(x_1,...,x_{\ell})$ ( see \cite{TF17} and \cite{Par}). We add  these intersection points as marked points to $S^{2}$ such that $S^2$ is stable. Denote the Riemann surface by $(S^{2},q,\mathbf{x})$. We may choose $\rho_{o}$ such that for any $(S^{2},q,v)\in O_{b_{oo}}(\rho_{o})$,  $v(S^2)$ and $Y$ transversally intersects, and  $v^{-1}(Y)$ has $\ell$ points.
Denote
$$
\widetilde O_{\hat b_{oo}}(1+\ell,\rho_{o})=\left\{(S^{2},q,\mathbf{x},v)|v(\mathbf{x})\in Y,\;v\in \widetilde O_{b_{oo}}(\rho_{o})\right\}.
$$
 Note that the additional marked points are unordered, so we consider the space
$$\widetilde O_{\hat b_{oo}}(1\mid\ell,\rho_{o}) £º=\widetilde O_{\hat b_{oo}}(1+\ell,\rho_{o})/Sy(\ell)$$
where $Sy(\ell)$ denotes the symmetric group of order $\ell$.  Denote $\hat b_{oo}:=(S^{2},q\mid \mathbf{x},u_2),$
where the points after $``\mid"$ are unordered.
Denote
$$G_{\hat{b}_{oo}}=\left\{\phi\in Diff^{+}(S^2)\mid \phi^{-1}(q)=q,\; u_2\circ\phi=u_2,\; \phi^{-1}\{x_1,...,x_{\ell}\}=\{x_1,...,x_{\ell}\}\right\}.$$
For any $\phi\in G_{b_{oo}}$, since $u_2\circ\phi=u_2$, we have $\phi^{-1}\{x_1,...,x_{\ell}\}=\{x_1,...,x_{\ell}\}$.
Then the following lemma holds.
\begin{lemma}\label{isotropy group}
$G_{b_{oo}}=G_{\hat b_{oo}}.$
\end{lemma}
\v
Let $\tilde{b}_{oo}:=(S^{2},q,\mathbf{x},u_2)$ be a representive of
$\hat b_{oo}:=(S^{2},q\mid\mathbf{x},u_2),$ where $\mathbf{x} =(x_1,...,x_{\ell})$ is an ordered set. We can construct a metric $\mathbf {g}$ on $(S^{2},q,\mathbf{x})$ as in section \S\ref{metric on surfaces} such that  $\mathbf g^{\diamond}$  is the standard cylinder  metric near marked points and nodal points. We choose cusp coordinates $z$ on $\Sigma^{P}$ and  $w$ on $S^{2}$ near $q$. Put $\Sigma_1=\Sigma^P$, $\Sigma_2=S^2$, $b_o=(b_{o1}, b_{o2})$. Let $G_{b_{oi}}$ be the isotropy group at $b_{oi}.$
Denote $\tilde{b}_o=(\tilde{b}_{o1},\tilde{b}_{o2})$, where $\tilde{b}_{o1}$ is a lift of $b_{o1}$ to the uniformization system, and $\tilde{b}_{o2}:=\tilde{b}_{oo}.$  Note that the cusp coordinates $z$ and $w$ are unique modulo rotations near nodal point $q$ and the metric $\mathbf{g}$ on $\Sigma^{P}$ is  $G_{b_{1}}$-invariant
   and  $\mathbf{g}$ on  $(S^{2},q,\mathbf{x})$ is $G_{b_{2}}$-invariant.  In the coordinates $z,w$   for any  $\phi_{i}\in G_{b_{i}}$,  $$\phi_{1}(z)=e^{-\sqrt{-1} \gamma_{1}}z,\;\;\;\phi_{2}(w)=e^{-\sqrt{-1} \gamma_{2}}w.$$ By the finitness of  $G_{b_{i}}$, we have $\gamma_{i}=\frac{2j_{i}\pi}{l_{i}}$ where $j_{i}<l_{i},j_{i},l_{i}\in \mathbb Z,i=1,2.$
\v
We choose $$\widetilde{K}_{b_{o}}=\bigoplus_{i=1}^{2}\widetilde{K}_{b_{oi}} \subset \widetilde{\E}|_{b_{o}}=\bigoplus_{i=1}^{2}\widetilde{\E}_{b_{oi}}$$ to be a finite dimensional subspace satisfying (1), (2) and (3) in \S\ref{without bubble tree}.
\v
Then we glue $\tilde{b}_{o1}$ and $\tilde{b}_{o2}$ at $q$ with gluing parameters $(r^*, \tau^*)$ in the coordinates $z$, $w$ to get representives of $\hat{p}^*:=(\Sigma_{(r^*)}, \mathbf{y}\mid \mathbf{x})$ and pregluing map $\hat{u}_{(r^*)}$. Let $\hat{b}^*_{o}=(\hat{p}^*, \hat{u}_{(r^*)})$, denote by $G_{\hat{b}^*_{o}}$ the  isotropy group at $\hat{b}^*_{o}$. Now we forget Y and the additional marked points $\mathbf{x}$. We get a element $\Sigma^*:=\Sigma_{(r^*)}$, which is a point $p^*=(\Sigma_{(r^*)},\mathbf{y})\in \overline{\mathcal{M}}_{g,n}$. Let $b^*_o=(p^*, u_{(r^*)})$, denote by $\mathbf{G}_{p^*}$ and $G_{b^*_o}$ the isotropy groups at $p^*$ and $b^*_o$ respectively.
The following lemma is obvious.
\begin{lemma}\label{isotropy group-1}
$G_{\hat{b}^*_{o}}=G_{b^*_{o}}.$
\end{lemma}

\v
We call this procedure a {\bf A-G-F} procedure ( Adding marked points-Gluing-Forgetting Y and marked points). This procedure can be extended to bubble tree and bubble chain in an obvious way.
\v
We use the same method as in \S\ref{without bubble tree} to construct the local regularization.
\v
\section{Local regularization for ${\mathcal{M}}_{A}({\mathbb{R}}\times \widetilde{M},g,m+\mu^{+}+\mu^{-},{\bf k}^+, {\bf k}^-,\nu)/\mathbb{C}^*$}\label{local regularization-2}

Let $[b_o]=[(p_o,u)]\in {\mathcal{M}}_{A}({\mathbb{R}}\times \widetilde{M},g,m+\mu^{+}+\mu^{-},{\bf k}^+, {\bf k}^-,\nu)/\mathbb{C}^*$.
We construct a local slice of the $\mathbb{C}^*$-action around $[b_o]$ as following. Let $p$ be one of $(p_1,...,p_{l})$. For example, $p=p_1$.
We fix a coordinate system $(a, \theta, \mathbf{w}_j)$ on ${\mathbb{R}}\times \widetilde{M}$, where $\mathbf{w}_j$, $j=1,2,...,l$, is a local coordinate on $Z$ near $\pi\circ\tilde{u}(p_j)$. Choose the local cusp cylinder coordinate $(s,t)$ near $p_j$. Suppose that
$$a( {s}, {t})-k_j {s}-\ell_j\rightarrow 0\;\;\;\;
\theta( {s}, {t})-k_j {t}-\theta_{j0}\rightarrow 0.$$
We can choose the coordinates  $(a,\theta)$ such that $\ell_1=0,\;\theta_{10}=0.$
Let
$$\widetilde{\mathbf{O}}^{*}_{b_{o}}(\delta,\rho):=\left\{(a,v)\in \mathbf{A}\times \widetilde{\mathcal{B}} \;|\; d_{\mathbf{A}}(a_o,a)<\delta, \|h+\hat{h}_0\|_{j_a,k,2}<\rho,\; a(h^1_{0})=0,\theta(h^1_{0})=0 \right\},
$$
where $v=\exp_u\{h+\hat{h}_0\}$, $h_0=(h^1_{0},...,h^l_{0})\in \ker H_{\infty}$, $a(h^1_{0}),\theta(h^1_{0})$ denote the components of $h^1_{0}$ with respect to $\{\frac{\p}{\p a}, \frac{\p}{\p \theta}\}$.
Denote by  $
\widetilde{\mathbf{O}}'_{b_{o}}(\delta,\rho)$ the image of $\mathbb C^{*}$-action on $\widetilde{\mathbf{O}}^{*}_{b_{o}}(\delta,\rho).
$ Then
$\widetilde{\mathbf{O}}^{*}_{b_{o}}(\delta,\rho)$ is a subspace of $\widetilde{\mathbf{O}}'_{b_{o}}(\delta,\rho)$.

We construct the local regularizations by the same method as in \S\ref{local regularization-top},
\S\ref{without bubble tree} and \S\ref{with bubble tree} such that $D\mathcal{S}_{(\kappa,b_o)}$ is surjective for any $b\in \widetilde{\mathbf{O}}^{*}_{b_{o}}(\delta,\rho)$.
 By $\mathbb C^{*}$-action we get the local  regularizations on $\widetilde{\mathbf{O}}'_{b_{o}}(\delta,\rho)$. Denote
$$
 \mathcal W_{*}^{k,2,\alpha}=\{h\in W^{k,2,\alpha}| a(h^1_{0})=0,\theta(h^1_{0})=0\}.
 $$
It is easy to see that  $D\mathcal S_{(\kappa,b_o)}(K\times \mathcal W_{*}^{k,2,\alpha})=D\mathcal S_{(\kappa,b_o)}(K\times \mathcal W^{k,2,\alpha}).$ Then  $D\mathcal S_{(\kappa,b_o)}|_{K\times \mathcal W_{*}^{k,2,\alpha}}$ is also surjective.
Applying the implicit theorem we conclude that both $\mathcal S^{-1}(0)|_{K\times \widetilde{\mathbf{O}}^{*}_{b_{o}}(\delta,\rho)}$ and $\mathcal S^{-1}(0)|_{K\times \widetilde{\mathbf{O}}'_{b_{o}}(\delta,\rho)}$  are smooth  manifolds.
There is a   $\mathbb C^{*}$-action on $\mathcal S^{-1}(0)|_{K\times \widetilde{\mathbf{O}}'_{b_{o}}(\delta,\rho)}.$
Obviously,  $\mathcal S^{-1}(0)|_{K\times \widetilde{\mathbf{O}}^{*}_{b_{o}}(\delta,\rho)}$ is a transverse to the $\mathbb C^{*}$-orbit through $(0,b_{o}).$
On the other hand, $\mathbb C^{*}$-action is proper on  $\mathcal S^{-1}(0)|_{K\times \widetilde{\mathbf{O}}'_{b_{o}}(\delta,\rho)}.$ Define
$$\eta:\mathbb C^{*}\times \mathcal S^{-1}(0)|_{K\times \widetilde{\mathbf{O}}^{*}_{b_{o}}(\delta,\rho)}\to \mathcal S^{-1}(0)|_{K\times \widetilde{\mathbf{O}}'_{b_{o}}(\delta,\rho)}$$ by
$\eta(g,v)=g\cdot v.$
We have proved
\begin{lemma}
For sufficiently small $\delta,\rho,$
$\eta:\mathbb C^{*}\times \mathcal S^{-1}(0)|_{K\times \widetilde{\mathbf{O}}^{*}_{b_{o}}(\delta,\rho)}\to \mathcal S^{-1}(0)|_{K\times \widetilde{\mathbf{O}}'_{b_{o}}(\delta,\rho)}$  maps
$\mathbb C^{*}\times \mathcal S^{-1}(0)|_{K\times \widetilde{\mathbf{O}}^{*}_{b_{o}}(\delta,\rho)}$ diffeomorphically onto a $\mathbb C^{*}$-invariant neighborhood $\mathcal S^{-1}(0)|_{K\times \widetilde{\mathbf{O}}'_{b_{o}}(\delta,\rho)}$ of the $\mathbb C^{*}$-orbit through $(0,b_{o}).$ So
$\mathcal S^{-1}(0)|_{K\times \widetilde{\mathbf{O}}'_{b_{o}}(\delta,\rho)}/ \mathbb C^{*}$ is a smooth manifold. The tangent space of $\mathcal S^{-1}(0)|_{K\times \widetilde{\mathbf{O}}^{*}_{b_{o}}(\delta,\rho)}$ at $(0,b_o)$ is
$$E^*:=\left\{(\kappa_{o},h+\hat{h}_0)\in ker D\mathcal S_{(\kappa,b_o)}\mid a(h^1_{0})=0,\theta(h^1_{0})=0.\right\}$$
\end{lemma}

\chapter{Gluing different holomorphic cascades}\label{gluing-pregluing}

We mainly discuss gluing a holomorphic cascade in $M^+$ and a holomorphic cascade in $(\mathbb {\mathbb{R}}\times \widetilde{M})$, other cases are similar.

\section{Pregluing}

\subsection{Gluing almost complex manifolds}\label{gluing almost complex manifolds}

Consider $\overline{M}^+\cup_Z \mathbb P(\mathcal N\oplus \mathbb C)$. We choose the coordinates $a_{1},\theta_{1}$, for $M^{+}$ and  $a_{2},\theta_{2}$, for $\mathbb {\mathbb{R}}\times \widetilde{M}$.
For any parameter $r>0,$ we can glue $M^+$ and ${\mathbb{R}}\times \widetilde{M}$ to get $M^+$ again as following.  We cut off the
part of $M^{+}$ with cylindrical
coordinate $|a_{1}|>\frac{3lr}{2}$ and glue the remainders along
the collars of length $lr$ of the cylinders with the gluing formulas:
 \begin{eqnarray} \label{gluing_local_relative_node_a}
&&a_1=a_2 + 2l r \\ \label{gluing_local_relative_node_theta}
&&\theta_1 = \theta_2   \;\;mod \;1.
\end{eqnarray}
In terms of the coordinates $(a_1, \theta_1)$ we write $M^+$ as $$M^+=M^{+}_{0}\bigcup\left\{[0,\infty)\times \widetilde{M}\right\}.$$
The line bundle $L$ over $M^+$ remains invariant.
Similarly, we can glue ${\mathbb{R}}\times \widetilde{M}$ and ${\mathbb{R}}\times \widetilde{M}$ to get ${\mathbb{R}}\times \widetilde{M}$ again, the line bundle $L$ over ${\mathbb{R}}\times \widetilde{M}$ remains invariant.

\subsection{Pregluing relative nodes}\label{relative_node_M_M-1}
Denote $\mathcal{M}^{0,1}:={\mathcal{M}}_{A}(\overline{M}^{+},g,m+\mu,{\bf k},\nu)\bigcup_Z
  {\mathcal{M}}_{A}(\mathbb P(\mathcal N\oplus \mathbb C) ,g,m+\mu^{+}+\mu^{-},{\bf k}^+, {\bf k}^-,\nu)/\mathbb{C}^*.$ Let $b=(b_1, b_2)\in \mathcal{M}^{0,1}$,
$b_1=(\Sigma_1,j_1,\bar u_1)$ and $b_2=(\Sigma_2,j_2,\bar u_2)$,
where  $(\Sigma_{1},j_{1})$ and $(\Sigma_{2},j_{2})$ are smooth Riemann surfaces of genus $g_{1}$ and $g_{2}$ joining at $q_{1},q_{2},...,q_{l}$ and
$\bar{u}_1: \Sigma_{1} \rightarrow \overline{M}^+$, $\bar{u}_2: \Sigma_{2} \rightarrow \mathbb P(\mathcal N\oplus \mathbb C)$ are $(j_{i},J)$-holomorphic maps such that
$\bar{u}_{i}(z)$ tangent to $Z$ at the  point $\bar{u}_{1}(q_{j})=\bar{u}_{2}(q_{j})\in Z $ with order $k_j$ as $z\rightarrow q_{j},j=1,2,...,l$. Suppose that both $(\Sigma_i,j_i,\mathbf{y}_i, \mathbf{q})$, $i=1,2$, are stable.

 We choose local Darboux coordinate systems ${\bf w}_{j}$
 near $u(q_j)\in Z$, with ${\bf w_j}(u(q_j))=0$. Choose the local cusp cylinder coordinates $( {s}_{ij}, {t}_{ij})$ on $\Sigma_i$ near $q_j$.
Suppose that
\begin{equation} \label{gluing_local_relative_node-1}
a_{i}({s}_{ij},{t}_{ij})-k_{j}{s}_{ij}-l_{ij}\rightarrow 0,\;\;\;\;
\theta_{i}({s}_{ij},{t}_{ij})-k_{j}{t}_{ij}-\theta_{ij0}\rightarrow 0,\;\;\;\;\;i=1,2,\;j=1,2,...,l
\end{equation}

\begin{figure}[ht]\label{Figure3}
\includegraphics[height=5cm]{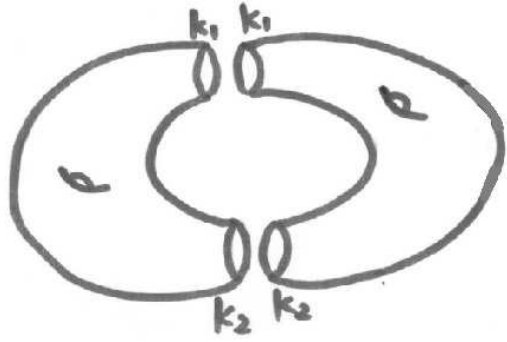}
\caption{Two Relative nodal points} 
\end{figure}
Since there exists a $\mathbb C^{*}$ action on ${\mathbb{R}}\times \widetilde{M},$ we can choose the coordinates  $(a_{2},\theta_{2})$ such that
$\ell_{11}=\ell_{21},\;\theta_{110}=\theta_{210}.$

\v
For any $(r, 0)$ we glue $M^+$ and ${\mathbb{R}}\times \widetilde{M}$ to get again $M^+$  as in \eqref{gluing_local_relative_node_a} and \eqref{gluing_local_relative_node_theta}.
Set
\begin{equation} \label{gluing_local_relative_node-2}
r_{j}=r+\frac{\ell_{2j}-\ell_{1j}}{2l},\;\;\;\; \tau_{j}=\theta_{2j0}-\theta_{1j0},\;\;\;\;j=1,2,...,l.
\end{equation}
We construct a surface
$\Sigma_{(\mathbf {r})} =\Sigma_1\#_{({\mathbf r})} \Sigma_2 $ with gluing formulas:
\begin{eqnarray} \label{glu_j_s}
&& {s}_{1j}= {s}_{2j} + \tfrac{2lr_{j}}{k_j} \\ \label{glu_j_t}
&& {t}_{1j}= {t}_{2j} + \tfrac{\tau_{j} + n_{j}}{k_j}
\end{eqnarray}
for some $n_{j} \in \mathbb Z_{k_{j}}$. Denote
$$
\mathbf w_{1j}=\mathbf w_{j}\circ u_{1},\;\;\;\mathbf w_{2j}=\mathbf w_{j}\circ u_{2}.
$$
 In terms of $(s_{ij}, t_{ij})$ we construct pre-gluing map $u_{(\mathbf r)}:\Sigma_{(\mathbf r)}\rightarrow M^{+}$ as follows: for every $j=1,2$
\begin{align*}
a_{(\mathbf r)}(s_{1j},t_{1j})= \;&k_{j}s_{1j} +l_{1j} +  \beta\left(3-\frac{4k_{j}s_{1j}}{lr_{j}}\right) (a_1(s_{1j},t_{1j})-
k_{j}s_{1j}- l_{1j} )  \\
&+ \beta\left(\frac{4k_{j}s_{1j}}{lr_{j}}-5\right) (a_2(s_{1j},t_{1j})- ks_{2j}
- l_{2j} ) ,
\\
 \theta_{(\mathbf r)}(s_{1j},t_{1j}) =\; & k_{j}t_{1j} +\theta_{1j0} +  \beta\left(3-\frac{4k_{j}s_{1j}}{lr_{j}}\right) (\theta_{1}(s_{1j},t_{1j})-
k_{j}t_{1j}-\theta_{1j0})  \\
&+ \beta\left(\frac{4k_{j}s_{1j}}{lr_{j}}-5\right) (\theta_{2}(s_{1j},t_{1j})- k_{j}t_{2j}-\theta_{2j0}) , \\
{\bf w}_{(\mathbf r)}(s_{1j},t_{1j})=\;&
 \beta\left(3-\frac{4k_{j}s_{1j}}{lr_{j}}\right)   {\bf w}_{1j}(s_{1j},t_{1j}) +
\beta\left(\frac{4k_{j}s_{1j}}{lr_{j}}-5\right)
 {\bf w}_{2j}(s_{1j},t_{1j})  .
 \end{align*}
We associate each $p_{j}$ with
$$\bar{\mathbf{t}}_j=\exp\{2l r_{j}+ 2\pi\sqrt{-1}\tau_{j}\}.$$
Put
$$\mathbb{D}^\circ :=\{\mathbf{t}^\circ_j \mid (\mathbf{t}^\circ _j)^{k_j}=\bar{\mathbf{t}}_j\}.$$
Set $z_{1j}=e^{-s_{1j}-2\pi\sqrt{-1}t_{1j}}$ and $z_{2j}=e^{s_{2j}+2\pi\sqrt{-1}t_{2j}}$. In term of $\mathbf t^{\circ}_{j}$, \eqref{glu_j_s} and \eqref{glu_j_t} can be written as
$$
z_{1j}\cdot z_{2j}=\mathbf t^{\circ}_{j}.
$$

\section{Gluing local regularizations}\label{gluing reg}

Let  ${K}_{b_{i}}\in C^{\infty}\left(\Sigma_{i}, \bar u^{*}_{i}TN\otimes\wedge^{0,1}_{j_{i}}T^{*} \Sigma_{i}\right)$ be the local regularization at $b_{i}$, which supports in the compact subset $\Sigma_{i}(R_0) $ of $\Sigma_{i}$. Then $({K}_{b_{1}}, {K}_{b_{2}})$ can be naturally identified with the subspace in $C^{\infty}\left(\Sigma_{i,0}, \bar u^{*}_{(\mathbf r)}TN\otimes\wedge^{0,1}_{j}T^{*} \Sigma_{(\mathbf r)}\right)$.
\v

We consider the case gluing one node. The general cases are similar.
Let $b=(b_{1}, b_{2})$ be as in \S\ref{relative_node_M_M-1} with one node $q$.
Put
$$E_1:=\{(\kappa_{10},h_1 + \hat{h}_{10})\mid D \mathcal{S}_{(\kappa_{o1},b_1)}(h_1 + \hat{h}_{10})=0,\;\;h_{10}\in \mathbb{H}_q\},$$
$$E_2:=\{(\kappa_{20},h_2 + \hat{h}_{20})\mid D \mathcal{S}_{(\kappa_{o2},b_2)}(h_2 + \hat{h}_{20})=0,\;\;h_{20}\in \mathbb{H}_q\}.$$
Put
$$Ker D \mathcal{S}_{(\kappa_{o},b)}:=E_1\bigoplus_{\bh}E_2=\left\{(\kappa_{10},(h_1,h_{10}), (\kappa_{20},h_2,h_{20}))\in E_1\oplus E_2 \mid h_{10}=h_{20}\in \bh\right\}.$$

\subsection{Estimates of right inverse}

For any $(\kappa,h,h_{0})\in Ker D \mathcal {S}_{(\kappa_{o},b)},$
 where $h\in W^{k,2,\alpha}(\Sigma;u^{*}TN),$  we define $$ \|(\kappa,h)\|_{k,2,\alpha}=\|\kappa\|_{k-1,2,\alpha}+\|h\|_{k,2,\alpha},\;\;\;\;
\|(\kappa,h,h_{0})\|=\|(\kappa,h)\|_{k,2,\alpha}+|h_{0}|.$$
 For any $(\kappa,h_{(r)})\in Ker D \mathcal {S}_{(\kappa_{o},b_{(r)})},$ we define
$$\|(\kappa,h_{(r)})\|=\|\kappa\|_{k-1,2,\alpha}+\|h_{(r)}\|_{k,2,\alpha,r}.
$$
By using the exponential decay of $u_i$ one can
easily prove that $u_{(r)}$ are a family of approximate
$(j,J)$-holomorphic map, precisely the following lemma holds ( for the proof see \cite{LS-1}). \vskip
0.1in \noindent
\begin{lemma} For any $r>R_{0},$  we have
\begin{equation}\label{estimate_p_u_r}
\|\bar{\partial}_{j,J}(u_{(r)})\|_{k-1,2,\alpha,r}\leq
Ce^{-(\mathfrak{c} -\alpha)r} .
\end{equation}  The constants C in the
above estimates are independent of $r$.
\end{lemma}

\begin{lemma}\label{right_inverse_after_gluing}
 Suppose that $D{\mathcal S}_{(\kappa_{o},b)}|_{K_{b}\times  {W}^{k,2,\alpha}}:K_{b}\times  {W}^{k,2,\alpha}\rightarrow L^{k-1,2,\alpha}$ is surjective.
Denote by $Q_{(\kappa_{o},b)}:L^{k-1,2,\alpha}\rightarrow K_{b}\times W^{k,2,\alpha}$ a bounded right inverse of $D{\mathcal
S}_{(\kappa_{o},b)} .$  Then $D{\mathcal
S}_{(\kappa_{o},b_{(r)})}$ is surjective for $r$ large enough. Moreover, there are a right inverses $Q_{(\kappa_{o}, b_{(r)})}$
such that
\begin{equation}
\label{right_verse}
D{\mathcal S}_{(\kappa_{o}, b_{(r)})}\circ Q_{(\kappa_{o}, b_{(r)})}=Id \end{equation}
\begin{equation}
\label{right_estimate}
\|Q_{(\kappa_{o},b_{(r)})}\|\leq  {C}
\end{equation}
 for some constant $C>0$ independent of $ r $.
  \end{lemma}
  \vskip
0.1in
\vskip
0.1in

\noindent {\bf Proof:} We first construct an approximate right inverse $Q'_{(\kappa_{o}, b_{(r)})}$ such that the following estimates holds
\begin{eqnarray}
\label{approximate_right_inverse_estimate_1}
\|Q'_{(\kappa_{o}, b_{(r)})}\|\leq C_1 \\
\label{approximate_right_inverse_estimate_2}
\|D{\mathcal S}_{(\kappa_{o}, b_{(r)})}\circ Q'_{(\kappa_{o}, b_{(r)})}-Id\|\leq \frac{1}{2}.
\end{eqnarray}
 Then the operator $D{\mathcal S}_{(\kappa_{o}, b_{(r)})}\circ Q'_{(\kappa_{o}, b_{(r)})}$ is invertible and a right inverse $Q_{(\kappa_{o}, b_{(r)})}$  of $D{\mathcal S}_{(\kappa_{o}, b_{(r)})}$ is given by
 \begin{equation}
 \label{express_right_inverse}
Q_{(\kappa_{o},b_{(r)})}=Q'_{(\kappa_{o},b_{(r)})}( D{\mathcal S}_{(\kappa_{o}, b_{(r)})}\circ Q'_{(\kappa_{o}, b_{(r)})})\inv
 \end{equation}
Denote $\beta_1=\beta(3/2-\frac{ks_1}{lr}).$ Let $\beta_{2}\geq 0$ be a smooth function such that $\beta_2^2=1-\beta_1^2.$ Given
$\eta\in L_{r}^{2,\alpha}$, we have a pair $(\eta_1,\eta_2)$, where
$$\eta_1=\beta_{1}\eta,\;\;\;\; \eta_2=\beta_{2}\eta.$$ Let
$Q_{(\kappa_{o},b)}(\eta_1,\eta_2)=(\kappa_b,h).$ We may write $h$ as
$(h_1,h_2)$, and define
\begin{equation}
\label{def_h_r}
h_{(r)}=
h_1\beta_1+ h_2\beta_{2}.
\end{equation}
Note that on $\{\frac{lr}{2k}\leq s_1\leq \frac{3lr}{2k}\},$ $\kappa=0$ and  we have $$u_{(r)}|_{\{ s_1 \leq \frac{lr}{2k}\}}=u_1|_{\{ s_1 \leq \frac{lr}{2k}\}},\;\;\;\;\;u_{(r)}|_{\{ |s_2|\leq \frac{lr}{2k}\}}=2lr+u_2|_{\{|s_2|\leq \frac{lr}{2k}\}},$$
so  along $u_{(r)}$ we have $ \kappa_{(r)} =\kappa_b$. Then we define
\begin{equation}
\label{def_approximate_right_inverse}
Q_{(\kappa_{o}, b_{(r)})}^{\prime}\eta = ( \kappa_{(r)} ,
h_{(r)})=(\kappa_b,
h_{(r)}).\end{equation}
Since $|\beta_1|\leq 1$ and $|\frac{\partial \beta_1}{\partial s_1}|\leq \frac{Ck}{lr},$ \eqref{approximate_right_inverse_estimate_1} follows from $\|Q_{(\kappa_{o},b)}\|\leq C_2$ for some constant $C_2>0$.
 We prove \eqref{approximate_right_inverse_estimate_2}. Since $\kappa_b+D_{u}h=\eta$ we have
\begin{equation}\label{app_DS_right_inverse}
 D{\mathcal S}_{(\kappa_{o}, b_{(r)})}\circ Q'_{(\kappa_{o}, b_{(r)})}\eta=\eta\;\;\;\;\;\;for \;\;|s_{i}|\leq \frac{lr}{2k}.
 \end{equation}
 It suffices to estimate the left hand side in the left annulus $\frac{lr}{2k}\leq |s_i|\leq \frac{3lr}{2k}.$
 Note that in this annulus
 $$D_{(\kappa_{o},u_{i})}= D_{u_i},\;\;\; \beta_1^2+\beta_2^2=1,\;\;\kappa_{b}=0,\;\;\; D_{u_{i}}h_{i}=\eta_{i},$$
 $$\beta_1 D_{u_{1}}h_{1}+\beta_2 D_{u_{2}}h_{2}=(\beta_1^2+\beta_2^2)\eta. $$
 Since near the periodic orbit $x(kt)$, $D_{u_{i}}=\bar{\partial}_{J_{0}}+F^{1}_{u_{i}}+F^{2}_{u_{i}}\frac{\p}{\p t}$, we have
\begin{eqnarray}
&&D \mathcal{S}_{(\kappa_{o}, b_{(r)})}\circ Q'_{(\kappa_{o}, b_{(r)})}\eta-(\beta_1^2+\beta_2^2)\eta = \kappa_{(\kappa_{o}, b_{(r)})} + D_{u_{(r)}}h_{(r)}-(\beta_1^2+\beta_2^2)\eta \nonumber\\\label{approximate_difference}
&=&\sum_{i=1}^{2}(\bar{\partial}\beta_{i}) h_{i} +\sum_{i=1}^{2}\beta_{i}(F^{1}_{u_{(r)}}-F^{1}_{u_{i}})h_{i}+\sum_{i=1}^{2}\beta_{i}(F^{2}_{u_{(r)}}-F^{2}_{u_{i}})\p_{t}h_{i}.
\end{eqnarray}
By the exponential decay of $F^{j}_{u_i},j=1,2$ and $\beta_1^2+\beta_2^2=1$ we get
\begin{align}
& \left\|D \mathcal{S}_{(\kappa_{o}, b_{(r)})}\circ Q'_{(\kappa_{o}, b_{(r)})}\eta-\eta\right\|_{k-1,2,\alpha,r} =  \left\|D \mathcal{S}_{(\kappa_{o}, b_{(r)})}\circ Q'_{(\kappa_{o}, b_{(r)})}\eta-(\beta_{1}^2+\beta_{2}^2)\eta\right\|_{k-1,2,\alpha,r} \nonumber \\
&\leq \frac{C_3}{r}\left(\|h_1\|_{k-1,2,\alpha}+ \|h_2\|_{k-1,2,\alpha} \right)  \leq  \frac{C_4}{r} \|\eta\|_{k-1,2,\alpha,r}
\end{align}
for some constant $C_4>0.$ In the last inequality we used that  $\|Q_{(\kappa_{o},b)}\|\leq C_2$ and $(h_{1},h_{2})=\pi_{2}\circ Q_{(\kappa_{o},b)}(\eta_{1},\eta_{2})$, where $\pi_{2}(\kappa_b,h)=h.$
 Then \eqref{approximate_right_inverse_estimate_2} follows by choosing $r$ big enough. The estimate \eqref{approximate_right_inverse_estimate_2} implies that
\begin{equation}\label{app_DS_right_inverse_bound}
\frac{1}{2}\leq \|D{\mathcal S}_{(\kappa_{o}, b_{(r)})}\circ Q'_{(\kappa_{o}, b_{(r)})}\|\leq \frac{3}{2}.
\end{equation}
Then \eqref{right_estimate} follows.  $\Box$
\v
By the same method we can prove (see also \S 5 in \cite{LS-1} for the proof)

\begin{lemma}\label{right_inverse_after_gluing.1}
Suppose that $D{\mathcal S}_{(\kappa_{o},b_{o})}|_{K_{b_{o}}\times  {W}^{k,2,\alpha}}:K_{b_{o}}\times  {W}^{k,2,\alpha}\rightarrow L^{k-1,2,\alpha}$ is surjective.
Denote by $Q_{(\kappa_{o},b_{o})}:L^{k-1,2,\alpha}\rightarrow K_{(\kappa_{o},b_{o})}\times W^{k,2,\alpha}$ a right inverse of $D{\mathcal
S}_{(\kappa_{o},b_{o})}$ with $\|Q_{(\kappa_{o},b_{o})}\|\leq C_{1}$.
Then there exist two constants $\delta_{0}>0$ and $\rho_{0}>0$  depending only on $C_{1}$ and $b_{o}$ such that  for any
 \begin{equation}
 \delta<\delta_{0},\;\;\;\rho<\rho_{0}
 \end{equation}
$D{\mathcal
S}_{(\kappa,b)}$ is surjective for any $(\kappa,b)\in K_{b_{o}}\times O_{b_{0 }}(\delta,\rho)$.
  \end{lemma}
\section{Isomorphism between $Ker D\mathcal{S}_{(\kappa_{o},b)}$ and $Ker D \mathcal{S}_{(\kappa_o,b_{(r)})}$}\label{isomorphism of Ker}
\v\n
For a fixed gluing parameter $(r)=(r,\tau)$
we define a map
$I_r: Ker D\mathcal{S}_{(\kappa_{o},b)}\longrightarrow Ker D \mathcal{S}_{(\kappa_o,b_{(r)})}$ as follows. For any $(\kappa,h,h_{0})\in Ker D \mathcal {S}_{(\kappa_{o},b)},$
 where $h\in W^{k,2,\alpha}(\Sigma;u^{*}TN),$  we write $h=(h_1,h_2  ),$ and define
\begin{equation}
\label{definition_h_ker}
h_{(r)}= \hat{h}_{0} + h_{1}\beta_1+h_{2}\beta_2,
\end{equation}
\begin{equation}\label{definition_I}
I_{r}(\kappa,h,h_{0})=(\kappa,h_{(r)})-Q_{(\kappa_{o},b_{(r)})}\circ D\mathcal{S}_{(\kappa_{o},b_{(r)})}(\kappa,h_{(r)}).
\end{equation}

\begin{lemma}$I_r$ is an isomorphisms for $r$ big enough.
\end{lemma}
 \vskip 0.1in
\noindent {\bf Proof:} The proof is basically a similar gluing
argument as in \cite{D}. The proof is
divided into 2 steps. \vskip 0.1in \noindent {\bf Step 1}.  We define a map $I'_{r}:Ker D\mathcal{S}_{(\kappa_{o},b_{(r)})}\longrightarrow Ker D\mathcal{S}_{(\kappa_{o},b)} $
 and show that $I'_{r}$ is injective for $r$ big enough. For any $(\kappa,h)\in Ker D\mathcal{S}_{(\kappa_{o},b_{(r)})}$ we denote by $h_{i}$ the restriction of $h$
 to the part $|s_{i}|\leq \frac{lr}{k} +\frac{1}{\alpha},$ we get a pair $(h_{1},h_{2}).$ Let
 \begin{equation}
 h_{0}=\int_{S^1}h\left(\frac{lr}{k},t\right)dt.
 \end{equation}
  We denote
  $$\beta [h] =\left((h_{1}- \hat h_{0})\beta\left(\frac{\alpha lr}{k} +1-\alpha s_1\right) + \hat h_{0},\;\;(h_{2}- \hat h_{0})\beta\left(\frac{\alpha lr}{k}+1+\alpha s_2\right) + \hat h_{0}\right)$$
and define
$I'_{r}:Ker D\mathcal{S}_{(\kappa_{o},b_{(r)})}\longrightarrow  Ker D\mathcal{S}_{(\kappa_{o},b)}$ by
\begin{equation}
\label{definition_I'}
I'_{r}(\kappa,h)= (\kappa,\beta [h])-Q_{(\kappa_{o},b)}\circ D\mathcal{S}_{(\kappa_{o},b)}(\kappa,\beta [h]),
\end{equation}
 where  $Q_{(\kappa_{o},b)}$ denotes the right inverse of $D\mathcal{S}_{(\kappa_{o},b)}|_{K_{b}\times W^{k,2,\alpha}}: {K}_{b } \times  W^{k,2,\alpha}\rightarrow L^{k-1,2,\alpha}.$  Since
$D\mathcal{S}_{(\kappa_{o},b)}\circ Q_{(\kappa_{o},b)}=D\mathcal{S}_{(\kappa_{o},b)}|_{K_{b}\times W^{k,2,\alpha}}\circ Q_{(\kappa_{o},b)}=I,$ we have
$I'_{r}(Ker D\mathcal{S}_{(\kappa_{o},b_{(r)})})\subset  Ker D\mathcal{S}_{(\kappa_{o},b)}.$
 \v
Since $\kappa$ and $D_{u} (\beta(h-\hat{h}_{0}))$ have compact support and $F^{i}_{u},i=1,2\in L^{k-1,2,\alpha}$, we have
$D\mathcal{S}_{(\kappa_{o},b)}(\kappa,\beta [h]) \in L^{k-1,2,\alpha}.$
Then $Q_{(\kappa_{o},b)}\circ D\mathcal{S}_{(\kappa_{o},b)}(\kappa,\beta [h])\in K_{b}\times W^{k,2,\alpha}.$

\v
Let $(\kappa,h)\in Ker D\mathcal{S}_{(\kappa_{o},b_{(r)})}$ such that $I'_{r}(\kappa,h)=0.$ Since $\beta(h-\hat{h}_{0})\in W^{k,2,\alpha}$ and $Q_{(\kappa_{o},b)}\circ D\mathcal{S}_{(\kappa_{o},b)}(\kappa,\beta [h])\in K_{b}\times W^{k,2,\alpha},$ then  $I'_{r}(\kappa,h)=0$ implies that $h_{0}=0.$ From \eqref{definition_I'} we have
$$\|I'_{r}(\kappa,h) -  (\kappa, \beta h) \|_{k,2,\alpha}
 \leq C_1  \| \kappa + D_u(\beta h)\|_{k-1,2,\alpha}  $$
$$= C_1 \left \| \kappa  + \beta \left(D_u h
  + D_{u_{(r)}} h+\kappa- D_{u_{(r)}} h- \kappa\right) + (\bar{\partial}\beta) h \right\|_{k-1,2,\alpha} $$
for some constant $C_1>0$. Since $(\kappa,h)\in Ker D\mathcal {S}_{(\kappa_{o},b_{(r)})}$, we have $\kappa+D_{u_{(r)}}h=0.$  We choose $\frac{lr}{2k}>R_{0}$.
As $\kappa|_{|s_i|\geq R_{0}}=0$ and $\beta |_{|s_{i}|\leq \frac{lr}{k}}=1$ we have $\kappa = \beta\kappa$.
Therefore
\begin{eqnarray*}
 \left\|I'_{r}(\kappa,h)-(\kappa,\beta h)\right\|_{k,2,\alpha}
 &\leq&    C  \|(\bar{\partial}\beta)  h \|_{k-1,2,\alpha} +C \sum_{i=1}^{2}\|
\beta_{i;2} (F^{1}_{u_{i}} - F^{1}_{u_{(r)}}) h  \|_{k-1,2,\alpha} \\
&& +C \sum_{i=1}^{2}\|
\beta_{i;2}(F^{2}_{u_{i}} - F^{2}_{u_{(r)}}) \p_{t}h  \|_{k-1,2,\alpha}
\end{eqnarray*}
Note that
$$F^{i}_u = F^{i}_{u(r)},\;i=1,2 \;\;\;\;if \; \; s_1\leq \; \frac{lr}{2k}, \; or \;\; s_2 \geq -\frac{lr}{2k}. $$
By exponential decay of $F^{i}_{u}$ we have
$$\sum_{i=1}^{2}\|
\beta_{i;2} (F^{1}_{u_{i}} - F^{1}_{u_{(r)}}) h  \|_{k-1,2,\alpha}  +\sum_{i=1}^{2}\|
\beta_{i;2}(F^{2}_{u_{i}} - F^{1}_{u_{(r)}}) \p_{t}h  \|_{k-1,2,\alpha}  \leq C_2e^{-\mathfrak{c}\frac{ lr}{2k}}\|\beta h\|_{k,2,\alpha}$$
for some constant $C_2>0$.
Since $(\bar{\partial}\beta(\frac{\alpha lr}{k} + 1 -\alpha s_1))h_{1}$ supports in
$\frac{lr}{k}\leq s_1 \leq \frac{lr}{k}+ \frac{1}{\alpha}$, and over this part
$$|\bar{\partial}\beta(\frac{\alpha lr}{k}+ 1 -\alpha s_1)|\leq 2|\alpha|$$
$$\beta(\frac{\alpha lr}{k}+ 1 + \alpha s_2)= 1,\;\;\;e^{2\alpha|s_1|}\leq e^4e^{2\alpha|s_2|},$$
we obtain
$$\|(\bar{\partial}\beta(\frac{\alpha lr}{k} + 1 -\alpha s_1))h_{1}\|_{k-1,2,\alpha}\leq2
|\alpha|e^4 \|h_2\|_{k-1,2,\alpha} \leq 2|\alpha|e^4 \|\beta h\|_{k-1,2,\alpha}.$$
Similar inequality for $(\bar{\partial}\beta(\frac{\alpha lr}{k} + 1 + \alpha s_2))h_{2}$
also holds. So we have
$$\|(\bar{\partial}\beta)h\|_{k-1,2,\alpha} \leq 4|\alpha|e^4 \|\beta h\|_{k,2,\alpha}.$$
Hence
\begin{equation}\label{delta_I'}
\|I'_{r}(\kappa,h) - (\kappa, \beta h)\|_{k,2,\alpha} \leq (4e^4|\alpha| + C_3
e^{- \frac{\mathfrak{c}   lr}{2k}})\|\beta h\|_{k,2,\alpha} \leq   1/2  \|\beta h\|_{k,2,\alpha}  \end{equation}
for some constant $C_3>0$, here we choosed $0<\alpha<\frac{1}{16e^4}$ and  $r$   big enough  such
that $\frac{lr}{k}>\frac{1}{\alpha}$ and  $C_3 e^{- \frac{\mathfrak{c}   lr}{2k}} < 1/4$.

Then $I'_{r}(\kappa,h)=0$ and
\eqref{delta_I'} gives us $$\|\kappa\|_{k-1,2,\alpha}=0 , \;\;  \|\beta h\|_{k,2,\alpha}=0 .$$
It follows that $\kappa = 0, \;\; h=0$. So $I'_r$ is injective.
\vskip 0.1in
\noindent

{\bf Step 2}.
Since $\|Q_{(\kappa_{o},b_{(r)})}\|$ is uniformly bounded, from \eqref{definition_I}  and \eqref{right_estimate}, we have
$$\|I_r((\kappa,h),h_{0}) - (\kappa, h_{(r)})\|_{1,p,\alpha,r}\leq C_4\|D\mathcal {S}_{(\kappa_{o},b_{(r)})} (\kappa, h_{(r)})\|$$
for some constant $C_4>0$. By a similar culculation as in the proof of Lemma \ref{right_inverse_after_gluing} we obtain
\begin{equation}\label{delta_I}
\|I_r((\kappa,h),h_{0}) - (\kappa, h_{(r)})\|_{1,p,\alpha,r}\leq \frac{C_5}{r}
(\|h\|_{k-1,2,\alpha} + |h_0|)\end{equation}
for some constant $C_5>0$.
In particular, it holds for $p=2$. It remains to show that
$\|h_{(r)}\|_{2,\alpha,r}$ is close to $\|h\|_{2,\alpha}$. Denote $\pi_2$ the
projection into the second component, that is, $\pi_2((\kappa,h),h_{0})=h$. Then
$\pi(ker D{\mathcal S}_{(\kappa_{o},b)})$ is a finite dimentional space. Let $f_i,\;i=1,..,d$
be an orthonormal basis. Then $F=\sum f_i^2e^{2\alpha|s|}$ is an integrable
function on $\Sigma$. For any $\epsilon' >0$, we may choose $R_0$ so big that
$$\int_{|s_i|\geq R_0}F \leq \epsilon'.$$
Then the restriction of $h$ to $|s_i|\geq R_0$ satisfies
$$\|h|_{|s_i|\geq R_0}\|_{2,\alpha}\leq \epsilon'\|h\|_{2,\alpha},$$
therefore
\begin{equation}\label{lower_bound_h_r}
\|h_{(r)}\|_{2,\alpha,r}\geq \|h|_{|s_i|\leq R_0}\|_{2,\alpha} + |h_0|
\geq (1 - \epsilon')\|h\|_{2,\alpha} + |h_0|,
\end{equation}
for $r>R_0$. Suppose that $I_r((\kappa,h),h_{0})=0$. Then \eqref{delta_I} and \eqref{lower_bound_h_r} give us $h=0$ and $h_{0}=0$,
and so $\kappa=0$. Hence $I_r$ is injective.
\vskip 0.1in
\noindent
The {\bf step 1} and {\bf step 2} together show that both $I_{r}$ and $I'_r$
are isomorphisms for $r$ big enough.  $\Box$
\v
The above lemmas can be immediately generalize to the case gluing several nodes. In particular we have
\begin{lemma}\label{isomor of ker}
For $|\mathbf{r}|>R_0$ there is an isomorphism
$$I_{(\mathbf r)}: \ker D \mathcal{S}_{( \kappa_{o}, b_{o})}\longrightarrow \ker D \mathcal{S}_{( \kappa_{o},b_{(\mathbf r)})}.$$
\end{lemma}
\v\n

\v
\chapter{Global regularization}\label{global_r}
\v

\v
\section{A finite rank orbi-bundle over ${\overline{\mathcal{M}}}_{A}(M^{+};g,m+\mu,{\bk},\nu)$}\label{finite rank orbi-bundle}
\v
By the compactness of ${\overline{\mathcal{M}}}_{A}(M^{+};g,m+\mu,{\bk},\nu)$
there exist finite points $[b_i]\in {\overline{\mathcal{M}}}_{A}(M^{+};g,m+\mu,{\bk},\nu)$, $1\leq i \leq \mathfrak{m}$, such that
\begin{itemize}
\item[(1)] The collection $\{\mathbf{O}_{[b_i]}(\delta_i/3,\rho_i/3)
\mid 1\leq i \leq \mathfrak{m}\}$ is an open cover of
${\overline{\mathcal{M}}}_{A}(M^{+};g,m+\mu,{\bk},\nu)$.
\item[(2)] Suppose that $\widetilde{\mathbf{O}}_{b_i}(\delta_i,\rho_i)
\cap \widetilde{\mathbf{O}}_{b_j}(\delta_j,\rho_j)
\neq\phi$. For any $b\in \widetilde{\mathbf{O}}_{b_i}(\delta_i,\rho_i)
\cap \widetilde{\mathbf{O}}_{b_j}(\delta_j,\rho_j)$, $G_b$ can be imbedded into both $G_{b_i}$ and $G_{b_j}$ as subgroups.
\end{itemize}
\begin{remark}\label{finite covering}
We may choose $[b_i]$, $1\leq i \leq \mathfrak{m}$, such that if $[b_i]$ lies in the top strata for some $i$, then
$\mathbf{O}_{[b_i]}(\delta_i,\rho_i)$ lies in the top strata.
\end{remark}
\v\n
Set
$$\mathcal{U}=\bigcup_{i=1}^{\mathfrak{m}}
\mathbf{O}_{[b_{i}]}(\delta_i/2,\rho_i/2).$$
There is a forget map
$$\mathscr{P}: \mathcal{U}\to \overline{\mathcal{M}}_{g,m+\mu},\;\;\;[(j,{\bf y,p}, u)]\longmapsto [(j,{\bf y,p})].$$
We construct a finite rank orbi-bundle $\mathbf{F}$ over
$\mathcal{U}$. The construction imitates Siebert's construction. We can slightly
deform $\omega$ to get a rational class $[\omega^*]$ on $\overline{M}^{+}$.
By taking multiple, we can
assume that $[\omega^*]$ is an integral class on $\overline{M}^{+}$.
\v
 Therefore, it is the Chern class of a complex line bundle $L$ over $\overline{M}^{+}$ ( see \S\ref{line bundle}). Let $i$ be the complex structure on $L$. We choose a Hermition metric $G^L$ and the associate unitary connection $\nabla^{L}$ on $L$. For $M^{+} \cup_{Z} \mathbb P(\mathcal N\oplus \mathbb C)$ we have line bundle $L \cup L'$ over $M^{+} \cup_{Z} \mathbb P(\mathcal N\oplus \mathbb C)$, where $L' = p^{*}(L|_{Z} )$, and $p:\mathbb P(\mathcal N\oplus \mathbb C)\to Z$ is the projection. To simplify notations we simply write $L \cup L'$ as $L$.

\v
Let $(\Sigma,j,{\bf y}, {\bf p},\nu)$ be a marked nodal Riemann surface
of genus $g$ with $m$ distinct marked points ${\bf y}=(y_1,...,y_m)$, $\mu$  distinct puncture points ${\bf p}=(p_1,...,p_{\mu})$, and $u:\Sigma\to M^+$
be a smooth map satisfying the nodal conditions.
We have complex line bundle $u^*L$ over $\Sigma$ with complex structure $u^*i$.
The unitary connection $u^{*}\nabla^{L}$ splits into $ u^{*}\nabla^{L}:=u^{*}\nabla^{L,(1,0)}\oplus u^{*}\nabla^{L,(0,1)}$.
Denote
$$D^L:=u^{*}\nabla^{L,(0,1)}:W^{k,2}(\Sigma,u^{*}L)\to W^{k-1,2}(\Sigma,u^{*}L\otimes\wedge_{j}^{(0,1)}T^{\star}\Sigma).$$
$D^{L}$ takes $s\in W^{k,2}(\Sigma,u^{*}L)$ to the $\mathbb C$-antilinear part of $\nabla^{L}$, where $s$ is a section of $L$. One can check that
$$
D^{L}(f\xi)=\bar{\p}_{\Sigma}f \otimes \xi +  f\cdot D^{L}\xi.
$$
 $D^{L}$ determines a holomorphic structure on $u^*L$, for
which $D^{L}$ is an associated Cauchy-Riemann operator (see \cite{HLS,IS}).
Then $u^*L$ is a holomorphic line bundle.
\v
Let $\lambda_{(\Sigma, j)}$ be
the dualizing sheaf of meromorphic 1-form with at worst simple pole at the nodal points and for each
nodal point p, say $\Sigma_1$ and $\Sigma_2$ intersects at p,
$$Res_p(\lambda_{(\Sigma_1, j_1)}) + Res_p(\lambda_{(\Sigma_2, j_2)})=0.$$
Let $\Pi:\overline{\mathscr{C}}_{g}\to \overline{\mathcal{M}}_{g}$ be the universal curve. Let $\lambda$ be the relative dualizing sheaf over $\overline{\mathscr{C}}_{g}$, the restriction of $\lambda$ to $(\Sigma, j)$  is $\lambda_{(\Sigma, j)}$. Set $\Lambda_{(\Sigma, j)}:=\lambda_{(\Sigma, j)}\left(\sum_{i=1}^{n} y_i + \sum_{j=1}^\mu p_j\right)$, $\mathbf{L}:=\Lambda\otimes u^*L $.
Then $\mathbf{L}\mid_{b}$ is a holomorphic line bundle over $\Sigma$.
We have a Cauchy-Riemann operator $  \bar{\p}_b,$
$H^0(\Sigma, \widetilde{\mathbf{L}}\mid_{b})$ is the $ker \bar{\p}_b$. Here the $\bar{\p}$-operator depends on the complex structure $j$ on $\Sigma$ and the bundle $u^*L$, so we denote it by $\bar{\p}_b$.
\v
If $\Sigma_\nu$
is not a ghost component, there exist a constant $\hbar_{o}>0$ such that
$$\int_{u(\Sigma_\nu)}\omega^*> \hbar_{o} .$$ Therefore, $c_1(u^*L)(\Sigma_\nu)>0$. For ghost component $\Sigma_\nu$, $\lambda_{\Sigma_\nu}\left(\sum_{i=1}^{n} y_i + \sum_{j=1}^\mu p_j\right)$ is positive. So
by taking the higher power of $\widetilde{\mathbf{L}}$, if necessary, we can assume that
$\widetilde{\mathbf{L}}\mid_{b}$ is very ample for any $b=(a,v)\in \widetilde{\mathbf{O}}_{b_{o}}(\delta,\rho)$. Hence, $H^1(\Sigma, \widetilde{\mathbf{L}}\mid_{b})= 0$. Therefore,
$H^0(\Sigma, \widetilde{\mathbf{L}}\mid_{b})$ is of constant rank ( independent of $b\in \widetilde{\mathbf{O}}_{b_{o}}(\delta,\rho)$). We have a finite rank bundle $\widetilde{F}$ over $\widetilde{\mathbf{O}}_{b_{o}}(\delta,\rho)$,
whose fiber at $b=(j,{\bf y,p},v)\in\widetilde{\mathbf{O}}_{b_{o}}(\delta,\rho)$
is $H^0(\Sigma, \widetilde{\mathbf{L}}\mid_{b})$.
The finite group $G_{b}$ acts on the bundle on $\widetilde{F}\mid_{b}$ in a natural way.
\v
\begin{lemma}\label{Transformation for F}
For any $\varphi\in Diff^+(\Sigma)$ denote
	$$ b'=(j',{\bf y',p'},u') =\varphi\cdot (j,{\bf y,p},u)=(\varphi^*j, \varphi^{-1}{\bf y,p}, \varphi^{*}u).$$
	Then the following hold
\v {\bf (a).} $\widetilde{\mathbf{L}}|_{b'}=\varphi^{*}\widetilde{\mathbf{L}}|_{b},$\;\;\; $(u')^{*}i=\varphi^{*}(u^{*}i)$
\v {\bf (b).}
$D^{\widetilde{\mathbf L}}|_{b'}(\varphi^{*}\xi)=\varphi^{*}(D^{\widetilde{\mathbf L}}|_{b}(\xi)).$
\end{lemma}
\v\n
It follows from {\bf (b)} above that
if we choose another coordinate system $\mathbf{A}'$ and another local model $\widetilde{\mathbf{O}}_{b'_{o}}(\delta',\rho')/G_{b'_o}$, we have
$$H^0(\Sigma, \widetilde{\mathbf{L}}\mid_{b})\cong H^0(\Sigma, \widetilde{\mathbf{L}}'\mid_{b'}).$$
But the coordinate transformation is continuous. So we get a continuous bundle
$F\rightarrow \mathcal{U}$. Moreover,  by (1) and (2) we conclude that $F$ has a ``orbi-vector bundle" structure over $\mathcal{U}$.
\v

\v
Both $\widetilde{K}_{b_i}$ and $\widetilde{F}\mid_{b_i}$ are representation spaces of $G_{b_i}$. Hence they can be decomposed as sum of irreducible representations.
There is a result in algebra saying that the irreducible factors of group
ring contain all the irreducible representations of finite group. Hence, it is
enough to find a copy of group ring in $\widetilde{F}(b_i)\mid_{b_i}$. This is done by algebraic geometry.
We can assume that $\widetilde{\mathbf L}$
induces an embedding of $\Sigma$
into $\mathbb{C}P^{N_i}$ for some $N_i$. Furthermore, since $\widetilde{\mathbf L}$ is invariant under
$G_{b_i}$ , $G_{b_i}$ also acts effectively naturally on $\mathbb{C}P^{N_i}$. Pick any point $x_0 \in
im(\Sigma)\subset \mathbb{C}P^{N_i}$ such that $\sigma_k(x_0)$ are mutually different for any $\sigma_k\in G_{b_i} $.
Then, we can find a homogeneous polynomial f of some degree, say $\mathsf{k}_i,$ such
that $f(x_0)\ne 0$, $f(\sigma_k(x_0)) = 0$ for $\sigma_k\ne I_d$. Note that $f\in H^0(\mathcal O(\mathsf{k}_i))$. By pull back over $\Sigma$, $f$ induces
a section $v \in H^0(\Sigma, \widetilde{\mathbf L}
^{\mathsf{k}_i})$. We replace $\widetilde{\mathbf L}$
by $\widetilde{\mathbf L}^{\mathsf{k}_i}$ and redefine $F_i\mid_{b_i} =   H^0(\Sigma, \widetilde{\mathbf L}^{\mathsf{k}_i}\mid_{b_i})$. Then $G_{b_i}\cdot v$ generates a group ring, denoted by $\ll G_{b_i}\cdot v\gg$. It is obvious that $\ll G_{b_i}\cdot v\gg$ is isomorphic to $\mathbb{R}[G_{b_i}]$, so $F_i\mid_{b_i}$ contains a copy of group ring. We denote the obtained bundle by $\mathbf{F}(\mathsf{k}_i)$.
\begin{lemma}\label{finite cov}
We have a continuous ``orbi-vector bundle" $\mathbf{F}(\mathsf{k}_i)\rightarrow \mathcal{U}$ such that
$\mathbf{F}(\mathsf{k}_i)\mid_{b_i}$ contains a copy of group ring
$\mathbb{R}[G_{b_i}].$
\end{lemma}
\v\n
In \cite{LS-3} we proved
\begin{lemma}\label{smoothness of bundle-1} For the top strata, in the local coordinate system $\mathbf{A}$ the bundle $\widetilde{\mathbf{F}}$ is smooth. Furthermore, for any base $\{e_{\alpha}\}$ of the fiber at $b_o$ we can get a smooth frame fields $\{e_{\alpha}(a,h)\}$ for the bundle $\widetilde{\mathbf{F}}$ over $\widetilde{\mathbf{O}}_{b_o}(\delta_{o},\rho_{o})$ .
\end{lemma}

\v
\begin{remark}\label{isotropy group}
Let $G_{b_o}$ be the isotropy group at $b_o$. $D^{\widetilde{\mathbf L}}$ is $G_{b_o}$-equivariant and $G_{b_o}$ acts on $\mbox{ker} D^{\widetilde{\mathbf L}}|_{b_{o}}$. We may choose a $G_{b_o}$-equivariant right inverse $Q^{\widetilde{\mathbf{L}}}_{b_{o}}$.
So we have a $G_{b_o}$-equivariant version of Lemma \ref{smoothness of bundle-1}. In particular,
for any base $\{e_{\alpha}\}$ of the fiber at $b_o$ we can get a smooth $G_{b_o}$-equivariant frame fields $\{e_{\alpha}(a,h)\}$ for the bundle $\widetilde{\mathbf{F}}$ over $\widetilde{\mathbf{O}}_{b_o}(\delta_{o},\rho_{o})$ ( see \cite{LS-3} ).
\end{remark}

\v

\v\n
Put $\mathbf{F}=\bigoplus_{i=1}^{\mathfrak{m}}\mathbf{F}(\mathsf{k}_i).$

\section{Gluing the finite rank bundle $\widetilde{\mathbf{F}}$}

We recall some results in \cite{LS-3}.
Let $(U,z)$ be a local coordinates on $\Sigma$ around a nodal point ( or a marked point) $q$ with  $z(q)=0$ . Let $b=(\mathbf{s},u)\in \widetilde{\mathbf{O}}_{b_o}(\delta_{o},\rho_{o})$ and $e$ be a local holomorphic section of $u^{*}L|_{U}$ with $\|e\|_{G^L}(q)\neq 0$ for $q\in U$.
Then for any $\phi\in   \widetilde{\mathbf F}|_{b}$ we can write
\begin{equation}\label{eqn_phi_loc_hol}
\phi|_{U} =f \left(\frac{dz}{z}\otimes e \right)^{\mathsf{k}},\;\;\mbox{ where } f\in \mathcal{O}(U).
\end{equation}
In terms of the holomorphic cylindrical coordinates $(s,t)$ defined by $z=e^{-s+2\pi\sqrt{-1}t}$ we can re-written \eqref{eqn_phi_loc_hol} as
$$
\phi(s,t)|_{U}=f(s,t) \left((ds+2\pi\sqrt{-1}dt)\otimes e \right)^{\mathsf{k}},
$$
where $f(z)\in \mathcal O(U)$. It is easy to see that $|f(s,t)-f(-\infty, t)|$ uniformly exponentially converges to 0 with respect to $t\in S^1$ as
$|s|\to \infty$.
\v
For any $\zeta\in C^{\infty}_{c}(\Sigma,\widetilde{\mathbf L}|_{b})$ and
any section $\eta \in
C^{\infty}_{c}(\Sigma, \widetilde{\mathbf L}|_{b}\otimes\wedge^{0,1}_{j}T^{*}\Sigma)$ we
define weighted norms $\|\zeta\|_{j,k,2,\alpha}$ and $\|\eta\|_{j,k-1,2,\alpha}$. Denote by $W^{k,2,\alpha}(\Sigma;\widetilde{\mathbf L}|_{b})$ and
$W^{k-1,2,\alpha}(\Sigma, \widetilde{\mathbf L}|_{b}\otimes \wedge^{0,1}_{j}T^{*}\Sigma)$ the complete spaces with respect to the norms respectively. We also define the space ${\mathcal W}^{k,2,\alpha}(\Sigma;\widetilde{\mathbf L}|_{b})$.
\v
Let $(\Sigma, j, {\bf y})$ be a marked nodal Riemann surface of genus $g$ with $n$ marked points. Suppose that $\Sigma$ has $\mathfrak{e}$ nodal points $\mathbf{p}=(p_{1},\cdots,p_{\mathfrak{e}})$ and $\iota$ smooth components. We fix a local coordinate system $\mathbf{s}\in\mathbf{A}$ for the strata of $\widetilde{\mathcal{M}}_{g,n}$, where
$\mathbf{A}=\mathbf{A}_1\times \mathbf{A}_2\times...\times \mathbf{A}_\iota$. Let $b_o=(\mathbf{s},u)$ where $u:\Sigma\to M$ be $(j,J)$-holomorphic map. For each node $p_i$ we can glue $\Sigma$ and $u$ at $p_i$ with gluing parameters $(\mathbf{r})=((r_1,\tau_1),...,(r_\mathfrak{e}, \tau_\mathfrak{e}))$ to get $\Sigma_{(\mathbf{r})}$ and $u_{(\mathbf{r})}$, then we glue $\widetilde{\mathbf{F}}\mid_b$ to get $\widetilde{\mathbf{F}}\mid_{b_{(\mathbf{r})}}$.
Denote $|\mathbf r|=\min_{i=1}^{\mathfrak{e}}|r_{i}|.$
\begin{lemma}\label{aright_inverse_after_gluing}
$D^{\widetilde{\mathbf L}}|_{b_{(\mathbf{r})}}$ is surjective for $|\mathbf{r}|$ large enough. Moreover, there is a $G_{b_{(\mathbf{r})}}$-equivariant right inverse $ Q^{\widetilde{\mathbf{L}}}_{b_{(\mathbf{r})}}$
such that
\begin{equation}
\label{right_estimate}
\| Q^{\widetilde{\mathbf{L}}}_{b_{(\mathbf{r})}}\|\leq  \mathsf{C}
\end{equation}
 for some constant $ \mathsf{C}>0$ independent of $(\mathbf{r})$.
\end{lemma}

\begin{lemma}\label{lem_est_I_r}
{\bf (1)} $I^{\widetilde{\mathbf{L}}}_{(r)}: \ker D^{\widetilde{\mathbf L}}|_{b_{o}}\longrightarrow \ker D^{\widetilde{\mathbf L}}|_{b_{(r)}}$ is a $\frac{|G_{b_{o}}|}{|G_{b_{(r)}}|}$-multiple covering map for $r_i$, $1\leq i \leq \mathfrak{e}$, large enough, and
$$\|I^{\widetilde{\mathbf L}}_{(\mathbf{r})}\|\leq \mathsf C ,$$
for some constant $\mathsf C>0$ independent of $(\mathbf{r})$.
\v
{\bf (2)} $I^{\widetilde{\mathbf{L}}}_{(r)}$ induces a isomorphism $I^{\mathbf{L}}_{(r)}: \ker D^{\mathbf L}|_{b_{o}}\longrightarrow \ker D^{\mathbf L}|_{b_{(r)}}$.
\end{lemma}

For fixed $(\mathbf{r})$ we consider the family of maps:
$$
\mathcal {F}_{(\mathbf{r})}: \mathbf{A}\times   W^{k,2,\alpha}(\Sigma_{(\mathbf{r})},u^{\star}_{(\mathbf{r})}TM)
\times {\mathcal W}^{k,2,\alpha}(\Sigma_{(\mathbf{r})},\widetilde{\mathbf L}|_{b_{(\mathbf{r})}})\to W^{k-1,2,\alpha}(\Sigma_{(\mathbf{r})},\wedge^{0,1}T\Sigma_{(\mathbf{r})}\otimes \widetilde{\mathbf L}|_{b_{(\mathbf{r})}})$$
defined by
\begin{equation}\label{def_F_r}
\mathcal {F}_{(\mathbf{r})}(\mathbf{s},h,\xi)= P^{\widetilde{\mathbf{L}}}_{b,b_{(\mathbf{r})}}
\circ D^{\widetilde{\mathbf L}}_{b}\circ (P^{\widetilde{\mathbf{L}}}_{b,b_{(\mathbf{r})}})^{-1}\xi,
\end{equation}
where $b=((\mathbf r),\mathbf{s},v_{\mathbf r})$ and $v_{\mathbf r}=\exp_{u_{(\mathbf r)}}h$.
By implicit function theorem we have

\v\v
\begin{lemma}\label{gluing}
There exist $\delta>0$, $\rho>0$ and a small neighborhood $\widetilde O_{(\mathbf{r})}$ of $0 \in \ker\;D^{\widetilde{\mathbf L}}|_{b_{(\mathbf{r})}}$ and a unique smooth map
$$f^{\widetilde{\mathbf L}}_{(\mathbf{r})}: \widetilde{\mathbf{O}}_{b_{(\mathbf{r})}}(\delta,\rho)
\times \widetilde{O}_{(\mathbf{r})}\rightarrow W^{k-1,2,\alpha}(\Sigma_{(\mathbf{r})},\wedge^{0,1}T\Sigma_{(\mathbf{r})}\otimes \widetilde{\mathbf L}|_{b_{(\mathbf{r})}})$$ such that for any $(b,\zeta)\in \widetilde{\mathbf{O}}_{b_{(\mathbf{r})}}(\delta,\rho)
\times \widetilde{O}_{(\mathbf{r})}$
\begin{equation*}
D^{\widetilde{\mathbf L}}_{b}\circ(P^{\widetilde{\mathbf{L}}}_{b,b_{(\mathbf{r})}})^{-1}\left(\zeta + Q^{\widetilde{\mathbf L}}_{b_{(\mathbf{r})}}\circ f^{\widetilde{\mathbf L}}_{\mathbf{s},h,(\mathbf{r})}(\zeta)\right)=0.
\end{equation*}
\end{lemma}
\v Together with  $I^{\mathbf L}_{(\mathbf{r})}$ we have gluing map
$$Glu^{\mathbf L}_{(\mathbf r)}:\mathbf{F}\mid_{[b_o]}\to \mathbf{F}\mid_{[b]}\;\;\;for\; any \;[b]\in \;\mathbf{O}_{[b_{(\mathbf{r})}]}(\delta,\rho)$$ defined\;by
$$Glu^{\mathbf L}_{(\mathbf r)}([\zeta]):=\left[(P^{\widetilde{\mathbf{L}}}_{b,b_{(\mathbf{r})}})^{-1}\left(I^{\widetilde{\mathbf L}}_{(\mathbf{r})}\zeta + Q^{\widetilde{\mathbf L}}_{b_{(\mathbf{r})}}\circ f^{\widetilde{\mathbf L}}_{\mathbf{s},h,(\mathbf{r})}I^{\widetilde{\mathbf L}}_{(\mathbf{r})}\zeta \right)\right],\;\;\;\;\forall [\zeta]\in \mathbf{F}\mid_{[b_o]}.$$

Given a frame $e_{\alpha}(z)$ on $\widetilde{\mathbf{F}}\mid_{b_o}$, $1\leq \alpha\leq rank\; \widetilde{\mathbf{F}},$ as Remark \ref{isotropy group} we have a $G_{b_o}$-equivariant frame field
\begin{equation} \label{glu_solu}
e_{\alpha}((\mathbf{r}),\mathbf{s}, h)(z)=(P^{\widetilde{\mathbf{L}}}_{b,b_{(\mathbf{r})}})^{-1}
\left(I^{\widetilde{\mathbf L}}_{(\mathbf{r})}e_{\alpha}  + Q^{\widetilde{\mathbf L}}_{b_{(\mathbf{r})}}\circ f^{\widetilde{\mathbf L}}_{\mathbf{s},h,(\mathbf{r})}I^{\widetilde{\mathbf L}}_{(\mathbf{r})}e_{\alpha} \right)(z)
\end{equation}
over $D_{R_{0}}^{*}(0)\times \widetilde{\mathbf{O}}_{b_o}(\delta_{o},\rho_{o})$,
where $z$ is the coordinate on $\Sigma$, and
$$D_{R_{0}}^{*}(0):=\bigoplus_{i=1}^{\mathfrak{e}}\left\{(r, \tau)\mid R_0< r<\infty, \;\tau\in S^1\right\}.$$
For any fixed $(\mathbf{r})$, $e_{\alpha}$ is smooth with respect to $\mathbf{s},h$ over $\widetilde{\mathbf{O}}_{b_o}(\delta_{o},\rho_{o})$.

\v

Let $\alpha_{(r_{i})}:[0,2r_{i}]\to [0,2R_{0}]$ be a smooth function satisfying
$$
  \alpha_{(r_{i})}(s)=\left\{
\begin{array}{ll}
s\;\;\;\;  &   [0,\frac{R_{0}}{2}-1] \\
 \frac{R_{0}}{2}+\frac{R_{0}}{2r_{i}-R_{0}}(s-R_{0}/2) \;\;\;\;\;\; &  [R_{0}/2,2r_{i}-R_{0}/2]   \\
 s-2r_{i}+2R_{0} & [2r_{i}-\frac{R_{0}}{2}+1,2r_{i}]
\end{array}
\right.
$$
Set $\alpha_{(r_{i})}:[-2r_{i},0]\to [-2R_{0},0]$ by $\alpha_{(r_{i})}(s)=-\alpha_{(r_{i})}(-s).$
Let $(s_{1}^{i},t_{1}^{i})$ and $(s_{2}^{i},t_{2}^{i})$ be cusp cylinder coordinates around $p_{i}$, thus $z_{i}=e^{-s_{1}^{i}-2\pi \sqrt{-1}t_{1}^{i}}$ and $w_{i} =e^{ s_{2}^{i}+2\pi \sqrt{-1}t_{2}^{i}}$.  Denote $$W_{i}(R)=\{ |s_{1}^{i}|>R \}\cup \{ |s_{2}^{i}|>R \}.$$ Obviously, $W(R)=\cup_{i=1}^{\mathfrak{e}}W_{i}(R).$
We can define a map $\varphi_{(\mathbf{r})}:\Sigma_{(\mathbf{r})}\to \Sigma_{(\mathbf{R}_{0})}$ as follows:
$$\varphi_{(\mathbf{r})}=\left\{
\begin{array}{ll}
p, & p\in \Sigma(R_{0}/4).\\
(\alpha_{(r_{i})}(s_{i}),t_{i})   \;\;\;\;\;\;\;& (s_{1}^{i},t_{1}^{i})\in W_{i}(R_{0}/4),\;i=1,\cdots,\mathfrak{e}.
\end{array}
\right.
$$
Then we obtain a family of Riemann surfaces $\left(\Sigma_{(\mathbf{R}_{0})},(\varphi_{(\mathbf{r})}^{-1})^{*}j_{\mathbf{r}} ,\varphi_{(\mathbf{r})}^{-1}(\mathbf{y})\right)$.
Denote $u^{\circ}_{(\mathbf{r})}:=u_{(\mathbf{r})}\circ \varphi_{\mathbf{r}}^{-1}.$

In \cite{LS-2} we have proved the following lemma.
\begin{lemma}\label{smoothness of Glu^L}
There exists positive  constants  $ \mathsf{d},R$ such that
for any  $h\in W^{k,2,\alpha}\left(\Sigma_{(R_{0})},(u_{(R_{0})})^*TM \right),$  $\zeta\in \ker D^{\widetilde{\mathbf L}}|_{b_{o}}$    with    $$\|\zeta\|_{\mathcal W,k,2,\alpha}\leq \mathsf{d},\;\;\;\;\;\|h-\hat h_{(\mathbf r)}\|< \mathsf{d},\;\;\;\;\;|\mathbf r|\geq R, $$
$( \varphi_{\mathbf r}^{-1})^{*}(Glu^{\widetilde{\mathbf L}}_{\fs,(\mathbf r),h'}(e_{\alpha}) )$ is smooth with respect to $(\mathbf{s}, (\mathbf{r}),h)$ \; for any $e_{\alpha}\in \ker D^{\widetilde{\mathbf L}}|_{b_o}$,  where  $h'=(\exp_{u_{(\mathbf r)}}^{-1}\circ (\exp_{u_{(\mathbf R)_{0}}}(h)\circ \varphi_{(\mathbf{r})})$. In particular $Glu^{\widetilde{\mathbf L}}_{\fs,(\mathbf r),h'}(e_{\alpha})\mid_{\Sigma(R_0)}$
 is smooth.
\end{lemma}

\section{Global regularization and virtual neighborhoods}\label{global_regu}

We are going to construct a bundle map $\mathfrak{i}:\mathbf{F}\rightarrow  \E$. We first define a bundle map
$\mathfrak{i}:  {\mathbf{F}}(\mathsf{k}_{i})\rightarrow \mathcal{E}$.
Consider two different cases:
\v
{\bf Case 1.} $[b_i]$ lies in the top strata $\mathcal{M}_{g,m+\mu}(A)$. Denote $b_o=b_i$. Choose a local coordinate system $(\psi, \Psi)$ for $\mathcal{Q}$ and a local model $
\widetilde{\mathbf{O}}_{b_o}(\delta_{b_{o}},\rho_{b_{o}})/G_{b_o}$ around $[b_o]$. We have an isomorphism
\begin{equation}\label{trivialization-1}
P_{b_o,b}=\Phi\circ\Psi_{j_o,j_a}:\widetilde{\E}_{b_o}\rightarrow \widetilde{\E}_{b},\;\;\forall \;\;b\in \widetilde{\mathbf{O}}_{b_{o}}(\delta_{b_{o}},\rho_{b_{o}}).
\end{equation}
To simplify notations we denote  $\widetilde{\mathbf{F}}(\mathsf{k}_i)=\widetilde{H}$, $P_{b_o,b}=P$ in this section.
\v
Choosing a base $\{e_{\alpha}\}$ of the fiber $\widetilde{H}\mid_{b_o}$, by Lemma \ref{smoothness of bundle-1} we can get a smooth frame fields $\{e_{\alpha}\}$ for the bundle $\widetilde{H}$ over $\widetilde{\mathbf{O}}_{b_o}(\delta_{o},\rho_{o}),$ which
induces another isomorphism
\begin{equation}\label{trivialization-2}
Q:\widetilde{H}\mid_{b_o}\to \widetilde{H}\mid_{b},\;\;\;\forall \;\;b\in \widetilde{\mathbf{O}}_{b_{o}}(\delta_{b_{o}},\rho_{b_{o}})
\end{equation}
\begin{equation}\label{trivialization-2}
\sum c_{\alpha}e_{\alpha}\mid_{b_o}\longmapsto \sum c_{\alpha}e_{\alpha}\mid_{b}.
\end{equation}
\v\n
Let $\rho_{\widetilde K_{b_o}}: G_{b_o}\rightarrow GL(\widetilde{K}_{b_o})$ be the natural linear representation, and let $\rho_{\mathbb{R}}: G_{b_o}\rightarrow GL(\mathbb{R}[G_{b_o}])$ be the standard representation. Both $\widetilde{K}_{b_o}$ and $\widetilde{H}\mid_{b_o}$ can be decomposed as sum of irreducible representations.
Without loss of generality we assume that $\rho_{\widetilde K_{b_o}}$ is an irreducible representation.
Let $\eta_1,...,\eta_l$ be a base of $\widetilde{K}_{b_o}$, let $\widetilde{H}\mid_{b_o}= \bigoplus_{i=1}^m E_i$ be the decomposition of irreducible representations such that $E_1$ has base $e_1,...,e_l.$ Define map $em(\eta_i)=e_i$, $i=1,...,l$. Thus we have map $p:\widetilde{H}\mid_{b_o}\rightarrow \widetilde{K}_{b_o}$ with $p\cdot em=id$.
\v
Let $\mathbb R^+=\{x\in \mathbb R|x\geq 0\}$ and $f_{\delta_{o},\rho_{o}}:\mathbb R^+\times \mathbb R^+\rightarrow \mathbb R^+$ be a smooth cut-off function such that
\[
f_{\delta_{o},\rho_{o}}(x,y)=\left\{
\begin{array}{ll}
1 \;\;\;\;\; on \;\;\{(x,y)|\;0\leq x\leq  \delta_o/3 ,\;0\leq y\leq \rho_o/3 \},    \\  \\
0 \;\;\;\;\;
on \;\;    \{(x,y)|\;x\geq  2\delta_o/3 \}\bigcup\{(x,y)|\; y\geq  2\rho_o/3 \}.
\end{array}
\right.
\]
We define a  cut-off function $\alpha_{b_o}: \widetilde{\mathbf{O}}_{b_o}(\delta_{b_{o}},\rho_{b_{o}}) \to [0, 1]$
by
\begin{equation}\label{cut-off-1}
\alpha_{b_o}(b)=f_{\delta_{b_{o}},\rho_{b_{o}}}(d_{\mathbf{A}}^{2}(a_o,a),\|h\|_{j_{a},k,2}^{2}).
\end{equation}
For any $\kappa\in \widetilde H\mid_{b}$ with $b\in \widetilde{\mathbf{O}}_{b_o}(\delta_{b_{o}},\rho_{b_{o}})$, in terms of
the local coordinate system $(\psi, \Psi)$, we define
\[
\mathfrak{i}(\kappa,b)_{b_o}=\left\{
\begin{array}{ll}
\alpha_{b_o}(b)P\circ p\circ Q^{-1} (\kappa) \;\;\;\;\; \mbox{ if }\;\|h\|_{j_{a},k,2} < \rho_{b_{o}}, \mbox{ and } d_{\mathbf{A}}^{2}(a_{o},a)<  \delta_{b_{o}}\\  \\
0 \;\;\;\;\;otherwise.
\end{array}
\right.
\]
\begin{lemma}\label{smooth of bundle map}
In the local coordinates $(\psi,\Psi)$ on $U$ and in  $\widetilde{\mathbf{O}}_{b_o}(\delta_o,\rho_o)$ the bundle map $\mathfrak{i}(\kappa,b)_{b_o}: \widetilde{\mathbf{F}}(\mathsf{k}_i)\rightarrow \widetilde{\mathcal{E}}$ is smooth with respect to $(\kappa, a, h)$.
\end{lemma}
\v\n
{\bf Proof.} By Lemma \ref{smoothness of norms} we immediately obtain that the cut-off function $\alpha_{b_o}(b)$ is a smooth function. Note that, in the local coordinates $(\psi,\Psi)$,
$P$, $p$ and $Q^{-1}$ are smooth. We conclude that $\mathfrak{i}(\kappa,b)_{b_o}$ is a smooth function of $(\kappa, a,h)$.
 $\Box$
\v\v
We can transfer the definition to other local coordinate system $(\psi', \Psi')$ and local model $\widetilde{\mathbf{O}}_{b_o}(\delta'_{b_{o}},\rho'_{b_{o}})$.
Suppose that in the coordinate system $(\psi, \Psi)$
$$b_o=(a_o, u_o),\;\;b=(a, v),\;\;v=\exp_{u_o}h,$$
and in the coordinate system $(\psi', \Psi')$
$$b'_o=(a_o', u_o'),\;\;b'=(a', v'),\;\;v'=\exp_{u'_o}h',\;\;where \;[b]=[b'].$$
We have
$$(\psi'\circ \psi^{-1}, \Psi'\circ \Psi^{-1})\cdot(a,v)=(a', v'),\;\;\;a'=\psi'\circ \psi^{-1}(a),\;\;v'=v\circ (\Psi'\circ \Psi^{-1})\mid_a.$$
$$(\psi'\circ \psi^{-1}, \Psi'\circ \Psi^{-1})\cdot(a_o,u_o)=(a_o', u_o'),\;\;\;a_o'=\psi'\circ \psi^{-1}(a),\;\;u_o'=u_o\circ (\Psi'\circ \Psi^{-1})\mid_{a_o}.$$
$(\psi'\circ \psi^{-1}, \Psi'\circ \Psi^{-1})$ send $e_{\alpha}$ to $e'_{\alpha}$. Then $(\Psi'\circ \Psi^{-1})\mid_{a}$
induces an isomorphism $\varphi_a:\widetilde{H}\mid_{(a,v)}\to \widetilde{H}'\mid_{(a',v')}$.
In $(\psi', \Psi')$ we have isomorphism $$Q':\widetilde{H}'\mid_{(a'_o,u'_o)}\to \widetilde{H}'\mid_{(a',v')},\;\;\forall b\in \widetilde{O}_{b_{o}}(\delta'_{b_{o}},\rho'_{b_{o}}),$$
$$Q'=\varphi_a \circ Q\circ \varphi_{a_o}^{-1}.$$
\v\n
We have chosen a finite dimensional subspace $\widetilde{K}_{(a,v)} \subset \widetilde{\E}|_{(a,v)}$ in $(\psi,\Psi)$. Denote $\vartheta_{a}=(\Psi'\circ \Psi^{-1})\mid_{a}.$ Define $\widetilde{K}'_{(a',v')}=\{\kappa\circ  d\vartheta_{a}^{-1}|\;\;\forall \;\kappa\in \widetilde{K}_{(a,v)}\}$.  Then $(\Psi'\circ \Psi^{-1})|_{a}$ induces a map
\begin{equation}\label{K isomorphic}
\phi_{a}:\widetilde{K}_{(a,v)}\to \widetilde{K}'_{(a',v')},\;\;\;
\phi_{a}(\kappa)= \kappa\circ  d\vartheta_{a}^{-1} ,\;\;\; \forall \kappa \in \widetilde{K}_{(a,v)}.
\end{equation}
Denote $\kappa'= \phi_{a}(\kappa)$.
Define
$$P':\widetilde{\E}'_{(a_o', u_o')}\rightarrow \widetilde{\E}'_{(a',v')},\;\;\mbox{ by }\; P'=\phi_a \circ P\circ \phi_{a_o}^{-1},$$
and
$$p':\widetilde{H}'\mid_{(a_o',u_o')}\rightarrow \widetilde{K}'_{(a'_o,u'_o)},\;\;\mbox{ by }\;p'=\phi_{a_o}\circ p\circ \varphi^{-1}_{a_o}.$$
$(\Psi'\circ \Psi^{-1})|_{a}$ also induces a map  $$\lambda_{a}: G_{(a_o,u_o)}\to G_{(a_o',u_o')}\;\;\;g\longmapsto g'=d\vartheta_{a}\circ g\circ (d\vartheta_{a})^{-1}.$$
It is easy to check that $\rho_{\widetilde K_{(a_o,u_o)}}: G_{(a_o,u_o)}\rightarrow GL(\widetilde{K}_{(a_o,u_o)})$ and
$\rho_{\widetilde K_{(a'_o,u'_o)}}: G_{(a'_o,u'_o)}\rightarrow GL(\widetilde{K}'_{(a'_o,u'_o)})$ are equivariant. Let
$$\eta'_i=\phi_{a}(\eta_i),\;\;e'_i=\varphi_{a}(e_i),\;\;em'(\eta'_i)=e'_i,\;\;i=1,2,...,l.$$
Then $em'(\widetilde{K}'_{(a'_o,u'_o)})=span\{e'_1,...e'_l\}\subset \widetilde{H}'\mid_{(a'_o,u'_o)}$. In the coordinate system $(\psi', \Psi')$
we define
\[
\mathfrak{i}(\kappa',b')_{b'_o}=\left\{
\begin{array}{ll}
\alpha_{b'_o}(b')P'\circ p'\circ (Q')^{-1} (\kappa') \;\;\;\;\; \mbox{ if }\;\|h\|_{j_{a'},k,2 }< \rho_{b'_{o}}, \mbox{ and } d_{\mathbf{A'}}^{2}(a'_{o},a')<  \delta_{b'_{o}}\\  \\
0 \;\;\;\;\;otherwise.
\end{array}
\right.
\]
We have
\begin{equation}\label{coord trans forbundle map}
\mathfrak{i}(\kappa',b')_{b'_o}=
\phi_a\circ\mathfrak{i}(\kappa,b)_{b_o}\circ\varphi^{-1}_a.
\end{equation}
If we choose three local coordinate systems $(\psi, \Psi)$,
$(\psi', \Psi')$ and $(\psi'', \Psi'')$,
since
$$
(\Psi\circ (\Psi'')^{-1})\circ(\Psi''\circ (\Psi')^{-1})\circ(\Psi'\circ \Psi^{-1})=Id,
$$
one can easily check that
\begin{equation}\label{coord trans equality}
\phi''_{a''}\phi'_{a'}\phi_a=Id,\;\;\; \varphi''_{a''}\varphi'_{a'}\varphi_a=Id.
\end{equation}
It follows from \eqref{coord trans forbundle map} and \eqref{coord trans equality} that the bundle map
$\mathfrak{i}:  {\mathbf{F}}(\mathsf{k}_i)\rightarrow \mathcal{E}$ is well defined. Obviously, $\mathfrak{i}([\kappa_{i},b])=[\mathfrak{i}(\kappa_{i},b)].$
\v\v
\begin{remark}\label{restri}
Let $(\psi', \Psi')$ be a local coordinate system in $\mathbf{O}_{[b'_o]}(\delta'_{[b'_{o}]},\rho'_{[b'_{o}]})\subset \mathbf{O}_{[b_o]}(\delta_{[b_{o}]},\rho_{[b_{o}]})$ such that $[b_o] \notin  \mathbf{O}_{[b'_o]}(\delta'_{[b'_{o}]},\rho'_{[b'_{o}]})$. The restriction of $[\mathfrak{i}(\kappa,b)_{b_o}]$ to $\mathbf{O}_{[b'_o]}(\delta'_{[b'_{o}]},\rho'_{[b'_{o}]})$ is a element in
$\E|_{\mathbf{O}_{[b'_o]}(\delta'_{[b'_{o}]},\rho'_{[b'_{o}]})}$. We can transfer it to $(\psi', \Psi')$ by \eqref{K isomorphic}.
\end{remark}

\v\v
\v\v
{\bf Case 2.} $[b_i]$ lies in a lower strata.
We choose $(\fs,\ft)$ coordinates. Put $\ft_{i}=e^{-2r_{i}-2\pi\tau_{i}}$, sometimes we use $(\fs, (\mathbf{r}))$ coordinates, where $(\mathbf{r})=((r_1,\tau_1),...,(r_{\mathbf{e}},\tau_{\mathbf{e}}))$.
Denote $b_{o}=b_i=(0,0,u)$, $\mathbf{F}(\mathsf{k}_i)=H(\fs,\ft)$, $\mathbf{F}(\mathsf{k}_i )\mid_{b_i}=H(0,0)$. We choose $|\fs|$, $|
\ft|$ small enough. In terms of $(\fs,\ft)$ we have an isomorphism

$$P:\bar{\E}_{b_o}\rightarrow \bar{\E}_{b},\;\;\forall \;\;b\in \widetilde{\mathbf{O}}_{b_{o}}(\delta_o,\rho_o).$$
Denote $\bar{H}=\{\zeta\mid_{\Sigma(R_0)} \mid \zeta \in \widetilde{H}\}$. Choosing a base $\{e_{\alpha}\}$ of the fiber $\bar{H}\mid_{b_o}$, by \eqref{glu_solu}
we can get a frame fields $\{e_{\alpha}((\mathbf{r}), a, h)\mid_{\Sigma(R_0)}\}$ for the bundle $\bar{H}$ over $\widetilde{\mathbf{O}}_{b_o}(\delta_{o},\rho_{o}).$
We have another isomorphism in the $(\fs,\ft)$ coordinates
$$ Q:\bar{H}(0,0)\to \bar{H}(\fs,\ft),\;\;\;\forall \;\;b\in \widetilde{\mathbf{O}}_{b_{o}}(\delta_o,\rho_o).$$
Denote $\mathbf{O}(\delta_o)=\{p\in \overline{\mathcal{M}}_{g,m+\mu}\mid d^2_{\mathsf{wp}}(0,p)< \delta_o\}$. Since $\overline{\mathcal{M}}_{g,m+\mu}$ has a natural effective orbifold structure, we can choose a smooth cut-off function in orbifold sense $\beta_{\delta_{o}}: \mathbf{O}(\delta_o) \to [0, 1]$ such that
$$\beta_{\delta_{o}}|_{\mathbf{O}(\delta_o/3)}=1,\;\;\;\;\beta_{\delta_{o}}|_{\mathbf{O}(\delta_o)\setminus\mathbf{O}(2\delta_o/3)}=0. $$

\v\n
We define a  cut-off function $\alpha_{b_o}: \widetilde{\mathbf{O}}_{b_o}(\delta_o,\rho_o) \to [0, 1]$
by
\begin{equation}\label{cut-off-2}
\alpha_{b_o}(b)=f_{\delta_o,\rho_o}(\beta_{\delta_{o}}(\fs,\ft) ,\|\beta_{R_{0}}h
\|_{j_{\fs,\ft},k,2}^{2}),\end{equation}
where $\beta_{R_{0}}$ is the function in \eqref{def_beta1}.
Using $\alpha_{b_o}(b)$ defined in \eqref{cut-off-2},
we can define the bundle map
$\mathfrak{i}:  {\mathbf{F}}(\mathsf{k}_i)\rightarrow \bar{\mathcal{E}}$ by
\[
\mathfrak{i}(\kappa,b)_{b_o}=\left\{
\begin{array}{ll}
\alpha_{b_o}(b)P\circ p\circ Q^{-1} (\kappa) \;\;\;\;\; \mbox{ if }\;\|h\|_{j_{a},k,2} < \rho_{b_{o}}, \mbox{ and }
\beta_{\delta_{o}}(\fs,\ft)<  \delta_{b_{o}}\\  \\
0 \;\;\;\;\;otherwise.
\end{array}
\right.
\]
For any fixed $(\mathbf{r})$, $\mathfrak{i}(\kappa,b)_{b_o}$
and $Q$ are smooth with respect to $(\mathbf{s},h)$ in the coordinates $(\mathbf{s},(\mathbf{r}))$. In order to study the smoothness with respect to $(\mathbf{r})$ we
note that $\mathfrak{i}(\kappa,b)_{b_o}$ is supported in $\Sigma(R_0)$. For any $v=\exp_{u_{(\mathbf{r})}}h$, we let   $$h^\circ=\left((h-\hat h_{0})(s_1,t_1)\beta_{1;2}(s_1),
(h-\hat h_{0})(s_2,t_2)\beta_{2;2}(s_{2})\right),$$
where
$$h_0=\int_{ {S}^1}h(r,t)dt.$$
Denote $v^\circ=\exp_u h^\circ$.
We can view $\bar{\E}\mid_v$ to be $\bar{\E}\mid_{v^\circ}$.
Then we view $P$ to be a family of operators in $\E$ over $W^{k,2}(\Sigma;u^{\ast}TM)$, where $\E\to W^{k,2}(\Sigma;u^{\ast}TM)$ is independent of $(\mathbf{r})$.
Consider the map
$$\mathfrak{i}(\kappa,b)_{b_o}\circ Q:  \bar{H}(0,0)\times \mathbf{A}\times D_{R_{0}}^{*}(0)\times W^{k,2}(\Sigma;u^{\ast}TM)
\rightarrow \E$$
$$\mathfrak{i}(\kappa,b)_{b_o}\circ Q(\kappa, \mathbf{s},(\mathbf{r}),h)  =\alpha_{b_o}(\mathbf{s},(\mathbf{r}),v)P\circ p(\kappa).
$$

 \v
\begin{lemma}\label{smooth of bundle map-1}
In the local coordinates $(\mathbf{s},(\mathbf{r}))$, the bundle map $\mathfrak{i}(\kappa,b)_{b_o}\circ Q$ is smooth with respect to $(\kappa, \mathbf{s}, (\mathbf{r}), h)$ in  $\widetilde{\mathbf{O}}_{b_o}(\delta_o,\rho_o)$.
\end{lemma}
\v\n
{\bf Proof.} $\alpha_{b_o}(\mathbf{s},(\mathbf{r}),v)$ is smooth
with respect to $(\mathbf{s}, (\mathbf{r}), h)$. For any $l\in Z^{+},$ denote $b_{\ft}=(\mathbf{s},\exp_{u}(h+\sum_{i=1}^{l}t_{l}h_{l} ))$
 and
$$
T^{l}(h;h_{1},\cdots,h_{l})=\nabla_{t_{1}}\cdots\nabla_{t_{l}}
\left.\left(P_{b_o,b_{\mathbf{t}}}\right)\right|_{\ft=\mathbf 0}
$$
By the same method as in the proof of Lemma 3.1 of \cite{LS-2} we can show that
$T^{l}(h;\cdot\cdot\cdot)$ is a bounded linear operator. The proof is complete.  \;\;$\Box$
\v
In \cite{LS-2} we proved
\begin{theorem}\label{smooth_line}
 Let $u:\Sigma\to M$ be a $(j,J)$-holomorphic map.	Let $\fc\in (0,1)$ be a fixed constant.    For any $0<\alpha<\frac{1}{100\fc}$, there exists positive  constants  $ \mathsf{d},R$ such that for any  $\zeta\in \ker D^{\widetilde{\mathbf L}}|_{b_{o}}$,   $(\kappa,\xi)\in \ker D \mathcal{S}_{(\kappa_{o},b_{o})}$ with    $$\|\zeta\|_{\mathcal W,k,2,\alpha}\leq \mathsf{d},\;\;\;\;\;\|(\kappa,\xi)\|< \mathsf{d},\;\;\;\;\;|\mathbf r|\geq R, $$ the following holds. Let $h_{(\mathbf r)}=\Pi_{2}(Glu_{\fs,(\mathbf r)}(\kappa,\xi))$ where $\Pi_{2}:  \widetilde{\mathbf F}_{b_{(\mathbf r)}}\times T_{u_{(\mathbf r)}} \widetilde{\mathcal B} \to T_{u_{(\mathbf r)}}\widetilde{ \mathcal B}$ denotes the projection.
 If $\exp_{u_{(\mathbf r)}}(h_{(\mathbf r)})$ is smooth with respect to $\fs$ and $(\mathbf r)$,	$Glu_{\mathbf{s},h_{(\mathbf r)},(\mathbf r)}(\zeta) $ is smooth with respect to $\fs,(\mathbf r)$ and $(\kappa,\xi)$.
\end{theorem}

\v
By {\bf Case 1}, {\bf Case 2} we have defined $\mathfrak{i}([\kappa_i,b])_i$ for all $i=1,...,\mathfrak{m}$. Set
$$\mathfrak{i}([\kappa,b])=\sum_{l=1}^{\mathfrak{m}} \mathfrak{i}([\kappa_l,b])_{l}\;\;
for\; any\; \kappa=(\kappa_1,...,\kappa_{\mathfrak{m}})\in \mathbf{F}\mid_{b}.$$ Then $\mathfrak{i}:\mathbf {F}\to \mathcal{E}$ is a bundle map.
We define a global regularization to be the bundle map $\mathcal{S}:\mathbf{F}\to  \E$
$$\mathcal{S}([\kappa,b])
=[\bar{\partial}_{j,J}v] + \mathfrak{i}([\kappa,b]).
$$
It is obvious that $D\mathcal{S}$ is surjective.
 Denote $\msp:\mathbf F\to \mathcal U $ by the projection of the bundle. Set
$$\mathbf{U}=\mathcal{S}^{-1}(0)|_{\msp^{-1}(\mathcal{U})}.$$
By restricting the bundle $\mathbf{F}$ to $\mathbf{U}$ we have a bundle $\msp:\mathbf{E}\to \mathbf U$ of finite rank with a canonical section $\sigma$ defined by
$$\sigma([(\kappa,b)])=( [((\kappa,b)  , \kappa)]),\;\;\;\;\forall \;[(\kappa,b)]\in \mathbf {U}.$$ We call
$$(\mathbf{U},\mathbf{E},\sigma),$$  a virtual neighborhood for $\overline{\mathcal{M}}_{g,m+\mu}(A)$.
\v

\section{Smoothness of the top strata}\label{top strata}

Denote by $\mathbf{U}^T$ the top strata of $\mathbf{U}$.
In this section we prove
\begin{theorem} \label{Smooth}
$\mathbf{U}^T$ is a smooth oriented, effective orbifold of dimension $\mathcal{N}= rank(\mathbf{F}) + ind \;D\mathcal{S}.$
\end{theorem}

The proof is divided into two steps, the subsections
\S\ref{smoothness-1} and \S\ref{smoothness-2}.
\v
\subsection{Smoothness}\label{smoothness-1}
\v
\v
Let $[(\kappa_{o},b_{o})]\in \mathbf{U}^T$. To simplify notations we consider the following case, for the general case the argument are the same. We assume that
$$[b_o]\in \mathbf{O}_{[b_{1}]}(2\delta_1/3,2\rho_1/3)\bigcap
\mathbf{O}_{[b_{2}]}(2\delta_2/3,2\rho_2/3)$$
and
$$[b_o]\notin \overline{\mathbf{O}}_{[b_{i}]}(2\delta_i/3,2\rho_i/3)\;\;\;\forall i=3,...,\mathfrak{m}.$$
We choose a local coordinate system  $(\psi, \Psi)$ for $\mathcal{Q}$ and local model $
\widetilde{\mathbf{O}}_{b_o}(\delta_o,\rho_o)/G_{b_o}$ around $b_o$. Let $b_o=(a_o,u)$, and let $\widetilde{\mathbf{U}}^T$ be the local expression of $\mathbf{U}^T$ in terms of $(\psi, \Psi)$. We choose $(\delta_o,\rho_o)$ so small that
$$\mathbf{O}_{[b_o]}(\delta_o,\rho_o)\notin \mathbf{O}_{[b_{i}]}(2\delta_i/3,2\rho_i/3)\;\;\;\forall i=3,...,\mathfrak{m}.$$
Then we only need to consider the bundles ${\mathbf{F}}(\mathsf{k}_{1})$ and ${\mathbf{F}}(\mathsf{k}_{2})$.
We consider two different cases.
\v
{\bf Case 1.} Both $[b_1]$ and $[b_2]$ lie in the top strata. By Remark \ref{finite covering} we may assume that
both $\mathbf{O}_{[b_{1}]}(2\delta_1/3,2\rho_1/3)$ and $\mathbf{O}_{[b_{2}]}(2\delta_2/3,2\rho_2/3)$ lie in the top strata. Let
$$b_1=(a_1, u_1)\;\;in\;(\psi_1, \Psi_1),\;\;\;\;b_2=(a_2, u_2)\;in \;(\psi_2, \Psi_2).$$
In terms of the coordinate system $(\psi, \Psi)$, let $b=(a,v)\in \widetilde{\mathbf{O}}_{b_o}(\delta_o,\rho_o)$. Suppose that, in the coordinate system $(\psi_1, \Psi_1)$,
$$[b']=[b],\; b'=(a', v'),\;\; v'=\exp_{u_1}h_1,$$
and in the coordinate system $(\psi_2, \Psi_2)$,
$$[b'']=[b],\;\;b''=(a'', v''),\;\; v''=\exp_{u_2}h_2.$$
The bundle maps are given respectively by
$$\mathfrak{i}(\kappa_1,b')_{b_1}=\alpha_{b_1}(b)P_1\circ p_1\circ Q_1^{-1}(\kappa_1): (\widetilde{H}_1)\mid_{b'}\to \widetilde{K}_1\mid_{b'}\;\;\;in\; (\psi_1, \Psi_1),$$
$$\mathfrak{i}(\kappa_2,b'')_{b_2}=\alpha_{b_2}(b)P_2\circ p_2\circ Q_2^{-1}(\kappa_2):(\widetilde{H}_2)\mid_{b''}\to \widetilde{K}_2\mid_{b''}\;in\;(\psi_2, \Psi_2),$$
where $P_1=P_{b_1,b'}$ in $(\psi_1, \Psi_1)$, $P_2=P_{b_2,b''}$ in $(\psi_2, \Psi_2)$. By Lemma \ref{smooth of bundle map}, $\mathfrak{i}(\kappa_1,b)_{b_1}$ in $(\psi_1, \Psi_1)$ ( resp. $\mathfrak{i}(\kappa_2,b)_{b_2}$ in $(\psi_2, \Psi_2)$ ) is smooth with respect to $(\kappa_1,b)$ ( resp. $(\kappa_2,b)$).
\v
We transfer from both the local coordinate systems $(\psi_1, \Psi_1)$ and $(\psi_2, \Psi_2)$ to the coordinates $(\psi, \Psi)$. We have
$$(\psi\circ \psi_1^{-1}, \Psi\circ \Psi_1^{-1})\cdot(a',v')=(a, v),\;\;\;a=\psi\circ \psi_1^{-1}(a'),\;\;v=v'\circ (\Psi\circ \Psi_1^{-1})\mid_{a'},$$
$$(\psi\circ \psi_2^{-1}, \Psi\circ \Psi_2^{-1})\cdot(a'',v'')=(a, v),\;\;\;a=\psi\circ \psi_2^{-1}(a''),\;\;v=v''\circ (\Psi\circ \Psi_2^{-1})\mid_{a''}.$$
The $(\psi\circ \psi_i^{-1}, \Psi\circ \Psi_i^{-1})$, $i=1,2$, induces maps
$$\phi^1_{a'}:\widetilde{K}_1\to \widetilde{K}^\diamond_1,\;\;
\phi^2_{a''}:\widetilde{K}_2\to \widetilde{K}^\diamond_2$$
$$\varphi^1_{a'}:\widetilde{H}_1\to \widetilde{H}^\diamond_1,\;\;\varphi^2_{a''}:\widetilde{H}_2\to \widetilde{H}^\diamond_2.$$
Put
$$\widetilde{H}^\diamond=(\widetilde{H}^\diamond_1)\mid_{b_o}\oplus (\widetilde{H}^\diamond_2)\mid_{b_o},\;\;\kappa=(\kappa_1,\kappa_2)\in \widetilde{H}^\diamond,\;\;(Q_1^\diamond \kappa_1, Q_2^\diamond \kappa_2):=Q^\diamond \kappa.$$
Here $\widetilde{H}^\diamond $, $\widetilde{K}^\diamond$ and $Q^\diamond$ denote the spaces and operator in $(\psi, \Psi)$.
By Remark \ref{restri} the bundle map in $(\psi, \Psi)$ becomes
$$\mathfrak{i}(\kappa,b)=\mathfrak{i}(\kappa_1,b')_{b_1}\circ  d\vartheta_{1}^{-1}
 +\mathfrak{i}(\kappa_2,b'')_{b_2}\circ  d\vartheta_{2}^{-1},
$$
where $\vartheta_{1}=(\Psi\circ \Psi_1^{-1})\mid_{a'}$, $\vartheta_{2}=(\Psi\circ \Psi_2^{-1})\mid_{a''}$. The key point is that $\Psi\circ \Psi_i^{-1}$, $i=1,2$, is a family of diffeomorphisms of $\Sigma$ depending on $a$.
For $v\in W^{k,2}$, $\frac{\p}{\p a}(v\circ\Psi\circ \Psi_i^{-1})$ is not in $W^{k,2}$. But for any fixed $a$ , $\mathfrak{i}(\kappa,b)$ and $Q_1^\diamond$, $Q_2^\diamond$ are smooth.
\v\n
Consider the map
$$ F_{(\kappa_o,b_o)}: \mathbf{A}\times \widetilde{H}^\diamond \times W^{k,2}(\Sigma;u^{\ast}TM)
\rightarrow W^{k-1,2}(u^{\ast}TM\otimes \wedge_{j_o}^{0,1})$$
$$F_{(\kappa_o,b_o)}(a,\kappa,h)=P_{b,b_o}\left(\bar{\partial}_{j_a,J}
v + \mathfrak{i}(Q^\diamond\kappa, b)
\right),$$ where  $b=(a,v)$, $v=\exp_{u}(h) $ for some
$h\in W^{k,2}(\Sigma, u^*TM)$. For any  $(a,\kappa ,h)\in F_{(\kappa_o,b_o)}^{-1}(0)$ we have
\begin{equation}\label{PDE}
 \bar{\partial}_{j_a,J}v + \mathfrak{i}(Q^\diamond\kappa,b)=0,
\end{equation}
 where $b=(a,v)$.
For any fixed $a$, it follows from the standard elliptic estimates and the smoothness of $\mathfrak{i}$ that $v\in C^{\infty}(\Sigma,M).$ Then by Lemma \ref{smooth of bundle map} and the smoothness of the frame field $e_{\alpha}$ we
conclude that $\mathfrak{i}\mid_v$ and $Q^\diamond\mid_v$ are smooth with respect to $(a,\kappa,h)$.
It is easy to see that
$F_{(\kappa_o,b_o)}(a, \kappa,h)$ is smooth with respect to $(a,\kappa,h)$. Then we use the implicity theorem with parameter $a$
to conclude that $v$ is smooth with respect to $(a,\kappa,h)$. It follows that $\widetilde{\mathbf{U}}^T\bigcap \pi^*\widetilde{\mathbf{O}}_{b_o}(\delta_o,\rho_o)$ is smooth, where $\pi:\widetilde{\mathbf{U}}^T\to \widetilde{ \mathcal B}$ is the projection.

\v
\v
{\bf Case 2.} $[b_2]$ lies in the top strata, $[b_1]$ lies in a lower strata. Without loss of generality we assume that $b_2=(\Sigma, j, {\mathbf{y}}, u)$, where $\Sigma$ has one node $q$, $\mathbf{s}_o\in \mathbf{A}=\mathbf{A_1\times \mathbf{A}}_2$. We glue $\Sigma$ at $q$ with gluing parameter $(r)$. We have bundle maps $\mathfrak{i}(\kappa_1,b)_{b_1}=\alpha_{b_1}(b)P_1\circ p_1\circ Q_1^{-1}(\kappa_1)$ and $\mathfrak{i}(\kappa_2,b')_{b_2}=\alpha_{b_2}(b')P_2\circ p_2\circ Q_2^{-1}(\kappa_2)$.
Then we transfer to the coordinates $(\psi, \Psi)$, we choose $(\mathbf{s}, \mathbf{t})$-coordinates. We use Lemma \ref{smooth of bundle map-1} and the same method as in {\bf Case 1} to prove that $v$ is smooth with respect to $(\mathbf{s}, (\mathbf{r}),\kappa,h)$. Then we use Lemma \ref{smooth_line} to prove that $Q_1$ is smooth with respect to $(\mathbf{s}, (\mathbf{r}),\kappa,h)$. Then we can prove
the smoothness of $\mathbf{U}^T\bigcap \pi^*\widetilde{\mathbf{O}}_{b_o}(\delta_o,\rho_o)$.
\v
The proof of the orientation of $\mathbf{U}^T$ is standard, we omit here.

\subsection{The oribifold structure}\label{smoothness-2}

We introduce a notation. For any $(\kappa_o, b_o)\in \mathbf{U}$ we choose a local coordinate system  $(\psi, \Psi)$ on $U\ni a_o$ and local model $\widetilde{\mathbf{O}}_{b_o}(\delta_o,\rho_o)/G_{b_o}$.
Set
$$\widetilde{\mathbf{U}}_{\kappa_o,b_o}(\varepsilon,\delta_o,\rho_o) =\left\{(\kappa, b)\in \widetilde{\mathbf{U}}\mid |\kappa-\kappa_{o}|_{\mathbf{h}}<\varepsilon,
b\in \widetilde{\mathbf{O}}_{b_o}(\delta_o,\rho_o)\right\},$$
$$\mathbf{U}_{\kappa_{o},b_{o}}(\varepsilon,\delta_o,\rho_o)
=\widetilde{\mathbf{U}}_{\kappa_{o},b_{o}}(\varepsilon,\delta_o,\rho_o)/G_{\kappa_o,b_o},$$
where $G_{\kappa_o,b_o}$ is the isotropy group at $(\kappa_o,b_o)$. For any $(\kappa,b)\in
\widetilde{\mathbf{U}}_{\kappa_{o},b_{o}}(\varepsilon,\delta_o,\rho_o)$
denote by $G_{\kappa,b}$ the isotropy group at $(\kappa,b)$.  Any element $\varphi\in G_{\kappa,b}$ satisfies
$\varphi^*(\kappa, b)=(\kappa, b).$
It follows that $G_{\kappa,b}$ is a subgroup of $\mathbf{G}_a$.

\v

 \begin{lemma}\label{orbi structure} Let $[(\kappa_{o},b_{o})]\in \mathbf{U}^T$.
Suppose that $\widetilde{\mathbf{U}}_{\kappa_{o},b_{o}}(\varepsilon,\delta_o,\rho_o)\subset \mathbf{U}^T.$ The following hold
 \begin{itemize}
\item[(1)] For any $p\in\widetilde{\mathbf{U}}_{\kappa_{o},b_{o}}(\varepsilon,\delta_o,\rho_o)$
let $G_{p}$ be the isotropy group at $p$, then $im(G_{p})$ is a subgroup of $G_{\kappa_o,b_o}$.
\item[(2)] Let $p\in\widetilde{\mathbf{U}}_{\kappa_{o},b_{o}}(\varepsilon,\delta_o,\rho_o)$
be an arbitrary point with isotropy group $G_p$, then there is a $G_p$-invariant neighborhood $O(p)\subset \widetilde{\mathbf{U}}_{\kappa_{o},b_{o}}(\varepsilon,\delta_o,\rho_o)$
such that for any $q\in O(p)$, $im(G_{q})$ is a subgroup of $G_{p}$, where $G_p$, $G_q$ denotes the isotropy groups at $p$ and $q$ respectively.
\end{itemize}
\end{lemma}
\v\n
{\bf Proof:} We only prove (1), the proof of (2) is similar. Denote $b_o=(a_o,u)$.
If the lemma not true, we can find a sequence $(\kappa_i,b_i)=(\kappa_i,a_i,u_i)\in \widetilde{\mathbf{U}}_{\kappa_{o},b_{o}}(\varepsilon,\delta_o,\rho_o)$
such that
\v
{(1)} $\delta_i\to 0$, $\rho_i\to 0$, $\kappa_i\to \kappa_o$,
\v
{(2)} $im(G_{\kappa_i,b_i})$ is not a subgroup of $G_{\kappa_o,b_o}$.
\v\n
It is obvious that $G_{\kappa_o,b_o}$ is a subgroup of $\mathbf{G}_{a_o}$, $G_{\kappa_i,b_i}$ is a subgroup of $\mathbf{G}_{a_i}$ and $\mathbf{G}_{a_i}$ can be imbedded into $\mathbf{G}_{a_o}$ as a subgroup for $i$ large enough. So we can view $im(G_{\kappa_i,b_i})$ as a subgroup of $\mathbf{G}_{a_o}$,  By choosing subsequence we may assume that $im(G_{\kappa_i,b_i})$ convergies to a subgroup $G_{\kappa,b}$ of $\mathbf{G}_{a_o}$ and $im(G_{\kappa_i,b_i})\cdot u_i$ converges to $im(G_{\kappa,b})\cdot u$ and $u_i$ converges to $u$ in $W^{k,2}$. By Sobolev imbedding theorem and elliptic estimates we have $im(G_{\kappa_i,b_i})\cdot (\kappa_i,u_i)$ converges to $im(G_{\kappa,b})\cdot (\kappa_o,u)$, $(\kappa_i,u_i)$ converges to $(\kappa_o,u)$ in $C^{\ell}$ for any $\ell>1$. It follows that $im(G_{\kappa,b})\subset G_{\kappa_o,b_o}$. Since there are only finite many subgroups of $G_{a_o}$, for $i$ large enough we have $im(G_{\kappa_i,b_i})=G_{\kappa,b}$. So $G_{\kappa_i,b_i}$ can be imbedded into $G_{\kappa_o,b_o}$ as a subgroup for $i$ large enough. We get a contradiction. $\Box$
\v
As corollary of Lemma \ref{orbi structure} we conclude that $\mathbf{U}^T$ is an orbifold. Since $(g,n)\ne (1,1), (2,0)$,
$\mathbf{U}^T$ has the structure of an effective orbifold.

\v
Combination of the subsections \S\ref{smoothness-1},
\S\ref{smoothness-2} give us the proof of Theorem \ref{Smooth}.
\v
\subsection{A metric on $\mathbf{E}$}\label{a metric}
In this section we construct a metric on $\mathbf{E}|_{\mathbf U_{\epsilon}}.$
 By the compactness of $\mathbf{U}_{2\varepsilon}$ we may find finite many points $(\kappa_1,b_1),...,(\kappa_{\mathbf{n}}, b_{\mathbf{n}})\in \mathbf{U}_{\varepsilon}$ such that
\begin{itemize}
\item $\{\mathbf{U}_{[(\kappa_{\mathbf{a}},b_{\mathbf{a}})]}(\varepsilon_{\mathbf{a}},\delta_{\mathbf{a}},\rho_{\mathbf{a}})
, \;1\leq \mathbf{a}\leq \mathbf{n}\}$ is a covering of $\mathbf{U}_{2\varepsilon}$.
\item For any $\mathbf{a}\in \{1,...,\mathbf{n}\}$ there is $i_{\mathbf{a}}\in \{1,...,\mathbf{m}\}$ such that
$$\mathsf{p}({\mathbf{U}}_{[(\kappa_{\mathbf{a}},b_{\mathbf{a}})]}(\varepsilon_{\mathbf{a}},\delta_{\mathbf{a}},\rho_{\mathbf{a}}))
\subset  {\mathbf{O}}_{b_{i_{\mathbf a}}}(\delta_{i_{\mathbf a}},\rho_{i_{\mathbf a}}),$$  where ${\mathbf{O}}_{b_{i_{\mathbf a}}}(\delta_{i_{\mathbf a}},\rho_{i_{\mathbf a}})$ is as in subsection \S\ref{finite rank orbi-bundle},
	\item $ \widetilde{{\mathbf{U}}}_{(\kappa_{\mathbf{a}},b_{\mathbf{a}})}(\varepsilon_{\mathbf{a}},\delta_{\mathbf{a}},\rho_{\mathbf{a}})\subset \widetilde{\mathbf U}^{T}$ for all $1\leq \mathbf{a}\leq \mathbf{n}_{t}$.
\end{itemize}

 Let $\{e^{i_{\mathbf a}}_{\alpha}\}_{1\leq\alpha\leq \mathsf r}$ be a local smooth frame field of $\mathbf F$ over $\mathbf{O}_{b_{i_{\mathbf{a}}}}(\delta_{i_{\mathbf{a}}},\rho_{i_{\mathbf{a}}})$ as in section \S\ref{global_regu}. Let $\mathsf{p}: \mathbf{U}\to \mathcal{U}$ denote the projection. Denote $e^{\mathbf a}_{\alpha}=\mathsf{p}^{*}e^{i_{\mathbf a}}_{\alpha}|_{{\mathbf{U}}_{[(\kappa_{\mathbf{a}},b_{\mathbf{a}})]}(\varepsilon_{\mathbf{a}},\delta_{\mathbf{a}},\rho_{\mathbf{a}}) }.$
Then we have a smooth frame field $\{e_{\alpha}^{\mathbf a}\}_{1\leq\alpha\leq \mathsf r}$ of $\mathbf{E}$ over ${\mathbf{U}}_{[(\kappa_{\mathbf{a}},b_{\mathbf{a}})]}(\varepsilon_{\mathbf{a}},\delta_{\mathbf{a}},\rho_{\mathbf{a}})$, where $\mathsf r$ denotes the rank of $\mathbf{E}$.
We define a local metric $h_{\mathbf{a}}$ on $\mathbf{E}|_
{\mathbf{U}_{[(\kappa_{\mathbf{a}},b_{\mathbf{a}})]}(\varepsilon_{\mathbf{a}},\delta_{\mathbf{a}},\rho_{\mathbf{a}})}$ by
$$
h_{\mathbf a}(e^{\mathbf a}_{\alpha},e^{\mathbf a}_{\beta})=\delta_{\alpha \beta}.
$$

\v
Now we choose smooth cutoff functions $\mathbf{\Gamma}'$ as follows.
Let $(\kappa_o, b_o)$ be one of $(\kappa_1,b_1),...,(\kappa_{\mathbf{n}}, b_{\mathbf{n}})$.
We consider two cases.
\v
{\bf (1).} \;$(\kappa_o,b_o)$ lies in $\widetilde{\mathbf{U}}^T$. We define a cut-off function $\alpha_{b_o}: \widetilde{\mathbf{O}}_{b_o}(\delta_{b_{o}},\rho_{b_{o}}) \to [0, 1]$
by \eqref{cut-off-1} and let $\mathbf{\Gamma}'_{o}=\mathsf{p}^*\alpha_{b_o}(b).$ \v
{\bf (2).} \;$(\kappa_o,b_o)$ lies in a lower strata.
We define a  cut-off function $\alpha_{b_o}: \widetilde{\mathbf{O}}_{b_o}(\delta_o,\rho_o) \to [0, 1]$
by \eqref{cut-off-2}
and let $\mathbf{\Gamma}'_{o}=\mathsf{p}^*\alpha_{b_o}(b).$

\v
Thus we have $\mathbf{\Gamma}'_{\mathbf{a}}$ for every $1\leq \mathbf{a}\leq \mathbf{n}.$ Set $$\mathbf{\Gamma}_{\mathbf{a}}=\frac{\mathbf{\Gamma}'_{\mathbf{a}}}{\sum_{l=1}^{\mathbf{n} } \mathbf{\Gamma}'_{l}}.$$ Then $\sum \mathbf{\Gamma}_{\mathbf{a}}=1$ and $\mathbf{\Gamma}_{\mathbf{a}}$ is smooth on $\mathbf{U}^T_{\epsilon}$ in orbifold sense. We define a metric $\mathbf{h}$ on $\mathbf E$ over $\mathbf{U}_{\varepsilon}$ by
$$\mathbf{h}=\sum_{\mathbf{a}=1}^{\mathbf{n}}\mathbf{\Gamma}_{\mathbf{a}}  h_{\mathbf{a}}.
$$
\v
 We define a connection on $\mathbf E$ as follows.
Let $\{e^{\mathbf a}_{\alpha}\}_{1\leq\alpha\leq \mathsf r}$ be a local smooth frame field of $\mathbf E$ over ${\mathbf{U}}_{[(\kappa_{\mathbf{a}},b_{\mathbf{a}})]}(\varepsilon_{\mathbf{a}},\delta_{\mathbf{a}},\rho_{\mathbf{a}})$ as above. Consider the  Gram-Schmidt process with respect to the metric $\mathbf h$ and denote by $\hat{e}^{{\mathbf a}}_{1}, ..., \hat{e}^{{\mathbf a}}_{\mathsf r}$ the Gram-Schmidt orthonormalization of $\{e^{\mathbf a}_{\alpha}\}$.  We define a local connection $\nabla^{\mathbf a}$ by
$$
\nabla^{\mathbf a}\hat{e}^{\mathbf a}_{\alpha}=0,\;\;\;\;\;\alpha=1,\cdots,\mathsf r.
$$
For any section $e\in \mathbf E|_{{\mathbf{U}}_{\epsilon}},$ we define
\begin{equation}\label{def_conn}
\nabla e=\sum \mathbf{\Gamma}_{{\mathbf a}}\nabla^{\mathbf a}(e|_{{\mathbf{U}}_{[(\kappa_{\mathbf{a}},b_{\mathbf{a}})]}(\varepsilon_{\mathbf{a}},\delta_{\mathbf{a}},\rho_{\mathbf{a}})}).
\end{equation}
It is easy to see that $\nabla$ is a compatible connection of the metric $\mathbf{h}$. Denote $$\nabla \hat e^{\mathbf a}_{\alpha}=\sum_{\beta} \omega_{\alpha \beta}^{\mathbf {a}}\hat e^{\mathbf a}_{\beta},\;\;\;\nabla^2\hat e^{\mathbf a}_{\alpha}=\sum  \Omega_{\alpha\beta}^{\mathbf a}\hat e^{\mathbf a}_{\beta}.$$
For any ${\mathbf{U}}_{[(\kappa_{\mathbf{a}},b_{\mathbf{a}})]}(\varepsilon_{\mathbf{a}},\delta_{\mathbf{a}},\rho_{\mathbf{a}})\bigcap  {\mathbf{U}}_{[(\kappa_{\mathbf{c}},b_{\mathbf{c}})]}(\varepsilon_{\mathbf{c}},\delta_{\mathbf{c}},\rho_{\mathbf{c}})\neq \emptyset,$ let $(\hat a^{\mathbf a\mathbf c}_{\alpha\beta})_{1\leq \alpha,\beta\leq \mathsf r}$ be   functions such that $$\hat e^{\mathbf a}_{\alpha}=\sum_{\beta=1}^{\mathsf r} \hat a^{\mathbf a\mathbf c}_{\alpha \beta}\hat e^{\mathbf c}_{\beta},\alpha=1,\cdots,\mathsf r.$$   It is easy to see that
\begin{equation}\label{local_omega}
\omega^{\mathbf a}_{\alpha\beta}=\sum_{\mathbf c}\sum_{\beta=1}^{\mathsf r} \alpha_{b_{\mathbf c}}d\hat a^{\mathbf a\mathbf c}_{\alpha \gamma}\hat a^{\mathbf c \mathbf a}_{\gamma \beta}.
\end{equation}
We get a  metric $\mathbf h$ and a connection $\nabla$ in $\mathbf E$ over $\mathbf{U}_{\varepsilon}$.
\v

\section{Gluing estimates}\label{estimates}
\v
\subsection{Gluing maps}\label{gluing map}

We have two kinds of gluing maps. \\
{\bf Case 1.}{ Gluing maps in a holomorphic cascade.}
Let $\Sigma$ be a marked nodal Riemann surfaces. Suppose that $\Sigma$
has nodes $p_{1},\cdots,p_{\mathbf{e}}$ and  marked points $y_{1},\cdots,y_{m}$. We choose local coordinate system $\mathbf{A}$. Let $u:\Sigma\to M$ be perturbed $J$-holomorphic map. We glue $\Sigma$ and $u$ at each node with gluing parameters $(\mathbf{r})$ to get $\Sigma_{(\mathbf{r})}$ and the pregluing map $u_{(\mathbf{r})}: \Sigma_{(\mathbf{r})}\to M$. Set
 $$\ft_{i}=e^{-2r_{i}-2\pi\tau_{i}}, \;\;\;|\mathbf{r}|=min\{r_1,...,r_{\mathbf{e}}\},\;\;\;b_{(\mathbf r)}:=(0,\mathbf{(r)}, u_{(\mathbf r)}).$$
The following lemma is proved in
\cite{LS-1}.
\begin{lemma}\label{isomor of ker}
For $|\mathbf{r}|>R_0$ there is an isomorphism
$$I_{(\mathbf r)}: \ker D \mathcal{S}_{( \kappa_{o}, b_{o})}\longrightarrow \ker D \mathcal{S}_{( \kappa_{o},b_{(\mathbf r)})}.$$
\end{lemma}
\v\n
Using Theorem 5.3 in \cite{LS-1} and the implicit function theorem with parameters we immediately obtain
\begin{lemma}\label{gluing map-1} There are constant $\varepsilon>0$, $R_0>0$ and a neighborhood $O_1\subset \mathbf{A}$ of $\mathbf{s}_o$ and a neighborhood $O$ of $0$ in $ker D{\mathcal S}_{(\kappa_{o},b_{o})}$
 such that
$$glu_{(\kappa_{o},b_{o})}: O_1\times(\mathbb{D}_{\mathbf{c}}^*(0))^{\mathbf{e}}\times O\rightarrow glu_{(\kappa_{o},b_{o})}(O_1\times(\mathbb{D}_{\mathbf{c}}^*)^{\mathbf{e}}\times O)\subset\mathbf{U}^T$$
 is an orientation preserving local diffeomorphisms, where $$\mathbb{D}^{*}_{\mathbf{c}}(0):=\{\mathbf{t}\mid 0<|\mathbf{t}|<\mathbf{c}\},\;\;\;\mathbf{c}=e^{-2R_0}.$$  \end{lemma}
\n
Denote
$$ Glu_{\fs,(\mathbf{r})}=I_{(\mathbf{r})} + Q_{(\kappa_{o},b_{(\mathbf{r})})}\circ f_{\fs,(\mathbf{r})}\circ I_{(\mathbf{r})}.$$

{\bf Case 2.} Gluing maps between different holomorphic cascades.
\v
Let $b=(b_1, b_2)\in \mathcal{M}^{0,1}$ be as in \S\ref{relative_node_M_M-1}. For every puncture point $q_j$ there are constants $(\ell_{1j}, \theta_{1j0})$ and
$(\ell_{2j}, \theta_{2j0})$ such that \eqref{gluing_local_relative_node-1} hold. Since there exists a $\mathbb C^{*}$ action on ${\mathbb{R}}\times \widetilde{M},$ we can choose the coordinates  $(a_{2},\theta_{2})$ such that
$\ell_{11}=\ell_{21},\;\theta_{110}=\theta_{210},$ that is we fix a slice for  ${\mathcal{M}}_{A}(\mathbb P(\mathcal N\oplus \mathbb C) ,g,m+\mu^{+}+\mu^{-},{\bf k}^+, {\bf k}^-,\nu)/\mathbb{C}^*.$ We associate a point
$$\bar{\mathbf{t}}_j=\exp\{(\ell_{1j}-\ell_{2j})+ 2\pi\sqrt{-1}(\theta_{2j0}-\theta_{1j0})\}.$$
Put
$$\mathbb{D}^*:=\{\mathbf{t}^*_j \mid (\mathbf{t}^*_j)^{k_j}=\bar{\mathbf{t}}_j\}.$$
For each puncture point $q_j$ we have $\mathbb{H}_{q_j}$.
Set $\mathds{H}:=\bigoplus_{j=1}^l \mathbb{H}_{q_j}$, $\mathds{H}^*:=\bigoplus_{j=1}^l \mathbb{H}^*_{q_j}$. Denote
$$E_2^*:=\left\{(\kappa_{20},h_2+\hat{h}_{20})\mid D \mathcal{S}_{u_2}(\kappa_{20},h_2 + \hat{h}_{20})=0,\;\;h_{20}\in \mathds{H}\;\;with\;\;
 a(h^1_{20})=0,\theta(h^1_{20})=0,\right\}$$
$$Ker \mathbb{D} \mathcal{S}^*_{(\kappa_{o},b)}:=E_1\bigoplus_{\bh^*}E^*_2=\left\{(\kappa_{10},(h_1,h_{10}), (\kappa_{20},h_2,h_{20}))\in E_1\oplus E^*_2 \mid \pi_{*}h_{10}=\pi_{*}h_{20}\in \mathds{H}^*\right\},$$
$$Ker \mathbb{D} \mathcal{S}_{(\kappa_{o},b)}:=\left\{((\kappa_{10},h_1,h_{10}), (\kappa_{20},h_2,h_{20}))\in E_1\oplus E_2\mid \pi_{*}h_{10}=\pi_{*}h_{20}\in \mathds{H}^*\right\}.$$
Obvioulsy, $Ker D \mathcal {S}_{(\kappa_{o},b)}$ is a subspace of $Ker \mathbb D \mathcal {S}_{(\kappa_{o},b)}$. We choose a Euclidean metric $<<\cdot>>$ on $Ker \mathbb D \mathcal {S}_{(\kappa_{o},b)}$.
Let $\mathbb{E}$ be a subspace of $Ker \mathbb D \mathcal {S}_{(\kappa_{o},b)}$ such that
$$Ker \mathbb D \mathcal S_{(\kappa_{o},b)}=\mathbb E\oplus Ker D \mathcal S_{(\kappa_{o},b)}.$$
It is easy to see that  $\dim \mathbb E=2l.$ For any
$((h_1,h_{10}), (h_2,h_{20}))\in \mathbb E\setminus \{0\}$ there is unique $\left((c^1_1, c^1_2),...,(c^l_1, c^l_2)\right)\in \mathds{H}$, $\left((c^1_1, c^1_2),...,(c^l_1, c^l_2)\right)\ne (0,...,0)$,
such that
$$\left(a(h^i_{10}-h^i_{20}), \theta(h^i_{10}-h^i_{20})\right)=(c^i_1,c^i_2),\;\; 1\leq i\leq l.$$
We fix a basis $\{\mathbf{f}_{1},\mathbf{l}_1\cdots,\mathbf{f}_{l},\mathbf{l}_l\}$ of $\mathbb E$ such that $\mathbf{f}_{i}$ corresponding to $(c^i_1,c^i_2)=(1,0)$, $\mathbf{l}_i$ corresponding to $(c^i_1,c^i_2)=(0,1)$ and
$(c^j_1,c^j_2)=(0,0)$ for all $j\ne i$. With respect to this base $\left(a(h^i_{10}-h^i_{20}), \theta(h^i_{10}-h^i_{20})\right)$ is the coordinate system of $\mathbb E$. We can also view  $\mathbf{t}^*=(\mathbf{t}^*_1,...,\mathbf{t}^*_l)$ as a coordinate system of $\mathbb E$. It is easy to see that there is an isomorphism
$$\psi:(\mathbb D^{*})^l\times Ker D\mathcal{S}_{(\kappa_{o},b)}\to Ker \mathbb{D} \mathcal{S}_{(\kappa_{o},b)}.$$
It is obvious that there is an isomorphism
$$\eta: \mathbb{C}^*\times Ker \mathbb{D}\mathcal{S}^*_{(\kappa_{o},b)} \to Ker \mathbb{D} \mathcal{S}_{(\kappa_{o},b)}.$$
We use the gluing parameters $(\mathbf{r})$ to glue
at $q_{1},q_{2},...,q_{l}$ to get $\Sigma_{(\mathbf{r})}$ and $u_{(\mathbf{r})}$
as in \S\ref{relative_node_M_M-1}.
Using Theorem 5.3 in \cite{LS-1} and the implicit function theorem with parameters we immediately obtain
\begin{lemma}\label{gluing map-2} There are constant $\varepsilon>0$, $R_0>0$ and a neighborhood $O_1\subset \mathbf{A}$ of $\mathbf{s}_o$ and a neighborhood $O$ of $0$ in $ker \mathbb D{\mathcal S}^*_{(\kappa_{o},b_{o})}$
 such that
$$glu: O_1\times \mathbb{C}^{*}_{\epsilon}(0)\times O\rightarrow glu(O_1 \times  \mathbb{C}^{*}_{\epsilon}(0)\times O)\subset\mathbf{U}^T$$
 is an orientation preserving local diffeomorphisms,
where $$\mathbb{C}^{*}_{\epsilon}(0):=\{z=e^{r+2\pi\sqrt{-1}\vartheta}\mid 0<|z|<\epsilon\}.$$
\end{lemma}
\v\n
This lemma can be immediately generalize to the case gluing several cascades.

\v\v
\begin{remark}\label{isotropy group}
Let $G_{(\kappa_{o},b_o)}$ be the isotropy group at $(\kappa_{o},b_o)$.
It is easy to check that  the operator  $D\mathcal{S}_{(\kappa_o,b_o)}$ is $G_{(\kappa_{o},b_{o})}$-equivariant. Then we may choose a $G_{(\kappa_{o},b_{o})}$-equivariant right inverse $Q_{(\kappa_o,b_o)}$.
$G_{(\kappa_{o},b_{o})}$ acts on $ker D\mathcal{S}_{(\kappa_o,b_o)}$ in a natural way. So we have  $G_{(\kappa_{o},b_o)}$-equivariant versions of Lemma \ref{isomor of ker}, Lemma \ref{gluing map-1}, Lemma \ref{gluing map-2} ( see \cite{LS-2}, \cite{LS-3}.)
\v
\end{remark}

\subsection{Exponential decay of gluing maps}

The following theorem is proved in \cite{LS-1}.
\begin{theorem}\label{coordinate_decay-2} Let $l\in \mathbb Z^+$ be a fixed integer.
There exists positive  constants  $\mathsf{C}_{l}, \mathsf{d}, R_{0}$ such that for any $(\kappa,\xi)\in \ker D \mathcal{S}_{(\kappa_{o},b_{o})}$ with $\|(\kappa,\xi)\|< \mathsf{d}$ and for any $X_{i}\in \{\frac{\p}{\p r_{i}},\frac{\p}{\p \tau_{i}}\},i=1,\cdots,\mathbf{e}$, restricting to the compact set $\Sigma(R_0)$, the following estimate hold
$$ \left\|X_{i} \left(Glu_{\fs,(\mathbf{r})}(\kappa,\xi) \right)  \right\|_{C^{l}(\Sigma(R_{0}))}
	\leq  \mathsf{C}_{l}e^{-(\fc-5\alpha)\tfrac{r_{i}}{4} },$$
$$ \left\|X_{i}X_{j}\left(Glu_{\fs,(\mathbf{r})}(\kappa,\xi) \right)  \right\|_{C^{l}(\Sigma(R_{0}))}
	\leq  \mathsf{C}_{l}e^{-(\fc-5\alpha)\tfrac{r_{i}+r_{j}}{4} },$$
$1\leq i\neq j\leq\mathbf{e},$
	for any $\mathbf{s}\in \bigotimes_{l=1}^{\iota} O_l$ when $|r|$ big enough.
\end{theorem}

\subsection{Estimates of exponential decay of the line bundle}

The following theorem is proved in \cite{LS-3}

\begin{theorem}\label{thm_est_mix_deri}
Let $l\in \mathbb Z^+$ be a fixed integer.
 Let $u:\Sigma\to M$ be a $(j,J)$-holomorphic map.	Let $\fc\in (0,1)$ be a fixed constant.    For any $0<\alpha<\frac{1}{100\fc}$, there exists positive  constants  $\mathsf C_{l}, \mathsf{d},R$ such that for any  $\zeta\in \ker D^{\widetilde{\mathbf L}}|_{b_{o}}$,   $(\kappa,\xi)\in \ker D \mathcal{S}_{(\kappa_{o},b_{o})}$ with    $$\|\zeta\|_{\mathcal W,k,2,\alpha}\leq \mathsf{d},\;\;\;\;\;\|(\kappa,\xi)\|< \mathsf{d},\;\;\;\;\;|\mathbf r|\geq R, $$ restricting to the compact set $\Sigma(R_0)$, the following estimate hold.
	\begin{equation}\label{eqn_est_L_1st}
	\left\|X_{i}  \left( Glu_{\mathbf{s},h_{(\mathbf r)},(\mathbf{r})}^{\widetilde{\mathbf L}}(\zeta) \right) 	\right\|_{C^{l}(\Sigma(R_{0}))}
	\leq \mathsf C_{l} e^{-(\fc-5\alpha)\tfrac{r_{i}}{4} } ,
	\end{equation}
	\begin{align} \label{eqn_est_L_2nd}
	&\left\|X_{i}X_{j} \left( Glu^{\widetilde{\mathbf L}}_{\mathbf{s},h_{(\mathbf r)},(\mathbf{r})}(\zeta) \right)
	\right\|_{C^{l}(\Sigma(R_{0}))} \leq   \mathsf C_{l} e^{-(\fc-5\alpha)\tfrac{r_{i}+r_{j}}{4} }
	\end{align}
	for any $X_{i}\in \{\frac{\p}{\p r_{i}},\frac{\p}{\p \tau_{i}}\},i=1,\cdots,\mathbf{e}$,    $\mathbf{s}\in \bigotimes_{l=1}^{\iota} O_l$ and any $1\leq i\neq j\leq \mathbf{e}$, where $h_{(\mathbf r)}=\Pi_{2}(Glu_{\fs,(\mathbf r)}(\kappa,\xi))$ and $\Pi_{2}:  \widetilde{\mathbf F}_{b_{(\mathbf r)}}\times T_{u_{(\mathbf r)}} \widetilde{\mathcal B} \to T_{u_{(\mathbf r)}}\widetilde{ \mathcal B}$ denotes the projection.
\end{theorem}

\v
\subsection{Estimates of Thom forms}\label{est_Thom_E}

\v
We estimate the derivatives  of the metric $\mathbf h$ near the boundary of $\mathbf F|_{\mathbf U^{T}}$. Let $(\kappa_{o},b_{o})$ be one of $\{(\kappa_{\mathbf{a}},b_{\mathbf{a}}), \mathbf{a}=\mathbf{n}_{t}+1,\cdots,\mathbf{n}_{c}\}$ and  $b_{o}=(a_{o},u).$ We use the notations in section \S\ref{gluing map}.
\v
Fix  a basis $\{\mathbf{e}_{1},\cdots,\mathbf{e}_{d}\}$   of $Ker\;D{\mathcal  S}_{(\kappa_o,b_o)}$ and let $\mathfrak{z}=(\mathfrak{z}_{1},\cdots,\mathfrak{z}_{d})$ be the corresponding coordinates. Set $\ft_{i}=e^{-2r_{i}-2\pi\tau_{i}}$, $1\leq i\leq \mathbf{e}$. Denote
$$
\mathcal L(\fs,(\mathbf r),\mathfrak{z}):= I_{(\mathbf r)}\left(\sum_{i=1}^{d} \mathfrak{z}_{i}\mathbf{e}_{i}\right)+Q_{( \kappa_{o},b_{(\mathbf r)})}\circ f_{\fs,(\mathbf r)} \circ I_{(\mathbf r)} \left(\sum_{i=1}^{d} \mathfrak{z}_{i}\mathbf{e}_{i}\right),$$
where    $b_{(\mathbf r)}=(0,(\mathbf r),u_{(\mathbf r)}).$ Then $(\fs,(\mathbf r),\fkz)$ is  a local coordinates of ${\mathbf{U}}_{(\kappa_{\mathbf{a}},b_{\mathbf{a}})}(\varepsilon_{\mathbf{a}},\delta_{\mathbf{a}},\rho_{\mathbf{a}})$.
We say that $f(\mathbf{s},(\mathbf r),\fkz)$ satisfies $(\mathbf r)$-exponential decay if
\begin{equation}
 \left(\left| \frac{\p f}{\p r_{i}}\right|+ \left| \frac{\p f}{\p \tau_{i}}\right|\right) \leq  C e^{-\delta r_{i}},\;\;\forall \; 1\leq i\leq \mathbf{e}
\end{equation}
\begin{equation}
\left| \frac{\p f}{\p s_{j}}\right|+ \left| \frac{\p f}{\p \fkz_{\alpha}}\right|\leq C,\;\;\forall \;1\leq j\leq \iota,\; 1\leq \alpha \leq d.
\end{equation}
 Let
$$
\Pi_{1}:\widetilde{\mathbf F}_{b_{(\mathbf r)}}\times T_{u_{(\mathbf r)}} \widetilde{\mathcal B} \to \widetilde{\mathbf F}_{b_{(\mathbf r)}},\;\;\;\;\Pi_{2}:  \widetilde{\mathbf F}_{b_{(\mathbf r)}}\times T_{u_{(\mathbf r)}} \widetilde{\mathcal B} \to T_{u_{(\mathbf r)}}\widetilde{ \mathcal B}
$$
be the projection.
By Theorem \ref{coordinate_decay-2},  the implicit function Theorem  and  \eqref{cut-off-2}, we conclude that $\mathbf{\Gamma}_{\mathbf a}$  satisfies $(\mathbf r)$-exponential decay,
  where $\mathbf{\Gamma}_{\mathbf a}$ is the cutoff function defined in section   \S\ref{a metric}.

\v

For any $ {\mathbf{U}}_{(\kappa_{\mathbf{a}},b_{\mathbf{a}})}(\varepsilon_{\mathbf{a}},\delta_{\mathbf{a}},\rho_{\mathbf{a}})\bigcap  {\mathbf{U}}_{(\kappa_{\mathbf{c}},b_{\mathbf{c}})}(\varepsilon_{\mathbf{c}},\delta_{\mathbf{c}},\rho_{\mathbf{c}})\neq \emptyset,$
 let $a^{\mathbf a\mathbf c}_{\alpha\beta},\alpha,\beta=1,\cdots,\mathsf r$ be   functions such that $e^{\mathbf a}_{\alpha}=\sum_{\beta=1}^{\mathsf r_{i}} a^{\mathbf a\mathbf c}_{\alpha \beta}e^{\mathbf c}_{\beta},\alpha=1,\cdots,\mathsf r.$
By  the implicit function theorem, Theorem \ref{thm_est_mix_deri} we have, for any $p\in \Sigma(R_{0}),$
 $e^{\mathbf a}_{\alpha}(p),e^{\mathbf c}_{\beta}(p)$ satisfies $(\mathbf r)$-exponential decay.
Since   $a^{\mathbf a\mathbf c}_{\alpha \beta}$ is a  function of $(\mathbf{s},(\mathbf r),\fkz)$, we have
\begin{equation}\label{equ_exp_a}
d (e^{\mathbf a}_{\alpha}(p))=\sum_{\beta=1}^{\mathsf r} e^{\mathbf c}_{\beta} (p) \cdot da^{\mathbf a\mathbf c}_{\alpha \beta}+\sum_{\beta=1}^{\mathsf r} a^{\mathbf a\mathbf c}_{\alpha \beta}\cdot d(e^{\mathbf c}_{\beta}(p)),\;\;\;\;\forall\; p\in \Sigma(R_{0}).
\end{equation}
Recall that $e^{\mathbf a}_{\alpha}=\left(I_{(\mathbf r)}^{\widetilde{\mathbf L}}+Q_{(\mathbf r)}^{\widetilde{\mathbf L}}f^{\widetilde{\mathbf L}}_{\fs,h_{(\mathbf r)},(\mathbf r)}I_{(\mathbf r)}^{\widetilde{\mathbf L}}\right)(e^{\mathbf a}_{\alpha}|_{(\kappa_{\mathbf a},b_{\mathbf a})}).$ Using the implicit function theorem we get $$\|Q_{(\mathbf r)}^{\widetilde{\mathbf L}}f^{\widetilde{\mathbf L}}_{\fs,h_{(\mathbf r)},(\mathbf r)}I_{(\mathbf r)}^{\widetilde{\mathbf L}}(e^{\mathbf a}_{\alpha}|_{(\kappa_{\mathbf a},b_{\mathbf a})})\|_{k,2,\alpha,\mathbf r}\leq 2C\left\|D^{\widetilde{\mathbf L}}_{b}\circ(P^{\widetilde{\mathbf{L}}}_{b,b_{(\mathbf{r})}})^{-1}(I_{(\mathbf r)}^{\widetilde{\mathbf L}}(e^{\mathbf a}_{\alpha}|_{(\kappa_{\mathbf a},b_{\mathbf a})}))\right\|.$$ Choosing $\delta_{\mathbf a}$ and $\rho_{\mathbf a}$ small enough, by the exponential estimates of $e^{\mathbf a}_{\alpha}|_{b_{\mathbf a}}$ we have
$$
\|(e^{\mathbf a}_{\alpha}|_{(\kappa,b)})|_{\Sigma(R_{0})}\|_{k,2,\alpha}\geq \frac{1}{4}\|e^{\mathbf a}_{\alpha}\|_{k,2,\alpha}.
$$
 So $\max\limits_{\Sigma(R_{0})}|e_{\alpha}^{\mathbf a}|$ has uniform lower bound. Then we obtain   the  $(\mathbf r)$-exponential decay of $a^{\mathbf a\mathbf c}_{\alpha \beta}$.
Denote $h^{\mathbf a}_{\alpha\beta}=\langle e_{\alpha}^{\mathbf a},e_{\beta}^{\mathbf a}  \rangle_{\mathbf h}.$ By the definition of $\mathbf h$ and the $(\mathbf r)$-exponential decay of $\mathbf{\Gamma}_{\mathbf a}$,  $a^{\mathbf a\mathbf c}_{\alpha \beta}$ we conclude that $h^{\mathbf a}_{\alpha\beta}$ satisfies   the  $(\mathbf r)$-exponential decay.
By   the Gram-Schmidt orthonormalization and the similar argument above
we obtain the $(\mathbf r)$-exponential decay of $\hat a^{\mathbf a\mathbf c}_{\alpha\beta}.$

Let $\Delta_r$ be the open disk in $\mathbb C$
with radius $r$, let  $\Delta_r^*=\Delta_r\setminus\{0\}$ and $\Delta ^*=\Delta
\setminus\{0\}$. Set $N=3g-3+m+\mu$. For each point $p\in \p \mathcal{M}_{g,m+\mu}$ we can find a coordinate
chart $(U,\fs_{1},\cdots, \fs_{N-\mathbf{e}},\ft_1,\cdots,\ft_{\mathbf{e}})$ around $p$ in $\overline{\mathcal{M}}_{g,m+\mu}$ such that
$U\cong\Delta^N$ and $V=U\cap \overline{\mathcal{M}}_{g,m+\mu}\cong \Delta^{N-\mathbf{e}}\times(\Delta^*)^{\mathbf{e}}$. We assume that $U\cap \Delta$ is
defined by the equation $\ft_1\cdots \ft_{\mathbf{e}}=0$.
Let $\{U_{\alpha}\}$ be the local chart of $\overline{\mathcal{M}}_{g,m+\mu}$.
On each chart $U_{\alpha}$ of $\overline{\mathcal{M}}_{g,m+\mu}$ we can define a local Poincare metric:
\begin{eqnarray}\label{localpo}
g^{\alpha}_{loc}= \sum_{i=1}^{\mathbf{e}} \frac{|d\ft_i|^2}{|\ft_i|^2   (\log|\ft_i|)^2}+
 \sum_{j=1}^{N-\mathbf{e}} |d\fs_j|^2.
\end{eqnarray}
We let $U_{\alpha}(r)\cong
\Delta_r^{N}$ for $0<r<1$ and let $V_{\alpha}(r)=U_{\alpha}(r)\cap \mathcal{M}_{g,m+\mu}$.

Let $\fs, (\mathbf r),\fkz$ be the local coordinates of ${\mathbf{U}}_{(\kappa_{\mathbf{a}},b_{\mathbf{a}})}(\varepsilon_{\mathbf{a}},\delta_{\mathbf{a}},\rho_{\mathbf{a}}).$ In the coordinates $(\fs,(\mathbf r),\fkz)$
the local Poincare metric $g_{loc}$ can be written as
\begin{equation}\label{omega_loc}
g_{loc}= \sum_{i=1}^{\mathbf{e}} \frac{4(d^2r_i+d^2\tau_i)}{   r_i^2}+
\sum_{j=1}^{3g-3+n-\mathbf{e}} |d\fs_j|^2+\sum_{i=1}^{d}d\fkz_{i}^2.
\end{equation}

\v

\begin{lemma}\label{lem_omega}
	There exists a constant $C>0$ such that
$$|\omega^{\mathbf a}_{\alpha\beta}(X_{1})|^2\leq g_{loc}(X_{1},X_{1}),\;\;\;\;
|\Omega^{\mathbf a}_{\alpha\beta}(X_{1},X_{2})|^2\leq \Pi_{i=1}^{2}g_{loc}(X_{i},X_{i})$$
$$|d\Omega^{\mathbf a}_{AB}(X_{1},X_{2},X_{3})|^2\leq \Pi_{i=1}^{3}g_{loc}(X_{i},X_{i})$$
for any $X_{i}\in T\mathbf U^{T},i=1,2,3.$
\end{lemma}
{\bf Proof.}
	The first inequality follows from \eqref{local_omega} and $(\mathbf r)$-exponential decay of $\hat a^{\mathbf a\mathbf c}_{\alpha\beta}.$ By $\Omega_{\alpha\beta}=d\omega_{\alpha \beta}+\sum_{\gamma}\omega_{\alpha\gamma}\wedge \omega_{\gamma\beta}$ and $(\mathbf r)$-exponential decay of $\hat a^{\mathbf a\mathbf c}_{\alpha\beta}$ and $\mathbf{\Gamma}_{\mathbf a}$,  we can get the second inequality. The last inequality follows from the   Bianchi identity.

Let $\mathsf{p}^*\mathbf{E}$ be the pull-back of the bundle $\mathbf E$ to a bundle over $\mathbf E$, where $\mathsf{p}:\mathbf E\to \mathbf U$ is the projection.  Then the bundle $\mathsf{p}^* \mathbf E$ has a metric $\mathsf{p}^*\mathbf{h}$ with compatible connection $\mathsf{p}^*\nabla$. To simply notation we write these as $\mathbf{h}$ and $\nabla$.
Let $\hat\sigma$ be  the tautological section of  $\mathsf{p}^*\mathbf E.$ Then the elements $|\hat \sigma|_{\mathbf{h}}^2\in  \mathcal{A}^{0}(\mathbf E,  \wedge^{0}(\mathsf{p}^*\mathbf E))$, and the covariant derivative $\nabla \sigma\in  \mathcal{A}^{1}(\mathbf E,  \wedge^{1}(\mathsf{p}^*\mathbf E)).$
The curvature $\mathsf{p}^*\Omega$ of the connection $\nabla$ on $\mathbf E$ can also seen as an element of $\mathcal{A}^2(\mathbf E, \wedge^2(\mathsf{p}^*\mathbf E))$.
The Mathai-Quillen type Thom form can be written as
\begin{equation}\label{MQ_thom_gauss}
\Theta_{MQ}=c(\mathsf{r})\int^{B}e^{-\frac{|  \hat \sigma|^2_{\mathbf{h}}}{2}-\nabla\hat \sigma-\mathsf{p}^* \Omega} \in \mathcal A^{\mathsf{r}}(\mathbf E)
\end{equation}
where $c(\mathsf{r})$ is a constant depending only $\mathsf{r},$ $\int^{B}$ denotes the Berezin integral on $\wedge^* (\mathsf{p}^*\mathbf E)$.
Here $\Theta_{\mathbf E}$ is Gaussian shaped Thom form. Let $B_{\epsilon}(0)$ denote the open $\epsilon$-ball in $R^{2\mathbf r}$ and consider the map $\rho_{\epsilon}:B_{\epsilon}(0)\to \mathbb R^{\mathbf r} $ defined by $\rho_{\epsilon}(v)=\frac{v}{\epsilon^2-|v|^2}.$
If we extend $\rho_{\epsilon}^*\Theta_{MQ}$ by setting it
equal to zero outside $B_{\epsilon}(0)$, still denoted by $\Theta_{\mathbf E}:=\rho_{\epsilon}^*\Theta_{MQ},$ we obtain a form $\Theta_{\mathbf E}$ of compact support.

\v
Finally, we have the following estimate for $\sigma^{*}\Theta_{\mathbf E}$.
\begin{lemma}\label{lem_Thom}
	There exists a constant $C>0$ such that
$$	|\sigma^{*}\Theta_{\mathbf E}(X_{1},\cdots,X_{\mathsf r})|^2\leq C \Pi_{i=1}^{\mathsf r}g_{loc}(X_{i},X_{i})$$for any $X_{i}\in T\mathbf U^{T},i=1,2,3.$
\end{lemma}
{\bf Proof.}
	One can easily check that
	$$\sigma^{*}\Theta_{\mathbf E} =\sigma^{*} \rho^{*}\Theta_{MQ}=c(\mathsf r)\int^{B}e^{-\frac{|\sigma|^2_{\mathbf{h}}}{(\epsilon^2-|\sigma|^2_{\mathbf{h}})^2}-\nabla (\frac{\sigma}{\epsilon^2-|\sigma|^2})- \Omega} \in \mathcal A^{\mathsf{r}}(\mathbf M).$$
Denote $\sigma=\sum_{\alpha} \sigma_{\alpha} \hat e_{\alpha}^{\mathbf a}.$  For any $p\in \Sigma(R_{0})$, by $d\sigma(p)=\sum_{\alpha} d \hat e_{\alpha}^{\mathbf a}(p)\sigma_{\alpha}^{\mathbf a}+\sum_{\alpha}  \hat e_{\alpha}^{\mathbf a}(p) d\sigma^{\mathbf a}_{\alpha}$, as above we obtain the $(\mathbf r)$ exponential decay $\sigma^{\mathbf a}_{\alpha}.$ 	Since $\nabla \sigma= \sum _{\alpha}d\sigma_{\alpha}\hat e_{\alpha}^{\mathbf a}+ \sum_{\alpha,\beta}\sigma_{\beta}\omega_{\alpha \beta}\hat e_{\alpha}^{\mathbf a} $, $\Omega=\sum \Omega_{\alpha\beta}\hat e_{\alpha}^{\mathbf a}\wedge\hat  e_{\beta}^{\mathbf a}$ and
$$
\int^{B}\hat e_{1}^{\mathbf a}\wedge \cdots\wedge\hat  e_{\mathsf r}^{\mathbf a}=1,
$$
the lemma follows from Lemma \ref{lem_omega} and a direct calculation.

\chapter{Relative GW-invariants}

Recall that we have two natural maps
$$ev_i: \mathbf{U} \longrightarrow
M^{+} $$
$$ (\kappa; j,{\bf y},{\bf p},u) \longrightarrow
u(y_{i}) $$ for $1\leq i\leq m$ defined by evaluating at marked points
and $$e_j: \mathbf{U}\longrightarrow Z$$
$$ (\kappa;j,{\bf y},{\bf p}, u) \longrightarrow
u(p_{j}) $$
for $1\leq j\leq \mu$ defined by projecting to its
periodic orbit.
\v
We have another map
$$\mathscr{P}: \mathbf{U}^T \longrightarrow \mathcal{M}_{g,m+\mu}\;\;\;\;\left(\Sigma,j,{\bf y},{\bf p},(\kappa, u)\right)\longmapsto (\Sigma,j,{\bf y},{\bf p}).$$
Choose a smooth metric $\mathbf{h}$ on the bundle $\mathbf{F}$. Using $\mathbf{h}$ we construct a Thom form $\Theta$ supported in a small $\varepsilon$-ball of the $0$-section of $\mathbf{E}$. The relative Gromov-Witten invariants are defined as
\begin{equation}\label{integral-0}
\Psi_{A,g,m+\mu}(K;\alpha_1,...,
\alpha_{m};\beta_1,...,\beta_{\mu})=\int_{\mathbf{U}^T}\mathscr{P}^*(K)\wedge\prod^m_{i=1} ev^*_i\alpha_i\wedge\prod^\mu_{j=1}ev^*_j\beta_j\wedge \sigma^*\Theta
\end{equation}
 for $\alpha_i\in H^*(\overline{M}^+, {\mathbb{R}})$, $\beta_j\in H^*(Z, {\mathbb{R}})$ represented by differential form and $K$ represented by a good differential form defined on $\mathcal{M}_{g,m+\mu}$ in Mumford's sense. Clearly, $\Psi=0$ if $deg(K)+\sum \deg(\alpha_i)\neq Index $.
\v

\v
The following theorem is obvious.
\begin{theorem} \label{thm_3.4}
 Restricting to $\mathbf{U}^{T}$, the following hold:
\v
(a) the forgetful map $\mathscr{P}$ is smooth,
\v
(b) the evaluation map $ev$ is smooth.
\end{theorem}

Denote
$$\mathbf {U}_{\epsilon}=\{(\kappa,b)\in \mathbf U\;| \;|\kappa|_{\mathbf h}\leq \epsilon\},\;\;\;\;\mathbf {U}_{\epsilon}^{T}=\{(\kappa,b)\in \mathbf U^{T}\;| \;|\kappa|_{\mathbf h}\leq \epsilon\}. $$
We choose open covering
$$\{\mathbf{U}_{\kappa_{\mathbf{a}},b_{\mathbf{a}}}(\varepsilon_{\mathbf{a}},\delta_{\mathbf{a}},\rho_{\mathbf{a}})
, \;1\leq \mathbf{a}\leq \mathbf{n}_c\}$$
of $\mathbf{U}_{2\varepsilon}$
and a family of cutoff  functions $\{\mathbf{\Gamma}_{\mathbf{a}},\;1\leq \mathbf{a}\leq \mathbf{n}_c \}$ as in \S\ref{a metric}.
Let $\Theta_{\mathbf E}$ be the Thom form of $\mathbf E$  supported in a small $\varepsilon$-ball of the $0$-section of $\mathbf{E}$. Denote $\Theta_{\mathbf E}$ by $\Theta$ to simplify notation.

\v

\v
\begin{remark}\label{partition of unity}
$\{\mathbf{\Gamma}_{\mathbf{a}}\}$ is not a partition of unity in the classical sense, since it is not smooth on lower stratum, and it is not compactly supported. But it is smooth on $\mathbf{U}^T$ and $\mathbf{\Gamma}_{\mathbf{a}}\sigma^{*}\Theta$ is compactly supported. This is enough to define the relatve Gromov-Witten invariants.
\end{remark}

\v
Denote
$$\mathbf{V}_{\kappa_{\mathbf{a}},b_{\mathbf{a}}}(\varepsilon_{\mathbf{a}},\delta_{\mathbf{a}},\rho_{\mathbf{a}})
:=\mathbf{U}_{\kappa_{\mathbf{a}},b_{\mathbf{a}}}(\varepsilon_{\mathbf{a}},\delta_{\mathbf{a}},\rho_{\mathbf{a}})\cap
\mathbf{U}^T,$$
$$\widetilde{\mathbf{V}}_{\kappa_{\mathbf{a}},b_{\mathbf{a}}}(\varepsilon_{\mathbf{a}},\delta_{\mathbf{a}},\rho_{\mathbf{a}})
:=\widetilde{\mathbf{U}}_{\kappa_{\mathbf{a}},b_{\mathbf{a}}}(\varepsilon_{\mathbf{a}},\delta_{\mathbf{a}},\rho_{\mathbf{a}})\cap
\widetilde{\mathbf{U}}^T.$$
Sometimes we write the above two sets by $\mathbf{V}_{\mathbf{a}}$ and $\widetilde{\mathbf{V}}_{\mathbf{a}}$ to simplify notations. Let $p:
\widetilde{\mathbf{V}}_{\mathbf{a}}\to \mathbf{V}_{\mathbf{a}}$,
let $\widetilde{\mathbf{\Gamma}}_{\mathbf{a}}$, $\widetilde{K}$ and $\widetilde{\Theta}$ be the lift of $\mathbf{\Gamma}_{\mathbf{a}}$, $K$ and $\Theta$ to $\widetilde{\mathbf{V}}_{\mathbf{a}}$.
We write the relative Gromov-Witten invariants as
\begin{align}\label{def_GRW}
\Psi_{A,g,m+\mu}(K;\alpha_1,...,
\alpha_{m};\beta_1,...,\beta_{\mu}) =\sum_{\mathbf{a}=1}^{\mathbf{n}_{c}} (\mathbf{I})_{\mathbf{a}}
\end{align}
 where
\begin{equation}\label{integrand}
(\mathbf{I})_{\mathbf{a}}:=\int_{\mathbf{V}_\mathbf{a}}
\mathbf{\Gamma}_{\mathbf{a}}
\cdot \mathscr{P}^*(K)\wedge\prod^n_j ev^*_j\alpha_j\wedge \sigma^*\Theta.
\end{equation}
\begin{theorem}\label{Conver}
The integral \eqref{integral-0} is convergent.
\end{theorem}

{\bf Proof.} Note that the integration region $\overline{\mathbf{U}}_\mathbf{a}$ for $1\leq \mathbf{a}\leq \mathbf{n}_{t}$ are compact set in $\mathbf{U}^T$ and the integrand in \eqref{integrand} are smooth we conclude that $\sum_{\mathbf{a}=1}^{\mathbf{n}_{t}}(I)_{\mathbf{a}}$ is bounded. So we only need to prove the convergence of
$(\mathbf{I})_{\mathbf{a}}$ for $\mathbf{a}=\mathbf{n}_{t}+1,\cdots,\mathbf{n}_{c}.$
Denote
$$(\mathbf{J})_{\mathbf{a}}=\int_{\widetilde{\mathbf{V}}_\mathbf{a}}
\widetilde{\mathbf{\Gamma}}_{\mathbf{a}}
\cdot \mathscr{P}^*(\widetilde{K})\wedge\prod^n_j \widetilde{ev}^*_j\alpha_j\wedge \widetilde{\sigma}^*\widetilde{\Theta}.$$
It suffices to prove the convergence of $(\mathbf{J})_{\mathbf{a}}$.
\v
Let $(\kappa_{o},b_{o})$ be one of $\{(\kappa_{\mathbf{a}},b_{\mathbf{a}}), \mathbf{a}=\mathbf{n}_{t}+1,\cdots,\mathbf{n}_{c}\}$  and  $b_{o}=(a_{o},u).$ We choose  coordinates $(\fs,\ft,\fkz)$ and use the notations in section \S\ref{gluing map}.
\v

To simplify notation we   denote
$$dV=\bigwedge_{i} \left(dr_{i}\wedge  d\tau_{i}\right)\wedge \left(\bigwedge_{j}(\frac{\sqrt{-1}}{2}d\fs_{j}\wedge d\bar \fs_{j})\right)\wedge d\fkz_1\wedge \cdots \wedge d\fkz_{d} $$ and
 $$\delta_{i}=glu(\ft)_*\left(\tfrac{\p}{\p r_{i}}\right),\;\; \eta_{i}=glu(\ft)_*\left(\tfrac{\p}{\p \tau_{i}}\right),\;\;1\leq i \leq \mathbf{e}$$

$$\delta_{\alpha}=glu(\ft)_*\left(\tfrac{\p}{\p \fs_{\alpha-\mathbf{e}}}\right),\;\;\;\eta_{\alpha}=glu(\ft)_*\left(\tfrac{\p}{\p \bar\fs_{\alpha-\mathbf{e}}}\right), \;\;\mathbf{e}+1\leq \alpha\leq 3g-3+n$$ $$\varrho_{ \mathsf{i}}=glu(\ft)_*\left(\tfrac{\p}{\p \fkz_\mathsf{i}}\right),\;\;\;1\leq \mathsf{i} \leq d.$$
We will denote by $(E_1, E_2,...,E_{6g-6+2n+d})$
the frame
$$\left(
\delta_{1},...,\delta_{\mathbf{e}},\eta_1,...,\eta_{\mathbf{e}},
\delta_{\mathbf{e}+1},...,\delta_{3g-3+n},\eta_{\mathbf{e}+1},...,\eta_{3g-3+n},
\varrho_1,...,\varrho_d \right).$$
Then, for $\mathbf{a}=\mathbf{n}_{t}+1,\cdots,\mathbf{n}_{c}$,
\begin{align*}
 (\mathbf{J})_{\mathbf{a}}=&\int_{\widetilde{\mathbf{V}}_\mathbf{a}}\widetilde{\mathbf{\Gamma}}_{\mathbf{a}}\cdot
 \left( \mathscr{P}^*\widetilde{K}\wedge\prod^n_i \widetilde{ev}^*_i\alpha_i\wedge \widetilde{\sigma}^{*}\widetilde{\Theta} (E_1,E_2,...,E_{6g-6+2n+d})\right)dV.
\end{align*}
{\bf 1. Estimates for $\mathscr{P}^*\widetilde{K}$ }
\v

We can choose $\fs,\ft,\mathfrak{z}_{1},\cdots,\mathfrak{z}_{d}$ as the local coordinates of $\mathbf{U}^T$.
In this coordinates  $\mathscr{P}:\mathbf U^{T}\to {\mathcal M}_{g,n}$ can be written as
$$\mathscr{P}(\fs,\ft,\mathfrak{z}_{1},\cdots,\mathfrak{z}_{d})=(\fs,\ft).$$   Noth that  $\mathscr{P}_{*}E_{i}=E_{i}$ for $i\leq 6g-6+2n$, $\mathscr{P}_{*}E_{i}=0$ for $i\geq 6g-6+2n+1.$ We assume that for any $1\leq j\leq \deg(K),$  $E_{i_j}\in \{E_{1},\cdots, E_{6g-6+2n}\}$. Since $K$ has Poincare growth we have
\begin{equation}\label{eqn_est_Lclass}
 |\mathscr{P}^* \widetilde{K}  (E_{i_{1}},\cdots,E_{i_{\deg(\widetilde{K} )}} )|= |\widetilde{K}   (E_{i_{1}},\cdots,E_{i_{\deg(\widetilde{K} )}} )|
\leq   C\left[\Pi_{j=1}^{\deg(\widetilde{K} )}  g_{loc} (E_{i_{j}},E_{i_{j}})\right]^{\frac{1}{2}}.\end{equation}
\v\n
{\bf 2. Estimates for $\prod^n_i \widetilde{ev}^*_i\alpha_i$ }
\v\n
  For any $p\in M$ and $\xi\in T_pM$ we denote
$D\exp_p(\xi): T_pM\rightarrow T_{\exp_{p}\xi}M,$ then
\begin{equation}\label{dExp}
D\exp_p(\xi)\xi':=\frac{d}{dt}\exp_{p}(\xi + t\xi')\mid_{t=0}.
\end{equation}
Obviously, $D\exp_p(\xi)$ is an isomorphism when $|\xi|$ small enough.  By a  direct calculation we have,  for any $X\in \{\frac{\p}{\p  s_{i}},\frac{\p}{\p \bar s_{i}},\frac{\p}{\p r_{l}}, \frac{\p}{\p \tau_{l}},\frac{\p}{\p \fkz_{j}},1\leq i\leq 3g-3+n-\mathbf{e},1\leq l\leq \mathbf{e},1\leq j\leq d\},$
\begin{align}\label{eqn_ev_glu}
 | (\widetilde{ev}_{i})_{*}(glu(\ft))_*X |=|\Pi_{2,u}(X(glu(\fs,\ft,\fkz))(y_{i}))|=\left| D\exp_{u}(\Pi_{2,u}\mathcal{L})(\Pi_{2,u}X (\mathcal{L}))(y_{i})) \right|.
\end{align}
By Theorem \ref{coordinate_decay-2} and \eqref{eqn_ev_glu}  we have
 $$
 \|\widetilde{ev}_{*}E_{i}\|_{G_{J}}+\|\widetilde{ev}_{*}E_{\mathfrak e+i}\|_{G_{J}}\leq Ce^{-\delta r_{i}} ,\;\;\;\;\; \|\widetilde{ev}_{*}E_{j}\|_{G_{J}}\leq C,
 $$
 $$\left[g_{loc} (E_{i},E_{i})\right]^{\frac{1}{2}}=\left[g_{loc} (E_{\mathfrak e+i},E_{\mathfrak e+i})\right]^{\frac{1}{2}}=\frac{2}{r_{i}},\;\;\;\;  \left[g_{loc} (E_{j},E_{j})\right]^{\frac{1}{2}}=1$$
for $1\leq i\leq \mathfrak{e}, \mathfrak{2e}+1\leq j\leq 6g-6+2n+d.$
It follows that
\begin{equation}\label{eqn_est_ev}
|\Pi \widetilde{ev}^*_i\alpha_i(E_{i_1},\cdots,E_{i_c})| \leq C \Pi_{i_{j}=i_{1}}^{i_{c}}\left[g_{loc} (E_{i_{j}},E_{i_{j}})\right]^{\frac{1}{2}},
\end{equation}
where $\{E_{i_{1}},\cdots,E_{i_{c}}\}\subset \{E_1, E_2,...,E_{6g-6+2n+d}\}.$
\v\n
{\bf 3. Estimates for the Thom form}
\v
By Lemma \ref{lem_Thom}  we have
\begin{equation}\label{eqn_est_sthom}
|\widetilde{\sigma}^*\widetilde{\Theta}(E_{i_{1}},\cdots,E_{i_{\mathsf r}})|\leq C \Pi_{i_1}^{i_{\mathsf r}}\left[g_{loc} (E_{i_{j}},E_{i_{j}})\right]^{\frac{1}{2}},
\end{equation}
where $\{E_{i_{1}},\cdots,E_{i_{\mathsf r}}\}\subset \{E_1, E_2,...,E_{6g-6+2n+d}\}.$
It follows from \eqref{eqn_est_Lclass}, \eqref{eqn_est_ev} and \eqref{eqn_est_sthom} that
$$
\left|\mathscr{P}^*\widetilde{K} \wedge\prod_i \widetilde{ev}^*_i\alpha_i\wedge \widetilde{\sigma}^{*}\widetilde{\Theta} (E_{1},\cdots,E_{6g-6+2n+d})\right|\leq C \Pi_{i= 1}^{ 6g-6+2n+d}\left[g_{loc} (E_{i},E_{i })\right]^{\frac{1}{2}}\leq \frac{C}{\Pi_{i=1}^{\mathbf{e}} r_{i}^2}.
$$
 Hence the integral $(\mathbf{J})_{\mathbf{a}}$ is convergence.
\v
It is easy to see that $\Psi_{(A,g, m+\mu)}(\alpha_1,..., \alpha_m; \beta_{l},...,
\beta_\mu)$ is multi-linear and skew symmetric, and is independent of the
choice of $\widetilde{J}$ and $J$. Moreover, the following hold
\begin{lemma}\label{well-defined}
The integral \eqref{integral-0} is independent of
\v
(1) the choices of the forms $\alpha_i$ in $[\alpha_i]$ and $\beta_A$ in $[\beta_A]$,
\v
(2) the choice of the Thom form $\Theta$,
\v
(3) the choice of the partition of unity $\{\Lambda_k\}$,
\v
(4) the choice of the regularization.
\end{lemma}
\vskip 0.1in
\noindent
The proofs are the same for the Gromov-Witten invariants ( see \cite{LS-2}), we omit here.

\v\v

\chapter{A gluing formula}

In this chapter we prove a general gluing formula relating
GW-invariants of a closed symplectic manifold in terms of the relative
GW-invariants of its symplectic cut. The proof has two steps. The first step
is to define an invariant for $M_{\infty}$  and prove that it is
the same as the invariant of $M_r$. Then, we write the invariant
of $M_{\infty}$ in terms of relative invariants of $M^{\pm}$.

\v

\v
\section{Moduli space of stable maps in $M_\infty$}
\v

\subsection{Line bundle over $M_{\infty}$ and $M_{\varrho}$}\label{line bundle_1*}
\v
Through a diffeomorphism as in \S\ref{cylin almost struc}  we consider
$M^{\pm}$ to be
$$M^+=M^{+}_{0}\bigcup\left\{[0,\infty)\times \widetilde{M}\right\},\;\;\;\;M^-=M^{-}_{0}\bigcup\left\{(-\infty ,0]\times \widetilde{M}\right\}.$$
Put
$$M_{\infty}=M^+\bigcup M^-.$$
By Lemma \ref{line bundle-1} we have line bundles $L^+$ over $M^+$ and $L^-$ over $M^-$ such that
$$L^+\mid_{\{\infty\}\times \widetilde{M}}=L^-\mid_{\{-\infty\}\times \widetilde{M}}=\pi^*L_{Z}.$$
We can assume that
\begin{equation}\label{equ_om_hpm}
\omega^{+*}|_{[R_{0}+1,\infty)\times \widetilde{M} }=\pi^{*}\tau_{0},\;\;\;\;\omega^{-*}|_{[R_{0}+1,\infty)\times \widetilde{M} }=\pi^{*}\tau_{0}.
\end{equation}
In fact, in $[R_{0},\infty)\times \widetilde{M} ,$ $\omega^{+*}$ can be written as
$$
\omega^{+*}=\pi^{*}\tau_{0}+d(y\lambda).
$$
Let $\beta_{R_{0}}$ be the cut-off function as before.
 Set
 $$
\omega'^{+} =\omega^{+*}-d((1-\beta_{R_{0}})y \lambda).
 $$
Since $d((1-\beta_{R_{0}})y \lambda)$ is exact form on $M^{+},$ we have
$[\omega'^{+}]=[\omega^{+}]$ and  $\omega'^{+}$ satisfies  \eqref{equ_om_hpm}.
Note that the almost complex structure $\tilde{J}^*$ on $Z$ is different from $\tilde{J}$ on $Z$ ( see \S\ref{cylin almost struc}). Both the GW-invariants and the relative GW-invariants are independent of the choice of the almost complex structures. To simplify notations we still use $J$ to denote the new almost complex structure.

\v
Given a $\varrho>0$ large enough, from $M^+\bigcup M^-$, we construct
a almost complex manifold $M_{\varrho}$ as follows.  We choose the coordinates  $(a^+,\theta^+)$ (resp. $(a^-,\theta^-)$) on the  cylinder end of $M^+$ (resp. $M^{-}$).
We cut off the
part of $M^{\pm}$   with cylindrical
coordinate $|a^{\pm}|>\frac{3l\varrho}{2}$  and glue the remainders along
the collars of length $l\varrho$ of the cylinders with the gluing formulas:
\begin{eqnarray}\label{gluing_M_pm}
&&a^+=a^- + 2l\varrho,\;\;\;\;\;\;\;\;\;\;\theta^+=\theta^-, \mod 1.
\end{eqnarray}
\v
We fix a large number $R>4R_{0}>0$ and construct $M_{R}$. We take an integral symplectic form $\omega^*_{R}$ and construct a line bundle $L^{R}$ over $M_R$. Let $\varphi_\varrho: (R_{0}+1,\varrho)\to (R_{0},R)$ be a $C^{\infty}$ diffeomorphism, which induces a diffeomorphism $\varphi_\varrho: M_\varrho\to M_R$. Define $L^\varrho:=\varphi_\varrho^*L^{R}$.
Then $$\omega^{*}|_{[R_{0}+1,2l\varrho -R_{0}-1]\times \widetilde{M}}=\pi^{*}\tau_{0},$$
where $\pi:[R_{0}+1,2l\varrho -R_{0}-1]\times \widetilde{M}\to Z$ denotes the projection.
Since the first Chern class classifies smooth complex line bundles on a space,    $L^\varrho$ can also be obtained by gluing.  Let  $\nabla^{\pm}$ be a Hermitian connection on line bundle $\pi_{\pm}:L^{\pm}\to M^{\pm}$ such that $\omega^{\pm*}$ is the curvature.

\subsection{Moduli space of stable maps in $M_\infty$}
\v

To define the Moduli space of stable maps in $M_\infty$ we need to extend holomorphic cascades system in $M^{\pm}$ to include several holomorphic cascades systems.
\v
 Let $\mathcal{M}_{G(\mathfrak{d})^{+}_{i^{+}}}$, $1\leq i^+\leq l^+$, ( resp. $\mathcal{M}_{G(\mathfrak{d})^{-}_{i^{-}}}$, $1\leq i^-\leq l^-$) be holomorphic cascades systems in $M^+$ ( resp. in $M^-$ ). We have data
$$\left\{\Sigma^{\pm}_{i^{\pm}},
A^{\pm}_{i^{\pm}},g^{\pm}_{i^{\pm}}, m^{\pm}_{i^{\pm}}+\mu^{\pm}_{i^{\pm}},{\bf k}^{\pm}_{i^{\pm}} ,i^{\pm}=1,...,l^{\pm}\right\}.$$
Assume that $\sum{\mu}^{+}_{i^+}= \sum{\mu}^{-}_{i^-}:=\mu$.
Put
$$G(\mathfrak{d})^{+}:=\bigcup G(\mathfrak{d})^{+}_{i^{+}},\;\;G(\mathfrak{d})^{-}:=\bigcup G(\mathfrak{d})^{-}_{i^{-}},\;\;
\Sigma^+=\cup \Sigma^{+}_{i^{+}},\;\;\Sigma^-=\cup \Sigma^{-}_{i^{-}}$$
$$\mathcal{M}_{G(\mathfrak{d})^{+}}:=\bigoplus_{i^+=1}^{l^+} \mathcal{M}_{G(\mathfrak{d})^{+}_{i^{+}}},\;\;
\mathcal{M}_{G(\mathfrak{d})^{-}}:=\bigoplus_{i^-=1}^{l^-} \mathcal{M}_{G(\mathfrak{d})^{-}_{i^{-}}}.$$


\noindent
\vskip 0.1in
\noindent
\begin{definition}\label{moduli space} Given integers $g$, $m$, and $A\in H_2(M,\mathbb{Z})$ and $(G(\mathfrak{d})^+, G(\mathfrak{d})^-, \rho)$. A stable $(j,J)$-holomorphic configuration
of genus $g$ and class $A$ in $M_{\infty}$ of type $(G(\mathfrak{d})^+, G(\mathfrak{d})^-, \rho)$ is a triple
$(\mathbf{b}^{+}, \mathbf{b}^{-}, \rho)$, where
$\mathbf{b}^{+}\in \mathcal{M}_{G(\mathfrak{d})^{+}}$, $\mathbf{b}^{-}\in \mathcal{M}_{G(\mathfrak{d})^{-}}$
and $\rho:\{p^+_1,...,p^+_{\mu}\}\rightarrow
\{p^-_1,...,p^-_{\mu}\} $ is a one-to-one map satisfying

\begin{itemize}
\item[(1)] If we identify $p^{+}_i$ and $\rho(p^{+}_i)$ then $\Sigma^{+}\bigcup
\Sigma^{-}$ forms a connected closed nodal Riemann surface of genus $g$ with $m=m^+ + m^-$ marked points;
\item[(2)] Put $u=(u^+,u^-)$.  Then $\hat u^{+}(p^{+}_{i})=\hat u^{-}(\rho(p^{+}_{i}))$. We fix the coordinates
$(a^+, \theta^+)$ on $M^+$ and $(a^-, \theta^-)$ on $M^-$. For each $p^+_i$ and $\rho(p^{+}_i)$ we choose local cusp cylinder coordinates $(s^+, t^+)$ and $(s^-,t^-)$ as in \S\ref{Deligne-Mumford moduli space} and choose local Darboux coordinate systems $\mathbf{w}_i$ on $Z$ near $u^+(p^+_i)$.
Then $\widetilde{u}^+(s^+,t^+)$ and $\widetilde{u}^-(s^-,t^-)$ converge to the same
periodic orbit when $(s^+,t^+)\rightarrow p^+_i$ and $(s^-,t^-) \rightarrow \rho(p^{+}_i)$
respectively;
\item[(3)] $(\mathbf{b}^{+}, \mathbf{b}^{-}, \rho)$ represents the homology class $A\circeq \sum_{i=1}^{l^+} A^+_i + \sum_{j=1}^{l^-}A_j^{-}$.
\end{itemize}
\end{definition}
Denote by $\mathcal{M}_{(G(\mathfrak{d})^+,G(\mathfrak{d})^-, \rho)}$ the moduli space of equivalence classes of all
$(j,J)$-holomorphic configuration
of genus $g$ and class $A$ in $M_{\infty}$ of type $(G(\mathfrak{d})^+,G(\mathfrak{d})^-, \rho)$.
Suppose that ${\mathcal  C}^{J,A}_{g,m}$ is the set of indices $(G(\mathfrak{d})^+,G(\mathfrak{d})^-, \rho)$.
Let $C \in
{\mathcal  C}^{J,A}_{g,m}$. Denote by ${\mathcal  M}_C$ the set of stable maps
corresponding to $C$.
\v
We introduce another moduli space ${\widehat{\mathcal  M}}_C$: in (2) of the definition \ref{moduli space} we use the condition
$$\bar{u}^+ \; and \;\bar{u}^- \;are \; tangent \; to \; Z\; at \; p\; with \;the \;same \; order$$
to instead the condition
$$\widetilde{u}^+(s^+,t^+)\; and \;\widetilde{u}^-(s^-,t^-)\; converge \; to \; the \;same\;periodic \; orbit\; $$$$when \;(s^+,t^+)\rightarrow p^+_i \; and \;(s^-,t^-) \rightarrow \rho(p^{+}_i).$$
\v
The following lemma is obvious.
\vskip 0.1in
\noindent
\begin{lemma} ${\mathcal  C}^{J,A}_{g,m}$ is a finite
set.
\end{lemma}
 \vskip 0.1in
We define
$$\overline{\mathcal {M}}_{A}(M_{\infty},g,m)=\bigcup_{C\in {\mathcal  C}^{J,A}_{g,m}} \mathcal{M}_C.$$
\v

\v\n
\v
Using Lemmas in subsection \S\ref{bubble_phenomenon} we immediately obtain the following compactness theorem and convergence
theorem:
\v\n
\begin{theorem}\label{compact_moduli_space}
${\overline{\mathcal{M}}}_{A}(M_{\infty},g,m)$ is compact.
\end{theorem}
\vskip 0.1in \noindent
\begin{theorem}\label{convergence} Let $b_{\varrho\upper} \in
{\overline{\mathcal{M}}}_{A} (M_{\varrho\upper},g,m)$ be a sequence with $\lim\limits_{i\to\infty}\varrho\upper=\infty$. Then there is
a subsequence, still denoted by $\varrho\upper$, which {\em weakly converges} to a stable
$(j,J)$-holomorphic configuration in
${\overline{\mathcal{M}}}_{A}(M_{\infty},g,m)$  as $i\rightarrow \infty$.
\end{theorem} \vskip 0.1in
\noindent
{\bf Proof.} For any $\varrho$ we write $M_{\varrho}=M_{0}^{+}\bigcup\{[0,2l\varrho]\times \widetilde M\}\bigcup M_{0}^{-}.$ Let ${M^{+}_{\varrho}}=M_{0}^{+}\bigcup \{[0,l\varrho]\times \widetilde M\}$ and  ${M^{-}_{\varrho}}=M_{0}^{-}\bigcup \{[-l\varrho,0]\times \widetilde M\}.$ Suppose that
$$b_{\varrho\upper}=(j\upper, {\bf y}\upper,\nu, u\upper),$$
where
$u\upper:\Sigma\upper\rightarrow M_{\varrho\upper}$.
We may assume that $(\Sigma\upper;j\upper, {\bf y}\upper,\nu)$ is stable and converges to $(\Sigma;j,{\bf y},\nu)$ in
${\overline{\mathcal{M}}}_{g,m}$. Denote by $P\subset \Sigma $ the set of  singular points for $u\upper$, marked points and
the double points. By Lemma \ref{lower_bound_of_singular_points} and \eqref{energy_bound}, $P$ is
a finite set. Then $|du\upper|$ is uniformly bounded on every
compact subset of $\Sigma - P$. By passing to a  subsequence, possible a $\mathcal{T}$-rescalling, we may
assume that $u\upper$ converges uniformly  with all derivatives on every compact
subset of $\Sigma - P $ to a $J$-holomorphic map $u:\Sigma - P
\rightarrow M_{\infty}.$ Obviously, $u$ is a finite energy $J$-holomorphic map.

\v Suppose that  $\Sigma-P=\bigcup(\Sigma_{l}-P_{l})$, each $\Sigma_{l}-P_{l}$ is a connected compentont of $\Sigma-P$.
For any compact set $K\subset \Sigma_{l}-P_{l}$, $K$ can be identify with a set of $\Sigma\upper$ as $i$ big enough.
 If there exists $K\subset \Sigma_{l}-P_{l},$ a point $z\in K$ and a subsequence, still denoted by $i$, such that $u\upper(z)\in {M^{+}_{\varrho\upper}}$ for all $i$, then $(j_{l},\mathbf{y}_{l},u_{l})$ belong
 to a holomorphic cascade of  $ {\overline{\mathcal{M}}}_{A^{+}}(M^{+};g^{+},m^{+}+\mu^{+},{\bk}^{+},\nu^{+})$, $\Sigma_{l}$ is a component of $\Sigma^{+},$
  otherwise, $(j_{l},\mathbf{y}_{l},u_{l})$ belong to a holomorphic cascade of  ${\overline{\mathcal{M}}}_{A^{-}}(M^{-};g^{-},m^{-}+\mu^{-},{\bk}^{-},\nu^{-})$,  $\Sigma_{l}$ is a component of $\Sigma^{-}$. Then we get $(u^{\pm},\Sigma^{\pm};j^\pm,\mathbf {y}^{\pm})$.

\v
We construct bubble tree as in subsection \S\ref{bubble_phenomenon}. Let $q$ be a singular point for $u^{(i)}$, suppose that $z\upper\rightarrow q$ such that
$|du\upper|(z\upper)\rightarrow \infty$.
 We may assume that $a\upper(z\upper)\ne l\varrho_{i}$. By choosing subsequence we assume that $a\upper(z\upper)$
are in one of ${M^{+}_{\varrho\upper}}$ and ${M^{-}_{\varrho\upper}}$ for all $i$.

If there a subsequence, still denoted by $i$, such that $u\upper(z\upper)\in {M^{+}_{\varrho\upper}}$ for all $i$, then we  constructe bubble tree to get ${\Sigma^+}'$ and $({u^{+}}',j^{+},\mathbf{y}^{+},\mathbf{p}^{+}),$  otherwise
we  constructe  bubble tree to   get ${\Sigma^-}'$ and $({u^{-}}',j^{-},\mathbf{y}^{-},\mathbf p^{-}).$

   Let $q^{\pm} $ be a relative node of $ {\Sigma^{\pm}}'$ with $q^{\pm}\in \Sigma_{1}^{\pm}\wedge \Sigma_{2}^{\pm}$.   By the proof of Lemma \ref{bubble_tree_a-1} we have   $u^{\pm}|_{ \Sigma_{1}^{\pm}}$ and $u^{\pm}|_{ \Sigma_{2}^{\pm}}$ converge to a same periodic orbit as the variable tend to $q^{\pm}.$
Then we conclude that
$b^{\pm}=({u^{\pm}}',j^{\pm},\mathbf{y}^{\pm},\mathbf{p}^{\pm})$  belong to ${\overline{\mathcal{M}}}_{A^{\pm}}(M^{\pm};g^{\pm},m^{\pm}+\mu^{\pm},{\bk}^{\pm},\nu^{\pm})$, and $b=(b^{+},b^{-})$ belongs to  ${\overline{\mathcal{M}}}_{A}(M_{\infty},g,m).$
  $\;\;\;\Box$
\v\n

\section{Construction of a virtual neighborhood for  ${\overline{\mathcal {M}}}_{A}(M_{\infty},g,m)$ }\label{virtual_infty}

\subsection{Construction of a virtual neighborhood for $\overline{\mathcal{M}}_C$}\label{virtual_infty-1}
Choose $3\sigma<\sigma_0$. We construct a virtual neighborhood $\mathbf{U}_{C}$ of $\overline{\mathcal{M}}_C$ for every $C \in
{\mathcal  C}^{J,A}_{g,m}$. The construction is divided into 2 steps:
\v
{\bf 1.} For each holomorphic block in $\overline{\mathcal{M}}_C$ we construct the local regularization as in Chapter \S\ref{regulization}.
\v
{\bf 2.} Construct global regularization. Since $\overline{\mathcal{M}}_C$ is compact, there exist finite points $b_{i_c}$, $1\leq i_c \leq \mathfrak{m}_C$, such that
\begin{itemize}
\item[(1)] The collection $\{\mathbf{O}_{[b_i]}(\delta_i/3,\rho_i/3)
\mid 1\leq i \leq \mathfrak{m}_C\}$ is an open cover of
$\overline{\mathcal{M}}_C$.
\item[(2)] Suppose that $\widetilde{\mathbf{O}}_{b_i}(\delta_i,\rho_i)
\cap \widetilde{\mathbf{O}}_{b_j}(\delta_j,\rho_j)
\neq\phi$. For any $b\in \widetilde{\mathbf{O}}_{b_i}(\delta_i,\rho_i)
\cap \widetilde{\mathbf{O}}_{b_j}(\delta_j,\rho_j)$, $G_b$ can be imbedded into both $G_{b_i}$ and $G_{b_j}$ as subgroups.
\end{itemize}
We may choose $[b_i]$, $1\leq i \leq \mathfrak{m}_C$, such that if $[b_i]$ lies in the top strata for some $i$, then
$\mathbf{O}_{[b_i]}(\delta_i,\rho_i)$ lies in the top strata.
\v\n
Set
$$\mathcal{U}^C=\bigcup_{i=1}^{\mathfrak{m}_C}
\mathbf{O}_{[b_{i}]}(\delta_i/2,\rho_i/2).$$

By Lemma \ref{finite cov} we have a continuous orbi-bundle $\mathbf{F}(\mathsf{k}_i)\rightarrow \mathcal{U}^C$ such that
$\mathbf{F}(\mathsf{k}_i)\mid_{b_i}$ contains a copy of group ring
$\mathbb{R}[G_{b_i}].$
Set
$$\mathbf{F}^C=\bigoplus_{i=1}^{\mathfrak{m}_C}\mathbf{F}(\mathsf{k}_i).
$$
We define a bundle map
$\mathfrak{i}^C:  {\mathbf{F}}^C\rightarrow \mathcal{E}^C$
and a global regularization $\mathcal{S}:\mathbf{F}^C\to \E^C$
as in \S\ref{global_regu}.

\v
For each $C\in \mathcal{C}$ we do this and put
$$\mathcal{U}_{\infty}=\bigcup_{C\in \mathcal{C}}\mathcal{U}^{C},\;\;\;\;\mathbf{F}_{\infty}=\bigoplus_{C\in \mathcal{C}}{\mathbf{F}}^C.$$
Define a bundle map
$\mathfrak{i}_{\infty}:  {\mathbf{F}}_{\infty}\rightarrow \mathcal{E}_{\infty}$ as in \S\ref{global_regu}.
We define a global regularization for $\overline{\mathcal {M}}_{A}(M_{\infty},g,m)$ to be the bundle map $\mathcal{S}_{\infty}:\mathbf{F}_{\infty}\to \E_{\infty}$
by
$$\mathcal{S}_{\infty}([\kappa,b])
=[\bar{\partial}_{j,J}v] + [\mathfrak{i}_{\infty}(\kappa,b)].
$$
The meaning of $\E^C$ and $\E_{\infty}$ above are obvious.
Denote
$$\mathbf{U}_{\infty}=\mathcal{S}_{\infty}^{-1}(0)|_{\mathcal{U}_{\infty}}.$$
There is a bundle of finite rank $\mathbf{E}_{\infty}$ over $\mathbf{U}_{\infty}$ with a canonical section $\sigma_{\infty}$. We have a virtual neighborhood for $\overline{\mathcal {M}}_{A}(M_{\infty},g,m)$:
$$(\mathbf{U}_{\infty},\mathbf{E}_{\infty},\sigma_{\infty}).$$
Denote by $\mathbf{U}_{\infty}^T$ the top strata of $\mathbf{U}_{\infty}$.
By the same method as in \cite{LS-2} we can prove
\begin{theorem} \label{Smooth}
$\mathbf{U}_{\infty}^T$ is a smooth oriented, effective orbifold.
\end{theorem}

We can define GW-invariants $\Psi_{C}$ for each component $C$ and define $\Psi_{(M_{\infty},A,g,m)}$ by
\begin{equation}
\Psi_{(M_{\infty},A,g,m)}=\sum_{{\mathcal
C}^{J,[A]}_{g,m}}\Psi_{C}.
\end{equation}
 \vskip 0.1in \noindent
\begin{remark} It is easy to see that
 \begin{itemize}
 \item[(i)] For $C= \{A^+, g^+, m^+\}$, we have
\begin{equation}\Psi_{C}(\alpha^+)=
\Psi^{(\overline{M}^+,Z)}_{(A^+,g^+,m^+)}(\alpha^+);
\end{equation}
\item[(ii)] For $C= \{A^-, g^-, m^-\}$, we have
\begin{equation}\Psi_{C}(\alpha^-)=
\Psi^{(\overline{M}^-,Z)}_{(A^-,g^-,m^-)}(\alpha^-).
\end{equation}
\end{itemize}
\end{remark}

\v

\section{Construction of a virtual neighborhood for  ${\overline{\mathcal {M}}}_{A}(M_{\varrho},g,m)$ }
\label{virtual_varrho}

We fix a $\varrho>0$ large enough and construct
almost complex manifold $M_{\varrho}$ as in \S\ref{line bundle}. There exists a smooth diffeomorphism  $\phi_{\varrho}:M_{\varrho}\to M_R$. Set $\omega_{\varrho}^{*}=\phi^{*}_{\varrho}\omega^*$. Then $(M,\omega^{*})$ and   $(M_{\varrho},\omega_{\varrho}^{*})$ are symplectic   diffeomorphism.
For $(M_{\varrho},\omega_{\varrho}^{*}),$
 we can construct a finite rank bundle $L_{\varrho}=\phi^{*}_{\varrho}L$.

\v
Let $C=(G(\mathfrak{d})^+, G(\mathfrak{d})^-,\rho)$ and $b=(b^+, b^-)\in \mathcal{M}_{C}$, where
$$b^+=(\Sigma^+,j^+,u^+)\in \mathcal{M}_{G(\mathfrak{d})^{+}},\;\;
b^-=(\Sigma^-,j^-,u^-)\in \mathcal{M}_{G(\mathfrak{d})^{-}}.$$
Suppose that
\begin{equation}\label{eqn_conv-pm}
a^{\pm}({s}^{\pm}_{j},{t}^{\pm}_{j})-k_{j}{s}^{\pm}_{j}-l^{\pm}_{j}\rightarrow 0,\;\;\;\;
\theta^{\pm}({s}^{\pm}_{j},{t}^{\pm}_{j})-k_{j}{t}^{\pm}_{j}-\theta^{\pm}_{j0}\rightarrow 0.
\end{equation}
Set
\begin{equation} \label{gluing_local_relative_node-3}
r_{j}=\varrho-\frac{\ell^+_{j}-\ell^-_{j}}{2l},\;\;\;\; \tau_{j}=\theta^+_{j0}-\theta^-_{j0}.
\end{equation}
We construct a surface
$\Sigma_{(\mathbf {r})} =\Sigma_1\#_{({\mathbf r})} \Sigma_2 $ with gluing formulas:
\begin{eqnarray}
&& {s}^+_{j}= {s}^-_{j} + \tfrac{2lr_{j}}{k_j} \\
&& {t}^-_{j}= {t}^-_{j} + \tfrac{\tau_{j} + n_{j}}{k_j}
\end{eqnarray}
for some $n_{j} \in \mathbb Z_{k_{j}}$. Denote
$$
\mathbf w^+_{j}=\mathbf w_{j}\circ u^+,\;\;\;\mathbf w^-_{j}=\mathbf w_{j}\circ u^-.
$$
 In terms of $(s^{\pm}_{j}, t^{\pm}_{j})$ we construct pre-gluing map $u_{(\mathbf r)}=(a_{(\mathbf r)}, \theta_{(\mathbf r)},{\bf w}_{(\mathbf r)} ):\Sigma_{(\mathbf r)}\rightarrow M_{\varrho}$ as follows:
  for every $j=1,...,\mu$
\begin{align}\label{pre-glu-a}
a_{(\mathbf r)}(s^+_{j},t^+_{j})= \;&k_{j}s^+_{j} +l^+_{j} +  \beta\left(3-\frac{4k_{j}s^+_{j}}{lr_{j}}\right) (a^{+}(s^+_{j},t^+_{j})-
k_{j}s^+_{j}- l^+_{j} )  \\
&+ \beta\left(\frac{4k_{j}s^{+}_{j}}{lr_{j}}-5\right) (a^{-}(s^{+}_{j},t^{+}_{j})- ks^{-}_{j}
- l^{-}_{j} ) ,\nonumber
\\ \label{pre-glu-theta}
 \theta_{(\mathbf r)}(s^{+}_{j},t^{+}_{j}) =\; & k_{j}t^+_{j} +\theta^{+}_{j0} +  \beta\left(3-\frac{4k_{j}s^+_{j}}{lr_{j}}\right) (\theta^{+}(s^+_{j},t^+_{j})-
k_{j}t^+_{j}-\theta^+_{j0})  \\
&+ \beta\left(\frac{4k_{j}s^{+}_{j}}{lr_{j}}-5\right) (\theta^{-}(s^{+}_{j},t^{+}_{j})- k_{j}t^{-}_{j}-\theta^{-}_{j0}) ,\nonumber\\ \label{pre-glu-w}
{\bf w}_{(\mathbf r)}(s^{+}_{j},t^{+}_{j})=\;&
 \beta\left(3-\frac{4k_{j}s^{+}_{j}}{lr_{j}}\right)   {\bf w}^{+}_{j}(s^{+}_{j},t^{+}_{j}) +
\beta\left(\frac{4k_{j}s^{+}_{j}}{lr_{j}}-5\right)
 {\bf w}^{-}_{j}(s^{+}_{j},t^{+}_{j})  .
\end{align}
\v
There is a constant $L>0$ such that the number of components of $\Sigma $ $< L$ for every $C\in {\mathcal  C}^{J,A}_{g,m}$. Let $b_{i_c}$, $1\leq i_c \leq \mathfrak{m}_C$, be as in \S\ref{virtual_infty-1}.
For each $b_{i_c}$ we do gluing at relative nodes as follows. For each relative node in $\Sigma^{\pm}$ we use $\frac{\varrho}{L}$ to glue $M^{\pm}$ and ${\mathbb{R}}\times \widetilde{M}$, for each relative node, at which $\Sigma^+$ and $\Sigma^-$ joint, we use $\varrho$ to glue $M^+$ and $M^-$. Then we
choose gluing parameters $(\mathbf{r})$ and construct $b_{i_c,\mathbf{(r)},\varrho}=(\Sigma_{i_c,(\mathbf r)},j_{i_c}, {\bf y}_{i_c},\nu, u_{i_c,(\mathbf r)})$ such that $u_{i_c,(\mathbf r)}:\Sigma_{i_c,(\mathbf r)}\to M_{\varrho}$.
We define the  $\{\mathbf{O}_{[b_{i_c,\varrho}]}(\delta_{i_c} ,\rho_{i_c})
\mid 1\leq i_c \leq \mathfrak{m}_C\}$ as before. For each $C\in {\mathcal  C}^{J,A}_{g,m}$ we do this, then we get $\mathfrak{m}$ points, denoted by $\{b_1,b_2,...,b_{\mathfrak{m}}\}$.
\begin{lemma} \label{virtual_varrho-1} There exist two constants $\varrho_{o}>0$ and $\epsilon>0,$ such that for any $\varrho>\varrho_{o},$
\begin{itemize}
\item[(1)] The collection $\{\mathbf{O}_{[b_{i,\varrho}]}((1+\epsilon)\delta_i/3,(1+\epsilon)\rho_i/3)
\mid 1\leq i \leq \mathfrak{m}\}$ is an open cover of
$\overline{\mathcal{M}}_{A}(M_{\varrho},g,m)$.
\item[(2)] Suppose that $\widetilde{\mathbf{O}}_{b_{i,\varrho}}((1-\epsilon)\delta_i,(1-\epsilon)\rho_i)
\cap \widetilde{\mathbf{O}}_{b_{j,\varrho}}((1-\epsilon)\delta_j,(1-\epsilon)\rho_j)
\neq\phi$. For any $b\in \widetilde{\mathbf{O}}_{b_{i,\varrho}}((1-\epsilon)\delta_i,(1-\epsilon)\rho_i)
\cap \widetilde{\mathbf{O}}_{b_{j,\varrho}}((1-\epsilon)\delta_j,(1-\epsilon)\rho_j)$, $G_b$ can be imbedded into both $G_{b_{i,\varrho}}$ and $G_{b_{j,\varrho}}$ as subgroups.
\end{itemize}
\end{lemma}
{\bf Proof.} If (1) does not hold, we can find a sequence
 $$b_k =(\Sigma_{k},j_{k},{\bf y}_{k},u_k)\in \mathcal{M}_{A}(M_{\varrho_{k}},g,m)\setminus \left( \bigcup_{i}\mathbf{O}_{[b_{i,\varrho_{k}}]}((1+\epsilon)\delta_i/3,(1+\epsilon)\rho_i/3)\right),$$
such that $\Sigma_{k}=\Sigma_{(\mathbf r_k)}$ with $\varrho_k\rightarrow \infty$, By the convergence theorem \ref{convergence} we conclude that $b_k$ weak converges to some $b\in \overline{\mathcal{M}}_{A}(M_{\infty},g,m+\mu)$.
  Then $b\in \mathbf{O}_{[b_{i}]}(\delta_i/3,\rho_i/3)$ for some $i.$ It follows that
$b_{k} \in  \mathbf{O}_{[b_{i,\varrho_{k}}]}((1+\epsilon)\delta_i/3,(1+\epsilon)\rho_i/3)$ as $k$ large enough. We get a contradiction.
\v The proof of (2) is standard ( see the proof of Lemma 4.3 in \cite{LS-2} ).

Obviously, the following lemma holds
\begin{lemma}
There exists a constant $\varrho_{o}>0$ such that for any $\varrho>\varrho_{0},$
$G_{b_{i,\varrho}}$  can be imbedded into $G_{b_{i}}$.
\end{lemma}
\v
Set
$$\mathcal{U}_{\varrho}=\bigcup_{i=1}^{\mathfrak{m}}
\mathbf{O}_{[b_{i,\varrho}]}(\delta_i/2,\rho_i/2).$$
By Lemma \ref{finite cov}, we have a continuous orbi-bundle $\mathbf{F}_{\varrho}(\mathsf{k}_i)\rightarrow \mathcal{U}_{\varrho}$ such that
$\mathbf{F}_{\varrho}(\mathsf{k}_i)\mid_{b_i}$ contains a copy of group ring
$\mathbb{R}[G_{b_i}].$
Set
$$\mathbf{F}_{\varrho}=\bigoplus_{i=1}^{\mathfrak{m}}\mathbf{F}_{\varrho}(\mathsf{k}_i).
$$
We define a bundle map
$\mathfrak{i}_{\varrho}:  {\mathbf{F}}_{\varrho}\rightarrow \mathcal{E}_{\varrho}$
and a global regularization $\mathcal{S}_{\varrho}:\mathbf{F}_{\varrho}\to \E_{\varrho}$
as in \S\ref{global_regu}. The meaning of   $\E_{\varrho}$ above are obvious.
Denote
$$\mathbf{U}_{\varrho}=\mathcal{S}_{\varrho}^{-1}(0)|_{\mathcal{U}_{\varrho}}.$$
There is a bundle of finite rank $\mathbf{E}_{\varrho}$ over $\mathbf{U}_{\varrho}$ with a canonical section $\sigma_{\varrho}$. We have a virtual neighborhood for $\overline{\mathcal {M}}_{A}(M_{\varrho},g,m)$:
$$(\mathbf{U}_{\varrho},\mathbf{E}_{\varrho},\sigma_{\varrho}).$$
Denote by $\mathbf{U}_{\varrho}^T$ the top strata of $\mathbf{U}_{\varrho}$.
By the same method as in \S\ref{top strata} we can prove
\begin{theorem} \label{Smooth}
$\mathbf{U}_{\varrho}^T$ is a smooth oriented, effective orbifold.
\end{theorem}

\v

 \section{ Equivariant gluing}\label{equ_glu}

We have several kinds of gluing maps. \\
{\bf (1).} Gluing maps in a holomorphic cascades as in {\bf Case 1} of \S\ref{gluing map}.
\v\n
{\bf (2).} Gluing maps at relative nodes in $\mathcal{M}_{G(\mathfrak{d})^{+}}$
or $\mathcal{M}_{G(\mathfrak{d})^{-}}$ as in {\bf Case 2} of \S\ref{gluing map}.
\v\n
{\bf (3).} Gluing maps at relative nodes between $\mathcal{M}_{G(\mathfrak{d})^{+}}$
and $\mathcal{M}_{G(\mathfrak{d})^{-}}$, which we discuss in this section.
\v

Let $C=(G(\mathfrak{d})^+, G(\mathfrak{d})^-,\rho),\kappa_{o}=(\kappa_{o}^{+},\kappa_{o}^{-}),$ and $b_{o}=(b_{o}^+, b_{o}^-)\in \mathcal{M}_{C}$, where
$$\;\;\;\;b_{o}^{\pm}=( a^{\pm},u^{\pm})\in \mathcal{M}_{G(\mathfrak{d})^{\pm}},\;\;\; a^{\pm}=(\Sigma^{\pm},j^{\pm},\mathbf{y}^{\pm}),$$ and $\Sigma^+$, $\Sigma^-$ denote marked nodal Riemann surfaces joining at $p_1,...,p_{\mu}$.  We may assume that $(\Sigma^{\pm},\mathbf y^{\pm})$ is stable.
We first discuss equivariant pregluing of Riemann surface. Denote by $\mathbf G_{a^{\pm}}$ the isotropy group at $a^{\pm}$.
Choose cusp coordinates near the nodes of $\Sigma^{\pm}$. Since the cusp coordinates are unique modulo rotations near nodes, each $g^{\pm}\in \mathbf G_{a^{\pm}}$ is a rotation in the cusp coordinate.
For any gluing parameter $(\mathbf r)$, we can obtain the gluing surface $ a_{(\mathbf r)}=(\Sigma_{(\mathbf r)},j,\mathbf y)$ as usual. Denote by $\mathbf G_{a_{(\mathbf r)}}$ the isotropy group at $a_{(\mathbf r)}$.
Obviously, $\mathbf G_{a_{(\mathbf r)}}$  is subgroup of  $ \mathbf G_{a^{+}}\times \mathbf G_{a^{-}}.$  Each element of $\mathbf G_{a_{(\mathbf r)}}$ is also a rotation in the domain of gluing.
Then the gluing map is the $\frac{|G_{\mathbf a^{+}}\times G_{a^{-}}|}{|G_{a_{(\mathbf r)}}|}$-multiple covering map of $(\Sigma_{(\mathbf r)},j,\mathbf y)$. We introduce some notations.  Put
$$\widetilde E^+:=\{(\kappa^+_0, h^++\hat{h}^+_{0})\mid D   \mathcal{S}_{u^+}(h^+ + \hat{h}^+_{0})=0,\;\;h^+_{0}\in \mathds{H}\},$$
$$\widetilde E^-:=\{(\kappa_{o}^{-},h^- + \hat{h}^-_{0})\mid D \mathcal{S}_{u^-}(h^- + \hat{h}^-_{0})=0,\;\;h^-_{0}\in \mathds{H}\}.$$
$$Ker D \mathcal S_{(\kappa_{o},b_{o})}:=\left\{\left((\kappa_{o}^{+},h^++\hat{h}^+_{0}),(\kappa_{o}^{-},h^-+\hat{h}^-_{0})\right)
\in \widetilde{E}^+ \oplus  \widetilde{E}^-\;|\; h^+_{0}=h^-_{0}\right\}.$$
The tangent space of $\widetilde{\mathbf{U}}_{\infty}$ at $(\kappa_{o},b_{o})$ can be defined by
 $$Ker \mathbb D \mathcal S_{(\kappa_{o},b_{o})}:=\left\{\left((\kappa_{o}^{+},h^++\hat{h}^+_{0}),(\kappa_{o}^{-},h^-+\hat{h}^-_{0})\right)
\in E^+ \oplus  E^-\;|\; \pi_{*}h^+_{0}=\pi_{*}h^-_{0}\right\}.$$
As in \S\ref{gluing map}, let $ {\mathbb{E}}$ be a subspace of $Ker \mathbb D \mathcal S_{(\kappa_{o},b_{o})}$ such that $Ker \mathbb D \mathcal S_{(\kappa_{o},b_{o})}= {\mathbb{E}}\oplus Ker  D \mathcal S_{(\kappa_{o},b_{o})}$.
Denote by $ G_{(\kappa_{o},b_{o})}$ the isotropy group at $(\kappa_{o},b_{o})$. For each $\phi\in  G_{(\kappa_{o},b_{o})},$ it induce a natural action on    $Ker \mathbb D \mathcal S_{(\kappa_{o},b_{o})}$.

 \v
 For every puncture point $p_j$ there are constants $(\ell^{\pm}_{j}, \theta^{\pm}_{j0})$ such that \eqref{eqn_conv-pm} hold.
We associate a point
$$\bar{\mathbf{t}}_j=\exp\{(\ell^+_{j}-\ell^-_{j})+ 2\pi\sqrt{-1}(\theta^-_{j0}-\theta^+_{j0})\}.$$
Put
$$\mathbb{D}^\circ :=\{\mathbf{t}^\circ_j \mid (\mathbf{t}^\circ _j)^{k_j}=\bar{\mathbf{t}}_j\}.$$
There is an $G_{(\kappa_{o},b_{o})}$-equivariant isomorphism
\begin{equation} \label{gluing_local_relative_node-4}
\psi:(\mathbb D^{\circ})^{\mu}\times Ker D \mathcal{S}_{(\kappa_{o},b_{o})}\to Ker \mathbb{D} \mathcal{S}_{(\kappa_{o},b_{o})}.\end{equation}
\v
For any  $\zeta=(\zeta_1,\zeta_2)\in \mathbb{E},$ we define a map $v=\exp_{u}(\zeta)$. Then $v$ satisfies \eqref{gluing_local_relative_node-1}.
For any $(\rho,0)$ we glue $M^+$ and $M^{-}$ to get $M_\varrho$  as in section \S\ref{relative_node_M_M-1}.
We construct a surface
$\Sigma_{(\mathbf {r})} $ with gluing parameter $(\mathbf r)$, where $(\mathbf r)$ satisfies \eqref{gluing_local_relative_node-3}.
 By pregluing as in \eqref{pre-glu-a}, \eqref{pre-glu-theta} and \eqref{pre-glu-w} we get $v_{(\mathbf r)}$. Let $b_{ (\mathbf r)}=(j,\mathbf{y},v_{(\mathbf r)})$.
 Denote by $G_{b_{(\mathbf{r})}}$ (resp. $G_{(\kappa_{o},b_{(\mathbf{r})})}$) the isotroy group at $b_{(r)}$ (resp. $(\kappa_{o},b_{(\mathbf{r})})$). It is easy to see that $G_{b_{(\mathbf{r})}}$ is a subgroup of $G_{b_{o}}.$ It follows that
$G_{(\kappa_{o},b_{(\mathbf{r})})}$ is a subgroup of $G_{(\kappa_{o},b_{o})}.$
Then $G_{(\kappa_{o},b_{(\mathbf{r})})}$ can be seen as rotation in the gluing part. The gluing map is the $\frac{|G_{(\kappa_{o},b_{o})}|}{|G_{(\kappa_{o},b_{(\mathbf{r})})}|}$-multiple covering map.
 Denote
$$
\ker D\mathcal{S}_{[\kappa_o,b_o]}=\ker D\mathcal{S}_{(\kappa_o,b_o)}/G_{(\kappa_{o},b_{o})},\;\;\;\ker D\mathcal{S}_{[\kappa_o,b_{(\mathbf{r})}]}=\ker D\mathcal{S}_{(\kappa_o,b_{(\mathbf{r})})}/G_{(\kappa_{o},b_{(\mathbf{r})})}.
$$
 By the same method as in \cite{LS-2}
 we can prove that
 \begin{lemma}\label{equi_glu_rho}
{\bf (1)} $I_{(\mathbf{r})}: \ker D\mathcal{S}_{(\kappa_o,b_o)} \longrightarrow \ker D\mathcal{S}_{(\kappa_o,b_{(\mathbf{r})})}$ is a $\frac{|G_{(\kappa_{o},b_{o})}|}{|G_{(\kappa_{o},b_{(\mathbf{r})})}|}$-multiple covering map.
\v
{\bf (2)} $I_{(\mathbf{r})}$ induces a isomorphism $I_{(\mathbf{r})}:  \ker D\mathcal{S}_{[\kappa_o,b_o]} \longrightarrow \ker D\mathcal{S}_{[\kappa_o,b_{(\mathbf{r})}]}$.
\end{lemma}

   \v\n
We introduce some notations. Denote
$$\mathbf {U}_{\infty,\epsilon}=\{(\kappa,b)\in \mathbf U_{\infty}\;| \;|\kappa|_{\mathbf h}\leq \epsilon\},\;\;\;\;  \mathbf {U}_{\varrho,\epsilon}=\{(\kappa,b)\in \mathbf U_{\varrho}\;| \;|\kappa|_{\mathbf h}\leq \epsilon\},$$
$$\mathbf{U}_{\infty;\kappa_o,b_o}(\varepsilon,\delta_o,\rho_o) =\left\{(\kappa, b)\in \mathbf{U}\mid |\kappa-\kappa_{o}|_{\mathbf{h}}<\varepsilon,
b\in \mathbf{O}_{b_o}(\delta_o,\rho_o)\right\}.$$
$$\mathbf{U}^{T}_{(\infty;\kappa,b)}(\varepsilon,\delta,\rho):=
\mathbf{U}_{(\infty;\kappa,b)}(\varepsilon,\delta,\rho)
\bigcap \mathbf {U}^T_{\infty,\epsilon}\;\;\;for \;(\kappa,b)\in \mathbf {U}_{\infty,\epsilon}.$$
We choose open covering
$$\{ \mathbf{U}^{T}_{\infty,(\kappa_{\mathbf{a}},b_{\mathbf{a}})}(\varepsilon_{\mathbf{a}},\delta_{\mathbf{a}},\rho_{\mathbf{a}})
, \;1\leq \mathbf{a}\leq \mathbf{n}_c\}$$
of $\mathbf{U}^{T}_{\infty,2\varepsilon}$
and a family of cutoff  functions $\{\mathbf{\Gamma}_{\infty,\mathbf{a}},\;1\leq \mathbf{a}\leq \mathbf{n}_c \}$ as in \S\ref{a metric}.
For each $1\leq \mathbf a\leq \mathbf n_{c},$ we fix a basis $\{\mathbf{f}_{1},\mathbf{l}_1\cdots,\mathbf{f}_{\mu},\mathbf{l}_{\mu}\}$ of $\mathbb E_{(\kappa_{\mathbf a},b_{\mathbf a})}$ and a basis $\mathbf{e}_{1},\cdots,\mathbf{e}_{d}$ of $Ker D\mathcal S_{(\kappa_{\mathbf a},b_{\mathbf a})}$. Let $\fkz=(\fkz_1,...,\fkz_d)$ be the coordinate system of $D\mathcal S_{(\kappa_{\mathbf a},b_{\mathbf a})}$, let $\mathbf{t}^*=(\mathbf{t}^*_1,...,\mathbf{t}^*_\mu)$ be the coordinate system of $\mathbb E_{(\kappa_{\mathbf a},b_{\mathbf a})}$. Recall that for bubble trees with nonstable domain we add some additional marked points ( see \S\ref{with bubble tree}).
By choosing $\varepsilon_{\mathbf{a}},\delta_{\mathbf{a}},\rho_{\mathbf{a}}$ small we can view  $(\mathbf{s}, \mathbf{t},\mathbf{t}^*, \fkz)$ with each $|\ft_{i}|\neq 0$ and each $ |\ft^{*}_{j}|\neq 0$ as a local coordinate system in $ \mathbf{U}^{T}_{\infty,(\kappa_{\mathbf{a}},b_{\mathbf{a}})}(2\varepsilon_{\mathbf{a}},2\delta_{\mathbf{a}},2\rho_{\mathbf{a}})$, where $(\mathbf{s}, \mathbf{t})=\left((\mathbf{s}^+, \mathbf{t}^+), (\mathbf{s}^-, \mathbf{t}^-)\right)$.
\v
Denote $\sigma(\varrho)=e^{-\frac{\varrho}{L}}.$ Let $D_{\mathbf a}=\{ (\mathbf{s}, \mathbf{t},\mathbf{t}^*, \fkz) \in \mathbf{U}_{\infty,(\kappa_{\mathbf{a}},b_{\mathbf{a}})}(\varepsilon_{\mathbf{a}},\delta_{\mathbf{a}},\rho_{\mathbf{a}})\;|\;  |\ft_{i}|\leq  \sigma(\varrho)\}.$ If $a\leq \mathbf n_{t}$, $D_{\mathbf a}=  \mathbf{U}_{\infty,(\kappa_{\mathbf{a}},b_{\mathbf{a}})}(\varepsilon_{\mathbf{a}},\delta_{\mathbf{a}},\rho_{\mathbf{a}})$. Denote $\mathbf U'_{\infty,\sigma(\varrho)}=\bigcup_{\mathbf a} D_{a}.$  We smooth the corner of $\mathbf U'_{\infty,\sigma(\varrho)}$ and denote the resulting by  $\mathbf U_{\infty,\sigma(\varrho)}.$
Then  we construct a family of smooth hypersurfaces $S_{\infty,\sigma(\varrho)}:=\p \mathbf U_{\infty,\sigma(\varrho)}\subset \mathbf U_{\infty}^{T},$   such that
  $$S_{\sigma(\varrho)}\cap \mathbf{U}^{T}_{\infty,(\kappa_{\mathbf{a}},b_{\mathbf{a}})}(\varepsilon_{\mathbf{a}},\delta_{\mathbf{a}},\rho_{\mathbf{a}})\subset \{(\mathbf{s}, \mathbf{t},\mathbf{t}^*, \fkz) \in \mathbf{U}_{\infty,(\kappa_{\mathbf{a}},b_{\mathbf{a}})}(\varepsilon_{\mathbf{a}},\delta_{\mathbf{a}},\rho_{\mathbf{a}})\;|\; \min |\ft_{i}|\geq  c\sigma(\varrho)\},$$
where $c\in (0,1)$ is  a constant independent of $\varrho.$
$S_{\sigma(\varrho)}$ divide $\mathbf{U}_{\infty}^{T}$ into two parts.  Denote  by $\mathbf {U}_{\infty,\epsilon}^{T,\sigma(\varrho)}$  the part which is relatively far away from lower strata. Set $$
\mathbf{U}^{T,\sigma(\varrho)}_{\infty,(\kappa_{\mathbf{a}},b_{\mathbf{a}})}(\varepsilon_{\mathbf{a}},\delta_{\mathbf{a}},\rho_{\mathbf{a}})=\mathbf {U}_{\infty,\epsilon}^{T,\sigma(\varrho)}\cap \mathbf{U}^{T}_{\infty,(\kappa_{\mathbf{a}},b_{\mathbf{a}})}(\varepsilon_{\mathbf{a}},\delta_{\mathbf{a}},\rho_{\mathbf{a}})
$$
By the Lemma \ref{gluing map-1}, Lemma \ref{gluing map-2} and the same argument for gluing nodes between $\Sigma^+$ and $\Sigma^-$ we immediately obtain
\begin{lemma}\label{gluing map_rho-1} For each $(\kappa_{\mathbf{a}},b_{\mathbf{a}})$, there exist  positive constants   $\varrho_0, \varepsilon,\delta$ and $\rho$ such that for any $\varrho>\varrho_{o}$
$$glu_{\varrho,(\kappa,b)}:\mathbf{U}^{T,\sigma(\varrho)}_{\infty,(\kappa_{\mathbf{a}},b_{\mathbf{a}})}(\varepsilon_{\mathbf{a}},\delta_{\mathbf{a}},\rho_{\mathbf{a}})
\to glu_{\varrho,(\kappa,b)}\left(\mathbf{U}^{T,\sigma(\varrho)}_{\infty,(\kappa_{\mathbf{a}},b_{\mathbf{a}})}(\varepsilon_{\mathbf{a}},\delta_{\mathbf{a}},\rho_{\mathbf{a}})\right)\subset \mathbf {U}_{\varrho,2\epsilon}^{T}$$
is an orientation preserving local diffeomorphisms in orbifold sense.
  \end{lemma}
\v
Then we can choose $(\fs,\ft,\mathbf{t}^*,\fkz)$ with each $|\ft_{i}|\neq 0$ and each $ |\ft^{*}_{j}|\neq 0$ as  a local coordinates of $glu_{\varrho}\left( \mathbf{U}^{T,\sigma(\varrho)}_{\infty,(\kappa_{\mathbf{a}},b_{\mathbf{a}})}
(2\varepsilon_{\mathbf{a}},2\delta_{\mathbf{a}},2\rho_{\mathbf{a}})\right).$
By the same argument of Lemma \ref{virtual_varrho-1} we can prove that
\begin{lemma}
$\{\mathbf U^{T}_{\varrho,(\kappa_{\mathbf a,\varrho},b_{\mathbf a,\varrho})}(\frac{3}{2}\varepsilon_{\mathbf{a}},\frac{3}{2}\delta_{\mathbf{a}},\frac{3}{2}\rho_{\mathbf{a}}), \;1\leq \mathbf{a}\leq \mathbf{n}_c\}$ is an open covering of $\mathbf U^{T}_{\varrho,\varepsilon}$ for $\varrho $ big enough.
\end{lemma}
\v

For any small $  \sigma',$ we construct a   smooth hypersurfaces $S_{\varrho,\sigma'}\subset \mathbf U_{\varrho}^{T}$ such that
  $$S_{ \varrho,\sigma'}\cap \mathbf U^{T}_{\varrho,(\kappa_{\mathbf a,\varrho},b_{\mathbf a,\varrho})}(\tfrac{3}{2}\varepsilon_{\mathbf{a}},\tfrac{3}{2}\delta_{\mathbf{a}},\tfrac{3}{2}\rho_{\mathbf{a}})\subset \{(\mathbf{s}, \mathbf{t},\mathbf{t}^*, \fkz_{\varrho}) \in \mathbf{U}_{\varrho,(\kappa_{\mathbf{a},\varrho},b_{\mathbf{a},\varrho})}(\varepsilon_{\mathbf{a}},\delta_{\mathbf{a}},\rho_{\mathbf{a}})\;|\; c\sigma'< \min|\ft_{i}|\},$$
where $c\in (0,1)$ is  a constant independent of $\varrho.$
$S_{\varrho,\sigma'}$ divide $\mathbf{U}_{\varrho}^{T}$ into two parts.  Denote  by $\mathbf {U}_{\varrho,\epsilon}^{T,\sigma'}$  the part, which relatively far away from the lower stratas.

\v

 \section{ Gluing estimates}

  As in \cite{LS-1} and \cite{LS-3} we can prove that
\begin{lemma}\label{gluing_rho_est}  Let $l\in \mathbb Z^+$ be a fixed integer.
There exists positive  constants  $\mathsf{C}_{l}, \mathsf{d}, R_{0}$ and a neighbor $O_{o}\subset \mathbb K$ of $0$ such that for any $(\kappa,\xi)\in \ker D \mathcal{S}_{(\kappa_{o},b_{o})}$ with $\|(\kappa,\xi)\|< \mathsf{d}$, restricting to the compact set $\Sigma(R_0)$, the following estimate hold
$$ \left\|\frac{\p}{\p \varrho}\left(X_{i}\left(glu_{\varrho}(\kappa,\xi) \right) \right) \right\|_{C^{l}(\Sigma(R_{0}))}
	\leq  \mathsf{C}_{l}e^{-\fc_{1}\varrho} e^{-(\fc-5\alpha)\tfrac{r_{i}}{4} },$$
$$ \left\|\frac{\p}{\p \varrho}\left(X_{i}X_{j}\left(glu_{\varrho}(\kappa,\xi) \right) \right) \right\|_{C^{l}(\Sigma(R_{0}))}
	\leq  \mathsf{C}_{l}e^{-(\fc-5\alpha)\tfrac{\varrho}{4} } e^{-(\fc-5\alpha)\tfrac{r_{i}+r_{j}}{4} },$$
$$ \left\|\frac{\p}{\p \varrho}\left(Y\left(glu_{\varrho}(\kappa,\xi) \right) \right) \right\|_{C^{l}(\Sigma(R_{0}))}+\left\|\frac{\p}{\p \varrho}\left(YZ\left(glu_{\varrho}(\kappa,\xi) \right) \right) \right\|_{C^{l}(\Sigma(R_{0}))}
	\leq  \mathsf{C}_{l}e^{-\fc_{1}\varrho},$$
for any $X_{i}\in \{\frac{\p}{\p r_{i}},\frac{\p}{\p \tau_{i}}\},$$1\leq i\neq j\leq\mathbf{e},$
 $Y,Z\in \{\frac{\p}{\p \mathbf{t}^*_{i}},\frac{\p}{\p \bar {\mathbf{t}}^*_{i}},\;\frac{\p }{\p s_{i}},\frac{\p }{\p \bar s_{i}}, \frac{\p }{\p \fkz_{l}}\},$ $\mathbf{s}\in \bigotimes_{l=1}^{\iota} O_l$ when $\varrho,|\mathbf r|$ big enough.
\end{lemma}

\begin{lemma}\label{gluing_rho_est_L}
Let $l\in \mathbb Z^+$ be a fixed integer.
 Let $u:\Sigma\to M$ be a $(j,J)$-holomorphic map.	Let $\fc\in (0,1)$ be a fixed constant.    For any $0<\alpha<\frac{1}{100\fc}$, there exists positive  constants  $\mathsf C_{l}, \mathsf{d},R$ such that for any  $\zeta\in \ker D^{\widetilde{\mathbf L}}|_{b_{o}}$,   $(\kappa,\xi)\in \ker D \mathcal{S}_{(\kappa_{o},b_{o})}$ with    $$\|\zeta\|_{\mathcal W,k,2,\alpha}\leq \mathsf{d},\;\;\;\;\;\|(\kappa,\xi)\|< \mathsf{d},\;\;\;\;\;|\mathbf r|\geq R, $$ restricting to the compact set $\Sigma(R_0)$, the following estimate hold.
 $$ \left\|\frac{\p}{\p \varrho}\left(X_{i}\left(Glu_{\mathbf{s},h_{(\mathbf r)},(\mathbf{r})}^{\widetilde{\mathbf L},\varrho}(\zeta) \right) \right) \right\|_{C^{l}(\Sigma(R_{0}))}
	\leq  \mathsf{C}_{l}e^{-\fc_{1}\varrho} e^{-(\fc-5\alpha)\tfrac{r_{i}}{4} },$$
$$ \left\|\frac{\p}{\p \varrho}\left(X_{i}X_{j}\left(Glu_{\mathbf{s},h_{(\mathbf r)},(\mathbf{r})}^{\widetilde{\mathbf L},\varrho}(\zeta)  \right) \right) \right\|_{C^{l}(\Sigma(R_{0}))}
	\leq  \mathsf{C}_{l}e^{-(\fc-5\alpha)\tfrac{\varrho}{4} } e^{-(\fc-5\alpha)\tfrac{r_{i}+r_{j}}{4} },$$
$$ \left\|\frac{\p}{\p \varrho}\left(Y\left(Glu_{\mathbf{s},h_{(\mathbf r)},(\mathbf{r})}^{\widetilde{\mathbf L},\varrho}(\zeta)  \right) \right) \right\|_{C^{l}(\Sigma(R_{0}))}+\left\|\frac{\p}{\p \varrho}\left(YZ\left(Glu_{\mathbf{s},h_{(\mathbf r)},(\mathbf{r})}^{\widetilde{\mathbf L},\varrho}(\zeta)  \right) \right) \right\|_{C^{l}(\Sigma(R_{0}))}
	\leq  \mathsf{C}_{l}e^{-\fc_{1}\varrho},$$
for any $X_{i}\in \{\frac{\p}{\p r_{i}},\frac{\p}{\p \tau_{i}}\},$$1\leq i\neq j\leq\mathbf{e},$
 $Y,Z\in \{\frac{\p}{\p \mathbf{t}^*_{i}},\frac{\p}{\p \bar {\mathbf{t}}^*_{i}},\;\frac{\p }{\p s_{i}},\frac{\p }{\p \bar s_{i}}, \frac{\p }{\p \fkz_{l}}\},$ $\mathbf{s}\in \bigotimes_{l=1}^{\iota} O_l$ when $\varrho,|\mathbf r|$ big enough.

\end{lemma}

\begin{remark}
$\varphi_{\varrho}^{*}\omega_{R}$ is a symplectic form of $M_{\varrho}.$ Recall that
$M_{\varrho}$ is gluing by $M^{+}$ and $M^{-}$ with
$$
a^{+}=a^{-}+2l\varrho,\;\;\;\theta^{+}=\theta^{-}.
$$
Since the  almost structure of $M^{+},M^{-}$ defined by \eqref{complex_structure_I} and \eqref{complex_structure_II},
we can choose the almost structure $J$ of $M_{\varrho}$ over gluing part  with
$$
J\frac{\p}{\p a^{+}} =\frac{\p}{\p \theta^{+}},\;\;\;J \frac{\p}{\p \theta^{+}}=\frac{\p}{\p a^{+}},\;\;\;\;J|_{\xi}=\widetilde{J}\;\;\;\; in  \;\;\;\{ lR \leq a^{+}\leq 2l\varrho-lR\},
$$
which is restrict of the  almost structure of $M^{+}$.
It is easy to see that  $J_{\varrho}=(\varphi_\varrho^{-1})^{*}J$ is a  family  smooth almost structure on $M_{R}$.
\v
Let $b_{o}=(b^{+}_{o},b^{-}_{o})$ as in section \ref{equ_glu}.
 By gluing we can obtain $u_{(\mathbf r)}$ with $(\mathbf r)$ satifying \eqref{gluing_local_relative_node-3}.
 For any $h\in W^{k,2,\alpha}(\Sigma_{(\mathbf r)},u^{*}_{(\mathbf r)}TM_{\varrho}),$ we can define
 $$
 h^{*}=(h^{+},h^{-})=\left(\beta\left(\tfrac{3}{2}-\tfrac{s_1^{+}}{r_{j}}\right)h,\; \left[1-\beta^2\left(\tfrac{3}{2}-\tfrac{s_j^{+}}{r_{j}}\right)\right]^{\frac{1}{2}}h\right),\;\;\;for \;|s_{j}^{\pm}|\geq R_{0}
 $$
 which can been seen as a element $W^{k,2,\alpha}(\Sigma,u^{*}TM_{\infty}).$
 Let $e_{u}$ (resp. $e_{u_{(\mathbf r)}}$) be the local frame of $\mathbf L|_{b_{o}}$ (resp. $\mathbf L|_{b_{(\mathbf r)}}$) near the nodal points.
 Then by cut-off function a element in  $W^{k,2,\alpha}(\Sigma_{(\mathbf r)},\mathbf L|_{b_{(\mathbf r)}})$ can be seen as $ {\mathcal W}^{k,2,\alpha}(\Sigma,\mathbf L|_{b_{o}})$
Similar  in \cite{LS-1} and \cite{LS-3} we can define two maps
$$
 glu_{\varrho}^{*}: Ker D\mathcal S_{(\kappa_{o}, b_{o})}\to  W^{k,2,\alpha}(\Sigma,u^*TM)$$
$$
 Glu_{\mathbf{s},h_{(\mathbf r)},(\mathbf{r})}^{\widetilde{\mathbf L},\varrho,*} : Ker D^{\mathbf{L}}_{b_{o}} \to  {\mathcal W}^{k,2,\alpha}(\Sigma,\mathbf L|_{b_{o}})
$$

 Then we can obtain the Theorem
 \ref{gluing_rho_est} and   \ref{gluing_rho_est_L}
 by repeating the argument of \cite{LS-1} and \cite{LS-3}.
\end{remark}

\subsection{A metric of $\mathbf E_{\varrho}$}
Using Lemma \ref{gluing_rho_est}  we can prove that
 \begin{lemma}\label{lemma_rho_9.8}
For any $\sigma'>0,$ there exists a positive constant  $\varrho_{o} $   such that for any $\varrho>\varrho_{o}$
\begin{itemize}
\item[(1)]
 $\{ glu_{\varrho}(\mathbf{U}^{T,\sigma(\varrho)}_{\infty,(\kappa_{\mathbf{a}},b_{\mathbf{a}})}(2\varepsilon_{\mathbf{a}},2\delta_{\mathbf{a}},2\rho_{\mathbf{a}}))
, \;1\leq \mathbf{a}\leq \mathbf{n}_c\}$ is an open covering of
 $\mathbf{U}^{T,\sigma'}_{\varrho,\varepsilon}$,
\item[(2)] Denote $\hat{\mathbf {\Gamma}}_{\varrho,\mathbf{a}}= (glu_{\varrho}^{-1})^{*}\mathbf{\Gamma}_{\infty,\mathbf{a}}.$ Then we have
    $$|\sum \hat{\mathbf \Gamma}_{\varrho,\mathbf{a}}-1|\leq Ce^{-\fc_{1} \varrho},\;\;\;\; \lim_{\varrho\to \infty}\sum \hat{\mathbf {\Gamma}}_{\varrho,\mathbf{a}}=1,$$
\end{itemize}
\end{lemma}
\v
{\bf Proof.} (1). If  (1) is not true, we can find   a sequence
 $$b^{(i)} =(\Sigma^{(i)},j^{(i)},{\bf y}^{(i)},u^{(i)})\in \mathbf{U}^{T,\sigma'}_{\varrho^{(i)},\varepsilon}\setminus   \bigcup_{C\in \mathcal C_{g,m}^{J,A}} \bigcup\limits_{1\leq  \mathbf{a}\leq \mathfrak{m}_{C} }  glu_{\varrho^{(i)}}(\mathbf{U}^{T,\sigma^{(i)}}_{\infty,(\kappa_{\mathbf{a}},b_{\mathbf{a}})}(\varepsilon_{\mathbf{a}},\delta_{\mathbf{a}},\rho_{\mathbf{a}})),$$
such that $$\lim_{i\to \infty}\varrho^{(i)}=\infty,\;\;\;\;\;\lim_{i\to \infty}\sigma^{(i)}=0.$$
By the convergence theorem \ref{convergence} and the definition of $S_{\varrho,\sigma'}$  we conclude that $b^{(i)}$ weak converges to some $b=(\Sigma,j,{\bf y},u)\in    \mathbf{U}^{T}_{(\kappa_{\mathbf{a}},b_{\mathbf{a}})}(\varepsilon_{\mathbf{a}},\delta_{\mathbf{a}},\rho_{\mathbf{a}}) $ for some $1\leq \mathbf{a}\leq \mathbf n_{c}$.   Let $b_{\mathbf{a}}=(\Sigma,j_{\mathbf a},{\bf y}_{\mathbf a},u_{\mathbf a}).$
We may assume that $(\Sigma,j,{\bf y})$ is stable.
 We choose   $(\mathbf{s}, \mathbf{t},\mathbf{t}^*, \fkz)$ as a local coordinate system in $ \mathbf{U}_{(\kappa_{\mathbf{a}},b_{\mathbf{a}})}(2\varepsilon_{\mathbf{a}},2\delta_{\mathbf{a}},2\rho_{\mathbf{a}})$  with $(\fs,\ft,\ft^{*})(j_{\mathbf{a}},\mathbf y_{\mathbf{a}})=0$, $|\ft_{i}(j,y)|\geq c\sigma'$  and $|\ft_{l}(j^{(i)},y^{(i)})|\geq c\sigma'.$
\v
 Let $q_{1},\cdots,q_{\mu}$ be relative nodes. Choose cusp cylinder  coordinates $(s^{\pm}_{j},t^{\pm}_{j})$ near each relative node $q_{i}$.   Suppose that $u$ satisfies \eqref{eqn_conv-pm}. Then there exist $h\in ker D\mathcal S_{(\kappa_{\mathbf a},b_{\mathbf a})}$ and $\zeta\in \mathbb E$ such that
$$v=\exp_{u}(\zeta),\;\;\;u=\exp_{v}(I_{b_{\mathbf{a}},b_{v}}(h)+Q_{(\kappa_{o},b_{v})}\cdot f\cdot  I_{b_{\mathbf{a}},b_{v}}(h)) $$ with $b_{v}=(\Sigma,j,{\bf y},v)$.
Then $\Sigma^{(i)}=\Sigma_{\fs^{(i)},\ft^{(i)},(\ft^{*})^{(i)})},$ such that  each $|\ft^{(i)}_{j}|>c\sigma'$ and each $ |(\ft^{*}_{j})^{(i)}|> 0$ and $\lim\limits_{i\to \infty}\fs^{(i)}=0.$  Let   $((r^{*}_{j})^{(i)},(\tau^{*}_{j})^{(i)})$ be is gluing parameter of the relative node $q_{i}$, with $2l(r^{*}_{j})^{(i)} +2\pi\sqrt{-1}(\tau^{*}_{j})^{(i)}=  (\ft^{*}_{j})^{(i)}+\frac{2l\varrho^{(i)}}{k}  $. Denote $$(\mathbf r^{*})^{(i)}=((r^{*}_{1})^{(i)},(\tau^{*}_{1})^{(i)},\cdots,((r^{*}_{\mathfrak{e}})^{(i)},(\tau^{*}_{\mathfrak{e}})^{(i)})).$$

  \v
  Next we prove that
  \begin{equation}\label{rho_diff}
\lim_{i\to \infty}(2l\varrho^{(i)}-2lk_{j}(r_{j}^{*})^{(i)})= \ell_{j}^{+}-\ell_{j}^{-},\;\;\;\ \lim_{i\to\infty}k_{j}(\tau^{*}_{j})^{(i)}=\theta_{j0}^{-}-  \theta_{j0}^{+} .
 \end{equation}
\v
 Denote $u^{(i)}=(a^{(i)},\theta^{(i)},\widetilde{u}^{(i)}).$ Since $u^{(i)}$ define in $(s_{j}^{+},t_{j}^{+})\in [-R_{0},2l(r_{j}^{*})^{(i)}-R_{0}]\times S^{1}\subset \Sigma^{(i)},$ by Theorem \ref{tube_L_decay}  and a direct integration we have for any $R_{1}>R_{0}>0$
  $$
   |a^{(i)}_{1}(R_{1},t^{+}_{j})-k_{j}R_{1}-a^{(i)}_{1}(l(r_{j}^{*})^{(i)},t_{j}^{+})+k_{j}l(r_{j}^{*})^{(i)}|\leq Ce^{-\fc(R_{1}-R_{0})},\;\;\;\;$$$$ |a^{(i)}_{2}(-R_{1},t_{j}^{-})+k_{j}R_{1}-a^{(i)}_{2}(-l(r_{j}^{*})^{(i)},t_{j}^{-})-k_{j}l(r_{j}^{*})^{(i)}|\leq Ce^{-\fc(R_{1}-R_{0})}
  $$
where $t_{j}^{+}=t_{j}^{-}$.
By $a^{(i)}_{1}(l(r_{j}^{*})^{(i)},t_{j}^{+})=a^{(i)}_{2}(lr^{(i)}_{j},t_{j}^{+})+2l\varrho^{(i)}$ and $s_{j}^{+}=s_{j}^{-}+2l(r_{j}^{*})^{(i)}$ we have
  $$
    \left|a^{(i)}_{1}(R_{1},t_{j}^{+})-2k_{j}R_{1}-2l\varrho^{(i)}+2lk_{j}(r_{j}^{*})^{(i)}-\a^{(i)}_{2}(-R_{1},t_{j}^{-})\right|\leq Ce^{-\fc(R_{1}-R_{0})}
  $$
By taking limits $i\to +\infty$
$$
\left|a_{1}(R_{1},t_{j}^{+})- k_{j}R_{1}-(\a_{2}(-R_{1},t_{j}^{+})+k_{j}R_{1})- \lim_{i\to\infty}2l(\varrho^{(i)}-k(r_{j}^{*})^{(i)})\right|\leq Ce^{-\fc(R_{1}-R_{0})}.
  $$
 We get the first equality of  \eqref{rho_diff} as $R_{1}\to \infty.$
 Similar we have the second equality of  \eqref{rho_diff}.
\v
By pre-gluing we get $v_{(\mathbf r^{*})^{(i)}}.$
Denote
$  u^{(i)}=\exp_{v_{( \mathbf r^{*})^{(i)}}}(\xi^{(i)})$ and $b_{v}^{(i)} =(\Sigma^{(i)},j^{(i)},{\bf y}^{(i)},v_{(\mathbf r^{*})^{(i)}}^{(i)}).$ By the convergence of $b^{(i)}$    we have for any $R_{1}>R_{0}$,
$$
\lim_{i\to \infty }\|\xi^{(i)}-I_{b_{\mathbf{a}},b_{v}}(h)+Q_{(\kappa_{o},b_{v})}\cdot f_{b_{v}}\cdot  I_{b_{\mathbf{a}},b_{v}} (h)|_{\mathfrak{D}(R_1)}\|=0.
$$
By Lemma \ref{coordinate_decay-2}  we can conclude that
$$
\left\|\left.I_{b_{\mathbf{a}},b_{v}}(h)+Q_{(\kappa_{o},b_{v})}\cdot f_{b_{v}}\cdot  I_{b_{\mathbf{a}},b_{v}} (h)-I_{b_{\mathbf{a}},b_{v^{(i)}}}(h)+Q_{(\kappa_{o},b_{v^{(i)}})}\cdot f_{b_{v^{(i)}}}\cdot  I_{b_{\mathbf{a}},b_{v^{(i)}}} (h)\right|_{|s^{\pm}_{j}|\geq R_{1}}\right\|_{C^{l}(\Sigma(R_{1}))}\leq  \epsilon
$$
Using Lemma \ref{exponential_estimates_theorem} and  Lemma \ref{tube_L_decay}   we have for any $\varepsilon>0,$
\begin{equation}\label{eqn_deg-diff}
\|\xi^{(i)}-I_{b_{\mathbf{a}},b_{v^{(i)}}}(h)+Q_{(\kappa_{o},b_{v^{(i)}})}\cdot f_{b_{v^{(i)}}}\cdot  I_{b_{\mathbf{a}},b_{v^{(i)}}} (h)\|_{k,2,\alpha,r}\leq  2\varepsilon
\end{equation}
 as $i$ big enough.
 Applying the implicit function theorem,  we get contradiction.
\v\n
  (2) We only to prove (2) in each $ glu_{\varrho}(\mathbf{U}^{T,\sigma(\varrho)}_{\infty,(\kappa_{\mathbf{a}},b_{\mathbf{a}})}(2\varepsilon_{\mathbf{a}},2\delta_{\mathbf{a}},2\rho_{\mathbf{a}})).$
  We choose the local coordinates $(\fs,\ft,\ft^{*},\fkz)$ as before.
   Denote   $\hat{\mathbf {\Gamma}}_{\varrho,\mathbf{a}}=(glu^{-1})^{*}\mathbf{\Gamma}_{\infty,\mathbf{a}}.$  It is easy to see that in this  coordinates,
  $$
  \frac{\p \hat{\mathbf {\Gamma}}_{\varrho,\mathbf{a}}}{\p \varrho}=0.
  $$
Assume that   $ glu_{\varrho}(\mathbf{U}^{T,\sigma(\varrho)}_{\infty,(\kappa_{\mathbf{a}},b_{\mathbf{a}})}(2\varepsilon_{\mathbf{a}},2\delta_{\mathbf{a}},2\rho_{\mathbf{a}})) \cap  glu_{\varrho}(\mathbf{U}^{T,\sigma(\varrho)}_{\infty,(\kappa_{\mathbf{a'}},b_{\mathbf{a'}})}(2\varepsilon_{\mathbf{a'}},2\delta_{\mathbf{a'}},2\rho_{\mathbf{a'}}))\neq \emptyset$.
Using Lemma  \ref{gluing_rho_est} and by a direct calculation we have
$$
\left|\frac{\p}{\p \varrho} \hat{\mathbf {\Gamma}}_{\varrho,\mathbf{a}'}\right|\leq e^{-\fc_{1}\varrho}
$$
where we use the smoothness of cut-off function.  Then (2) follows.\;\;\;$\Box$

\v

As in section \S\ref{a metric} we choose local metric $h_{\mathbf{a}}$ on
$\mathbf{U}_{\infty,(\kappa_{\mathbf{a}},b_{\mathbf{a}})}(\varepsilon_{\mathbf{a}},\delta_{\mathbf{a}},\rho_{\mathbf{a}})$.
We can define the metric $\mathbf h_{\infty}$ and a connection $\nabla^{L}_{\infty}$ of $E_{\infty}$ as  in section \S\ref{a metric} .
Denote $\mathbf{\Gamma}_{\varrho,\mathbf{a}}=\frac{\hat{\mathbf {\Gamma}}_{\varrho,\mathbf{a}}}{\sum_{\mathbf{a}} \hat{\mathbf {\Gamma}}_{\varrho,\mathbf{a}}}.$ Then by Lemma \ref{lemma_rho_9.8} $\mathbf{\Gamma}_{\varrho,\mathbf{a}}$ is partition of unity of $\mathbf{U}_{\varrho,\epsilon}^{T,\sigma'}$.

 We define a metric $\mathbf{h}^{\varrho}$ on $\mathbf E_{\varrho}|_{\mathbf U^{T,\sigma'}_{\varrho,\varepsilon}}$  by
$$\mathbf{h}^{\varrho}=\sum_{\mathbf{a}=1}^{\mathbf{n}_c} {\mathbf{\Gamma}_{\varrho,\mathbf{a}} } ((Glu_{\mathbf{s},h_{(\mathbf r)},(\mathbf{r})}^{\widetilde{\mathbf L},\varrho})^{-1})^{*}\mathbf h_{\mathbf{\mathbf a}}.
$$

 Next we define  a family  connection on $\mathbf E_{\varrho}$.
Let $\{e^{\mathbf a}_{\alpha}\}_{1\leq\alpha\leq \mathsf r_{j}}$ be the smooth orthonormalization frame field of $\mathbf E_{\infty}$ in $\widetilde{\mathbf{U}}^{T}_{\infty,(\kappa_{\mathbf{a}},b_{\mathbf{a}})}(\varepsilon_{\mathbf{a}},\delta_{\mathbf{a}},\rho_{\mathbf{a}})$.   Denote $e^{\mathbf a,\varrho}_{\alpha}= (Glu_{\mathbf{s},h_{(\mathbf r)},(\mathbf{r})}^{\widetilde{\mathbf L},\varrho})_{*}{e}^{\mathbf a}_{\alpha}$.
 Consider the  Gram-Schmidt process with respect to the metric $\mathbf h^{\varrho}$ and denote by   $\hat{e}^{\mathbf a,\varrho}_{1}, ..., \hat{e}^{\mathbf a,\varrho}_{\mathsf r_{i}}$ the Gram-Schmidt orthonormalization of $\{e^{\mathbf a,\varrho}_{\alpha}\}$.  We define a local connection $\nabla^{\mathbf a}$ by
$$
\nabla^{\mathbf a}\hat{e}^{{\mathbf a}}_{\alpha}=0,\;\;\;\;\;\alpha=1,\cdots,\mathsf r_{i}.
$$
Then we can define $\nabla^{\varrho}$ as before.
It is easy to see that $\nabla^{\varrho}$ is a compatible connection of the metric $\mathbf{h}^{\varrho}$. Denote $$\nabla^{\varrho} \hat e^{\mathbf a,\varrho}_{\alpha}=\sum_{\beta} \omega_{\alpha \beta}^{\mathbf {a},\varrho}\hat e^{\mathbf a,\varrho}_{\beta},\;\;\;(\nabla^{\varrho})^2\hat e^{\mathbf a}_{\alpha}=\sum  \Omega_{\alpha\beta}^{\mathsf a,\varrho}\hat e^{\mathbf a,\varrho}_{\beta}.$$
For any $ glu_{\varrho} ({\mathbf{U}}^{T,\sigma(\varrho)}_{\infty,(\kappa_{\mathbf{a}},b_{\mathbf{a}})}(\varepsilon_{\mathbf{a}},\delta_{\mathbf{a}},\rho_{\mathbf{a}}))\bigcap  glu_{\varrho}({\mathbf{U}}^{T,\sigma(\varrho)}_{\infty,(\kappa_{\mathbf{c}},b_{\mathbf{c}})}(\varepsilon_{\mathbf{c}},\delta_{\mathbf{c}},\rho_{\mathbf{c}}))\neq \emptyset,$ let $(\hat a^{\mathbf a\mathbf c}_{\alpha\beta})_{1\leq \alpha,\beta\leq \mathsf r}$ be   functions such that $\hat e^{\mathbf a,\varrho}_{\alpha}=\sum_{\beta=1}^{\mathsf r_{i}} \hat a^{\mathbf a\mathbf c}_{\alpha \beta}\hat e^{\mathbf c,\varrho}_{\beta},\alpha=1,\cdots,\mathsf r.$   It is easy to see that
\begin{equation}\label{local_omega}
\omega^{\mathbf a}_{\alpha\beta}=\sum_{\mathbf c}\sum_{\beta=1}^{\mathsf r_{i}} \Gamma_{\varrho,\mathbf c}d\hat a^{\mathbf a\mathbf c}_{\alpha \gamma}\hat a^{\mathbf c\mathbf a}_{\gamma \beta}.
\end{equation}
\n

Using Lemma \ref{gluing_rho_est}  and Lemma \ref{gluing_rho_est_L}, by the  same argument of section \S\ref{est_Thom_E} we can prove that
\begin{lemma}\label{est_Thom_varrho_L}
	There exists a constant $C>0$ such that in each ${\mathbf{U}}^{T,\sigma(\varrho)}_{(\kappa_{\mathbf{a}},b_{\mathbf{a}})}(\varepsilon_{\mathbf{a}},\delta_{\mathbf{a}},\rho_{\mathbf{a}}))$
$$	|\sigma^{*}(glu^{*}_{\varrho}\Theta_{\varrho}-\Theta_{\infty}) (X_{1},\cdots,X_{\mathsf r})|^2\leq C e^{-\fc_{1}\varrho}\Pi_{i=1}^{\mathsf r}g_{loc}(X_{i},X_{i})$$for any $X_{i}\in T\mathbf U_{\infty}^{T},i=1,2,3.$
\end{lemma}

Let $\Theta_{\infty,\mathbf E}$ (resp. $\Theta_{\varrho,\mathbf E}$) be the Thom form of $\mathbf E|_{\mathbf {U}_{\infty,\epsilon}}$ (resp.  $\mathbf E|_{\mathbf {U}_{\varrho,\epsilon}}$)  supported in a small $\varepsilon$-ball of the $0$-section of $\mathbf{E}$.

\section{The relation of $\Psi_{(M_{\infty},A,g,m)}$ and $\Psi_{(M_{R},A,g,m)}$ }

The following result is well-known:

\begin{lemma} \label{lem_cob}
For any fixed $R$ with  $R_0\leq R<\infty$ and any $\varrho > R$ we have
$$
\Psi_{(M_{R},A,g,m)}=\Psi_{(M_{\varrho},A,g,m)}.
$$
\end{lemma}
In this section we prove
\begin{theorem}\label{thm_re}
 For any R, $R_0<R<\infty$, we have
$$
\Psi_{(M_{\infty},A,g,m)}=\Psi_{(M_{R},A,g,m)}.
$$
\end{theorem}
\v
To prove this theorem we first introduce some notations.
\v
Let $C\in {\mathcal  C}^{J,A}_{g,m}$, we consider $\overline{\mathcal{M}}_C$. Note that the data $C$ gives

\begin{itemize}
\item[(1)] a fixed partition of the index set $\{1,\cdots, m \}=S^+\cup S^-$,
\item[(2)] a fixed partition of the index set $\{1,\cdots, \mu \}$, a map $\rho:\{p^+_1,...,p^+_{\mu}\}\rightarrow
\{p^-_1,...,p^{-}_{\mu}\}$,
\item[(3)] a fixed partition of $A$..
\end{itemize}
If we forget the data of the partition of $A$ we get a data denoted by $C'$.
Denote by $\overline{\mathcal{M}}_{C'}$ the moduli space of Riemann surfaces corresponding to $C'$. Let $\theta_{C'}:\overline{\mathcal{M}}_{C'}\to \overline{\mathcal{M}}_{g,m+\mu}$ be the embbeding submanifold.
We define the GW-invarians $\Psi_{(M_{\infty},A, g,m)}(K^+\times K^-; \{\alpha_i\})$  as
\begin{equation}\label{def_GW_sp}
\Psi_{(M_{\infty},A, g,m)}(K^+\times K^-; \{\alpha_i\})=\sum_{C\in {\mathcal
C}^{J,[A]}}
\int_{ \mathbf U^T_{c,\varepsilon} }\mathscr{P}^*(K_{C'})\wedge \prod_j ev'^*_j\alpha_j  \wedge \sigma_{c}^*\Theta_{c},
\end{equation}
where   $\Theta_{c}$  is the Thom form of   $\pi:\mathbf E_{c}\to \mathbf U_{c},$  $ev'_j$ denote the evaluation map
  $ev'_j: {\bf U'}_{c,\varepsilon} \longrightarrow
M^{\pm}$  at $j$-th marked point. We have
$$\Psi_{(M_\varrho,A,g,n)}(K^+\times K^-;\{\alpha_i\}) =\int_{\mathbf{U}^T_{\varrho,\varepsilon}}
\mathscr{P}^* \left(  \pi^{*}K_{\overline{\M} }\right) \wedge e_{j}^* (\prod_i
\alpha_i)\wedge \Theta.$$

\v

\v

\v\n
{\bf Proof of Theorem \ref{thm_re}.}
Let $\alpha_i\in H^*(M , \mathbb R)$ with $\sum \deg(\alpha_i)=ind$.
Denote
$$\mathbb{F}_{\varrho}=\mathscr{P}^* \left(  \pi^{*}K_{\overline{\M} }\right) \wedge \prod_j e^*_j\alpha_j\wedge \sigma^{*}\Theta_{\varrho}.$$
To simplify notations we denote $\mathbf{U}^{T,\sigma(\varrho_{0})}_{\mathbf{a}}:=\mathbf{U}^{T,\sigma(\varrho_{0})}_{\infty,(\kappa_{\mathbf{a}},b_{\mathbf{a}})}(\varepsilon_{\mathbf{a}},\delta_{\mathbf{a}},\rho_{\mathbf{a}}),$
$\mathbf{U}^{T}_{\mathbf{a}}:=\mathbf{U}^{T,}_{\infty,(\kappa_{\mathbf{a}},b_{\mathbf{a}})}(\varepsilon_{\mathbf{a}},\delta_{\mathbf{a}},\rho_{\mathbf{a}}).$
Here  $e_j$ denote the evaluation map
$e_j: {\bf U}(\varrho,\varepsilon) \longrightarrow
M $ at $j$-th marked point. For any $\varrho\geq \varrho_0,$ we calculate
$$
\Psi_{(M_{\varrho_0},A,g,m)}(\pi^{*}K_{\mathcal M},\{\alpha_i\})= \int_{ \mathbf U^{T}(\varrho_0,\varepsilon) }\mathbb{F}_{\varrho_0}
= I + II ,$$
where
$$I= \int_{ \mathbf U^{T}(\varrho_0,\varepsilon)\setminus \bigcup_{\mathbf a} glu_{\varrho_0}(\mathbf U_{\mathbf a}^{T,\sigma(\varrho_{0})})  }\mathbb{F}_{\varrho_0} +  \int_{\bigcup_{\mathbf a}   glu_{\varrho_0}(\mathbf U_{\mathbf a}^{T,\sigma(\varrho_{0})})  }(1-\sum \Gamma_{\varrho_0,\mathbf a})\mathbb{F}_{\varrho_0},$$
$$II= \sum \int_{glu_{\varrho_0}(U^{T,\sigma(\varrho_{0})}_{\mathbf{a}})}\Gamma_{\varrho_0,\mathbf a}\mathbb{F}_{\varrho_0}= \sum \int_{U^{T,\sigma(\varrho_{0})}_{\mathbf{a}}}\Gamma_{\infty,\mathbf a}glu_{\varrho_0}^{*}(\mathbb{F}_{\varrho_{0}})
.$$
 It follows from (1) of Lemma \ref{lemma_rho_9.8} that $$\mathbf U^{T}(\varrho_0,\varepsilon)\setminus \bigcup_{\mathbf{a}} glu_{\varrho_0}(\mathbf U_{\mathbf a}^{T})\subset  \mathbf U^{T}(\varrho_0,\varepsilon)\setminus \mathbf U^{T,\sigma'}_{\varrho_0,\varepsilon}.$$ Then  by (2) of Lemma \ref{lemma_rho_9.8} and the proof of Theorem \ref{Conver} we have
$$\left|I\right|\leq C\left(e^{-\frac{\alpha l \varrho_0}{32}}+\frac{1}{|\log \sigma'|} \right).$$
Note  $ K$ and $\theta^* {K}_{\mathcal M}$ are   in the same cohomology. We have
$$
\Psi_{(M_{\infty}, A,g,n)}(K,\{\alpha _i\})=\left(\sum
\int_{ \mathbf U^T_{\mathbf a}\setminus  U^{T,\sigma(\varrho_{0})}_{\mathbf a}}+\sum \int_{  U^{T,\sigma(\varrho_{0})}_{\mathbf a}}\right)\Gamma_{\infty,\mathbf a}\mathscr{P}^* \left(\theta^* {K}_{\mathcal M}\right) \wedge\prod_j e'^*_j\alpha_j  \wedge \sigma^{*}\Theta.
$$
By the proof of Theorem \ref{Conver} we have
$$
\left|\sum
\int_{ \mathbf U^T_{\mathbf a}\setminus  U^{T,\sigma(\varrho_{0})}_{\mathbf a}}
\Gamma_{\infty,\mathbf a}\mathscr{P}^* \left(\theta^* {K}_{\mathcal M}\right) \wedge\prod_j e'^*_j\alpha_j  \wedge \sigma^{*}\Theta
 \right|\leq \frac{C}{\varrho_{0}}
$$
Next we estimates
\begin{align*}
&\int_{  U^{T,\sigma(\varrho_{0})}_{\mathbf a}}\Gamma_{\infty,\mathbf a}\left(glu_{\varrho_{0}}^{*}(\mathscr{P}_{\varrho_{0}}^* \left(  \pi^{*}K_{\overline{\M} }\right) \wedge \prod_j e^*_j\alpha_j\wedge \sigma^{*}\Theta_{\varrho_{0}})-\mathscr{P}^* \left(\theta^* {K}_{\mathcal M}\right) \wedge\prod_j e'^*_j\alpha_j  \wedge \sigma^{*}\Theta\right)\\
=&(III)+(IV)+(V)
\end{align*}
where
$$
(III)=\int_{ U^{T,\sigma(\varrho_{0})}_{\mathbf a}}\Gamma_{\infty,\mathbf a}\left( glu_{\varrho_{0}}^{*}(\mathscr{P}_{\varrho_{0}}^* \left(  \pi^{*}K_{\overline{\M} }\right))-  \mathscr{P}^* \left(\theta^* {K}_{\mathcal M}\right) \right) \wedge glu_{\varrho_{0}}^{*}\left(\prod_j e^*_j\alpha_j\wedge \sigma^{*}\Theta_{\varrho_{0}}\right)
$$
$$
(IV)=\int_{ U^{T,\sigma(\varrho_{0})}_{\mathbf a}}\Gamma_{\infty,\mathbf a}  \mathscr{P}^* \left(\theta^* {K}_{\mathcal M}\right)  \wedge \left(glu_{\varrho_{0}}^{*}(\prod_j e^*_j\alpha_j) -\prod_j e'^*_j\alpha_j  \right)\wedge glu_{\varrho_{0}}^{*}\sigma^{*}\Theta_{\varrho_{0}})
$$
$$
(V)=\int_{ U^{T,\sigma(\varrho_{0})}_{\mathbf a}}\Gamma_{\infty,\mathbf a}  \mathscr{P}^* \left(\theta^* {K}_{\mathcal M}\right)  \wedge \prod_j e'^*_j\alpha_j  \wedge\left(glu_{\varrho_{0}}^{*}(\sigma^{*}\Theta_{\varrho_{0}})-\sigma^{*}\Theta
\right)
$$
Denote $b_{\mathbf a}=(\Sigma,j,\mathbf{y},u).$ Suppose that $\Sigma$ has $\mathbf{e}$ relative nodes.
Let $(\fs,\ft)$ be the plumbing coordinates of $\mathscr{P}(U^{T,\sigma(\varrho_{0})}_{\mathbf a}).$
We can choose
$(\fs,\ft,\ft^{*},\fkz)$ as the local coordinates of $U^{T,\sigma(\varrho_{0})}_{\mathbf a}.$
By the definition we have
$$\mathscr{P}_{\varrho_{0}} glu_{\varrho_{0}}(\fs,\ft,\ft^*,\fkz)=(\fs,\ft,\ft_{o})  .$$
where
$$\ft_{o}=  e^{-2l\rho/k-\ft^*} $$  denotes gluing parameters near relative  nodes.
Then $(\fs,\ft,\ft_{o})$ is a local coordinates of $\mathscr{P}_{\varrho_{0}} glu_{\varrho_{0}}\left(U^{T,\sigma(\varrho_{0})}_{\mathbf a}\right)$.   On the other hand, since the bundle $\mathcal{N}$ has a Riemannian structure, we can
choose a smooth orthonormal frame field. This defines a coordinate $\mathfrak{y}$ over fiber.
   Then $(\hat{\fs},\hat{\ft},\mathfrak{y})$ is also a local coordinates of $\mathscr{P}_{\varrho_{0}} glu_{\varrho_{0}}\left(U^{T,\sigma(\varrho_{0})}_{\mathbf a}\right)$. Denote the Jacobi matrix  by $(a_{ij})=\frac{\p(\fs,\ft,\ft_{o})}{\p (\hat{\fs},\hat{\ft},\mathfrak{y})}$.
Since $\overline{\mathcal M}^{red}_{g,m+\nu}$ is a smooth orbifold, $(a_{ij})$ and the inverse matrix $(a^{-1}_{ij})$ are uniform bounded in the coordinates.
In particular
$$
\left|\pi_{*}(\frac{\p}{\p \ft^*_{i}})\right|+\left|\pi_{*}(\frac{\p}{\p \bar \ft^*_{i}})\right|\leq Ce^{-\fc \varrho_{0}},$$
$$\left|\pi_{*}(\frac{\p}{\p \fs_{i}}|_{(\fs,\ft,\ft_{0})})- \frac{\p}{\p \fs_{i}}|_{(\hat {\fs},\hat{\ft},0)}\right|+\left|\pi_{*}(\frac{\p}{\p \ft_{i}}|_{(\fs,\ft,\ft_{0})})- \frac{\p}{\p \ft_{i}}|_{(\hat {\fs},\hat{\ft},0)}\right|\leq Ce^{-\fc \varrho_{0}} ,\;\;\;
$$
On the other hand, $|\ft^*|<C$, where $C$ is a constant independent of $\varrho_{0}.$ As in the proof of Theorem  \ref{Conver}, we have
\begin{equation}
|(III)|\leq Ce^{-\fc \varrho_{0}}.
\end{equation}

By Lemma  \ref{gluing_rho_est} we have for any $X_{1},\cdots,X_{c}$
\begin{equation}
|(glu_{\varrho_{0}}^{*}(e^*_j\alpha_j) -e'^*_j\alpha_j)(X_{1},\cdots,X_{c})|^2\leq Ce^{-\fc_{1} \varrho_{0}}\Pi_{i=1}^{\mathsf r}g_{loc}(X_{i},X_{i})
\end{equation}
By Lemma \ref{est_Thom_varrho_L}, we have
$$\left| (glu_{\varrho_{0}}^{*}(\sigma_{\varrho_{0}}^{*}\Theta_{\varrho_{0}})-\sigma^{*}\Theta_{\infty})(X_{1},\cdots,X_{a})\right|^2 \leq C_1 e^{-\fc_{1}\varrho_{0}}.$$
It follows that
$$
\left|\Psi_{(M_{\infty}, A,g,n)}(\{\alpha _i\})
-\Psi_{(M_{\varrho_0},A, g,n)}(\{\alpha _i\})\right|\leq  C_{5} \left(\frac{1}{\varrho_{0}}+\frac{1}{|\log\sigma'|}\right)
$$
Then by letting $\varrho_{0}\rightarrow \infty$ and $\sigma'\to 0$, using Lemma \ref{lem_cob}
we get the Theorem.

\section{Some calculations for $\Psi_{C}$}
 \vskip 0.1in

We derive a gluing formula for the component $C=\{A^+, g^+, m^+, k; A^-, g^-,
m^-, k\}$. For any component $C$ we can use this
formula repeatedly. Choose a homology basis $\{\beta_b\}$ of $H^{\ast}(Z,
{\mathbb{R}})$. Let $(\delta_{ab} )$ be its intersection matrix. \vskip
0.1in \noindent
\begin{theorem}
 Let $\alpha^{\pm}_i$ be
differential forms with $deg \alpha^{+}_i=deg \alpha^{-}_i$ even.
Suppose that $\alpha^+_i|_Z=\alpha^-_i|_Z$ and hence
$\alpha^+_i\cup_Z \alpha^-_i\in H^*(\overline{M}^+ \cup_Z
\overline{M}^-, \mathbb R)$. Let $\alpha_i=\pi^*(\alpha^+_i\cup_Z
\alpha^-_i)$. For $C=\{A^+, g^+, m^+, k; A^-, g^-, m^-,
k\}$, we have the gluing formula
\begin{equation}\Psi_{C}(\alpha_1,...,\alpha_{m^++m^-}) =k\sum
\delta^{ab}\Psi^{(\overline{M}^+,Z)}_{(A^+,g^+,m^+,k)}
(\alpha^{+}_1,...,\alpha^{+}_{m^+}, \beta_a)
\Psi^{(\overline{M}^-,Z)}_{(A^-,g^-,m^-,k)}(\alpha^{-}_{m^++1},...,
\alpha^{-}_{m^++m^-},\beta_b).
\end{equation}
 where we use
$\Psi^{(\overline{M}^{\pm},Z)}_{(A^{\pm},g^{\pm},m^{\pm},k)}$ to
denote
$\Psi^{(\overline{M}^{\pm},Z)}_{(A^{\pm},g^{\pm},T_{m^{\pm}})}$.
\end{theorem}
\noindent {\bf Proof: }  We denote by $\mathbf{U}_{C}$ the virtual neighborhood for ${\mathcal  M}_C$, and by $\widehat{\mathbf{U}}_{C}$ the virtual neighborhood for ${\widehat{\mathcal  M}}_C$.
There is a natural map of degree k
$$Q:\mathbf{U}_{C}\rightarrow \widehat{\mathbf{U}}_C.$$
Let $\Delta \subset Z\times Z $ be the diagonal. The Poincar\'e dual
$\Delta^{\ast}$ of $\Delta $ is $$\Delta^{\ast}=\Sigma
\delta^{ab}\beta_a \wedge \beta_b .$$ Choose a Thom form $\Theta=
\Theta^{+}\wedge \Theta^{-}$, where $\Theta^{\pm}$ are Thom forms
in $\mathbf{E}^{\pm}$ supported in a neighborhood of the zero section .
Then
$$\Psi_{C}(\alpha_1,...,\alpha_{m^++m^-})=\int_{U_{C}}
\prod_{1}^{m^+}\alpha_i\wedge \prod_{m^+}^{m^++m^-}\alpha_j\wedge
I^{\ast}\Theta $$ $$= k
\int_{\widehat{\mathbf{U}}_{C}}\prod_{1}^{m^+}\alpha_i\wedge
\prod_{m^+}^{m^++m^-}\alpha_j\wedge I^{\ast}\Theta $$
  $$=k\sum
\delta^{ab}\Psi^{(\overline{M}^+,Z)}_{(A^+,g^+,m^+,k)}(\alpha^{+}_1,
...,\alpha^{+}_{m^+},\beta_a)
\Psi^{(\overline{M}^-,Z)}_{(A^-,g^-,m^-,k)}(\alpha^{-}_{m^++1},...,
\alpha^{-}_{m^++m^-},\beta_b). \;\;\;\;\Box $$ \vskip 0.1in
\noindent For general $C=\{A^+, g^+, m^+, {\bf k}; A^-, g^-,
m^-, {\bf k}\}$, where ${\bf k}=(k_1,...,k_{\mu})$
we may easily obtain \vskip 0.1in \noindent
\begin{theorem}\label{glu formula}
\begin{equation}\Psi_C(\alpha)=\frac{|{\bf k}|}{\mu!}\sum_{I,J}\Psi^{(\overline{M}^+,Z)}_
{(A^+,g^+,m^+,{\bf k})}(\alpha^{+}, \beta_I)\delta^{I,J}
\Psi^{(\overline{M}^-,Z)}_{(A^-,g^-,m^-,{\bf
k})}(\alpha^{-},\beta_J),
\end{equation}
 where we associate
$\beta_i\delta^{i,j}\beta_j$ to every periodic orbit, and put $|{\bf k}|=k_1...k_{\mu}$,
$\delta^{I,J}=\delta^{i_1,j_1}...\delta^{i_{\mu},j_{\mu}}$, and
denote by
$\Psi^{(\overline{M}^{\pm},Z)}_{(A^{\pm},g^{\pm},m^{\pm},{\bf
k})}(\alpha^{\pm},\beta_J)$ the product of relative invariants
cooresponding to each component.
\end{theorem}
 \v  For
example, for $C=\{A^+, g^+, m^+, k_1,k_2; ,A_1^-, g_1^-, m_1^-,
k_1, A_2^-, g_2^-, m_2^-, k_2\}$, our formula \eqref{glu formula} reads:
$$\Psi_C(\alpha)=\frac{1}{2}k_1k_2\sum_{i_1,i_2,j_1,j_2}\Psi^{(\overline{M}^+,Z)}_
{(A^+,g^+,m^+,k_1,k_2)}(\alpha^{+},
\beta_{i_1},\beta_{i_2})\delta^{i_1,j_1} \delta^{i_2,j_2}$$
$$\Psi^{(\overline{M}^-,Z)}_{(A_1^-, g_1^-, m_1^-,
k_1)}(\alpha_1^{-},\beta_{j_1}) \Psi^{(\overline{M}^-,Z)}_{(A_2^-,
g_2^-, m_2^-, k_2)}(\alpha_2^{-},\beta_{j_2}).$$ \vskip 0.1in
\v
The $\mu!$ comes from the fact that there is isotropy group $Aut(b)$, which exchange the puncture point $p_i$.
\v\v

\backmatter

%

\bibliographystyle{amsalpha}

\end{document}